\documentclass[]{article}

\usepackage{subcaption}
\usepackage{ stmaryrd }
\usepackage{graphicx}
\usepackage{amsmath}
\usepackage{amssymb}
\usepackage{amsthm}
\usepackage{pxfonts}
\usepackage{enumerate}
\usepackage{color}
\usepackage{mathdots}
\usepackage{sectsty}
\usepackage{tikz}
\usepackage{adjustbox}
\usepackage{enumitem}
\usepackage{caption}
\usepackage[hidelinks]{hyperref}
\usepackage{mathrsfs}

\allowdisplaybreaks
\usepackage{fancyhdr}
\usepackage[nottoc,notlot,notlof]{tocbibind}

\usetikzlibrary{decorations.markings,decorations.pathreplacing}
\usetikzlibrary{arrows.meta}

\newcommand\shorttitle{Non-colliding inhomogeneous Markov chains}
\newcommand\authors{Theodoros Assiotis and Miles Foster Nyamundanda}

\newtheoremstyle{claimstyle}
  {\topsep}       
  {\topsep}       
  {\itshape}      
  {}              
  {\bfseries}     
  {:}             
  {.5em}          
  {}              
\theoremstyle{plain}
\newtheorem{theorem}{Theorem}[section]

\newtheorem{lemma}[theorem]{Lemma}
\newtheorem{proposition}[theorem]{Proposition}
\newtheorem{definition}[theorem]{Definition}
\theoremstyle{claimstyle}

\newtheorem{assumption}[theorem]{Assumption}

\newtheorem{example}[theorem]{Example}
\newtheorem{cor}[theorem]{Corollary}
\newtheorem{remark}[theorem]{Remark}
\newtheorem*{theorem*}{Theorem}

\title{\large \bf NON-COLLIDING SPACE-TIME INHOMOGENEOUS MARKOV CHAINS}
\author{\small THEODOROS ASSIOTIS\thanks{School of Mathematics, University of Edinburgh, James Clerk Maxwell Building, Peter Guthrie Tait Road, Edinburgh EH9 3FD, United Kingdom. Email: \texttt{theo.assiotis@ed.ac.uk}.} \;\;\textnormal{AND}\;\; MILES FOSTER NYAMUNDANDA\thanks{School of Mathematics, University of Edinburgh, James Clerk Maxwell Building, Peter Guthrie Tait Road, Edinburgh EH9 3FD, United Kingdom. Email: \texttt{s2432435@ed.ac.uk}.}}
\date{}

\begin{document}
\maketitle
\begin{abstract}
\noindent We establish the explicit leading order asymptotics, with a quantitative error bound, of tail probabilities of collision times for a class of integrable space-time inhomogeneous Markov chains, in discrete and continuous time. The corresponding process conditioned not to intersect arises in interacting particle systems with local push-block interactions thereby confirming a prediction from \cite{assiotis2023integrablemodelsinhomogeneousspace}. The generic discrete nature of the spatial inhomogeneities rules out powerful coupling-with-Brownian-motion techniques, so our proof strategy proceeds instead via a novel steepest-descent analysis combined with a Karlin--McGregor semigroup expansion in terms of dominant-index contributions.
\end{abstract}

\tableofcontents

\section{Introduction}

\subsection{Setting of the problem}
The purpose of this paper is to compute explicit asymptotics in Theorem \ref{collisionThm} for tail probabilities of collision times for certain integrable space-time inhomogeneous Markov chains which include general rate pure-birth chains and discrete-time chains with one-step inhomogeneous Bernoulli and geometric updates. This allows us to construct the process conditioned never to intersect and identify it with an explicit Doob $h$-transform \cite{Doob,RevuzYor} of the original independent chains in Corollary \ref{Thm1}. To solve this problem we need to develop some techniques to deal with asymptotics of determinants with entries involving a diverging number of parameters, which represent the space-time inhomogeneous environment, having little global structure. We begin with some motivation for this problem.

Conditioning a stochastic process to stay in a subset of its state space is a topic with a long and storied history in probability theory \cite{Doob}. Some of the most classical examples are conditioning a one-dimensional random walk or Brownian motion to stay non-negative \cite{McKean,Spitzer,Iglehart,Pitman,Bolthausen,BertoinDoney}. In higher dimensions things become more intricate. Arguably one of the most important domains is the Weyl chamber with $N$ ordered coordinates where the corresponding one-dimensional coordinates give rise to non-intersecting paths.

Random models of non-intersecting paths have been extensively studied in the past three decades in relation to random matrices \cite{Dyson1962,Grabiner,KonigOConnell,KatoriTanemuraMatrix,KatoriTanemuraFunctional}, random tiling models \cite{OkounkovReshetikhin,Johansson2002,BorodinGorin,BorodinGorinRains,DuitsKuijlaars,CharlierHexagon,RailYard} interacting particle systems \cite{BorodinFerrari,BorodinFerrariTilings,WarrenWindridge,Nordenstam,NikosTASEP,IwaoMotegiScrimshaw}, representation theory \cite{BianeBougerolOConnell,DoumercOConnell}, the Robinson-Schensted-Knuth correspondence \cite{OConnellRSK,OConnellYor} and KPZ universality class \cite{CorwinKPZ}. Moreover, lattice paths constrained to stay in domains have received a lot of attention in the enumerative combinatorics literature, see for example \cite{BousquetMelouMishna,KurkovaRaschel,RaschelJEMS,RaschelSPA,FayolleRaschel,WeightedWalks}. Nevertheless, despite this attention the asymptotic behaviour of weighted lattice paths with general weights depending on the individual vertices of the lattice is still poorly understood. Our discrete-time results when rephrased in these terms make advances in this direction for special models.

In the case of multidimensional Brownian motion detailed asymptotics of the distribution of exit times from general Euclidean space domains were obtained in \cite{BanuelosSmits,DeBlassie}. The connection to random matrices and Dyson Brownian motion was developed in \cite{Grabiner} using the celebrated Karlin-McGregor formula for non-coincidence probabilities \cite{karlin1959coincidence}. For conditioning special random walks (i.e. sum of independent identically distributed random increments) to stay in Weyl chambers early studies were based on the Karlin-McGregor and Lindstr\"{o}m–Gessel–Viennot formula \cite{gessel1985binomial,lindstrom1973vector} and generalisations, see \cite{OConnellYor,OConnellRSK,OrderedRW,KonigOConnellRoch,KonigSchmid}. Then, a very general approach was developed in \cite{DenisovWachtel1} which allowed to deal with both general increment distributions and general cones. This is based on a powerful coupling with Brownian motion technique which essentially allows to transfer the explicit results obtained using analysis of the Laplacian for the Brownian case \cite{BanuelosSmits,DeBlassie}. This led to many developments, see \cite{DenisovWachtel3,DenisovFitzGerald,DenisovFitzGeraldWachtel2024,RaschelTarrago,RaschelTarragoSurvey,JetlirRaschelTarragoWachtel}  for related subsequent works. More recently, by refining this approach, the asymptotics of tail probabilities of collision times for a quite general class of Markov chains, essentially ones which under diffusive scaling converge to Brownian motion (modulo some additional simplifying technical conditions) have been established in \cite{DenisovZhang}. We stress that, as far as we can tell, the generic space-inhomogeneous chains we will consider here do not converge to Brownian motion under some scaling and interestingly, due to the various assumptions in \cite{DenisovZhang} the only overlap with \cite{DenisovZhang} in terms of results turns out to be the homogeneous case only.

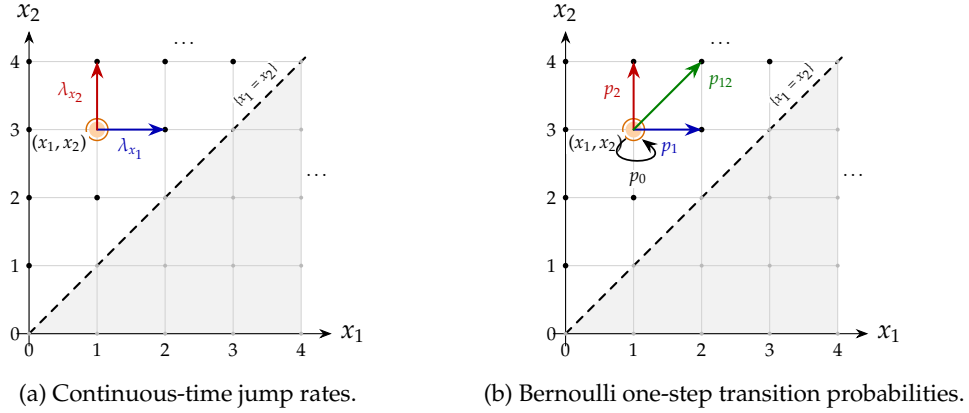
\begin{figure}[ht!]
        \centering
        \begin{subfigure}[t]{0.48\textwidth}
        \centering
        \begin{tikzpicture}[x=0.90cm,y=0.90cm,>=Stealth,line cap=round,line join=round]
            \tikzset{
                xonejump/.style={->,draw=blue!70!black,line width=0.75pt},
                xtwojump/.style={->,draw=red!75!black,line width=0.75pt},
                gridline/.style={draw=gray!35,line width=0.25pt},
                statelabel/.style={font=\scriptsize,inner sep=0.7pt,fill=white},
                ratelabel/.style={font=\scriptsize,inner sep=0.8pt,fill=white}
            }
            \fill[gray!10] (0,0)--(4,0)--(4,4)--cycle;
            \draw[step=1,gridline] (0,0) grid (4,4);
            \draw[->] (-0.18,0) -- (4.45,0) node[right] {$x_1$};
            \draw[->] (0,-0.18) -- (0,4.45) node[above] {$x_2$};
            \foreach \k in {0,...,4}{
                \node[below] at (\k,0) {\scriptsize$\k$};
                \node[left] at (0,\k) {\scriptsize$\k$};
            }
            \draw[dashed, thick] (0,0) -- (4.12,4.12);
            \node[rotate=45,fill=white,inner sep=0.6pt] at (3.34,3.65) {\tiny$\{x_1=x_2\}$};
            \foreach \i in {0,...,4}{
                \foreach \j in {0,...,4}{
                    \ifnum\j>\i
                        \fill (\i,\j) circle (1.05pt);
                    \else
                        \fill[gray!55] (\i,\j) circle (0.75pt);
                    \fi
                }
            }
            \coordinate (q) at (1,3);
            \fill[orange!35] (q) circle (3.0pt);
            \draw[orange!85!black,line width=0.55pt] (q) circle (4.1pt);
            \node[statelabel,anchor=east] at (0.88,2.80) {$(x_1,x_2)$};
            \draw[xonejump] (q) -- (2,3);
            \draw[xtwojump] (q) -- (1,4);
            \node[ratelabel,text=blue!70!black] at (1.50,2.72) {$\lambda_{x_1}$};
            \node[ratelabel,text=red!75!black,anchor=east] at (0.86,3.55) {$\lambda_{x_2}$};
            \node at (4.25,2.3) {\scriptsize$\cdots$};
            \node at (2.3,4.25) {\scriptsize$\cdots$};
        \end{tikzpicture}
        \caption{Continuous-time jump rates.}
        \end{subfigure}
        \hfill
        \begin{subfigure}[t]{0.48\textwidth}
        \centering
        \begin{tikzpicture}[x=0.90cm,y=0.90cm,>=Stealth,line cap=round,line join=round]
            \tikzset{
                xonejump/.style={->,draw=blue!70!black,line width=0.75pt},
                xtwojump/.style={->,draw=red!75!black,line width=0.75pt},
                bothjump/.style={->,draw=green!50!black,line width=0.8pt},
                stayput/.style={->,draw=black,line width=0.65pt},
                gridline/.style={draw=gray!35,line width=0.25pt},
                statelabel/.style={font=\scriptsize,inner sep=0.7pt,fill=white},
                ratelabel/.style={font=\scriptsize,inner sep=0.7pt,fill=white,align=center}
            }
            \fill[gray!10] (0,0)--(4,0)--(4,4)--cycle;
            \draw[step=1,gridline] (0,0) grid (4,4);
            \draw[->] (-0.18,0) -- (4.45,0) node[right] {$x_1$};
            \draw[->] (0,-0.18) -- (0,4.45) node[above] {$x_2$};
            \foreach \k in {0,...,4}{
                \node[below] at (\k,0) {\scriptsize$\k$};
                \node[left] at (0,\k) {\scriptsize$\k$};
            }
            \draw[dashed, thick] (0,0) -- (4.12,4.12);
            \node[rotate=45,fill=white,inner sep=0.6pt] at (3.34,3.65) {\tiny$\{x_1=x_2\}$};
            \foreach \i in {0,...,4}{
                \foreach \j in {0,...,4}{
                    \ifnum\j>\i
                        \fill (\i,\j) circle (1.05pt);
                    \else
                        \fill[gray!55] (\i,\j) circle (0.75pt);
                    \fi
                }
            }
            \coordinate (q) at (1,3);
            \fill[orange!35] (q) circle (3.0pt);
            \draw[orange!85!black,line width=0.55pt] (q) circle (4.1pt);
            \node[statelabel,anchor=east] at (0.88,2.80) {$(x_1,x_2)$};
            \draw[xonejump] (q) -- (2,3);
            \draw[xtwojump] (q) -- (1,4);
            \draw[bothjump] (q) -- (2,4);
            \draw[stayput] (0.88,2.86) .. controls (0.35,2.45) and (1.65,2.45) .. (1.12,2.86);
            \node[ratelabel,text=blue!70!black] at (1.54,2.72) {$p_1$};
            \node[ratelabel,text=red!75!black,anchor=east] at (0.86,3.55) {$p_2$};
            \node[ratelabel,text=green!45!black,anchor=west] at (2.08,3.70) {$p_{12}$};
            \node[ratelabel,text=black] at (1.07,2.28) {$p_0$};
            \node at (4.25,2.3) {\scriptsize$\cdots$};
            \node at (2.3,4.25) {\scriptsize$\cdots$};
        \end{tikzpicture}
        \caption{Bernoulli one-step transition probabilities.}
        \end{subfigure}
        \caption{Finite windows of the two-particle state graphs for continuous time pure-birth and discrete time Bernoulli dynamics. In both panels, the dashed diagonal $\{x_1=x_2\}$ is the collision set and the admissible states satisfy $x_1<x_2$; the shaded region is the non-admissible side of the diagonal. The orange point is a representative admissible state $(x_1,x_2)$, and the labelled arrows display the outgoing local transition rates or probabilities from that state. Blue denotes an $x_1$-jump, red denotes an $x_2$-jump, green denotes a simultaneous jump of both coordinates, and black denotes remaining at the same state. In the Bernoulli panel, $p_1=\lambda_{x_1}\alpha_t(1-\lambda_{x_2}\alpha_t)$, $p_2=(1-\lambda_{x_1}\alpha_t)\lambda_{x_2}\alpha_t$, $p_{12}=\lambda_{x_1}\lambda_{x_2}\alpha_t^2$, and $p_0=(1-\lambda_{x_1}\alpha_t)(1-\lambda_{x_2}\alpha_t)$.}
        \label{twoparticlestategraph}
        \end{figure}

Returning to our setting, the space-inhomogeneous Markov chains we consider are of course not completely arbitrary. See Figure \ref{twoparticlestategraph} for an illustration of the simplest case $N=2$ for pure-birth and Bernoulli dynamics; the geometric case is slightly more complicated, see Section \ref{noncolliding}. They enjoy some integrable structure in the sense of having explicit, but complicated, contour integral expressions for their transition probabilities. A consequence of this integrability will be the fact that the function governing the dependence on the initial condition in the tail asymptotics is completely explicit. It generalises the classical Vandermonde determinant. We note that, for general random walks \cite{DenisovWachtel1} the analogous function is not explicit but rather has a probabilistic expression which is unlikely to be explicitly computable.

These special Markov chains also arise in some integrable space-inhomogeneous generalisations \cite{assiotis2023integrablemodelsinhomogeneousspace} of (2+1)-dimensional (two space \& one  time dimension) growth models \cite{BorodinFerrari,BorodinFerrariTilings,WarrenWindridge} that couple (1+1)-dimensional systems such as the totally asymmetric simple exclusion process (TASEP) and PushTASEP \cite{BorodinFerrari,BorodinFerrariTilings} with non-colliding processes, see Section \ref{SectionInterlacingParticleSystem}. In fact, a consequence of our main result, in the form of Corollary \ref{Thm1}, was conjectured in \cite{assiotis2023integrablemodelsinhomogeneousspace}, see Remark \ref{Rmk:Prediction}. Models with inhomogeneities being more realistic physical toy models give rise to novel asymptotic behaviours, most famously in periodic tiling models, and have been extensively studied in the past ten years or so, see \cite{ChhitaYoung,ChhitaJohansson,JohanssonRahman,BorodinDuits,BerggrenBorodin,BerggrenDuits,DuitsKuijlaars,CharlierHexagon,Kuijlaars,bufetov2025dominotilingsaztecdiamond,zografos2025quenchedannealedcltsoneperiodic,duits2025gammadisorderedaztecdiamond,NicosTASEP,AssiotisKarlinMcGregor,NikosTASEP,IwaoMotegiScrimshaw,AggarwalBorodinPetroWheeler}. It is also worth mentioning that detailed analyses of space-inhomogeneous Markov chains have been performed in \cite{MustaphaGaussianEstimates,MustaphaSifiDiscreteHarmonic,EssifiMustaphaOrthant,DolgopyatSarig,DolgopyatHafouta}
. These works, whose setting is a bit different from ours, are very general and do not use any integrable structure at all but, as far as we can tell, do not give results to the asymptotic precision that we establish in our paper. 

Since, as far as we can tell, the generic inhomogeneous parameter chains we consider do not converge to Brownian motion we cannot directly use general results but rather rely on the integrable structure in the form of explicit formulae. Nevertheless, we are able to find some asymptotic Gaussian structure, for a single time, under a highly non-obvious implicit time-dependent transformation, see Theorem \ref{thm:stationarypoints}.

To the best of our knowledge, the present paper gives the first systematic quantitative study of asymptotics of intersection times of a class of generically discrete-space-inhomogeneous Markov chains (i.e. there is no natural continuous diffusion process behind them and we overall assume minimal global structure for the inhomogeneities). There are a number of allied models whose asymptotic analysis involves dealing with similar difficulties. It is plausible that aspects of the strategy we develop here could be adapted, after non-trivial effort, to those settings as well. We say more in Section \ref{SectionExtensions}.

In the rest of introduction we introduce the models precisely, present our main results, give an outline of the proof and discuss examples and extensions.

\subsection{Space-time inhomogeneous Markov chains}
\label{typeofbirth}

    We use the notation $\mathbb{N}\overset{\mathrm{def}}{=}\{1,2, 3,\dots\}$, $\mathbb{Z}_+\overset{\mathrm{def}}{=}\{0,1,2, 3,\dots\}$, and $\mathbb{R}_+\overset{\mathrm{def}}{=}\{z\in\mathbb{R}: z\geq0\}$. We define the basic one-dimensional Markov chains we will be studying.

    \begin{definition}[Continuous-time pure birth process]\label{contprocess}
        Let $(\lambda_n)_{n\in\mathbb{Z}_{+}}\subset(0,\infty)$. A non-explosive process $(\mathscr{X}(t):t\in\mathbb{R}_+)$ on $\mathbb{Z}_{+}$ is called a continuous-time pure birth process with rates $(\lambda_n)_{n\in \mathbb{Z}_+}$ if it is a time-homogeneous Markov process with generator
        \begin{equation*}
            (\mathsf{L}f)(y)\overset{\mathrm{def}}{=}\lambda_y\big(f(y+1)-f(y)\big),\qquad y\in\mathbb{Z}_{+},
        \end{equation*}
        for bounded test functions $f:\mathbb{Z}_{+}\to\mathbb{R}$. Equivalently, conditional on $\mathscr{X}(t)=y$, the process waits an exponential random time with parameter $\lambda_y$, $\mathsf{Exp}(\lambda_y)$, then jumps to $y+1$.
    \end{definition}
    \begin{figure}[ht]
        \centering
        \begin{tikzpicture}[scale=0.5]

	\def\nx{14}
	\def\ny{7}

	\foreach \x in {0,...,\numexpr\nx-1} {
	  \foreach \y in {0,...,\numexpr\ny} {
	            \draw[dotted] (\x,\y) -- (\x+1,\y);
	    
  }
}

\foreach \y in {0,...,\numexpr\ny} {
\node at (-0.5,\y) {\small$\y$};
}
  
\draw[thick,->] (0,0) -- (0.75,0);
\draw[thick] (0.75,0) -- (1.5,0);

\draw[thick,->] (1.5,0) -- (1.5,0.5);
\draw[thick] (1.5,0.5) -- (1.5,1);

\draw[thick,->] (1.5,1) -- (3,1);
\draw[thick] (3,1) -- (4.5,1);

\draw[thick,->] (4.5,1) -- (4.5,1.5);
\draw[thick] (4.5,1.5) -- (4.5,2);

\draw[thick,->] (4.5,2) -- (4.75,2);
\draw[thick] (4.75,2) -- (5,2);

\draw[thick,->] (5,2) -- (5,2.5);
\draw[thick] (5,2.5) -- (5,3);

\draw[thick,->] (5,3) -- (6,3);
\draw[thick] (6,3) -- (7,3);

\draw[thick,->] (7,3) -- (7,3.5);
\draw[thick] (7,3.5) -- (7,4);

\draw[thick,->] (7,4) -- (8,4);
\draw[thick] (8,4) -- (9,4);

\draw[thick,->] (9,4) -- (9,4.5);
\draw[thick] (9,4.5) -- (9,5);

\draw[thick,->] (9,5) -- (9.75,5);
\draw[thick] (9.75,5) -- (10.5,5);

\draw[thick,->] (10.5,5) -- (10.5,5.5);
\draw[thick] (10.5,5.5) -- (10.5,6);

\draw[thick,->] (10.5,6) -- (12,6);
\draw[thick] (12,6) -- (13.5,6);

\draw[thick,->] (13.5,6) -- (13.5,6.5);
\draw[thick] (13.5,6.5) -- (13.5,7);

	\draw[thick,->] (13.5,7) -- (14,7);
	\draw[->] (1.5,-0.5) -- (0,-0.5) node[midway, below] {\tiny{$E_0$}};
	\draw[->] (0,-0.5) -- (1.5,-0.5);
	\draw[->] (1.5,-0.5) -- (4.5,-0.5) node[midway, below] {\tiny{$E_1$}};
	\draw[->] (4.5,-0.5) -- (1.5,-0.5);
	\draw[->] (4.75,-0.5) -- (5,-0.5) node[midway, below] {\tiny{$E_2$}};
	\draw[->] (5,-0.5) -- (4.5,-0.5);
	\draw[->] (5,-0.5) -- (7,-0.5) node[midway, below] {\tiny{$E_3$}};
	\draw[->] (7,-0.5) -- (5,-0.5);
	\draw[->] (7,-0.5) -- (9,-0.5) node[midway, below] {\tiny{$E_4$}};
	\draw[->] (9,-0.5) -- (7,-0.5);
	\draw[->] (9,-0.5) -- (10.5,-0.5) node[midway, below] {\tiny{$E_5$}};
	\draw[->] (10.5,-0.5) -- (9,-0.5);
	\draw[->] (10.5,-0.5) -- (13.5,-0.5) node[midway, below] {\tiny{$E_6$}};
	\draw[->] (13.5,-0.5) -- (10.5,-0.5);

\end{tikzpicture}
        \caption{Sample path of a continuous-time pure-birth process started from $0$. The lower labels $E_j$ denote the successive independent random holding times, with $E_j$ distributed as $\mathsf{Exp}(\lambda_j)$.}
        \label{fig:continuous-birth-path}
    \end{figure}
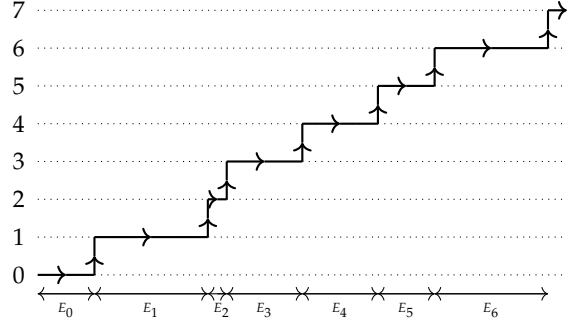
    \begin{definition}[Mixed Bernoulli-and-geometric jump process]\label{mixprocess}
    Let $(\lambda_n)_{n\in\mathbb{Z}_{+}}\subset(0,\infty)$, $(\alpha_t)_{t\in\mathbb{Z}_{+}}\subset[0,\infty)$, $(\beta_t)_{t\in\mathbb{Z}_{+}}\subset[0,\infty)$, and let $\tau=(\tau_t)_{t\in\mathbb{Z}_{+}}\in\{0,1\}^{\mathbb{Z}_{+}}$.
    A process $(\mathscr{X}(t):t\in\mathbb{Z}_{+})$ on $\mathbb{Z}_{+}$ is called a mixed Bernoulli-and-geometric jump process if, under the stochasticity conditions in Remark~\ref{stochasticconditions}, it is a time-inhomogeneous Markov chain with one-step kernel
    \begin{equation*}
        p_t(y,z)\overset{\mathrm{def}}{=}\mathbb{P}\!\left(\mathscr{X}(t+1)=z\mid \mathscr{X}(t)=y\right),
        \qquad y,z,t\in\mathbb{Z}_{+},
    \end{equation*}
    given as follows.
    \begin{enumerate}
        \item If $\tau_t=0$ (Bernoulli jump), then
        \begin{equation*}
            p_t(y,y+1)=\lambda_y\alpha_t,\qquad
            p_t(y,y)=1-\lambda_y\alpha_t,\qquad
            p_t(y,z)=0\ \text{for }z\notin\{y,y+1\}.
        \end{equation*}
        \item If $\tau_t=1$ (geometric jump), then for $k\in\mathbb{Z}_{+}$,
        \begin{equation*}
            p_t(y,y+k)
            =\frac{1}{1+\beta_t\lambda_{y+k}}
            \prod_{i=y}^{y+k-1}\left(\frac{\beta_t\lambda_i}{1+\beta_t\lambda_i}\right),
        \end{equation*}
        where an empty product is interpreted as $1$, and $p_t(y,z)=0$ for $z<y$.
    \end{enumerate}
    \end{definition}
    \begin{remark}\label{stochasticconditions}
        In order that the above transition mechanism define a stochastic kernel, we require
        \begin{equation*}
            0\leq \lambda_n\alpha_t\leq 1,\quad\text{for all }n,t\in\mathbb{Z}_+,\qquad\text{and}\qquad\lim_{m\to\infty}\prod_{i=y}^{y+m}\left(\frac{\beta_t\lambda_i}{1+\beta_t\lambda_i}\right)=0,\quad \text{for all }y,t\in\mathbb{Z}_+.
        \end{equation*}
        The first condition is the admissibility condition for the Bernoulli step, while the second is the non-defectiveness condition for the inhomogeneous geometric step. They ensure that the one-step transition probabilities have total mass one. These conditions follow from \textbf{(A1)}, introduced later in Section~\ref{generalassumption}, so we do not treat them separately here.
    \end{remark}

     For the one-particle models in Definition~\ref{contprocess} and Definition~\ref{mixprocess}, we write
    \begin{equation*}
        \mathscr{T}\overset{\textnormal{def}}{=}
        \begin{cases}
            \mathbb{R}_+,&\text{for Definition~\ref{contprocess}},\\
            \mathbb{Z}_+,&\text{for Definition~\ref{mixprocess}}.
        \end{cases}
    \end{equation*}
    With the parameter sequences fixed, for each initial condition $x\in\mathbb{Z}_+$, we denote by $\mathbb{P}_{x}$ the law of $(\mathscr{X}(t):t\in\mathscr{T})$ with $\mathscr{X}(0)=x$, and by $\mathbb{E}_{x}$ the corresponding expectation. For $\mathbf{x}=(x_1,\ldots,x_N)\in\mathbb{Z}_+^N$, we denote by $\mathbb{P}_{\mathbf{x}}^{N}$ the joint law of $N$ independent copies started from $\mathbf{x}$.
    
    \begin{figure}
        \centering
        \begin{tikzpicture}[scale=0.5]

	\def\nx{14}
	\def\ny{10}
\fill[blue, opacity=0.1] (0,0) rectangle (6,10);
\fill[red, opacity=0.1] (6,0) rectangle (10,10);
\fill[blue, opacity=0.1] (10,0) rectangle (12,10);
\fill[red, opacity=0.1] (12,0) rectangle (13,10);
\fill[blue, opacity=0.1] (13,0) rectangle (14,10);
\foreach \x in {0,1,2,3,4,5,10,11,13} {
  \foreach \y in {0,...,\numexpr\ny-1} {
	            \draw[dotted] (\x,\y) -- (\x+1,\y);
	    
	            \draw[dotted] (\x,\y) -- (\x+1,\y+1);
	  }
	}
	
	\foreach \x in {6,7,8,9,12} {
	  \foreach \y in {0,...,\numexpr\ny-1} {
	            \draw[dotted] (\x,\y) -- (\x+1,\y);
	    
	            \draw[dotted] (\x,\y) -- (\x,\y+1);
	  }
}

\foreach \y in {0,...,\numexpr\ny} {
\node at (-0.5,\y) {\small$\y$};
}
\foreach \x in {0,...,\numexpr\nx} {
\node at (\x,-0.5) {\small$\x$};
}
  
	\draw[thick,->] (0,0) -- (0.5,0);
	\draw[thick] (0.5,0) -- (1,0);

	\draw[thick,->] (1,0) -- (1.5,0.5);
	\draw[thick] (1.5,0.5) -- (2,1);

	\draw[thick,->] (2,1) -- (2.5,1.5);
	\draw[thick] (2.5,1.5) -- (3,2);

	\draw[thick,->] (3,2) -- (3.5,2);
	\draw[thick] (3.5,2) -- (4,2);

	\draw[thick,->] (4,2) -- (4.5,2);
	\draw[thick] (4.5,2) -- (5,2);

	\draw[thick,->] (5,2) -- (5.5,2);
	\draw[thick] (5.5,2) -- (6,2);

\draw[thick,->] (6,2) -- (6,2.5);
\draw[thick] (6,2.5) -- (6,3);
\draw[thick,->] (6,3) -- (6,3.5);
\draw[thick] (6,3.5) -- (6,4);
	\draw[thick,->] (6,4) -- (6,4.5);
	\draw[thick] (6,4.5) -- (6,5);

\draw[thick,->] (6,5) -- (6.5,5);
\draw[thick] (6.5,5) -- (7,5);
\draw[thick,->] (7,5) -- (7.5,5);
\draw[thick] (7.5,5) -- (8,5);
\draw[thick,->] (8,5) -- (8,5.5);
\draw[thick] (8,5.5) -- (8,6);
\draw[thick,->] (8,6) -- (8.5,6);
\draw[thick] (8.5,6) -- (9,6);
\draw[thick,->] (9,6) -- (9,6.5);
\draw[thick] (9,6.5) -- (9,7);
\draw[thick,->] (9,7) -- (9,7.5);
\draw[thick] (9,7.5) -- (9,8);

\draw[thick,->] (9,8) -- (9.5,8);
\draw[thick] (9.5,8) -- (10,8);
\draw[thick,->] (10,8) -- (10.5,8);
\draw[thick] (10.5,8) -- (11,8);
\draw[thick,->] (11,8) -- (11.5,8.5);
\draw[thick] (11.5,8.5) -- (12,9);
\draw[thick,->] (12,9) -- (12,9.5);
\draw[thick] (12,9.5) -- (12,10);
\draw[thick,->] (12,10) -- (12.5,10);
\draw[thick] (12.5,10) -- (13,10);
\draw[thick,->] (13,10) -- (13.5,10);
\draw[thick] (13.5,10) -- (14,10);

\end{tikzpicture}
        \caption{Sample path of a mixed Bernoulli-and-geometric jump process. Blue regions indicate Bernoulli updates and red regions indicate geometric updates.}
        \label{mixbernoulligeofig}
    \end{figure}
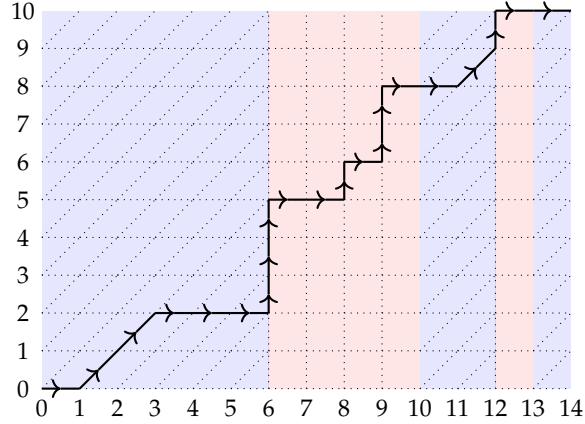
    \begin{figure}[ht!]
    \centering
    \begin{subfigure}{0.45\textwidth}
        \centering
        \begin{tikzpicture}[scale=0.4]
            \node at (-0.5,0) {\small{$x$}};
            \node at (-0.5,4) {\small{$x+1$}};

            \draw[thick,->] (0,0)--(4,0) node[midway, below] {\small{$1-\alpha_{t}\lambda_{x}$}};
            \draw[thick,->] (0,0)--(4,4) node[midway, above left] {\small{$\alpha_{t}\lambda_{x}$}};

        \end{tikzpicture}
        \caption{One-step Bernoulli transition probabilities from $x\in\mathbb{Z}_+$.}
        \label{fig:single}
    \end{subfigure}
    \hfill
    \begin{subfigure}{0.48\textwidth}
        \centering
        \begin{tikzpicture}[scale=0.6]
            \node at (-0.5,0) {\tiny{$x$}};
            \node at (-0.5,1) {\tiny{$x+1$}};
            \node at (-0.5,2) {\tiny{$x+2$}};
            \node at (-0.5,3) {\tiny{$x+3$}};
            \node at (-0.5,5) {\tiny{$x+n$}};

            \draw[thick,->] (0,0)--(0,3);
            \draw[thick,->] (0,0)--(1,0);
            \draw[thick,->] (0,1)--(1,1);
            \draw[thick,->] (0,2)--(1,2);
            \draw[thick,->] (0,3)--(1,3);
            \draw[dashed,->] (0,3)--(0,5);
            \draw[thick,->] (0,5)--(1,5);

            \node[right] at (1.2,0) {\tiny$\frac{1}{1+\beta_t\lambda_{x}}$};
            \node[right] at (1.2,1) {\tiny$\frac{1}{1+\beta_t\lambda_{x+1}}\left(\frac{\beta_t\lambda_x}{1+\beta_t\lambda_x}\right)$};
            \node[right] at (1.2,2) {\tiny$\frac{1}{1+\beta_t\lambda_{x+2}}\prod_{i=x}^{x+1}\left(\frac{\beta_t\lambda_i}{1+\beta_t\lambda_i}\right)$};
            \node[right] at (1.2,3) {\tiny$\frac{1}{1+\beta_t\lambda_{x+3}}\prod_{i=x}^{x+2}\left(\frac{\beta_t\lambda_i}{1+\beta_t\lambda_i}\right)$};
            \node[right] at (1.2,5) {\tiny$\frac{1}{1+\beta_t\lambda_{x+n}}\prod_{i=x}^{x+n-1}\left(\frac{\beta_t\lambda_i}{1+\beta_t\lambda_i}\right)$};
        \end{tikzpicture}
        \caption{One-step geometric transition probabilities from $x\in\mathbb{Z}_+$.}
        \label{fig:multi}
    \end{subfigure}
\end{figure}
\paragraph*{Assumptions}\label{generalassumption}
    Beyond the basic well-posedness input, our analysis relies on the following uniform control assumptions on the parameter sequences.
    \begin{assumption}\label{A1}
\textbf{(A1)} There exist constants $\mathfrak{L},\mathfrak{B}_1,\mathfrak{B}_2,\mathfrak{U}>0$ such that
\begin{equation*}
    \mathfrak{L}<\mathfrak{B}_1<\mathfrak{B}_2<\mathfrak{U}.
\end{equation*}
\textbf{(i) Continuous-time regime:} For the continuous-time pure birth process of Definition~\ref{contprocess},
the spatial rates satisfy
\begin{equation}\label{A1:continuous}
    \mathfrak{L}<\lambda_n< \mathfrak{B}_1,
    \qquad n\in\mathbb{Z}_+.
\end{equation}
\textbf{(ii) Discrete-time regime:} For the mixed Bernoulli-and-geometric jump process of
Definition~\ref{mixprocess},
the spatial parameters satisfy \eqref{A1:continuous}, and the temporal parameters obey
\begin{equation}\label{A1:discrete}
    \mathfrak{B}_2<\alpha_t^{-1}< \mathfrak{U},
    \qquad
    \mathfrak{L}<\beta_t<\mathfrak{U},
    \qquad t\in\mathbb{Z}_+.
\end{equation}
\end{assumption}

Enlarging $\mathfrak{U}$ if necessary we henceforth assume $\mathfrak{B}_2^{-1} \le \mathfrak{U}$. Then, $\alpha_t \le \mathfrak{U}$, $\beta_t \le \mathfrak{U}$ for all $t$.

 By \textbf{(A1)}, the stochastic-kernel conditions from Remark~\ref{stochasticconditions} hold, in particular by noting,
	\begin{equation*}
	    0\leq\prod_{i=y}^{y+m}\left(\frac{\beta_t\lambda_i}{1+\beta_t\lambda_i}\right)
	    \leq\left(\frac{\mathfrak{U}\mathfrak{B}_1}{1+\mathfrak{U}\mathfrak{B}_1}\right)^{m+1},
	    \qquad y,t\in\mathbb{Z}_+,
	\end{equation*} so the standing assumptions in Definition~\ref{contprocess} and Definition~\ref{mixprocess} are satisfied. Finally, the following assumption is imposed to identify the spatial reciprocal-rate centering used below.
        
    \begin{assumption}\label{A2}
    \textbf{(A2)} We assume that
    \begin{equation*}
        \quad\Lambda\overset{\textnormal{def}}{=}\lim_{N\rightarrow\infty}\frac{1}{N+1}\sum_{n=0}^{N}\frac{1}{\lambda_{n}}\in(0,\infty).
    \end{equation*}
    \end{assumption}
    
    \subsection{Non-colliding process}\label{noncolliding}
    
    We now define precisely collision times and conditioning on infinite-time horizon non-collision.
        \begin{definition}[Weyl chamber]
            For $N\in\mathbb{N}$, we define the strict ordered discrete Weyl chamber:
        \begin{equation*}
            \mathbf{W}_{\mathrm d}^{N}\overset{\textnormal{def}}{=}\{\mathbf{y}=(y_1,\ldots,y_N)\in\mathbb{Z}_+^{N}: y_1<\cdots<y_N\}.
        \end{equation*}
        \end{definition}

        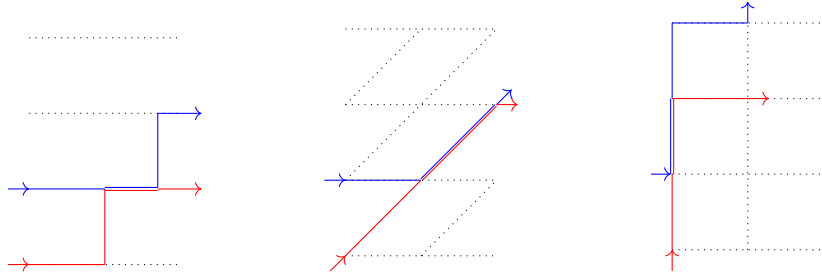
\begin{figure}[ht!]
        \centering
            \begin{subfigure}{0.3\textwidth}
            \centering
                \begin{tikzpicture}[scale=1]
                    \draw[dotted] (0,0)--(2,0);
                    \draw[dotted] (0,2)--(2,2);
                    \draw[dotted] (0,3)--(2,3);
                
                    \draw[->, red] (-0.28,0)--(0,0);
                    \draw[-, red] (0,0)--(1,0);
                    \draw[-, red] (1,0)--(1,1);
                    \draw[-, red] (1,0.98)--(1.7,0.98);
                    \draw[->, red] (1.7,1)--(2.28,1);

                    \draw[->, blue] (-0.28,1)--(0,1);
                    \draw[-,blue] (0,1)--(1,1);
                    \draw[-, blue] (1,1.02)--(1.7,1.02);
                    \draw[-, blue] (1.7,1.02)--(1.7,2);
                    \draw[->, blue] (1.7,2)--(2.28,2);
                \end{tikzpicture}
            \end{subfigure}
            \begin{subfigure}{0.3\textwidth}
            \centering
                \begin{tikzpicture}[scale=1]
                    \draw[dotted] (0,0)--(2,0);
                    \draw[dotted] (0,1)--(2,1);
                    \draw[dotted] (0,2)--(2,2);
                    \draw[dotted] (0,1)--(2,3);
                    \draw[dotted] (1,0)--(2,1);
                    \draw[dotted] (0,2)--(1,3);
                    \draw[dotted] (0,3)--(2,3);

                    \draw[->, red] (-0.2,-0.2)--(0,0);
                    \draw[-, red] (0,0)--(1,1);
                    \draw[-,red] (1.02,1)--(2,1.98);
                    \draw[->, red] (2,2)--(2.28,2);

                    \draw[->,blue] (-0.28,1)--(0,1);
                    \draw[-,blue] (0,1)--(1,1);
                    \draw[-, blue] (1,1.02)--(1.98,2);
                    \draw[->, blue] (2,2)--(2.2,2.2);
                \end{tikzpicture}
            \end{subfigure}
            \begin{subfigure}{0.3\textwidth}
            \centering
                \begin{tikzpicture}[scale=1]
                    \draw[dotted] (0,0)--(2,0);
                    \draw[dotted] (0,1)--(2,1);
                    \draw[dotted] (0,2)--(2,2);
                    \draw[dotted] (0,3)--(2,3);
                    \draw[dotted] (1,0)--(1,3);
                
                    \draw[->,red] (0,-0.28)--(0,0);
                    \draw[-, red] (0,0)--(0,1);
                    \draw[-,red] (0.02,1)--(0.02,2);
                    \draw[->,red] (0,2)--(1.28,2);

                    \draw[->, blue] (-0.28,1)--(-0.02,1);
                    \draw[-,blue] (-0.02,1)--(-0.02,2);
                    \draw[-, blue] (0,2)--(0,3);
                    \draw[-, blue] (0,3)--(1,3);
                    \draw[->, blue] (1,3)--(1,3.28);
                \end{tikzpicture}
            \end{subfigure}
            \caption{Three ways in which the embedded coordinate paths of $\mathbf{X}$ can intersect. Left: a continuous-time intersection as in Definition~\ref{contprocess}. Middle: an intersection created by a Bernoulli jump in Definition~\ref{mixprocess}. Right: an intersection created by a geometric jump in Definition~\ref{mixprocess}, where the embedded paths can intersect without the terminal configuration leaving $\mathbf{W}_{\mathrm d}^{N}$.}
            \label{breakinginterlacement}
        \end{figure}
        To treat all regimes uniformly, we embed the trajectories into piecewise-continuous planar paths, as in Figure~\ref{breakinginterlacement}, according to the following conventions:
        \begin{itemize}
            \item For Definition~\ref{contprocess}, use the usual cadlag embedding: horizontal motion during holding times, followed by a vertical unit jump at each jump time.
            \item For Bernoulli steps in Definition~\ref{mixprocess}, use linear interpolation between $(t,\mathscr{X}(t))$ and $(t+1,\mathscr{X}(t+1))$.
            \item For geometric steps in Definition~\ref{mixprocess}, use staircase interpolation between $(t,\mathscr{X}(t))$ and $(t+1,\mathscr{X}(t+1))$: first move vertically in the spatial coordinate, then horizontally in time.
        \end{itemize}
        With this convention, collision events are identified with intersections of the embedded paths.
        \begin{definition}[First collision time]\label{def:Tcollision}
            Let $(\mathbf{X}(t):t\in\mathscr{T})=(\mathscr{X}_1(t),\ldots,\mathscr{X}_N(t):t\in\mathscr{T})$ be an $N$-dimensional process whose coordinates are identical, independent copies of either Definition~\ref{contprocess} or Definition~\ref{mixprocess}.
            For each initial time $r\in\mathscr{T}$ and each $\mathbf{x}\in\mathbf{W}_{\mathrm d}^{N}$, let $\mathbb{P}_{r,\mathbf{x}}^{N}$ denote the law of the process started from $\mathbf{X}(r)=\mathbf{x}$ at absolute time $r$, and let $\mathbb{E}_{r,\mathbf{x}}^{N}$ denote the corresponding expectation.
            The canonical filtration  is
        \begin{equation*}
            \mathcal{F}_{r,t}\overset{\textnormal{def}}{=}\sigma\!\left(\mathbf{X}(u):u\in\mathscr{T}\cap[r,t]\right),\qquad t\geq r.
        \end{equation*}
        The first collision time of $\mathbf{X}$ relative to the initial time $r$ is then defined to be
        \begin{equation}\label{thestoppingtime}
            \mathcal{T}^{(r)}\overset{\textnormal{def}}{=}\inf\left\{t\in\mathscr{T}\cap[r,\infty): \textnormal{two embedded coordinate paths intersect by time }t\right\},
        \end{equation}
        where the embedding is the one specified above. 
        \end{definition}

In the continuous-time setting this dependence on $r$ is redundant, since the dynamics are time-homogeneous; in the discrete-time setting it is essential because the future evolution depends on the shifted environment $(\alpha_n,\beta_n,\tau_n)_{n\geq r}$ (although as will be seen in the proofs from an asymptotic perspective under \textbf{(A1)} and \textbf{(A2)} it is negligible) . For the case $r=0$ we write $\mathbb{P}_{\mathbf{x}}^{N}=\mathbb{P}_{0,\mathbf{x}}^{N}$, $\mathbb{E}_{\mathbf{x}}^{N}=\mathbb{E}_{0,\mathbf{x}}^{N}$, and $\mathcal{T}=\mathcal{T}^{(0)}$.

We now give an equivalent description of the stopping time $\mathcal{T}^{(r)}$.  In continuous time, the only possible collision mechanism is the nearest neighbour one, so the collision time $\mathcal{T}^{(r)}$ is simply the first exit time from the Weyl chamber:
    \begin{equation}\label{TWeylcont}
            \mathcal{T}_{\mathrm{W},\mathrm{c}}^{(r)}\overset{\textnormal{def}}{=}\inf\left\{t\in\mathbb{R}_+\cap[r,\infty):\mathbf{X}(t)\notin\mathbf{W}_{\mathrm d}^{N}\right\}.
        \end{equation}
            In discrete time, the update from time $t-1$ to time $t$ is governed by $\tau_{t-1}$. For a Bernoulli step, collision records the first Bernoulli endpoint at which the configuration leaves the Weyl chamber, the relevant stopping time being
        \begin{equation}\label{TWeyl}
            \mathcal{T}_{\mathrm{W},\mathrm{d}}^{(r)}\overset{\textnormal{def}}{=}\inf\left\{t\in\mathbb{Z}_+\cap(r,\infty):\mathbf{X}(t)\notin\mathbf{W}_{\mathrm d}^{N}\text{ and }\tau_{t-1}=0\right\}.
        \end{equation}
For a geometric step $\tau_{t-1}=1$ path collision is slightly more involved.
We say that  $\mathbf{x}=(x_1,\ldots,x_N) \in \mathbf{W}_{\mathrm{d}}^N$ interlaces $\mathbf{y}=(y_1,\ldots,y_N)\in\mathbf{W}_{\mathrm{d}}^N$ and denote this by $\mathbf{x}\preceq_N^N\mathbf{y}$ if
        \begin{equation*}
    x_1\leq y_1<x_2\leq y_2<\cdots<x_N\leq y_N.
        \end{equation*}
Then, in this case,  see Figure~\ref{breakinginterlacement} for an illustration, path collision is when configurations stop interlacing and the relevant stopping time is,
\begin{equation}\label{Tintl}
\mathcal{T}_{\mathrm{intl}}^{(r)}\overset{\textnormal{def}}{=}\inf\left\{t\in\mathbb{Z}_+\cap(r,\infty):\mathbf{X}(t-1)\not\preceq_N^N\mathbf{X}(t)\text{ and }\tau_{t-1}=1\right\}.
        \end{equation}
Hence, putting things together we see that in discrete time $\mathcal{T}^{(r)}=\min\left(\mathcal{T}_{\mathrm{W},\mathrm{d}}^{(r)},\mathcal{T}_{\mathrm{intl}}^{(r)}\right)$.

     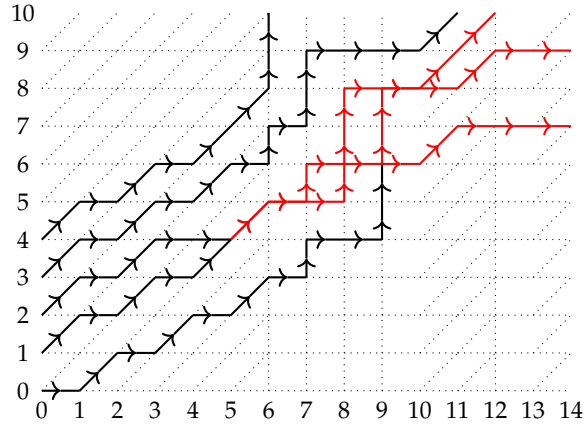
\begin{figure}[ht!]
    \centering
\begin{tikzpicture}[scale=0.5]

	\def\nx{14}
	\def\ny{10}

\foreach \x in {0,1,2,3,4,5,10,11,13} {
  \foreach \y in {0,...,\numexpr\ny-1} {
	            \draw[dotted] (\x,\y) -- (\x+1,\y);
	    
	            \draw[dotted] (\x,\y) -- (\x+1,\y+1);
	  }
	}
	
	\foreach \x in {6,7,8,9,12} {
	  \foreach \y in {0,...,\numexpr\ny-1} {
	            \draw[dotted] (\x,\y) -- (\x+1,\y);
	    
	            \draw[dotted] (\x,\y) -- (\x,\y+1);
	  }
}

\foreach \y in {0,...,\numexpr\ny} {
\node at (-0.5,\y) {\small$\y$};
}
\foreach \x in {0,...,\numexpr\nx} {
\node at (\x,-0.5) {\small$\x$};
}
  
	\draw[thick,->] (0,0) -- (0.5,0);
\draw[thick] (0.5,0) -- (1,0);
\draw[thick,->] (1,0) -- (1.5,0.5);
\draw[thick] (1.5,0.5) -- (2,1);
\draw[thick,->] (2,1) -- (2.5,1);
\draw[thick] (2.5,1) -- (3,1);
\draw[thick,->] (3,1) -- (3.5,1.5);
\draw[thick] (3.5,1.5) -- (4,2);
\draw[thick,->] (4,2) -- (4.5,2);
\draw[thick] (4.5,2) -- (5,2);
\draw[thick,->] (5,2) -- (5.5,2.5);
\draw[thick] (5.5,2.5) -- (6,3);
\draw[thick,->] (6,3) -- (6.5,3);
\draw[thick] (6.5,3) -- (7,3);
\draw[thick,->] (7,3) -- (7,3.5);
\draw[thick] (7,3.5) -- (7,4);
\draw[thick,->] (7,4) -- (7.5,4);
\draw[thick] (7.5,4) -- (8,4);
\draw[thick,->] (8,4) -- (8.5,4);
\draw[thick] (8.5,4) -- (9,4);
\draw[thick,->] (9,4) -- (9,4.5);
\draw[thick] (9,4.5) -- (9,5);

\draw[thick,->] (9,5) -- (9,5.5);
\draw[thick] (9,5.5) -- (9,6);
\draw[red,thick,->] (9,6) -- (9,6.5);
\draw[red,thick] (9,6.5) -- (9,7);
\draw[red,thick,->] (9,7) -- (9,7.5);
\draw[red,thick] (9,7.5) -- (9,8);
\draw[red,thick,->] (9,8) -- (9.5,8);
\draw[red,thick] (9.5,8) -- (10,8);
\draw[red,thick,->] (10,8) -- (10.5,8);
\draw[red,thick] (10.5,8) -- (11,8);
\draw[red,thick,->] (11,8) -- (11.5,8.5);
\draw[red,thick] (11.5,8.5) -- (12,9);

\draw[red,thick,->] (12,9) -- (12.5,9);
\draw[red,thick] (12.5,9) -- (13,9);
\draw[red,thick,->] (13,9) -- (13.5,9);
\draw[red,thick] (13.5,9) -- (14,9);

	\draw[thick,->] (0,1) -- (0.5,1.5);
\draw[thick] (0.5,1.5) -- (1,2);
\draw[thick,->] (1,2) -- (1.5,2);
\draw[thick] (1.5,2) -- (2,2);
\draw[thick,->] (2,2) -- (2.5,2.5);
\draw[thick] (2.5,2.5) -- (3,3);
\draw[thick,->] (3,3) -- (3.5,3);
\draw[thick] (3.5,3) -- (4,3);
\draw[thick,->] (4,3) -- (4.5,3.5);
\draw[thick] (4.5,3.5) -- (5,4);
\draw[red,thick,->] (5,4) -- (5.5,4.5);
\draw[red,thick] (5.5,4.5) -- (6,5);
\draw[red,thick,->] (6,5) -- (6.5,5);
\draw[red,thick] (6.5,5) -- (7,5);
\draw[red,thick,->] (7,5) -- (7,5.5);
\draw[red,thick] (7,5.5) -- (7,6);
\draw[red,thick,->] (7,6) -- (7.5,6);
\draw[red,thick] (7.5,6) -- (8,6);
\draw[red,thick,->] (9,6) -- (9.5,6);
\draw[red,thick] (9.5,6) -- (10,6);

\draw[red,thick,->] (10,6) -- (10.5,6.5);
\draw[red,thick] (10.5,6.5) -- (11,7);

\draw[red,thick,->] (11,7) -- (11.5,7);
\draw[red,thick] (11.5,7) -- (12,7);
\draw[red,thick,->] (12,7) -- (12.5,7);
\draw[red,thick] (12.5,7) -- (13,7);
\draw[red,thick,->] (13,7) -- (13.5,7);
\draw[red,thick] (13.5,7) -- (14,7);
\draw[red,thick,->] (8,6) -- (8.5,6);
\draw[red,thick] (8.5,6) -- (9,6);

	\draw[thick,->] (0,2) -- (0.5,2.5);
\draw[thick] (0.5,2.5) -- (1,3);
\draw[thick,->] (1,3) -- (1.5,3);
\draw[thick] (1.5,3) -- (2,3);
\draw[thick,->] (2,3) -- (2.5,3.5);
\draw[thick] (2.5,3.5) -- (3,4);
\draw[thick,->] (3,4) -- (3.5,4);
\draw[thick] (3.5,4) -- (4,4);
\draw[thick,->] (4,4) -- (4.5,4);
\draw[thick] (4.5,4) -- (5,4);
\draw[red,thick,->] (5,4) -- (5.5,4.5);
\draw[red,thick] (5.5,4.5) -- (6,5);

	\draw[red,thick,->] (6,5) -- (6.5,5);
\draw[red,thick] (6.5,5) -- (7,5);
\draw[red,thick,->] (7,5) -- (7.5,5);
\draw[red,thick] (7.5,5) -- (8,5);
\draw[red,thick,->] (8,5) -- (8,5.5);
\draw[red,thick] (8,5.5) -- (8,6);

\draw[red,thick,->] (8,6) -- (8,6.5);
\draw[red,thick] (8,6.5) -- (8,7);
\draw[red,thick,->] (8,7) -- (8,7.5);
\draw[red,thick] (8,7.5) -- (8,8);

\draw[red,thick,->] (8,8) -- (8.5,8);
\draw[red,thick] (8.5,8) -- (9,8);
\draw[red,thick,->] (9,8) -- (9.5,8);
\draw[red,thick] (9.5,8) -- (10,8);

\draw[red,thick,->] (10,8) -- (10.5,8.5);
\draw[red,thick] (10.5,8.5) -- (11,9);
\draw[red,thick,->] (11,9) -- (11.5,9.5);
\draw[red,thick] (11.5,9.5) -- (12,10);

	\draw[thick,->] (0,3) -- (0.5,3.5);
\draw[thick] (0.5,3.5) -- (1,4);
\draw[thick,->] (1,4) -- (1.5,4);
\draw[thick] (1.5,4) -- (2,4);
\draw[thick,->] (2,4) -- (2.5,4.5);
\draw[thick] (2.5,4.5) -- (3,5);
\draw[thick,->] (3,5) -- (3.5,5);
\draw[thick] (3.5,5) -- (4,5);
\draw[thick,->] (4,5) -- (4.5,5.5);
\draw[thick] (4.5,5.5) -- (5,6);
\draw[thick,->] (5,6) -- (5.5,6);
\draw[thick] (5.5,6) -- (6,6);
\draw[thick,->] (6,6) -- (6,6.5);
\draw[thick] (6,6.5) -- (6,7);
\draw[thick,->] (6,7) -- (6.5,7);
\draw[thick] (6.5,7) -- (7,7);
\draw[thick,->] (7,7) -- (7,7.5);
\draw[thick] (7,7.5) -- (7,8);
\draw[thick,->] (7,8) -- (7,8.5);
\draw[thick] (7,8.5) -- (7,9);
\draw[thick,->] (7,9) -- (7.5,9);
\draw[thick] (7.5,9) -- (8,9);

\draw[thick,->] (8,9) -- (8.5,9);
\draw[thick] (8.5,9) -- (9,9);
\draw[thick,->] (9,9) -- (9.5,9);
\draw[thick] (9.5,9) -- (10,9);
\draw[thick,->] (10,9) -- (10.5,9.5);
\draw[thick] (10.5,9.5) -- (11,10);

	\draw[thick,->] (0,4) -- (0.5,4.5);
\draw[thick] (0.5,4.5) -- (1,5);
\draw[thick,->] (1,5) -- (1.5,5);
\draw[thick] (1.5,5) -- (2,5);
\draw[thick,->] (2,5) -- (2.5,5.5);
\draw[thick] (2.5,5.5) -- (3,6);
\draw[thick,->] (3,6) -- (3.5,6);
\draw[thick] (3.5,6) -- (4,6);
\draw[thick,->] (4,6) -- (4.5,6.5);
\draw[thick] (4.5,6.5) -- (5,7);
\draw[thick,->] (5,7) -- (5.5,7.5);
\draw[thick] (5.5,7.5) -- (6,8);
\draw[thick,->] (6,8) -- (6,8.5);
\draw[thick] (6,8.5) -- (6,9);
\draw[thick,->] (6,9) -- (6,9.5);
\draw[thick] (6,9.5) -- (6,10);

\end{tikzpicture}
    \caption{Five independent mixed Bernoulli-and-geometric jump processes started from $(0,0),(0,1),\ldots,(0,4)$. The red segments indicate the portions of the trajectories after the first collision. In this example the second and third paths first intersect at time $5$, so $\mathcal{T}=5$.}
    \label{noncollidingfig}
\end{figure} 

 We now define the conditioned non-colliding law. Under \textbf{(A1)} and \textbf{(A2)}, Theorem~\ref{collisionThm} will give a quantitative collision time tail estimate with a vanishing right-hand side for $N\geq2$, so the eventual non-colliding law is a limiting conditioning on a zero-probability event. The rigorous construction therefore begins with finite-horizon conditioning.

Fix $r\leq T<s$ in $\mathscr{T}$. For each event $\mathscr{A}\in\mathcal{F}_{r,T}$, define
    \begin{equation*}
        \mathsf{P}_{r,\mathbf{x}}^{N,(s)}(\mathscr{A})\overset{\textnormal{def}}{=}\mathbb{P}_{r,\mathbf{x}}^{N}\left(\mathscr{A}\mid \mathcal{T}^{(r)}>s\right),
    \end{equation*}
    whenever $\mathbb{P}_{r,\mathbf{x}}^{N}[\mathcal{T}^{(r)}>s]>0$, this condition is automatic for every $\mathbf{x}\in\mathbf{W}_{\mathrm d}^{N}$ under \hyperref[A1]{\textbf{(A1)}}.  Define the space-time survival function
    \begin{equation*}
        H_{t,s}(\mathbf{y})\overset{\textnormal{def}}{=}\mathbb{P}_{t,\mathbf{y}}^{N}\left(\mathcal{T}^{(t)}>s\right),\qquad r\leq t\leq s,\quad \mathbf{y}\in\mathbf{W}_{\mathrm d}^{N}.
    \end{equation*}
    Then, the Markov property at the deterministic time $T$ gives
    \begin{equation}\label{construction}
        \mathsf{P}_{r,\mathbf{x}}^{N,(s)}(\mathscr{A})
        =
        \mathbb{E}_{r,\mathbf{x}}^{N}\left[
            \frac{H_{T,s}(\mathbf{X}(T))}{H_{r,s}(\mathbf{x})}
            \mathbf{1}_{\mathscr{A}}\mathbf{1}_{\{\mathcal{T}^{(r)}>T\}}
        \right],
        \qquad \mathscr{A}\in\mathcal{F}_{r,T}.
    \end{equation}
We will mainly be concerned with the case $r=0$. If the finite-horizon conditioned laws admit a consistent weak limit as $s\to\infty$, we denote the resulting limiting probability law by $\mathsf{P}_{\mathbf{x}}^{N}$:
\begin{equation*}
\mathsf{P}_\mathbf{x}^N(\cdot)\overset{\mathrm{def}}{=}\lim_{s\to \infty}\mathsf{P}_{0,\mathbf{x}}^{N,(s)}(\cdot).
\end{equation*}
We write $\mathsf{E}^N_{\mathbf{x}}$ for the corresponding expectation and denote its canonical process by $(\hat{\mathbf{X}}(t):t\in\mathscr{T})$. 
If moreover there exists strictly positive $(\mathfrak{H}_t(\cdot))_{t\in \mathscr{T}}$ on $ \mathbf{W}_{\mathrm{d}}^N$ such that for each $T\in \mathscr{T}$,
   \begin{equation*}
\mathsf{E}^N_{\mathbf{x}}\left[F\left(\hat{\mathbf{X}}(t): t\leq T\right)\right]
        =
        \mathbb{E}_{\mathbf{x}}^{N}\left[
            \frac{\mathfrak{H}_{T}(\mathbf{X}(T))}{\mathfrak{H}_{0}(\mathbf{x})}
            F\left(\mathbf{X}(t): t\leq T\right)\mathbf{1}_{\{\mathcal{T}>T\}}
        \right].
    \end{equation*}
then $\mathsf{E}_\mathbf{x}^N/\hat{\mathbf{X}}$ is the Doob $h$-transform \cite{Doob,RevuzYor} of $\mathbb{E}_\mathbf{x}^N/\mathbf{X}$  by $(\mathfrak{H}_t)_{t\in \mathscr{T}}$ and the transition probability of $\hat{\mathbf{X}}$ is given by
\begin{equation}\label{startofproofs}
        \mathsf{P}_{\mathbf{x}}^{N}\left[\hat{\mathbf{X}}(t)=\mathbf{y}\right]
        =
        \frac{\mathfrak{H}_{t}(\mathbf{y})}{\mathfrak{H}_{0}(\mathbf{x})}\mathbb{P}_{\mathbf{x}}^{N}\left[\mathbf{X}(t)=\mathbf{y},\,\mathcal{T}>t\right].
    \end{equation}

\subsection{Main results}
    We now state main results of the paper for the non-colliding construction from \eqref{construction}. We assume throughout this subsection that \textbf{(A1)} and \textbf{(A2)} are in force.  We need a few  definitions. We write $\llbracket N\rrbracket\overset{\textnormal{def}}{=}\{1,\ldots,N\}$ for integer coordinate indices. 
    \begin{definition}[Continuous Weyl chamber]\label{def:continuous-weyl}
        For $N\in\mathbb{N}$, set
        \begin{equation*}
            \mathbf{W}_{\mathrm c}^{N}\overset{\textnormal{def}}{=}\{\mathbf{u}=(u_1,\ldots,u_N)\in\mathbb{R}^{N}:u_1<\cdots<u_N\}.
        \end{equation*}
    \end{definition}

    \begin{definition}[Characteristic polynomials]\label{def:characteristic-polynomials}
        For $x\in\mathbb{Z}_+$, define
        the ``characteristic polynomials" associated to $(\lambda_x)_{x\in \mathbb{Z}_+}$ by
\begin{equation}\label{characteristic}
            \rho_0(w)=1,\qquad \rho_x(w)=\prod_{k=0}^{x-1}\left(1-\frac{w}{\lambda_k}\right),\qquad x\geq1.
        \end{equation}
    
    \end{definition}
    
The following function will play a key role in what follows.

\begin{definition}\label{harmonicfunction}
	        For $N\in\mathbb{N}$ and $\mathbf{x}=(x_1,\ldots,x_N)\in\mathbf{W}_{\mathrm d}^{N}$, define $\mathfrak{h}_N:\mathbf{W}_{\mathrm d}^{N}\rightarrow\mathbb{R}$ by
	        \begin{equation*}
	            \mathfrak{h}_N(\mathbf{x})=\det\!\left((-1)^{j-1}\frac{\rho_{x_i}^{(j-1)}(0)}{(j-1)!}\right)_{i,j=1}^{N}.
	        \end{equation*}
	    \end{definition}
Observe that, $\mathfrak{h}_N$ can also be written as,
\begin{equation*}
\mathfrak{h}_N(\mathbf{x})=
    \det\!\left(\mathfrak{e}_{j-1}\left(\lambda_0^{-1},\ldots,\lambda_{x_i-1}^{-1}\right)\right)_{i,j=1}^{N},
\end{equation*}
where $\mathfrak{e}_k$ is the $k$-th elementary symmetric polynomial, with $\mathfrak{e}_0 \equiv 1$. We note that $\mathfrak{h}_N$ is strictly positive in $\mathbf{W}_{\mathrm d}^{N}$, see Lemma \ref{prop:h-positive}.

	 As we will see shortly, the function $\mathfrak{h}_N$ records the dependence of the collision time tail estimate on the initial condition and later becomes the Doob $h$-function defining the conditioned dynamics. In the homogeneous case, when $\lambda_n\equiv1$, the determinant reduces to a multiple of the Vandermonde determinant $\Delta$, namely with $\mathbf{y}=(y_1,\ldots,y_N)\in\mathbb{R}^{N}$, 
        \begin{equation*}
            \Delta(\mathbf{y})\overset{\textnormal{def}}{=}\det\left(y_i^{j-1}\right)_{i,j=1}^N\equiv\prod_{1\leq j<i\leq N}\left(y_i-y_j\right),
        \end{equation*}
        recovering classical results, see \cite{OConnellRSK,OConnellYor,OrderedRW,DenisovWachtel1,KonigOConnellRoch}. However, unlike the Vandermonde determinant, $\mathfrak{h}_N$ cannot in general be extended to real-valued vectors.
    \begin{definition}\label{themeananddeviation}
        For all $t\in\mathscr{T}$, we set
        \begin{equation*}
            \mathrm{z}(t)\overset{\textnormal{def}}{=}\mathbf{1}_{\{\mathscr{T}=\mathbb{R}_+\}}\frac{t}{\Lambda}+\mathbf{1}_{\{\mathscr{T}=\mathbb{Z}_+\}}\frac{1}{\Lambda}\sum_{n=0}^{t-1}\left(\alpha_n(1-\tau_n)+\beta_n\tau_n\right),\quad\text{and}\quad\mathscr{Z}(t)\overset{\textnormal{def}}{=}\lfloor\mathrm{z}(t)\rceil,
        \end{equation*}
        where $\lfloor\cdot\rceil$ gives the closest integer to a real-valued input, with half-integers rounded upward. Moreover, we define for $t$ large enough,
        \begin{equation*}
            \Sigma(t)\overset{\textnormal{def}}{=}\sqrt{\sum_{n=0}^{\mathscr{Z}(t)}\frac{1}{\lambda_n^2}-\mathbf{1}_{\{\mathscr{T}=\mathbb{Z}_+\}}\sum_{n=0}^{t-1}\big(\alpha_n^2(1-\tau_n)-\beta_n^2\tau_n\big)}.
        \end{equation*}
    \end{definition}
We note that at this point the positivity of the quantity under the square root in Definition~\ref{themeananddeviation} has not yet been shown. This is established below in \textbf{(iv)} of Proposition~\ref{specificcasesub}: there is a deterministic lower time after which $\Sigma(t)$ is positive. We tacitly assume we are in this eventual-positive range in all statements below. The quantity $\mathscr{Z}(t)$ plays the role of the deterministic first-order location of $\mathscr{X}(t)$, and $\Sigma(t)$ is the natural fluctuation scale.

    \begin{definition}\label{gammadef}
        On the eventual-positive time range for $\Sigma(t)$ described above, define
        \begin{equation*}
            \gamma(t)\overset{\textnormal{def}}{=}\max\left(\sup_{\substack{s\in\mathscr{T}\\s\geq t}}\left|\frac{1}{\mathscr{Z}(s)+1}\sum_{n=0}^{\mathscr{Z}(s)}\frac{1}{\lambda_n}-\Lambda\right|,\sup_{\substack{s\in\mathscr{T}\\s\geq t}}\Sigma(s)^{-\frac{1}{2}}\right),
        \end{equation*}
        and set
        \begin{equation}\label{eq:psiG-def}
            \psi_{\mathrm{G}}(t)\overset{\textnormal{def}}{=}\left(1+\log(\gamma(t)^{-1})\right)^{-1/4},
        \end{equation}
        whenever $t$ has been chosen large enough so that $\gamma(t)<1$.
    \end{definition}
Observe that, $\psi_{\mathrm{G}}(t) \longrightarrow 0$ as $t \longrightarrow \infty$. The scale $\gamma(t)$ records the residual reciprocal-rate inhomogeneity together with the minimal fluctuation scale, while $\psi_{\mathrm{G}}(t)$ is the common logarithmic error scale used in the quantitative limits below. We also set,
	    \begin{equation*}
	        \mathrm{d}_{N}\overset{\textnormal{def}}{=}\frac{N}{2}(N-1),\qquad \mathrm{A}_{N}\overset{\textnormal{def}}{=}\frac{1}{N!}\prod_{i=1}^{N}\frac{\Gamma\!\left(1+\frac{i}{2}\right)}{\Gamma\!\left(\frac{3}{2}\right)}.
	    \end{equation*}

    The following theorem is our main result. It gives the asymptotics of $\mathbb{P}_{\mathbf{x}}^N[\mathcal{T}>t]$ as $t\longrightarrow \infty$ with an explicit error bound.
    \begin{theorem}[Collision tail]\label{collisionThm}
        Assume \textbf{(A1)} and \textbf{(A2)}. For every $\mathbf{x}\in\mathbf{W}_{\mathrm d}^{N}$ there exist constants $\mathrm{T}_{\mathrm{col},\mathbf{x}}\in\mathscr{T}$ and $\mathrm{C}_{\mathrm{col},\mathbf{x}}>0$ such that, for all $t>\mathrm{T}_{\mathrm{col},\mathbf{x}}$,
        \begin{equation*}
            \left|\mathbb{P}_{\mathbf{x}}^N[\mathcal{T}>t]-\mathrm{A}_{N}\mathfrak{h}_N(\mathbf{x})\Sigma(t)^{-\mathrm{d}_{N}}\right|\leq\mathrm{C}_{\mathrm{col},\mathbf{x}}\Sigma(t)^{-\mathrm{d}_{N}}\psi_{\mathrm{G}}(t).
        \end{equation*}
    \end{theorem}

This confirms in strong form a prediction by one of us from \cite{assiotis2023integrablemodelsinhomogeneousspace}. Namely, a relatively straightforward consequence of Theorem \ref{collisionThm} is the following corollary which identifies the explicit Doob $h$-transform structure of the conditioned process $\hat{\mathbf{X}}$.

        \begin{cor}\label{Thm1}
        Assume \textbf{(A1)} and \textbf{(A2)}. For all $N\in\mathbb{N}$, $\mathbf{x}\in\mathbf{W}_{\mathrm d}^{N}$ and $T\in\mathscr{T}$, the non-colliding process $(\hat{\mathbf{X}}(t):t\leq T)$ is the Doob $h$-transform of $(\mathbf{X}(t):t\leq T)$ by $\mathfrak{h}_N(\cdot)$: for every bounded $\mathcal{F}_{0,T}$-measurable functional $F$,
        \begin{equation*}
            \mathsf{E}^N_{\mathbf{x}}\!\left[F\!\left(\hat{\mathbf{X}}(t):t\leq T\right)\right]=\mathbb{E}_{\mathbf{x}}^N\!\left[\frac{\mathfrak{h}_N\!\left(\mathbf{X}(T)\right)}{\mathfrak{h}_N(\mathbf{x})}F\!\left(\mathbf{X}(t):t\leq T\right)\mathbf{1}_{\{\mathcal{T}>T\}}\right].
        \end{equation*}
    \end{cor}
    
 \begin{remark} \label{Rmk:Prediction}
 The way Corollary \ref{Thm1} was predicted in \cite{assiotis2023integrablemodelsinhomogeneousspace} was as follows. One adds a certain drift $\boldsymbol{\delta}_i$ to each coordinate $\mathscr{X}_i$ (the way it manifests in each of the three types of dynamics is different, see \cite{assiotis2023integrablemodelsinhomogeneousspace}). The drifts are assumed strictly ordered and in fact separated by some sufficiently large constant, namely $\boldsymbol{\delta}_{i+1}-\boldsymbol{\delta}_i>\mathsf{C}_{((\lambda_n),(\alpha_t),(\beta_t))}$ that depends on the environment. Then, it can be shown that the probability of non-collision for all times (i.e. taking $t=\infty$ in Theorem \ref{collisionThm} above) is strictly positive and in fact explicit. The strong ordered drifts make the individual chains asymptotically ordered and this renders the $t\to \infty$ limit essentially trivial. Then, looking at the ratio of these probabilities for different starting positions $\mathbf{y}$ and $\mathbf{x}$ and formally taking the singular limit $\boldsymbol{\delta}_i \to 0$ in the formulae, in particular violating the threshold $\boldsymbol{\delta}_{i+1}-\boldsymbol{\delta}_i>\mathsf{C}_{((\lambda_n),(\alpha_t),(\beta_t))}$ for which the argument is valid, one obtains the ratio $\mathfrak{h}_N(\mathbf{y})/\mathfrak{h}_N(\mathbf{x})$. The approach we take in this paper to actually prove the desired statement is completely different and does not introduce drifts.
 \end{remark}

Moving towards our second main result we define the path-space.

    \begin{definition}
        For $N\in\mathbb{N}$, we define the space of $N$, $\mathscr{T}$ indexed paths in $\mathbb{Z}_+$ by
        \begin{equation*}
            \mathscr{P}_N\overset{\textnormal{def}}{=}\{\mathbf{z}=(z_1,\ldots,z_N):\mathscr{T}\rightarrow\mathbb{Z}_+^{N}\}.
        \end{equation*}
        For ease of notation, we write $\mathscr{P}\overset{\textnormal{def}}{=}\mathscr{P}_1$.
    \end{definition}
    We next define the phase function.
    \begin{definition}[Phase function]\label{form}
        For $y\in\mathscr{P}$ and any $t\in\mathscr{T}$, we define the phase function $w\mapsto \mathfrak{F}_y(w,t)$ locally, away from the displayed logarithmic singularities, by
        \begin{equation}
        \begin{aligned}
        \mathfrak{F}_{y}(w,t)&\overset{\textnormal{def}}{=}-\mathbf{1}_{\{\mathscr{T}=\mathbb{R}_+\}}tw-\sum_{n=0}^{y(t)}\log\left(1-\frac{w}{\lambda_n}\right)\\&+\mathbf{1}_{\{\mathscr{T}=\mathbb{Z}_+\}}\sum_{n=0}^{t-1}\left[(1-\tau_{n})\log\left(1-\alpha_nw\right)-\tau_n\log\left(1+\beta_nw\right)\right].
        \end{aligned}
        \end{equation}
        The logarithms are taken on a common simply connected domain, with the branch normalized by $\log(1)=0$ near the origin.
    \end{definition}
    \noindent In this paper, $\mathfrak{F}_y$ is the analytic object that simultaneously records the spatial and temporal inhomogeneities, and the endpoint $y(t)$. Section \ref{glo} amounts to a systematic analysis of this phase function and of the stationary point it generates defined next. Below and throughout $\mathsf{f}^{(k)}$ denotes the $k$-th derivative of a function $\mathsf{f}$.

    \begin{definition}\label{definitionstationary}
        For each $t\in\mathscr{T}$ and $y\in\mathscr{P}$, define
        \begin{equation*}
            \mathscr{U}_{y}(t)
            \overset{\textnormal{def}}{=}
            \begin{cases}
            -\Sigma(t)\operatorname*{arg\,min}\left\{|u|:u\in\mathbb{C},\ \mathfrak{F}^{(1)}_{y}(u,t)=0\right\},&\textnormal{if this minimizer is unique},\\
            0,&\textnormal{otherwise}.
            \end{cases}
        \end{equation*}
    \end{definition}
    
    Proposition~\ref{axpro} shows that, after a deterministic lower time and uniformly for endpoints satisfying $y(t)\in\mathscr{C}_{t}$, where $\mathscr{C}_t$ is given in Definition \ref{def:endpoint-window}, the minimizer is unique. On this range we write the minimizer as $\mathrm{u}_{y}(t)$, so that $\mathscr{U}_{y}(t)=-\Sigma(t)\mathrm{u}_{y}(t)$. Proposition~\ref{moments}, applied to each coordinate, shows that almost surely there is a finite random time after which the particle endpoints lie in this regime.

Write $\mathcal{N}(0,1)^{\otimes N}$ for the law of an $\mathbb{R}^{N}$-valued random vector whose coordinates are independent standard Gaussian random variables. Moreover, define the following classical random matrix law \cite{AndersonGuionnetZeitouni2010,ForresterBook}.

    \begin{definition}[Gaussian Orthogonal Ensemble] \label{GOEeigenvalue}
        For $N\in\mathbb{N}$, we consider the Gaussian Orthogonal Ensemble, the law on $\mathrm{Sym}_N(\mathbb{R})$ ($N \times N$ real symmetric matrices) given by
        \begin{equation*}
2^{-N/2}\pi^{-N(N+1)/4}\exp\!\left(-\frac{1}{2}\mathrm{Tr}(\mathbf{H}^2)\right)\mathrm{d}\mathbf{H},\qquad \mathbf{H}\in\mathrm{Sym}_N(\mathbb{R}),
        \end{equation*}
        where $\mathrm{d}\mathbf{H}$ denotes Lebesgue measure on $\mathrm{Sym}_N(\mathbb{R})$. We denote by $\mathrm{GOE}_N$ the induced law of its ordered eigenvalues; that is,
        \begin{equation*}
\mathrm{GOE}_N(\mathrm{d}\mathbf{u})\overset{\textnormal{def}}{=}\frac{N!}{(2\pi)^{N/2}}\prod_{i=1}^{N}\frac{\Gamma(3/2)}{\Gamma(1+i/2)}\prod_{1\leq j<i\leq N}(u_i-u_j)\exp\!\left(-\frac12|\mathbf{u}|^2\right)\mathrm{d}\mathbf{u},\qquad \mathbf{u}\in\mathbf{W}_{\mathrm c}^{N}.
        \end{equation*}
    \end{definition}

    Then, the following is our second main result. The stochastic process $((\mathscr{X}_i(t))_{i=1}^N:t\in \mathscr{T})$ under the highly non-trivial transformation that maps it to the stationary point of the corresponding phase function becomes Gaussian in the large $t$ limit! Moreover, this stationary point stochastic process conditioned on $\mathcal{T}>t$ converges to the $\mathrm{GOE}_N$ law.

    \begin{theorem}[Stationary-point limits]\label{thm:stationarypoints}
        Assume \textbf{(A1)} and \textbf{(A2)}. Fix $\mathbf{x}\in\mathbf{W}_{\mathrm d}^{N}$. Let $(\mathbf{X}(t):t\in\mathscr{T})$ denote the $N$-particle system started from $\mathbf{x}$, and define
        \begin{equation*}
            \mathbf{U}(t)=(\mathscr{U}_1(t),\ldots,\mathscr{U}_N(t)),\qquad \mathscr{U}_i(t)\overset{\textnormal{def}}{=}\mathscr{U}_{\mathscr{X}_i}(t),\qquad i\in\llbracket N\rrbracket.
        \end{equation*}
        For each $i\in\llbracket N\rrbracket$, there exists a $\mathbb{P}_{\mathbf{x}}^N$-almost surely finite random time $\mathrm{T}_{i}$ such that the minimizer in Definition~\ref{definitionstationary}, with $y=\mathscr{X}_i$, is unique for every $t>\mathrm{T}_{i}$. Moreover, if $F:\mathbb{R}^{N}\to\mathbb{R}$ is bounded and Lipschitz, then there are constants $\mathrm{T}_{\mathrm{free},F,\mathbf{x}}\in\mathscr{T}$ and $\mathrm{C}_{\mathrm{free},F,\mathbf{x}}>0$ such that, for all $t>\mathrm{T}_{\mathrm{free},F,\mathbf{x}}$,
        \begin{equation*}
            \left|\mathbb{E}_{\mathbf{x}}^N\left[F\left(\mathbf{U}(t)\right)\right]-\frac{1}{(2\pi)^{N/2}}\int_{\mathbb{R}^N}F(\mathbf{u})\exp\!\left(-\frac12|\mathbf{u}|^2\right)\mathrm{d}\mathbf{u}\right|\leq\mathrm{C}_{\mathrm{free},F,\mathbf{x}}\psi_{\mathrm{G}}(t).
        \end{equation*}
        If $F:\mathbb{R}^{N}\to\mathbb{R}$ is bounded and Lipschitz, then there are constants $\mathrm{T}_{\mathrm{GOE},F,\mathbf{x}}\in\mathscr{T}$ and $\mathrm{C}_{\mathrm{GOE},F,\mathbf{x}}>0$ such that, for all $t>\mathrm{T}_{\mathrm{GOE},F,\mathbf{x}}$,
        \begin{equation*}
            \left|\mathbb{E}_{\mathbf{x}}^N\!\left[F\!\left(\mathbf{U}(t)\right)\middle|\mathcal{T}>t\right]-\frac{N!}{(2\pi)^{N/2}}\prod_{i=1}^{N}\frac{\Gamma(3/2)}{\Gamma(1+i/2)}\int_{\mathbf{W}_{\mathrm c}^{N}}F(\mathbf{u})\Delta(\mathbf{u})\exp\!\left(-\frac12|\mathbf{u}|^2\right)\mathrm{d}\mathbf{u}\right|\leq\mathrm{C}_{\mathrm{GOE},F,\mathbf{x}}\psi_{\mathrm{G}}(t).
        \end{equation*}
    \end{theorem}

Finally, the following result will be a relatively straightforward consequence of Theorem \ref{thm:stationarypoints} by essentially approximating the stationary point stochastic process $\mathbf{U}$ by the simpler process $\mathbf{V}$.
    \begin{cor}\label{cor:particlecoordinates}
        Assume \textbf{(A1)} and \textbf{(A2)}, and fix $N\in\mathbb{N}$ and $\mathbf{x}\in\mathbf{W}_{\mathrm d}^{N}$. Let $(\mathbf{X}(t):t\in\mathscr{T})$ denote the $N$-particle system started from $\mathbf{x}$. For $\mathbf{y}=(y_1,\ldots,y_N)\in\mathbb{Z}_+^N$ define
        \begin{equation*}
            \Lambda_{N}(\mathbf{y})\overset{\textnormal{def}}{=}\left(\sum_{n=0}^{y_1}\frac{1}{\lambda_n},\sum_{n=0}^{y_2}\frac{1}{\lambda_n},\ldots,\sum_{n=0}^{y_N}\frac{1}{\lambda_n}\right).
        \end{equation*}
        Let $\underline{1}=(1,\ldots,1)\in\mathbb{R}^{N}$ and set $\mathbf{V}(t)\overset{\textnormal{def}}{=}\Sigma(t)^{-1}\left(\Lambda_N(\mathbf{X}(t))-\Lambda\mathscr{Z}(t)\underline{1}\right)$. For every bounded Lipschitz $F:\mathbb{R}^{N}\to\mathbb{R}$ there are constants $\mathrm{T}_{\mathrm{coord},F,\mathbf{x}},\mathrm{T}_{\mathrm{coordGOE},F,\mathbf{x}}\in\mathscr{T}$ and $\mathrm{C}_{\mathrm{coord},F,\mathbf{x}},\mathrm{C}_{\mathrm{coordGOE},F,\mathbf{x}}>0$ such that, for every $t>\mathrm{T}_{\mathrm{coord},F,\mathbf{x}}$,
        \begin{equation*}
            \left|\mathbb{E}_{\mathbf{x}}^N\left[F\left(\mathbf{V}(t)\right)\right]-\frac{1}{(2\pi)^{N/2}}\int_{\mathbb{R}^{N}}F(\mathbf{u})\exp\!\left(-\frac12|\mathbf{u}|^2\right)\mathrm{d}\mathbf{u}\right|\leq\mathrm{C}_{\mathrm{coord},F,\mathbf{x}}\psi_{\mathrm{G}}(t),
        \end{equation*}
        and, for every $t>\mathrm{T}_{\mathrm{coordGOE},F,\mathbf{x}}$,
        \begin{align*}
            &\left|\mathbb{E}_{\mathbf{x}}^N\!\left[F\!\left(\mathbf{V}(t)\right)\middle|\mathcal{T}>t\right]-\frac{N!}{(2\pi)^{N/2}}\prod_{i=1}^{N}\frac{\Gamma(3/2)}{\Gamma(1+i/2)}\int_{\mathbf{W}_{\mathrm c}^{N}}F(\mathbf{u})\Delta(\mathbf{u})\exp\!\left(-\frac12|\mathbf{u}|^2\right)\mathrm{d}\mathbf{u}\right|\\
            &\hspace{2cm}\leq\mathrm{C}_{\mathrm{coordGOE},F,\mathbf{x}}\psi_{\mathrm{G}}(t).
        \end{align*}
    \end{cor}

\subsection{Examples}

    We give some representative examples in continuous time. The first four satisfy \textbf{(A1)}--\textbf{(A2)} and illustrate how the deterministic inhomogeneity scale $\gamma(t)$ behaves in concrete environments. In each admissible example, on the time range where Definition~\ref{definitionstationary} selects the real stationary point near the origin, the one-particle stationary-point stochastic process is obtained by solving the equation
    \begin{equation}\label{intro-stationary-equation}
        \sum_{n=0}^{\mathscr{X}(t)}\frac{1}{\lambda_n-\mathrm{u}_{\mathscr{X}}(t)}=t
    \end{equation}
    for the real stationary point $\mathrm{u}_{\mathscr{X}}(t)$ near the origin and then setting $\mathscr{U}(t)=-\Sigma(t)\mathrm{u}_{\mathscr{X}}(t)$.
    The fifth example lies outside the hypotheses, exhibiting different asymptotic behaviour, thus showing why our assumptions are structurally important.

    \begin{example}[Constant rates]\label{constantcase}
        Suppose $\lambda_n\equiv\lambda$ for some $\lambda>0$. Then, obviously \textbf{(A1)} and \textbf{(A2)} hold with $\Lambda=\lambda^{-1}$. Hence the inhomogeneity part of $\gamma(t)$ vanishes, and
        \begin{equation*}
            \gamma(t)=\sup_{\substack{s\in\mathscr{T}\\s\geq t}}\Sigma(s)^{-1/2}.
        \end{equation*}
        Since $\mathrm{z}(t)=\lambda t$, $\mathscr{Z}(t)=\lfloor\lambda t\rceil$, and $\Sigma(t)^2=(\mathscr{Z}(t)+1)/\lambda^2$, the stationary equation \eqref{intro-stationary-equation} can be solved explicitly:
        \begin{equation*}
            \mathrm{u}_{\mathscr{X}}(t)=\lambda-\frac{\mathscr{X}(t)+1}{t},\qquad\mathscr{U}(t)=\Sigma(t)\left(\frac{\mathscr{X}(t)+1}{t}-\lambda\right).
        \end{equation*}
        Thus, in the homogeneous case, the stationary-point coordinate is an explicit saddle-scale version of the usual centered particle coordinate.
    \end{example}

    \begin{example}[Periodic rates]\label{periodiccase}
        Fix $M\in\mathbb{N}$ and positive numbers $\lambda^{(0)},\ldots,\lambda^{(M-1)}$, and define
        \begin{equation*}
            \lambda_n=\lambda^{(r)},\qquad\text{whenever }n\equiv r\ (\mathrm{mod}\ M).
        \end{equation*}
Clearly \textbf{(A1)} holds. Moreover, if we write $a_j=(\lambda^{(j)})^{-1}$, then \textbf{(A2)} holds with
        \begin{equation*}
            \Lambda=\frac{1}{M}\sum_{j=0}^{M-1}a_j.
        \end{equation*}
        Moreover, writing $N+1=qM+r$ with $0\leq r<M$, the contribution of the complete periods cancels and the remaining incomplete period gives
        \begin{equation*}
            \left|\frac{1}{N+1}\sum_{n=0}^{N}\frac{1}{\lambda_n}-\Lambda\right|\leq\frac{2M\max_{0\leq j<M}a_j}{N+1}.
        \end{equation*}
        Hence,  the inhomogeneity component of $\gamma(t)$ is controlled explicitly by a constant multiple of $(\mathscr{Z}(s)+1)^{-1}$ inside the future supremum in $s\ge t$. 
    \end{example}

    \begin{example}[Algebraically convergent rates]\label{algebraiccase}
        Suppose that the positive rates satisfy $\lambda_n\to\lambda>0$ and, for some $\alpha\in(0,1)$ and $\mathrm{C}_{\mathrm{alg}}>0$,
        \begin{equation*}
            |\lambda_n-\lambda|\leq\mathrm{C}_{\mathrm{alg}}(n+1)^{-\alpha},\qquad n\in\mathbb{Z}_+.
        \end{equation*}
        Clearly \textbf{(A1)} holds. If $\mathfrak{L}_{0}\overset{\textnormal{def}}{=}\inf_{n\geq0}\lambda_n>0$, then
        \begin{equation*}
            \left|\frac{1}{\lambda_n}-\frac{1}{\lambda}\right|\leq\frac{\mathrm{C}_{\mathrm{alg}}}{\lambda\mathfrak{L}_{0}}(n+1)^{-\alpha}.
        \end{equation*}
        Averaging this bound gives, for every $N\in\mathbb{Z}_+$,
        \begin{equation*}
            \left|\frac{1}{N+1}\sum_{n=0}^{N}\frac{1}{\lambda_n}-\lambda^{-1}\right|\leq\frac{\mathrm{C}_{\mathrm{alg}}}{\lambda\mathfrak{L}_{0}(N+1)}\left(1+\frac{(N+1)^{1-\alpha}}{1-\alpha}\right).
        \end{equation*}
        Thus, again \textbf{(A2)} holds with $\Lambda=\lambda^{-1}$ and in particular, for every $s\in\mathscr{T}$,
        \begin{equation*}
            \left|\frac{1}{\mathscr{Z}(s)+1}\sum_{n=0}^{\mathscr{Z}(s)}\frac{1}{\lambda_n}-\lambda^{-1}\right|\leq\frac{\mathrm{C}_{\mathrm{alg}}}{\lambda\mathfrak{L}_{0}(\mathscr{Z}(s)+1)}\left(1+\frac{(\mathscr{Z}(s)+1)^{1-\alpha}}{1-\alpha}\right),
        \end{equation*}
        so the inhomogeneity component of $\gamma(t)$ is controlled by the future supremum of the above explicit bound over $s\geq t$. 
    \end{example}

    \begin{example}[Qualitatively convergent rates]\label{qualitativecase}
        Suppose now only that the  rates satisfy $\lambda_n\to\lambda>0$. Then obviously \textbf{(A1)} holds and Cesàro's theorem gives
        \begin{equation*}
            \lim_{N\to\infty}\frac{1}{N+1}\sum_{n=0}^{N}\frac{1}{\lambda_n}=\lambda^{-1},
        \end{equation*}
        so that \textbf{(A2)} holds with $\Lambda=\lambda^{-1}$. In this situation $\gamma(t)$ still tends to zero by Definition~\ref{gammadef}; its inhomogeneity component is governed by the rate at which the reciprocal-rate empirical average approaches $\lambda^{-1}$ along the deterministic centres $\mathscr{Z}(s)$, while the fluctuation-scale term in $\gamma(t)$ may also contribute. 
    \end{example}

    \begin{example}[Yule pure-birth case]\label{Yule}
        Fix $\lambda>0$ and set
        \begin{equation*}
            \lambda_n=\lambda n,\qquad n \ge 1.
        \end{equation*}
  Clearly \textbf{(A1)} fails and it is easy to see \textbf{(A2)} fails as well. It is well known,  see \cite{Yule1925,AthreyaNey1972}, that starting from $\mathscr{X}(0)=x \ge 1$,
        \begin{equation*}
        \lim_{t\to\infty}\mathrm{e}^{-\lambda t}\mathscr{X}(t)=\mathsf{\Gamma}_{x,1}\quad\textnormal{almost surely},
        \end{equation*}
        where $\mathsf{\Gamma}_{x,1}$ is a Gamma random variable with density $y^{x-1}\mathrm{e}^{-y}/\Gamma(x)$ for $y>0$. Using this fact and the Karlin-McGregor formula it is not hard to show that the infinite time horizon probability of non-collision for independent Yule chains is strictly positive (and in fact explicit). In particular, a statement like Theorem \ref{collisionThm} cannot always be true if one drops assumptions \textbf{(A1)} and \textbf{(A2)}. It would be interesting to extend our asymptotic analysis to cover cases of unbounded rates or rates decaying to zero in the future.
    \end{example}

\subsection{An inhomogeneous-space interacting particle system}\label{SectionInterlacingParticleSystem}

\begin{figure}[ht!]
    \centering
    \resizebox{\textwidth}{!}{\begin{tikzpicture}[x=0.92cm,y=0.92cm,>=Triangle,line cap=round,line join=round]
    \tikzset{
        particle/.style={circle,fill=black,inner sep=1.75pt},
        particlehi/.style={circle,fill=black,inner sep=2.05pt},
        lowerhi/.style={circle,fill=black,inner sep=1.9pt},
        futureblue/.style={circle,draw=blue!70!black,fill=blue!18,inner sep=1.65pt,line width=0.55pt},
        futuregreen/.style={circle,draw=green!45!black,fill=green!18,inner sep=1.65pt,line width=0.55pt},
        lattice/.style={draw=black!9,line width=0.26pt},
        guide/.style={draw=black!13,line width=0.38pt},
        rowline/.style={draw=black!18,line width=0.38pt},
        breakline/.style={draw=black!36,line width=0.42pt,line join=round,line cap=round,preaction={draw=white,line width=2.0pt}},
        freejump/.style={-{Triangle[length=4.4pt,width=5.2pt]},draw=blue!70!black,line width=1.05pt},
        blockjump/.style={-{Triangle[length=4.4pt,width=5.2pt]},draw=red!75!black,dashed,line width=1.0pt},
        pushjump/.style={-{Triangle[length=4.4pt,width=5.2pt]},draw=green!45!black,line width=1.05pt},
        title/.style={font=\small\bfseries,align=left},
        subtitle/.style={font=\scriptsize,align=left,text=black!70},
        level/.style={font=\scriptsize,inner sep=1pt,text=black!78},
        tinylabel/.style={font=\tiny,inner sep=0.8pt},
        particlelabel/.style={font=\tiny,inner sep=0.6pt},
        dotslabel/.style={font=\scriptsize,inner sep=0.6pt,text=black!60},
        paneltitle/.style={font=\scriptsize\bfseries,align=center},
        panelnote/.style={font=\scriptsize,align=center,text=black!70,text width=3.45cm},
        markblue/.style={circle,draw=blue!55!black,fill=blue!7,inner sep=1.2pt,font=\tiny,text=blue!55!black},
        markred/.style={circle,draw=red!65!black,fill=red!7,inner sep=1.2pt,font=\tiny,text=red!65!black},
        markgreen/.style={circle,draw=green!45!black,fill=green!8,inner sep=1.2pt,font=\tiny,text=green!35!black}
    }

    \begin{scope}[shift={(3.00,0)}]
        \fill[blue!4,rounded corners=2pt] (-0.25,3.94) rectangle (8.25,4.48);
        \fill[green!4,rounded corners=2pt] (0.30,3.10) rectangle (7.55,3.64);

        \foreach \x in {0.00,0.60,...,7.80} {\draw[lattice] (\x,-0.16)--(\x,4.48);}
        \foreach \y in {-0.08,0.62,1.34,2.50,3.35,4.20} {\draw[lattice] (-0.35,\y)--(7.90,\y);}

        \draw[rowline] (-0.35,4.20)--(1.50,4.20);
        \draw[rowline] (2.80,4.20)--(4.50,4.20);
        \draw[rowline] (5.80,4.20)--(7.55,4.20);
        \draw[rowline] (-0.35,3.35)--(2.10,3.35);
        \draw[rowline] (3.40,3.35)--(5.05,3.35);
        \draw[rowline] (6.40,3.35)--(7.55,3.35);
        \draw[rowline] (-0.35,2.50)--(2.65,2.50);
        \draw[rowline] (4.00,2.50)--(5.65,2.50);
        \draw[rowline] (2.05,1.34)--(5.25,1.34);
        \draw[rowline] (2.70,0.62)--(4.50,0.62);
        \draw[rowline] (3.20,-0.08)--(4.00,-0.08);
        \node[level,left] at (-0.42,4.20) {Level $N$};
        \node[level,left] at (-0.42,3.35) {Level $N-1$};
        \node[level,left] at (-0.42,2.50) {Level $N-2$};
        \node[dotslabel,left] at (-0.42,1.92) {$\vdots$};
        \node[level,left] at (-0.42,1.34) {Level $3$};
        \node[level,left] at (-0.42,0.62) {Level $2$};
        \node[level,left] at (-0.42,-0.08) {Level $1$};

        \foreach \a/\b in {0.00/0.60,1.20/0.60,1.20/1.80,3.00/3.60,4.20/3.60,4.20/4.80,6.00/6.60,7.20/6.60} {\draw[guide] (\a,4.20)--(\b,3.35);}
        \foreach \a/\b in {0.60/1.20,1.80/1.20,1.80/2.40,3.60/4.20,4.80/4.20,4.80/5.40} {\draw[guide] (\a,3.35)--(\b,2.50);}
        \draw[guide] (2.40,2.50)--(2.70,2.10);
        \draw[guide] (4.20,2.50)--(3.90,2.10);
        \draw[guide] (4.20,2.50)--(4.20,2.10);
        \draw[guide] (5.40,2.50)--(5.70,2.10);
        \draw[guide] (2.05,1.72)--(2.40,1.34);
        \draw[guide] (3.15,1.72)--(2.40,1.34);
        \draw[guide] (3.15,1.72)--(3.60,1.34);
        \draw[guide] (4.65,1.72)--(3.60,1.34);
        \draw[guide] (4.65,1.72)--(4.80,1.34);
        \draw[guide] (5.40,1.72)--(4.80,1.34);
        \foreach \a/\b in {2.40/3.00,3.60/3.00,3.60/4.20,4.80/4.20} {\draw[guide] (\a,1.34)--(\b,0.62);}
        \foreach \a/\b in {3.00/3.60,4.20/3.60} {\draw[guide] (\a,0.62)--(\b,-0.08);}

        \draw[breakline] (1.95,4.45)--(2.19,4.24)--(2.16,3.98)--(2.47,3.72)--(2.44,3.46)--(2.75,3.20)--(2.72,2.94)--(3.03,2.68)--(3.00,2.42)--(3.31,2.16)--(3.15,1.92);
        \draw[breakline] (5.05,4.45)--(5.29,4.24)--(5.26,3.98)--(5.57,3.72)--(5.54,3.46)--(5.85,3.20)--(5.82,2.94)--(6.13,2.68)--(6.10,2.42)--(6.41,2.16)--(6.35,1.92);
        \draw[breakline] (-0.05,1.84)--(0.35,1.76)--(0.75,1.92)--(1.15,1.76)--(1.55,1.92)--(1.95,1.76)--(2.35,1.92)--(2.75,1.76)--(3.15,1.92)--(3.55,1.76)--(3.95,1.92)--(4.35,1.76)--(4.75,1.92)--(5.15,1.76)--(5.55,1.92)--(5.95,1.76)--(6.35,1.92)--(6.75,1.76)--(7.15,1.92)--(7.55,1.76)--(7.95,1.92)--(8.15,1.84);

        \node[particlehi] at (0.00,4.20) {};
        \node[particlehi] at (1.20,4.20) {};
        \node[particlehi] at (3.00,4.20) {};
        \node[particlehi] at (4.20,4.20) {};
        \node[particlehi] at (6.00,4.20) {};
        \node[particlehi] at (7.20,4.20) {};
        \node[particlelabel,above=5pt] at (0.00,4.20) {$\mathscr X_1^{(N)}$};
        \node[particlelabel,above=5pt] at (1.20,4.20) {$\mathscr X_2^{(N)}$};
        \node[particlelabel,above=5pt] at (3.00,4.20) {$\mathscr X_i^{(N)}$};
        \node[particlelabel,above=5pt] at (4.20,4.20) {$\mathscr X_{i+1}^{(N)}$};
        \node[particlelabel,above=5pt] at (6.00,4.20) {$\mathscr X_{N-1}^{(N)}$};
        \node[particlelabel,above=5pt] at (7.20,4.20) {$\mathscr X_N^{(N)}$};

        \node[lowerhi] at (0.60,3.35) {};
        \node[lowerhi] at (1.80,3.35) {};
        \node[lowerhi] at (3.60,3.35) {};
        \node[lowerhi] at (4.80,3.35) {};
        \node[lowerhi] at (6.60,3.35) {};
        \node[particlelabel,below=5pt] at (0.60,3.35) {$\mathscr X_1^{(N-1)}$};
        \node[particlelabel,below=5pt] at (1.80,3.35) {$\mathscr X_2^{(N-1)}$};
        \node[particlelabel,below=5pt] at (3.60,3.35) {$\mathscr X_i^{(N-1)}$};
        \node[particlelabel,below=5pt] at (4.80,3.35) {$\mathscr X_{i+1}^{(N-1)}$};
        \node[particlelabel,below=5pt] at (6.60,3.35) {$\mathscr X_{N-1}^{(N-1)}$};

        \node[particle] at (1.20,2.50) {};
        \node[particle] at (2.40,2.50) {};
        \node[particle] at (4.20,2.50) {};
        \node[particle] at (5.40,2.50) {};
        \node[particlelabel,below=2pt] at (1.20,2.50) {$\mathscr X_1^{(N-2)}$};
        \node[particlelabel,below=2pt] at (2.40,2.50) {$\mathscr X_2^{(N-2)}$};
        \node[particlelabel,below=4pt] at (4.08,2.50) {$\mathscr X_i^{(N-2)}$};
        \node[particlelabel,below=4pt] at (5.56,2.50) {$\mathscr X_{i+1}^{(N-2)}$};
        \node[particle] at (2.40,1.34) {};
        \node[particle] at (3.60,1.34) {};
        \node[particle] at (4.80,1.34) {};
        \node[particlelabel,below=4pt] at (2.40,1.34) {$\mathscr X_1^{(3)}$};
        \node[particlelabel,below=4pt] at (3.60,1.34) {$\mathscr X_2^{(3)}$};
        \node[particlelabel,below=4pt] at (4.80,1.34) {$\mathscr X_3^{(3)}$};
        \node[particle] at (3.00,0.62) {};
        \node[particle] at (4.20,0.62) {};
        \node[particlelabel,below=5pt] at (3.00,0.62) {$\mathscr X_1^{(2)}$};
        \node[particlelabel,below=5pt] at (4.20,0.62) {$\mathscr X_2^{(2)}$};
        \node[particle] at (3.60,-0.08) {};
        \node[particlelabel,below=5pt] at (3.60,-0.08) {$\mathscr X_1^{(1)}$};
    \end{scope}

    \begin{scope}[shift={(0,-3.00)}]
        \node[paneltitle] at (2.00,1.95) {(1) Free jump};
        \draw[black!12,rounded corners=2pt,line width=0.35pt] (-0.30,-0.66) rectangle (4.30,1.32);
        \fill[blue!4,rounded corners=1.5pt] (-0.18,0.78) rectangle (4.20,1.22);
        \fill[green!4,rounded corners=1.5pt] (-0.18,-0.22) rectangle (4.20,0.22);
        \foreach \x in {0,1,2,3,4} {\draw[lattice] (\x,-0.28)--(\x,1.26);}
        \foreach \y in {0,1} {\draw[lattice] (-0.18,\y)--(4.20,\y);}
        \draw[rowline] (-0.18,1)--(4.20,1);
        \draw[rowline] (-0.18,0)--(4.20,0);
        \node[particle] at (1,1) {};
        \node[futureblue] at (2,1) {};
        \node[particle] at (4,1) {};
        \node[lowerhi] at (3,0) {};
        \draw[guide] (1,1)--(3,0);
        \draw[guide] (3,0)--(4,1);
        \node[particlelabel,above=6pt] at (1,1) {$\mathscr X_i^{(N)}$};
        \node[particlelabel,above=6pt] at (4,1) {$\mathscr X_{i+1}^{(N)}$};
        \node[particlelabel] at (3,-0.34) {$\mathscr X_i^{(N-1)}$};
        \draw[freejump] (1.12,1.08) to[out=62,in=118,looseness=1.45] node[midway,below,tinylabel] {$\lambda_x$} (1.88,1.08);
    \end{scope}

    \begin{scope}[shift={(4.95,-3.00)}]
        \node[paneltitle] at (2.00,1.95) {(2) Blocked attempt};
        \draw[black!12,rounded corners=2pt,line width=0.35pt] (-0.30,-0.66) rectangle (4.30,1.32);
        \fill[blue!4,rounded corners=1.5pt] (-0.18,0.78) rectangle (4.20,1.22);
        \fill[green!4,rounded corners=1.5pt] (-0.18,-0.22) rectangle (4.20,0.22);
        \foreach \x in {0,1,2,3,4} {\draw[lattice] (\x,-0.28)--(\x,1.26);}
        \foreach \y in {0,1} {\draw[lattice] (-0.18,\y)--(4.20,\y);}
        \draw[rowline] (-0.18,1)--(4.20,1);
        \draw[rowline] (-0.18,0)--(4.20,0);
        \node[particle] at (1,1) {};
        \node[particle] at (3,1) {};
        \node[lowerhi] at (1,0) {};
        \node[lowerhi] at (4,0) {};
        \draw[guide] (1,1)--(1,0);
        \draw[guide] (1,0)--(3,1);
        \draw[guide] (3,1)--(4,0);
        \node[particlelabel,above=6pt] at (1,1) {$\mathscr X_i^{(N)}$};
        \node[particlelabel,above=6pt] at (3,1) {$\mathscr X_{i+1}^{(N)}$};
        \node[particlelabel] at (0.92,-0.34) {$\mathscr X_i^{(N-1)}$};
        \node[particlelabel] at (4,-0.34) {$\mathscr X_{i+1}^{(N-1)}$};
        \draw[blockjump] (1.12,1.08) to[out=62,in=118,looseness=1.45] node[midway,below,tinylabel] {$\lambda_x$} (1.88,1.08);
        \draw[red!75!black,line width=0.9pt] (1.86,0.86)--(2.14,1.14);
        \draw[red!75!black,line width=0.9pt] (1.86,1.14)--(2.14,0.86);
    \end{scope}

    \begin{scope}[shift={(9.90,-3.00)}]
        \node[paneltitle] at (2.00,1.95) {(3) Push};
        \draw[black!12,rounded corners=2pt,line width=0.35pt] (-0.30,-0.66) rectangle (4.30,1.32);
        \fill[blue!4,rounded corners=1.5pt] (-0.18,0.78) rectangle (4.20,1.22);
        \fill[green!4,rounded corners=1.5pt] (-0.18,-0.22) rectangle (4.20,0.22);
        \foreach \x in {0,1,2,3,4} {\draw[lattice] (\x,-0.28)--(\x,1.26);}
        \foreach \y in {0,1} {\draw[lattice] (-0.18,\y)--(4.20,\y);}
        \draw[rowline] (-0.18,1)--(4.20,1);
        \draw[rowline] (-0.18,0)--(4.20,0);
        \node[particle] at (0,1) {};
        \node[particle] at (2,1) {};
        \node[futuregreen] at (3,1) {};
        \node[lowerhi] at (1,0) {};
        \node[futureblue] at (2,0) {};
        \draw[guide] (0,1)--(1,0);
        \draw[guide] (1,0)--(2,1);
        \node[particlelabel,above=6pt] at (0,1) {$\mathscr X_i^{(N)}$};
        \node[particlelabel,above=6pt] at (2,1) {$\mathscr X_{i+1}^{(N)}$};
        \node[particlelabel] at (1,-0.34) {$\mathscr X_i^{(N-1)}$};
        \draw[freejump] (1.12,0.08) to[out=62,in=118,looseness=1.45] node[midway,below,tinylabel] {$\lambda_x$} (1.88,0.08);
        \draw[pushjump] (2.12,1.08) to[out=62,in=118,looseness=1.45] (2.88,1.08);
    \end{scope}
\end{tikzpicture}
}
    \caption{Continuous-time pure-birth push-block dynamics on an interlacing array, following \cite{assiotis2023integrablemodelsinhomogeneousspace}. The upper part shows a generic triangular interlacing configuration; jagged breaks indicate omitted particles and omitted intermediate levels. The top level is highlighted in blue and the next level is highlighted in green. The lower panels isolate the three local mechanisms on these two levels, for the particle attempting to move being at site $x$: a free right jump at rate $\lambda_x$, a blocked attempted jump when the same-rank lower-level particle is at the jumper's current site, and a push caused by a lower-level jump into the site of the upper-right particle. Future or proposed particle locations are coloured to match the jump arrow that creates them.}
    \label{interlacingpushblockfig}
\end{figure}
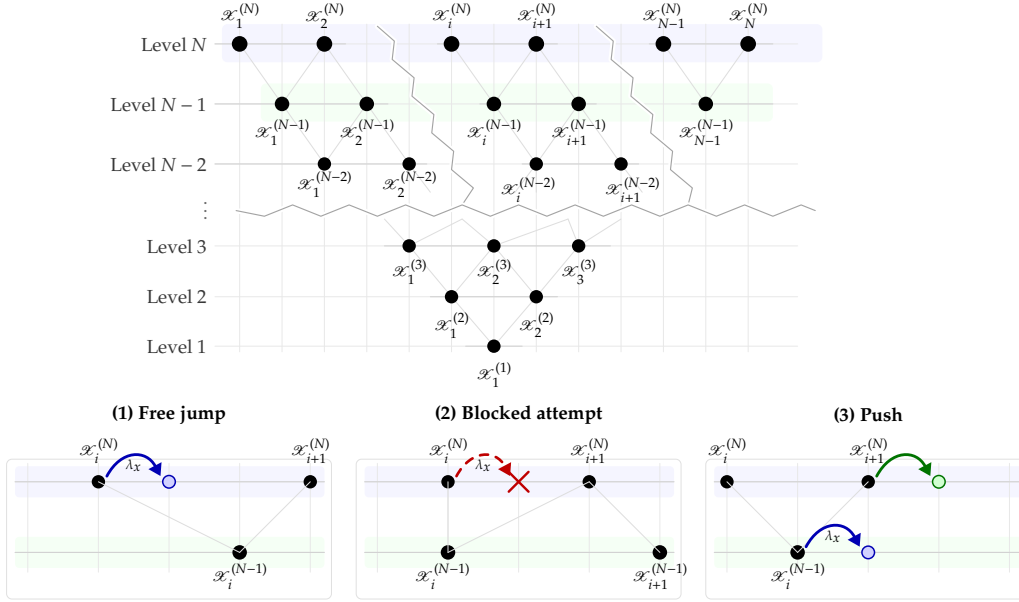

In this section we describe an interacting particle system in space inhomogeneous environment from \cite{AssiotisDeterminantal,assiotis2023integrablemodelsinhomogeneousspace} which partly motivated the present study. This model couples inhomogeneous space TASEP-like and Push-TASEP like particle systems to non-colliding processes. We will describe the continuous time model as it is the simplest. Discrete time Bernoulli and geometric jump models also exist but are more complicated, see \cite{assiotis2023integrablemodelsinhomogeneousspace}. For the special case of homogeneous rates $\lambda_x\equiv 1$ this model and its discrete time variants specialises to seminal works of Borodin-Ferrari \cite{BorodinFerrari}, Warren-Windridge \cite{WarrenWindridge} and Nordenstam \cite{Nordenstam}, see also \cite{BorodinFerrariTilings,MauriceInterlacing}.

We say $\mathbf{x}\in \mathbf{W}_{\mathrm{d}}^N$ interlaces with $\mathbf{y}\in \mathbf{W}_{\mathrm{d}}^{N+1}$ and denote this by $\mathbf{x}\prec_N^{N+1}\mathbf{y}$ if the following inequalities hold:
\begin{equation*}
y_1\le x_1 < y_2 \le x_2 <\cdots < y_N \le x_N <y_{N+1}.
\end{equation*}
We define the interlacing array $\mathsf{IA}_N$ with $N$-levels by,
\begin{equation}
\mathsf{IA}_N=\left\{\left(\mathbf{x}^{(1)},\dots,\mathbf{x}^{(N)}\right)\in \prod_{j=1}^N \mathbf{W}_\mathrm{d}^j: \mathbf{x}^{(i)}\prec_{i}^{i+1} \mathbf{x}^{(i+1)} \right\}.
\end{equation}

We denote the $\mathsf{IA}_N$-valued interacting process by $(\mathscr{X}^{(1)},\dots,\mathscr{X}^{(N)})$ as in Figure \ref{interlacingpushblockfig}. It is described as follows. Individual coordinates $\mathscr{X}_i^{(j)}$ evolve as independent pure birth chains with rates $(\lambda_x)_{x\in \mathbb{Z}_+}$ except when interlacing is about to break, i.e. at the boundary of $\mathsf{IA}_N$, at which point the coordinates interact according to  the following push-block dynamics. Lower level (lower superscript) particles \textbf{(a)} push (induce instantaneous jumps) of the upper level particle(s) or \textbf{(b)} block the jumps of upper level particles for interlacing to remain true. See Figure \ref{interlacingpushblockfig} for a precise illustration of the rules.

Observe that, the projections on the left and right edge particle systems of the array, $(\mathscr{X}^{(1)}_1,\mathscr{X}_1^{(2)},\dots,\mathscr{X}_1^{(N)})$ and $(\mathscr{X}_1^{(1)},\mathscr{X}_2^{(2)},\dots,\mathscr{X}_N^{(N)})$ respectively, are autonomous. The right edge process is inhomogeneous space PushTASEP, see \cite{PushTASEPetrov,AssiotisDeterminantal,assiotis2023integrablemodelsinhomogeneousspace}, while the left edge process under a shift in coordinates $x_1^{(i)}\to x_1^{(i)}-i+1$ becomes a variation of inhomogeneous TASEP.

The connection to the non-colliding process $\hat{\mathbf{X}}$ of this paper is much less obvious.  The following is a non-trivial theorem from \cite{AssiotisDeterminantal,assiotis2023integrablemodelsinhomogeneousspace} which showed this connection to the Doob transformed by $\mathfrak{h}_N$ pure-birth chains. By Corollary \ref{Thm1} this stochastic process indeed matches the conditioned to never intersect process $\hat{\mathbf{X}}$ confirming the prediction from \cite{assiotis2023integrablemodelsinhomogeneousspace}. First, for each fixed $\mathbf{x}\in \mathbf{W}_{\mathrm{d}}^{N}$ define the positive kernel $(\mathbf{y}^{(1)},\dots,\mathbf{y}^{(N-1)})\mapsto \mathcal{L}_N\left(\mathbf{y}^{(1)},\dots,\mathbf{y}^{(N-1)}\big|\mathbf{x}\right)$ on $\prod_{j=1}^{N-1}\mathbf{W}_\mathrm{d}^j$ by the expression:
\begin{equation*}
\mathcal{L}_N\left(\mathbf{y}^{(1)},\dots,\mathbf{y}^{(N-1)}\big|\mathbf{x}\right)=\frac{1}{\mathfrak{h}_N\left(\mathbf{x}\right)} \prod_{j=1}^{N-1} \prod_{i=1}^j \frac{1}{\lambda_{y_i^{(j)}}}\mathbf{1}_{\mathbf{y}^{(1)}\prec_{1}^{2}\mathbf{y}^{(2)}\prec_{2}^{3} \cdots \prec_{N-2}^{N-1} \mathbf{y}^{(N-1)}\prec_{N-1}^N\mathbf{x}}
\end{equation*}
It is known, see \cite{AssiotisDeterminantal,assiotis2023integrablemodelsinhomogeneousspace}, that for any $\mathbf{x}\in \mathbf{W}_{\mathrm{d}}^N$, $\mathcal{L}_N(\cdot|\mathbf{x})$ is actually a stochastic kernel. Let $\mu_N$ be a probability measure on $\mathbf{W}_\mathrm{d}^N$. Suppose the process $(\mathscr{X}^{(i)})_{i=1}^N$ in $\mathsf{IA}_N$ described above is distributed at time $t=0$ according to the probability measure on $\mathsf{IA}_N$ with probability mass function,
\begin{equation*}
(\mathbf{x}^{(1)},\dots,\mathbf{x}^{(N)})\mapsto \mu_N(\mathbf{x}^{(N)})\mathcal{L}_N\left(\mathbf{x}^{(1)},\dots,\mathbf{x}^{(N-1)}\big|\mathbf{x}^{(N)}\right).
\end{equation*}
Then, the projection on $\mathscr{X}^{(N)}$ is a Markov process (in its own filtration) and it coincides in distribution with the conditioned process $\hat{\mathbf{X}}$ with initial distribution $\mu_N$, see \cite{AssiotisDeterminantal,assiotis2023integrablemodelsinhomogeneousspace} for more details.

\subsection{Strategy of proof}
Since the proof is rather long and technical we give an outline. There are certain technical parts of the strategy one could try to condense or circumvent altogether using shorter model specific arguments. We have chosen to develop and present them in rather general form so that it will be easier to re-use or adapt them for future work. Broadly the proof proceeds in seven steps.

\textbf{Step 1: deterministic phase control.}
Section~\ref{glo} establishes the basic analytic properties of the phase $\mathfrak{F}_y(\cdot,t)$ from \hyperref[form]{Definition~\ref{form}}. The first difficulty is that the phase depends on the terminal endpoint $y(t)$ and is not centered a priori. Assumption \hyperref[A2]{\textbf{(A2)}} identifies the correct centering through the constant $\Lambda$ and therefore the deterministic location $\mathrm{z}(t)$, together with its nearest-integer version $\mathscr{Z}(t)$, from \hyperref[themeananddeviation]{Definition~\ref{themeananddeviation}}. In particular,
\begin{equation*}
    \mathfrak{F}_{y}^{(1)}(0,t)=\sum_{n=0}^{y(t)}\frac{1}{\lambda_n}-\Lambda\mathrm{z}(t),
\end{equation*}
so the first derivative at the origin measures the centering error, up to the bounded rounding error incurred when $\mathrm{z}(t)$ is replaced by $\mathscr{Z}(t)$. The second difficulty is to show that, after centering, the phase has a stable quadratic approximation with controlled higher derivatives. This is resolved in \hyperref[specificcasesub]{Proposition~\ref{specificcasesub}}, which proves holomorphicity, identifies the fluctuation scale, and bounds the higher derivatives in a fixed neighbourhood of the origin. 
Moreover, on the finite-time endpoint window $\mathscr{C}_{t}$ from Definition \ref{def:endpoint-window}, the same proposition gives an explicit comparison between $\mathfrak{F}_{y}^{(2)}(w,t)$ and $\Sigma(t)^2$ for $w$ in a fixed neighbourhood of the origin. This is the deterministic input used later to solve the stationary-point equation and to justify the local contour expansion.

\textbf{Step 2: stationary point and contour geometry.}
Sections~\ref{ASP} and \ref{CM} determine where the contour integral should be expanded and along which path it should be evaluated. The first difficulty is to prove that the stationary-point equation has a unique relevant solution near the origin and that this is the solution selected by the minimisation rule in Definition~\ref{definitionstationary}. Section~\ref{ASP} proves this directly on the endpoint window $\mathscr{C}_{t}$: \hyperref[axpro]{Proposition~\ref{axpro}} shows that the phase is a controlled perturbation of a quadratic function and that the selected stationary point $\mathrm{u}_y(t)$ satisfies
\begin{equation*}
    \left|\mathrm{u}_y(t)+\frac{\mathfrak{F}_y^{(1)}(0,t)}{\mathfrak{F}_y^{(2)}(0,t)}\right|
    \leq
    \mathrm{C}_{\mathrm{A}}\!\left(\frac{\mathfrak{F}_y^{(1)}(0,t)}{\mathfrak{F}_y^{(2)}(0,t)}\right)^{\!2};
\end{equation*}
see \eqref{secondconditiontosat}. Combined with the endpoint estimates in \hyperref[specificcasesub]{Proposition~\ref{specificcasesub}}, this quadratic comparison is the quantitative reason the contour geometry can be chosen with constants independent of the terminal endpoint. The second difficulty is geometric: one needs a contour on which the real part of the phase decreases away from this stationary point. This is resolved in Section~\ref{CM}. The local and global pieces of the contour are constructed in \hyperref[contour1]{Proposition~\ref{contour1}} and \hyperref[descentcurve]{Proposition~\ref{descentcurve}}, and then combined in \hyperref[steepestdescentpath]{Corollary~\ref{steepestdescentpath}}. The outcome is a contour along which the integral separates into a small neighborhood of the stationary point and a remainder on which the exponential term is strictly smaller.

\textbf{Step 3: one-particle decoupling on the endpoint window.}
Section~\ref{decouplingsection} is the culmination of the deterministic one-particle contour analysis. At this stage the endpoint is not yet random; rather, we fix an arbitrary terminal value in the window $\mathscr{C}_{t}$ and prove an expansion that is uniform over that whole window. The main difficulty is to separate the universal contribution of the exponential phase from the dependence on the starting point $x$. The contour decomposition from Step~2 reduces the transition probability to a local integral near $\mathrm{u}_y(t)$ and a complementary contour integral. The complementary part is bounded exponentially, while the local part is expanded by Taylor expanding the prefactor $\rho_x$ and evaluating the resulting Gaussian-type coefficients. In this way \hyperref[Thm2]{Proposition~\ref{Thm2}} proves that, for constants $\mathrm{C}_{x},\mathrm{c}_{\varphi}>0$ independent of the endpoint,
\begin{equation*}
    \left|\frac{\mathbb{P}_{x}[\mathscr{X}(t)=y(t)]}{\mathbb{P}_{0}[\mathscr{X}(t)=y(t)]}
    -
    \sum_{n=0}^{x}\frac{\rho_x^{(n)}(\mathrm{u}_y(t))}{n!}\eta_{n,y}(t)\right|
    \leq \mathrm{C}_{x}\mathrm{e}^{-\mathrm{c}_{\varphi}\Sigma(t)^{2/3}},
\end{equation*}
uniformly for endpoints $y(t)\in\mathscr{C}_{t}$. This formula separates the dependence on the starting point, which appears only through derivatives of $\rho_x$, from the phase dependence, which is contained in the coefficients $\eta_{n,y}(t)$. The coefficient estimates in the same proposition then identify the Gaussian orders of the even and odd terms. This is the one-particle expansion later inserted into the Karlin--McGregor determinant \cite{karlin1959coincidence,gessel1985binomial}.

\textbf{Step 4: probabilistic localisation to the endpoint window.}
The deterministic expansion from Step~3 is useful for random endpoints only after proving that the process enters the endpoint window with high probability. This is the role of Section~\ref{MOM}. The main difficulty is that the process is inhomogeneous, so concentration must be proved in the scale determined by the reciprocal rates rather than in a homogeneous diffusive scale. Using the eigenfunction identity in \hyperref[eigenfunction]{Proposition~\ref{eigenfunction}}, Section~\ref{MOM} first proves the quantitative one-time bound in \hyperref[localizationtail]{Proposition~\ref{localizationtail}} for deviations of the form
\begin{equation*}
    |\mathscr{X}(t)-\mathscr{Z}(t)|\geq \kappa\gamma(t)\Sigma(t)^2,
\end{equation*}
and then upgrades this to path tail localization in \hyperref[moments]{Proposition~\ref{moments}}. Applied coordinatewise, this shows that every particle endpoint lies in the corresponding endpoint window after an almost surely finite random time. The quantity that governs the quality of this localization is
\begin{equation*}
    \left|\frac{1}{\mathscr{Z}(t)+1}\sum_{n=0}^{\mathscr{Z}(t)}\frac{1}{\lambda_n}-\Lambda\right|,
\end{equation*}
namely the discrepancy between the local reciprocal-rate average and its limiting value $\Lambda$. The definition of $\gamma(t)$ in \hyperref[gammadef]{Definition~\ref{gammadef}} is chosen precisely so that this discrepancy is absorbed into the fluctuation window. Once Step~4 is complete, the later determinant expansions may restrict their endpoint sums to $\mathscr{C}_{t}$, with the complement controlled by the tail localization estimate.

\textbf{Step 5: stationary point stochastic process and Gaussian limit.}
Section~\ref{stationarypointprocess} next turns the stationary point itself into a random observable through Definition~\ref{definitionstationary}. In the single-valued regime supplied by \hyperref[axpro]{Proposition~\ref{axpro}},
\begin{equation*}
    \mathscr{U}_i(t)=-\Sigma(t)\,\mathrm{u}_{\mathscr{X}_i(t)}(t),
    \qquad i\in\llbracket N\rrbracket.
\end{equation*}
The difficulty is to show that this quantity has a non-degenerate limit and therefore is the correct fluctuation coordinate. Using \hyperref[axpro]{Proposition~\ref{axpro}} together with \hyperref[specificcasesub]{Proposition~\ref{specificcasesub}}\,\hyperref[phaseprop:iv]{(iv)}, Section~\ref{stationarypointprocess} proves on the endpoint-window event that
\begin{equation*}
    \left|\mathscr{U}_i(t)-\frac{\mathfrak{F}^{(1)}_{\mathscr{X}_i}(0,t)}{\Sigma(t)}\right|
    \leq
    \frac{4\mathrm{C}_{\mathrm{A}}}{\Sigma(t)^3}\bigl|\mathfrak{F}^{(1)}_{\mathscr{X}_i}(0,t)\bigr|^2
    +
    \frac{2\mathrm{C}_{\mathscr{C}}\gamma(t)^{1/2}}{\Sigma(t)}\bigl|\mathfrak{F}^{(1)}_{\mathscr{X}_i}(0,t)\bigr|;
\end{equation*}
see \eqref{eq:Uidentify}. Thus $\mathscr{U}_i(t)$ differs from the normalized first derivative $\mathfrak{F}^{(1)}_{\mathscr{X}_i}(0,t)/\Sigma(t)$ by an explicitly bounded error on the endpoint-window event. Moment bounds are obtained in \hyperref[boundedmoments]{Proposition~\ref{boundedmoments}}, and the Gaussian limit is proved in \hyperref[NormalisedProcess]{Proposition~\ref{NormalisedProcess}} by evaluating characteristic functions through the eigenfunction identity. The endpoint-ordering statement in \hyperref[twobadprobs-fixed]{Proposition~\ref{twobadprobs-fixed}} then shows that, on $\mathscr{C}_{t}^{N}$, ordering of the terminal particles is equivalent to ordering of the stationary-point coordinates. This step is needed because the final determinant asymptotics are expressed in terms of the random vector $\mathbf{U}(t)$ rather than directly in the particle coordinates.

\textbf{Step 6: determinant asymptotics and dominant-index analysis.}
We next lift the one-particle asymptotic expansion from Step~3 to the $N$-particle Karlin--McGregor determinant. Using this determinant formula we can write
\begin{equation*}
    \mathbb{P}_{\mathbf{x}}^N[\mathcal{T}>t]
    =
    \sum_{\mathbf{y}\in\mathbf{W}_{\mathrm d}^{N}}
    \det\left(\mathbb{P}_{x_i}\left[\mathscr{X}(t)=y_j\right]\right)_{i,j=1}^{N}.
\end{equation*}
On the endpoint window, we insert the expansion from \hyperref[Thm2]{Proposition~\ref{Thm2}}, use \hyperref[twobadprobs-fixed]{Proposition~\ref{twobadprobs-fixed}} to replace particle ordering by stationary-point ordering, and then reorganize the resulting determinant by means of \hyperref[Gcauchybinet]{Lemma~\ref{Gcauchybinet}}. This reduces the problem to a sum over index families $(\mathbf n,\mathbf k)$, and the remaining task is to determine which of these families contributes at leading order. The two key mechanisms are cancellation and parity. Cancellation comes from antisymmetry of the determinant, which excludes repeated values of $n_i+k_i$. Parity comes from the explicit coefficient estimates in Section~\ref{decouplingsection}: for every $t>\mathrm{T}_{\mathrm{dec}}$ and every path $y\in\mathscr{P}$ satisfying $y(t)\in\mathscr{C}_{t}$,
\begin{equation*}
    \begin{aligned}
    &\left|\eta_{n,y}(t)-(-1)^{\frac{n}{2}}(n-1)!!\,\Sigma(t)^{-n}\right|
    \leq
    \mathrm{C}_n\gamma(t)^{1/2}\Sigma(t)^{-n},
    \quad n\textnormal{ even},\\
    &|\eta_{n,y}(t)|
    \leq \mathrm{C}_n\Sigma(t)^{-n-1},
    \quad n\textnormal{ odd}.
    \end{aligned}
\end{equation*}
The extra loss in the odd case determines the leading index family $\mathcal{M}_N$. On this family, the leading coefficients are deterministic, and the determinant factors into
\begin{equation*}
    \mathfrak{h}_N(\mathbf{x})\,\Delta(\mathbf{U}(t))\,\Sigma(t)^{-\frac{N}{2}(N-1)}.
\end{equation*}
This is proved in \hyperref[ThmKM]{Proposition~\ref{ThmKM}}. The proposition reduces the collision asymptotics and the conditioned limit problem to the expectation in \eqref{expectform}, whose main observable is $\Delta(\mathbf{U}(t))F(\mathbf{U}(t))$ on the ordered endpoint window. That expectation is then evaluated by combining Step~5 with the moment bounds from \hyperref[boundedmoments]{Proposition~\ref{boundedmoments}}.

\textbf{Step 7: conditioning and main limits.}
The final section combines the previous steps to prove the statements from the introduction. Theorem~\ref{collisionThm} is obtained by taking $F\equiv 1$ in the asymptotic formula of \hyperref[ThmKM]{Proposition~\ref{ThmKM}}. Corollary~\ref{Thm1} then follows by inserting this collision-tail asymptotic into the finite-horizon conditioning identity \eqref{construction}. The shifted survival tail $H_{t,s}$ is compared with the initial-time survival tail $H_{0,s}$ using \hyperref[lem:shifted-deterministic-scales]{Lemma~\ref{lem:shifted-deterministic-scales}}; this is the step that constructs the infinite-horizon limit and converts the space-time conditioning into the spatial Doob $h$-transform with $\mathfrak h_N$. The survival-conditioned stationary-point limit is then isolated in Proposition~\ref{lastcorcor}, and Theorem~\ref{thm:stationarypoints} follows by combining that result with the free Gaussian limit from \hyperref[NormalisedProcess]{Proposition~\ref{NormalisedProcess}}. Finally, Corollary~\ref{cor:particlecoordinates} transfers the limit from $\mathbf{U}(t)$ to the centered inverse-rate coordinates by using the comparison \eqref{eq:Uidentify}. In this way the one-particle contour analysis from Subsections~\ref{SubsectionTransitionKernel}--\ref{decouplingsection} is converted into the multi-particle conditioned limit laws stated in the introduction.

\subsection{Extensions}\label{SectionExtensions}

We believe that a few of the technical innovations developed in this paper may be adapted to tackle qualitatively similar problems. Examples include allied problems concerned with asymptotics of determinants with entries contour integrals involving a diverging number of parameters (having minimal global structure as they do here) such as the following:

\begin{itemize}
	    \item \textbf{Non-colliding birth-and-death chains.}  
    It is natural to extend the method from pure-birth chains to birth-and-death dynamics with rates $\mathfrak{u}(x)$ to jump from $x$ to $x+1$ and $\mathfrak{d}(x)$ to jump from $x$ to $x-1$. By the seminal works of Karlin-McGregor \cite{KarlinMcGregorB&D1,KarlinMcGregorB&D2}, the transition probabilities $\mathfrak{p}_t^{\mathfrak{u},\mathfrak{d}}(x,y)$ of this chain are given in terms of integrals of orthogonal polynomials on $\mathbb{R}_+$. This type of expression for  $\mathfrak{p}_t^{\mathfrak{u},\mathfrak{d}}(x,y)$  does not immediately fit our framework. Nevertheless, an alternative expression in terms of a contour integral of the resolvent of the birth-and-death chain's generator exists which is likely amenable to our approach. 
	    \item \textbf{Large-$N$ limits of non-colliding chains.}  
It would be interesting to consider large $N$ limits of the non-colliding process $\hat{\mathbf{X}}$ we constructed here. For example, in the homogeneous case, at an edge scaling, it is known \cite{DauvergneNicaVirag} that one obtains the Airy line ensemble \cite{CorwinHammond} from KPZ universality class \cite{CorwinKPZ} while in the bulk one obtains the extended discrete sine kernel, see \cite{LocalStatRW}. For quite general random walks (but still space-homogeneous) Airy line ensemble limits have recently been proven in  \cite{DenisovFitzGeraldWachtel2024}, see also \cite{Baik,BaikSuidan} for earlier work. In our setting, for space-inhomogeneous perturbations of the homogeneous case it is very likely one would still get the Airy line ensemble. For complicated space-inhomogeneous environments it is less clear what types of local asymptotic behaviour are possible.
    
	    \item \textbf{Asymptotics for the interlacing interacting particle system.}  
    The present steepest-descent framework could possibly be used to obtain local and global asymptotics for the (2+1)-dimensional interacting particle system  in $\mathsf{IA}_N$ from Section \ref{SectionInterlacingParticleSystem} using the formulae for its space-time correlations obtained in \cite{AssiotisDeterminantal,assiotis2023integrablemodelsinhomogeneousspace}. For asymptotics in related works see \cite{BorodinFerrari,MauriceInterlacing,KnizelPetrovSaenz,expjumpmodel,PushTASEPetrov}.
    
	    \item \textbf{Non-colliding Markov chains in random environment.}  
    A further direction is to take both the space and time environments $(\lambda_n,\alpha_n,\beta_n)$  to be random. For recent works on the closely connected random tiling models with random weights see  \cite{duits2025gammadisorderedaztecdiamond,moulard2025dimerslayereddisorder,bufetov2025dominotilingsaztecdiamond,zografos2025quenchedannealedcltsoneperiodic}.
    
    \item \textbf{Boundary of the inhomogeneous Gelfand-Tsetlin graph.} The Gelfand-Tsetlin graph is a combinatorial model describing the branching of irreducible representations of the chain of unitary groups $\mathbb{U}(1)\subset \mathbb{U}(2)\subset \cdots$. The explicit description of the boundary of this graph, see \cite{BorodinOlshanskiBouquet} for background, is a foundational result in integrable probability and asymptotic algebraic combinatorics, see \cite{VershikKerov,BorodinOlshanskiBoundary,PetrovBoundary,GorinPanova} for different proofs. There exists a natural inhomogeneous generalisation of the Gelfand-Tsetlin graph, closely related to the model from Section \ref{SectionInterlacingParticleSystem}, where edge weights depend in an inhomogeneous way on the vertices they connect via parameters $(\lambda_x)_{x\in \mathbb{Z}_+}$ (for $\lambda_x \equiv 1$ one gets back the Gelfand-Tsetlin graph). For certain classes of inhomogeneous parameters this graph's boundary may also be explicit, see the discussion  in \cite{assiotis2023integrablemodelsinhomogeneousspace}, and the asymptotic analysis problem one needs to solve to obtain this is somewhat connected to the asymptotic considerations in this paper.
 
    \end{itemize}

\paragraph*{Organisation of the paper.}
Section~\ref{glo} develops the analysis of the phase function $\mathfrak{F}_y$, including centering and scaling estimates.  
Section~\ref{ASP} constructs and controls the relevant stationary point.  
Section~\ref{CM} builds suitable steepest-descent contours.  
Section~\ref{decouplingsection} derives the one-particle asymptotic expansion.  
Section~\ref{MOM} proves the large-deviation/concentration inputs that allow this deterministic endpoint-window expansion to be applied to random terminal states.  
Section~\ref{stationarypointprocess} identifies the limiting stationary point stochastic process and performs the Karlin--McGregor determinant asymptotics.  
The final section proves the main theorems and corollaries from these ingredients.

\paragraph*{Acknowledgements.} TA is grateful to Kilian Raschel, Denis Denisov and Vitali Wachtel for useful discussions on asymptotics of collision times of stochastic processes. MFN is grateful to the School of Mathematics of University of Edinburgh for funding through an EPSRC PhD studentship.

\paragraph*{Tool and computational resources disclosure.} In the final stages of preparation of this paper AI tools (ChatGPT 5.5 Pro) were used to assist proofreading and editing the manuscript. The majority of figures were generated using AI tools. The computations in Example \ref{ex:homogeneous-unit-rates} were obtained using AI tools. All ideas and proofs were initially generated by the human authors.

    \section{Analysis of the Markov Kernel}\label{glo}
    This section develops deterministic control of the phase function  $\mathfrak{F}_y$ from Definition~\ref{form}.
    
    \subsection{Contour representation for the Markov kernel}\label{SubsectionTransitionKernel}
    
     We begin by constructing the exact one-particle contour representation for the transition probabilities, from which the phase $\mathfrak{F}_y$ is extracted. The following proposition extends, with a somewhat different proof, results of \cite{assiotis2023integrablemodelsinhomogeneousspace} which imposed an additional technical assumption on the parameters.

    \begin{proposition}[Contour representation in phase form; cf.\ \cite{assiotis2023integrablemodelsinhomogeneousspace}]\label{contourrepresentation}
        Under Assumption \textbf{(A1)}, for all $x\in\mathbb{Z}_+$, $y\in\mathscr{P}$ and $t\in\mathscr{T}$, the transition probabilities of the Markov processes from Definition~\ref{contprocess} and Definition~\ref{mixprocess} are given by
        \begin{equation}\label{cir}
            \mathbb{P}_{x}[\mathscr{X}(t)=y(t)]=-\frac{1}{\lambda_{y(t)}2\pi \mathrm{i}}\oint_{\Gamma_{\lambda}}\rho_x(w)\mathrm{e}^{\mathfrak{F}_y(w,t)}\mathrm{d}w.
        \end{equation}
        Here $(\rho_x)_{x\in\mathbb{Z}_+}$ are the characteristic polynomials from \eqref{characteristic}, $\mathrm{e}^{\mathfrak{F}_y(w,t)}$ is understood through the branch-free product form in \eqref{phaseproductform}, and $\Gamma_{\lambda}$ is a simple positively oriented contour such that $[\mathfrak{L},\mathfrak{B}_1]\subset\Gamma_{\lambda}^{\circ}$ and $\overline{\Gamma_{\lambda}^{\circ}}\subset\{w\in\mathbb{C}:\mathfrak{Re}(w)>-\mathfrak{U}^{-1}\}$, where $\Gamma_{\lambda}^{\circ}$ denotes the interior of $\Gamma_{\lambda}$. The right-hand side depends on the path $y$ only through its terminal value $y(t)$. In particular, in the discrete-time regime, $\Gamma_{\lambda}$ encloses all spatial poles $\lambda_n$ and avoids all geometric poles $-\beta_s^{-1}$, $s\in\mathbb{Z}_+$.
    \end{proposition}

    \begin{proof}
        \medskip

        \noindent\textbf{Construction of the contour kernel.}
        For $r\in\mathbb{R}$ and an open set $U\subset\mathbb{C}$, we define the $r$-shifted right half-plane and the set of holomorphic functions on $U$ by
        \begin{equation*}
            \mathbb{H}_{r}\overset{\textnormal{def}}{=}\{u\in\mathbb{C}:\mathfrak{Re}(u)>r\},\quad\text{and}\quad\mathrm{Hol}(U)\overset{\textnormal{def}}{=}\{f:U\rightarrow\mathbb{C}:\text{$f$ is holomorphic on $U$}\}.
        \end{equation*}
        For $f\in\mathrm{Hol}(\mathbb{H}_{-1/\mathfrak{U}})$ and $x,m\in\mathbb{Z}_+$, we define the kernel
        \begin{equation}\label{defofTf}
            \mathrm{T}^{f}(x,m)\overset{\textnormal{def}}{=}-\frac{1}{\lambda_{m}2\pi\mathrm{i}}\oint_{\Gamma_{\lambda}}\frac{\rho_x(w)}{\rho_{m+1}(w)}f(w)\mathrm{d}w,\qquad x,m\in\mathbb{Z}_{+},
        \end{equation}
        where $\Gamma_{\lambda}$ is a positively oriented contour enclosing the spatial interval $[\mathfrak{L},\mathfrak{B}_1]$ while remaining inside $\mathbb{H}_{-1/\mathfrak{U}}$. By \textbf{(A1)}, the spatial poles satisfy $\mathfrak{L}<\lambda_n<\mathfrak{B}_1$, for all $n\in\mathbb{Z}_+$. Consequently, there is a strictly positive gap between the spatial-pole region $(\mathfrak{L},\mathfrak{B}_1)$ and the boundary of $\mathbb{H}_{-1/\mathfrak{U}}$. A family of admissible choices is given, for $\mu\in(0,1)$, by the positively oriented circles
        \begin{equation}\label{admissiblecontours}
            \Gamma_{\lambda}(\mu)=\left\{w\in\mathbb{C}:\left|w-\frac{\mu\left(\mathfrak{B}_2-\frac{1}{\mathfrak{U}}\right)+(1-\mu)\left(\mathfrak{B}_1+\mathfrak{L}\right)}{2}\right|=\frac{\mu\left(\mathfrak{B}_2+\frac{1}{\mathfrak{U}}\right)+(1-\mu)\left(\mathfrak{B}_1-\mathfrak{L}\right)}{2}\right\},
        \end{equation}
        where $\mathfrak{B}_2$ is the constant from \textbf{(A1)} satisfying $\mathfrak{B}_1<\mathfrak{B}_2$. Indeed, the leftmost point of $\Gamma_{\lambda}(\mu)$ is $(1-\mu)\mathfrak{L}-\mu\mathfrak{U}^{-1}>-\mathfrak{U}^{-1}$, while its rightmost point is $(1-\mu)\mathfrak{B}_1+\mu\mathfrak{B}_2>\mathfrak{B}_1$. Hence $\Gamma_{\lambda}(\mu)$ lies inside $\mathbb{H}_{-1/\mathfrak{U}}$ and encloses every spatial pole $\lambda_n\in(\mathfrak{L},\mathfrak{B}_1)$. This is sketched in Figure~\ref{contourseparationfig}. Consequently, at least one admissible contour exists, and the kernel $\mathrm{T}^f$ is well-defined for every $f\in\mathrm{Hol}(\mathbb{H}_{-1/\mathfrak{U}})$.
        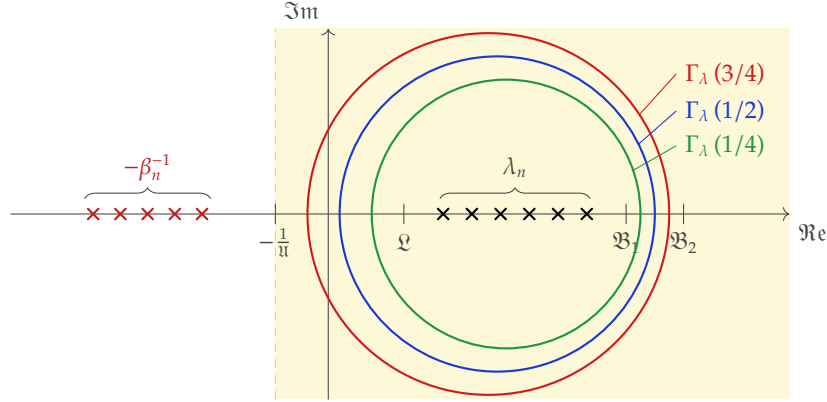
\begin{figure}
            \centering
            \captionsetup{type=figure}
            \begin{tikzpicture}[x=1cm,y=1cm,line cap=round,line join=round]
                \definecolor{poleRed}{RGB}{190,30,30}
                \definecolor{contourBlue}{RGB}{25,55,210}
                \definecolor{contourRed}{RGB}{210,25,25}
                \definecolor{contourGreen}{RGB}{25,150,70}
                \definecolor{contourYellow}{RGB}{245,215,40}
                \definecolor{axisGray}{RGB}{70,70,70}

                \fill[contourYellow!18] (-0.70,-2.45) rectangle (6.10,2.45);
                \draw[contourYellow, dashed] (-0.70,-2.45) -- (-0.70,2.45);

                \draw[->,axisGray] (-4.2,0) -- (6.1,0) node[below right] {\small $\mathfrak{Re}$};
                \draw[->,axisGray] (0,-2.45) -- (0,2.45) node[above left] {\small $\mathfrak{Im}$};

                \draw[poleRed,thick] (-3.18,-0.07) -- (-3.04,0.07);
                \draw[poleRed,thick] (-3.18,0.07) -- (-3.04,-0.07);
                \draw[poleRed,thick] (-2.82,-0.07) -- (-2.68,0.07);
                \draw[poleRed,thick] (-2.82,0.07) -- (-2.68,-0.07);
                \draw[poleRed,thick] (-2.46,-0.07) -- (-2.32,0.07);
                \draw[poleRed,thick] (-2.46,0.07) -- (-2.32,-0.07);
                \draw[poleRed,thick] (-2.10,-0.07) -- (-1.96,0.07);
                \draw[poleRed,thick] (-2.10,0.07) -- (-1.96,-0.07);
                \draw[poleRed,thick] (-1.74,-0.07) -- (-1.60,0.07);
                \draw[poleRed,thick] (-1.74,0.07) -- (-1.60,-0.07);
                \draw[axisGray,decorate,decoration={brace,amplitude=4pt,raise=3pt}] (-3.22,0.10) -- (-1.56,0.10);
                \node[poleRed] at (-2.39,0.62) {\small $-\beta_n^{-1}$};

                \draw[axisGray] (-0.70,-0.13) -- (-0.70,0.13);
                \node[axisGray] at (-0.70,-0.38) {\small $-\frac{1}{\mathfrak{U}}$};

                \draw[axisGray] (1,-0.13) -- (1,0.13);
                \node[axisGray] at (1,-0.38) {\small $\mathfrak{L}$};

                \draw[axisGray] (3.94,-0.13) -- (3.94,0.13);
                \node[axisGray] at (3.94,-0.38) {\small $\mathfrak{B}_1$};

                \draw[axisGray] (4.70,-0.13) -- (4.70,0.13);
                \node[axisGray] at (4.70,-0.38) {\small $\mathfrak{B}_2$};

                \draw[black,thick] (1.45,-0.07) -- (1.59,0.07);
                \draw[black,thick] (1.45,0.07) -- (1.59,-0.07);
                \draw[black,thick] (1.83,-0.07) -- (1.97,0.07);
                \draw[black,thick] (1.83,0.07) -- (1.97,-0.07);
                \draw[black,thick] (2.21,-0.07) -- (2.35,0.07);
                \draw[black,thick] (2.21,0.07) -- (2.35,-0.07);
                \draw[black,thick] (2.59,-0.07) -- (2.73,0.07);
                \draw[black,thick] (2.59,0.07) -- (2.73,-0.07);
                \draw[black,thick] (2.97,-0.07) -- (3.11,0.07);
                \draw[black,thick] (2.97,0.07) -- (3.11,-0.07);
                \draw[black,thick] (3.35,-0.07) -- (3.49,0.07);
                \draw[black,thick] (3.35,0.07) -- (3.49,-0.07);
                \draw[axisGray,decorate,decoration={brace,amplitude=4pt,raise=3pt}] (1.48,0.10) -- (3.46,0.10);
                \node[axisGray] at (2.47,0.62) {\small $\lambda_n$};

                \draw[thick,contourBlue] (0.25*4.70+0.25*3.94-0.25*0.70+0.25,0) circle (0.25*4.70+0.25*3.94+0.25*0.70-0.25);
                \draw[contourBlue,thin] (4.60,1.36) -- (4.10,0.92);
                \node[contourBlue,anchor=west] at (4.60,1.36) {\small $\Gamma_{\lambda}\left(1/2\right)$};

                \draw[thick,contourRed] (0.5*0.75*4.70+0.5*0.25*3.94-0.5*0.75*0.70+0.5*0.25,0) circle (0.5*0.75*4.70+0.5*0.25*3.94+0.5*0.75*0.70-0.5*0.25);
                \draw[contourRed,thin] (4.60,1.85) -- (4.12,1.32);
                \node[contourRed,anchor=west] at (4.60,1.85) {\small $\Gamma_{\lambda}\left(3/4\right)$};

                \draw[thick,contourGreen] (0.5*0.25*4.70+0.5*0.75*3.94-0.5*0.25*0.70+0.5*0.75,0) circle (0.5*0.25*4.70+0.5*0.75*3.94+0.5*0.25*0.70-0.5*0.75);
                \draw[contourGreen,thin] (4.60,0.90) -- (4.02,0.57);
                \node[contourGreen,anchor=west] at (4.60,0.90) {\small $\Gamma_{\lambda}\left(1/4\right)$};

            \end{tikzpicture}
            \caption{Schematic location of the real singularities of the integrand in the contour representation \eqref{cir}, for a sample choice of parameters satisfying \textbf{(A1)}. The dashed vertical line marks the boundary $\mathfrak{Re}(w)=-\mathfrak{U}^{-1}$ of the half-plane $\mathbb{H}_{-1/\mathfrak{U}}$, and the shaded yellow region indicates the domain on which the input functions $f\in\mathrm{Hol}(\mathbb{H}_{-1/\mathfrak{U}})$ are holomorphic. The geometric poles $-\beta_t^{-1}$ lie to the left of this boundary, while the spatial poles $\lambda_n$ lie in $(\mathfrak{L},\mathfrak{B}_1)$; hence there is a positive separation between the geometric and spatial singularities. The displayed contours $\Gamma_{\lambda}(1/4)$, $\Gamma_{\lambda}(1/2)$, and $\Gamma_{\lambda}(3/4)$ are nested admissible contours lying inside $\mathbb{H}_{-1/\mathfrak{U}}$ and enclosing all spatial poles. The Bernoulli factors contribute only zeros and are therefore omitted.}
            \label{contourseparationfig}
        \end{figure}
        Now by the explicit construction in \eqref{admissiblecontours}, we have the following embedding of contours:
        \begin{equation}\label{nestedcontours}
            [\mathfrak{L},\mathfrak{B}_1]\subset\Gamma_{\lambda}\left(1/4\right)^{\circ}\subset \Gamma_{\lambda}\left(1/2\right)^{\circ}\subset \Gamma_{\lambda}\left(3/4\right)^{\circ}\subset\mathbb{H}_{-1/\mathfrak{U}}.
        \end{equation}
        Here, for any closed contour $\Gamma$, $\Gamma^{\circ}$ denotes its interior. We shall use the following contour-independence repeatedly: if two admissible contours are deformed into one another inside $\mathbb{H}_{-1/\mathfrak{U}}$ without crossing the spatial poles enclosed by both contours, then Cauchy's theorem gives the same value of the corresponding integral in \eqref{defofTf}. For any $h\in\mathrm{Hol}(\mathbb{H}_{-1/\mathfrak{U}})$, we also have the upper-triangularity relation
        \begin{equation}\label{uppertriangularT}
            \mathrm{T}^{h}(x,m)=0,\qquad x>m.
        \end{equation}
        Indeed, in this regime, by definition of the characteristic polynomials, the rational function
        \begin{equation*}
            \frac{\rho_x(w)}{\rho_{m+1}(w)}=\prod_{k=m+1}^{x-1}\left(1-\frac{w}{\lambda_k}\right)
        \end{equation*}
        is a polynomial, with the empty product interpreted as one, and is thus holomorphic inside $\Gamma_{\lambda}$. This also proves the composition identity in the case $x>m$, since for each $k\in\mathbb{Z}_+$ either $x>k$ or $k>m$, and therefore \eqref{uppertriangularT} gives $\mathrm{T}^{f}(x,k)\mathrm{T}^{g}(k,m)=0$, while $\mathrm{T}^{fg}(x,m)=0$. Consequently, for $f,g\in\mathrm{Hol}(\mathbb{H}_{-1/\mathfrak{U}})$, it remains to consider $x\le m$, where the upper-triangularity relation \eqref{uppertriangularT} reduces the convolution of the two kernels $\mathrm{T}^{f}$ and $\mathrm{T}^{g}$ to the finite sum
        \begin{equation}\label{finiteconvolutionT}
            (\mathrm{T}^{f}\mathrm{T}^{g})(x,m)=\sum_{k=x}^{m} \mathrm{T}^{f}(x,k)\mathrm{T}^{g}(k,m).
        \end{equation}
       
        \textbf{Multiplicative composition.}
        By the contour-independence following \eqref{nestedcontours}, the contours in the definition of $\mathrm{T}^f$ and $\mathrm{T}^g$ may be deformed to the nested admissible contours used below, because the deformation takes place inside $\mathbb{H}_{-1/\mathfrak{U}}$, keeps the same spatial singularities enclosed, and crosses no singularity of the corresponding integrands. We choose $\Gamma_{\lambda}(1/2)$ for the $u$-integral and $\Gamma_{\lambda}(3/4)$ for the $w$-integral. Therefore, the finite convolution \eqref{finiteconvolutionT} can be rewritten as
        \begin{equation*}
            (\mathrm{T}^{f}\mathrm{T}^{g})(x,m)=\sum_{k=x}^{m}\left[-\frac{1}{\lambda_k2\pi\mathrm{i}}\oint_{\Gamma_{\lambda}\left(1/2\right)}\frac{\rho_x(u)}{\rho_{k+1}(u)}f(u)\mathrm{d}u\right]\left[-\frac{1}{\lambda_m2\pi\mathrm{i}}\oint_{\Gamma_{\lambda}\left(3/4\right)}\frac{\rho_k(w)}{\rho_{m+1}(w)}g(w)\mathrm{d}w\right],
        \end{equation*}
        and subsequently as the double contour integral
        \begin{equation}\label{convolT}
            (\mathrm{T}^{f}\mathrm{T}^{g})(x,m)=\frac{1}{\lambda_m(2\pi\mathrm{i})^2}\oint_{\Gamma_{\lambda}\left(1/2\right)}\oint_{\Gamma_{\lambda}\left(3/4\right)}
            f(u)g(w)\frac{\rho_x(u)}{\rho_{m+1}(w)}
            \left[\sum_{k=x}^{m}\frac{\rho_k(w)}{\lambda_k\rho_{k+1}(u)}\right]\mathrm{d}w\mathrm{d}u.
        \end{equation}
        The preceding step is justified by parametrising the two compact contours and by using the finiteness of the sum over $k$. We now compute the sum in the square bracket explicitly. For $k\in\mathbb{Z}_+$, again by the definition of the characteristic polynomials, we have
        \begin{equation*}
            \frac{\rho_k(w)}{\rho_k(u)}-\frac{\rho_{k+1}(w)}{\rho_{k+1}(u)}=\frac{\rho_k(w)}{\rho_{k+1}(u)}\left(\frac{\rho_{k+1}(u)}{\rho_{k}(u)}-\frac{\rho_{k+1}(w)}{\rho_k(w)}\right)=
            (w-u)\frac{\rho_k(w)}{\lambda_k\rho_{k+1}(u)}.
        \end{equation*}
        Since $\Gamma_{\lambda}(1/2)$ and $\Gamma_{\lambda}(3/4)$ are disjoint, we have $w\neq u$ throughout the double integral in \eqref{convolT}. Hence we may divide the preceding identity by $w-u$, and subsequently telescope the resulting finite sum to obtain
        \begin{equation}\label{telescopingrho}
            \sum_{k=x}^{m}\frac{\rho_k(w)}{\lambda_k\rho_{k+1}(u)}
            =\frac{1}{w-u}\sum_{k=x}^{m}\left(\frac{\rho_k(w)}{\rho_k(u)}-\frac{\rho_{k+1}(w)}{\rho_{k+1}(u)}\right)=\frac{1}{w-u}\left(\frac{\rho_x(w)}{\rho_x(u)}-\frac{\rho_{m+1}(w)}{\rho_{m+1}(u)}\right).
        \end{equation}
        Substituting the telescoping identity \eqref{telescopingrho} back into the square bracket in the double contour integral \eqref{convolT}, we obtain
        \begin{equation*}
            (\mathrm{T}^{f}\mathrm{T}^{g})(x,m)=\frac{1}{\lambda_m(2\pi\mathrm{i})^2}\oint_{\Gamma_{\lambda}\left(1/2\right)}\oint_{\Gamma_{\lambda}\left(3/4\right)}
            f(u)g(w)\left[
            \frac{\rho_x(w)}{\rho_{m+1}(w)}\frac{1}{w-u}
            -\frac{\rho_x(u)}{\rho_{m+1}(u)}\frac{1}{w-u}
            \right]\mathrm{d}w\mathrm{d}u.
        \end{equation*}
        Splitting the preceding display into two double contour integrals, we denote
        \begin{equation}\label{IoneItwo}
            \begin{aligned}
                &I_1(x,m)=\frac{1}{\lambda_m(2\pi\mathrm{i})^2}\oint_{\Gamma_{\lambda}\left(1/2\right)}\oint_{\Gamma_{\lambda}\left(3/4\right)}f(u)g(w)\frac{\rho_x(w)}{\rho_{m+1}(w)}\frac{1}{w-u}\mathrm{d}w\mathrm{d}u,\\&I_2(x,m)=-\frac{1}{\lambda_m(2\pi\mathrm{i})^2}\oint_{\Gamma_{\lambda}\left(1/2\right)}\oint_{\Gamma_{\lambda}\left(3/4\right)}f(u)g(w)\frac{\rho_x(u)}{\rho_{m+1}(u)}\frac{1}{w-u}\mathrm{d}w\mathrm{d}u.
            \end{aligned}
        \end{equation}
        We now deform the contour for the inner integral in $I_1(x,m)$ from $\Gamma_{\lambda}\left(3/4\right)$ to $\Gamma_{\lambda}\left(1/4\right)$. This deformation does not add or subtract any residue contribution from the spatial singularities, since by construction all of them lie inside both contours. The only pole crossed during the deformation is the simple pole at $w=u$, because $g$ is holomorphic on $\mathbb{H}_{-1/\mathfrak{U}}$ and $u\in\Gamma_{\lambda}(1/2)$ lies between the two contours. Hence, by Cauchy's residue theorem, we have
	        \begin{equation}\label{IoneDeform}
	            I_1(x,m)=\frac{1}{\lambda_m(2\pi\mathrm{i})^2}\oint_{\Gamma_{\lambda}\left(1/2\right)}f(u)\left(\oint_{\Gamma_{\lambda}\left(1/4\right)}g(w)\frac{\rho_x(w)}{\rho_{m+1}(w)}\frac{1}{w-u}\mathrm{d}w+2\pi\mathrm{i}g(u)\frac{\rho_x(u)}{\rho_{m+1}(u)}\right)\mathrm{d}u.
	        \end{equation}
	        The sign of the residue term in \eqref{IoneDeform} is positive because the integral over the outer contour $\Gamma_{\lambda}(3/4)$ equals the integral over the inner contour $\Gamma_{\lambda}(1/4)$ plus $2\pi\mathrm{i}$ times the residue of the pole enclosed by the annular region between the two contours.
        Similarly, we can deform the contour for the inner integral in $I_2(x,m)$ from $\Gamma_{\lambda}\left(3/4\right)$ to $\Gamma_{\lambda}\left(1/4\right)$. In this case the $w$-dependent part of the integrand is $g(w)/(w-u)$, so the only pole crossed is again the simple pole at $w=u$. Hence, by Cauchy's residue theorem, we have
        \begin{equation}\label{ItwoDeform}
            I_2(x,m)=-\frac{1}{\lambda_m(2\pi\mathrm{i})^2}\oint_{\Gamma_{\lambda}\left(1/2\right)}f(u)\frac{\rho_x(u)}{\rho_{m+1}(u)}\left(\oint_{\Gamma_{\lambda}\left(1/4\right)}g(w)\frac{1}{w-u}\mathrm{d}w+2\pi\mathrm{i}g(u)\right)\mathrm{d}u.
        \end{equation}
        The same annular orientation gives the positive residue term inside the parentheses in \eqref{ItwoDeform}; the overall negative sign is the one already present in the definition of $I_2$ in \eqref{IoneItwo}.
        For any $u\in\Gamma_{\lambda}\left(1/2\right)$, by construction $u\not\in\Gamma_{\lambda}\left(1/4\right)^{\circ}$. Since $g$ is holomorphic on $\mathbb{H}_{-1/\mathfrak{U}}$, the function $w\mapsto g(w)/(w-u)$ is holomorphic inside and on $\Gamma_{\lambda}(1/4)$, and hence by Cauchy's theorem,
        \begin{equation}\label{innerCauchyZero}
            \oint_{\Gamma_{\lambda}\left(1/4\right)}g(w)\frac{1}{w-u}\mathrm{d}w=0,
        \end{equation}
        and consequently, $I_2(x,m)$ simplifies to
        \begin{equation}\label{ItwoSimplified}
            I_2(x,m)=-\frac{1}{\lambda_m2\pi\mathrm{i}}\oint_{\Gamma_{\lambda}\left(1/2\right)}f(u)g(u)\frac{\rho_x(u)}{\rho_{m+1}(u)}\mathrm{d}u.
        \end{equation}
        The expression in \eqref{ItwoSimplified} is exactly the negative of the residue contribution in \eqref{IoneDeform}, and therefore
        \begin{equation}\label{IoneItwoCancel}
                I_1(x,m)=-I_2(x,m)+\frac{1}{\lambda_m(2\pi\mathrm{i})^2}\oint_{\Gamma_{\lambda}\left(1/2\right)}\oint_{\Gamma_{\lambda}\left(1/4\right)}f(u)g(w)\frac{\rho_x(w)}{\rho_{m+1}(w)}\frac{1}{w-u}\mathrm{d}w\mathrm{d}u.
        \end{equation}
        Therefore, recombining the two contributions, we obtain
        \begin{equation}\label{compositionReduced}
            (\mathrm{T}^{f}\mathrm{T}^{g})(x,m)=I_1(x,m)+I_2(x,m)=\frac{1}{\lambda_m(2\pi\mathrm{i})^2}
            \oint_{\Gamma_{\lambda}\left(1/2\right)}\oint_{\Gamma_{\lambda}\left(1/4\right)}
            f(u)g(w)\frac{\rho_x(w)}{\rho_{m+1}(w)}\frac{1}{w-u}\,\mathrm{d}w\mathrm{d}u.
        \end{equation}
        Since both contours are compact and the integrand is continuous on their product, we can apply Fubini's theorem to change the order of integration. For fixed $w\in\Gamma_{\lambda}\left(1/4\right)$, the point $w$ lies inside $\Gamma_{\lambda}\left(1/2\right)$, and applying Cauchy's formula on the resulting inner integral gives
        \begin{equation}\label{CauchyInnerComposition}
            \oint_{\Gamma_{\lambda}\left(1/2\right)}\frac{f(u)}{w-u}\mathrm{d}u=-2\pi\mathrm{i}f(w).
        \end{equation}
        Finally, substituting \eqref{CauchyInnerComposition} back into the reduced double contour integral \eqref{compositionReduced} and simplifying gives
        \begin{equation}\label{multiplicativeT}
            (\mathrm{T}^{f}\mathrm{T}^{g})(x,m)
            =-\frac{1}{\lambda_m2\pi\mathrm{i}}
            \oint_{\Gamma_{\lambda}\left(1/4\right)}\frac{\rho_x(w)}{\rho_{m+1}(w)}f(w)g(w)\,\mathrm{d}w
            =\mathrm{T}^{fg}(x,m),
        \end{equation}
        where the last equality uses the contour-independence already established by Cauchy's theorem.
        Therefore, the kernels $\mathrm{T}^f$ satisfy the multiplicative composition rule \eqref{multiplicativeT}. We now identify the finite-time continuous kernel and the one-step discrete kernels from Definition~\ref{contprocess} and Definition~\ref{mixprocess} in terms of our kernel $\mathrm{T}^{f(\cdot,t)}$ for suitable classes of functions $\{f(\cdot,t)\}_{t\in\mathscr{T}}\subset\mathrm{Hol}(\mathbb{H}_{-1/\mathfrak{U}})$.
       
        \textbf{Continuous-time kernel.}
        For Definition~\ref{contprocess}, for each $t\in\mathbb{R}_+$, $w\in\mathbb{C}$ and $x,m\in\mathbb{Z}_+$, we set
        \begin{equation}\label{continuousKernelDefinition}
            \psi^{\textnormal{C}}(w,t)\overset{\textnormal{def}}{=}\mathrm{e}^{-tw}\quad\text{and then}\quad \mathrm{T}^{\psi^{\textnormal{C}}(\cdot,t)}(x,m)=-\frac{1}{\lambda_{m}2\pi\mathrm{i}}\oint_{\Gamma_{\lambda}}\frac{\rho_x(w)}{\rho_{m+1}(w)}\mathrm{e}^{-tw}\mathrm{d}w.
        \end{equation}
        Since $\psi^{\textnormal{C}}(w,0)=1$, the contour definition gives the correct initial condition
        \begin{equation}\label{continuousinitialT}
            \mathrm{T}^{\psi^{\textnormal{C}}(\cdot,0)}(x,m)=-\frac{1}{\lambda_{m}2\pi\mathrm{i}}\oint_{\Gamma_{\lambda}}\frac{\rho_x(w)}{\rho_{m+1}(w)}\mathrm{d}w=\mathbf{1}_{\{m=x\}}.
        \end{equation}
        Indeed, when $x>m$ the integrand in $\mathrm{T}^{1}(x,m)$ is holomorphic inside $\Gamma_{\lambda}$ and the contour integral is zero, while when $x=m$ the residue at $w=\lambda_x$, together with the prefactor in \eqref{defofTf}, gives one. It remains to justify the case $x<m$. Let $\mathcal{C}_R$ be the positively oriented circle $|w|=R$, where $R$ is chosen so that $\Gamma_\lambda$ lies inside $\mathcal{C}_R$ and $R>2\max(\lambda_x,\ldots,\lambda_m,1)$. Since the integrand has no poles in the annular region between $\Gamma_\lambda$ and $\mathcal{C}_R$, Cauchy's theorem allows us to replace $\Gamma_\lambda$ by $\mathcal{C}_R$. Moreover, by \textbf{(A1)}, for all $w\in\mathcal{C}_R$,
        \begin{equation}\label{largecircleTone}
            \left|\frac{\rho_x(w)}{\rho_{m+1}(w)}\right|=\prod_{k=x}^{m}\left|1-\frac{w}{\lambda_k}\right|^{-1}\leq\left(\frac{2\mathfrak{B}_1}{R}\right)^{m-x+1}.
        \end{equation}
        Consequently,
        \begin{equation}\label{initialOffDiagonalBound}
            \left|\mathrm{T}^{\psi^{\textnormal{C}}(\cdot,0)}(x,m)\right|\leq\frac{1}{\lambda_m}R\left(\frac{2\mathfrak{B}_1}{R}\right)^{m-x+1}\leq\frac{(2\mathfrak{B}_1)^{m-x+1}}{\mathfrak{L}R^{m-x}}.
        \end{equation}
	        Since the value of the contour integral is independent of the auxiliary radius $R$, while the right-hand side of \eqref{initialOffDiagonalBound} can be made arbitrarily small by increasing $R$, the original contour integral is zero for all $x<m$. Next, for each bounded time interval, the compactness of $\Gamma_{\lambda}$ and the positive distance between $\Gamma_{\lambda}$ and the spatial poles give a uniform bound on the differentiated integrand in \eqref{continuousKernelDefinition}. Hence we may differentiate under the integral sign to obtain
        \begin{equation}\label{continuousderivativeT}
            \frac{\mathrm{d}}{\mathrm{d}t}\mathrm{T}^{\psi^{\textnormal{C}}(\cdot,t)}(x,m)=-\mathrm{T}^{(\cdot)\psi^{\textnormal{C}}(\cdot,t)}(x,m)=-\left(\mathrm{T}^{(\cdot)}\mathrm{T}^{\psi^{\textnormal{C}}(\cdot,t)}\right)(x,m),
        \end{equation}
        where the last equality uses the multiplicative composition rule \eqref{multiplicativeT} and where $(\cdot)$ denotes the identity function $w\mapsto w$. We now compute $\mathrm{T}^{(\cdot)}$ explicitly:
        \begin{equation}\label{Tidentitygenerator}
            \mathrm{T}^{(\cdot)}(x,m)=-\frac{1}{\lambda_{m}2\pi\mathrm{i}}\oint_{\Gamma_{\lambda}}\frac{\rho_x(w)}{\rho_{m+1}(w)}w\mathrm{d}w=\lambda_x\left(\mathbf{1}_{\{m=x\}}-\mathbf{1}_{\{m=x+1\}}\right).
        \end{equation}
        For $x>m$, the integrand is holomorphic inside $\Gamma_{\lambda}$. For $m>x+1$, the same large-circle argument applies: for $w\in\mathcal{C}_R$,
        \begin{equation}\label{largecircleTidentity}
            \left|\frac{\rho_x(w)}{\rho_{m+1}(w)}w\right|\leq R\left(\frac{2\mathfrak{B}_1}{R}\right)^{m-x+1},
        \end{equation}
        and hence
	        \begin{equation}\label{TidentityOffDiagonalBound}
	            \left|\mathrm{T}^{(\cdot)}(x,m)\right|\leq\frac{1}{\lambda_m}R^2\left(\frac{2\mathfrak{B}_1}{R}\right)^{m-x+1}\leq\frac{(2\mathfrak{B}_1)^{m-x+1}}{\mathfrak{L}R^{m-x-1}}.
	        \end{equation}
	        Since $m>x+1$, the exponent $m-x-1$ is positive. The value of the contour integral is independent of the auxiliary radius $R$, whereas the explicit upper bound in \eqref{TidentityOffDiagonalBound} can be made arbitrarily small by increasing $R$. Hence $\mathrm{T}^{(\cdot)}(x,m)=0$ in this case.
        The remaining cases are $m=x$ and $m=x+1$. For $m=x$, the integrand has a single simple pole and
        \begin{equation}\label{TidentityDiagonal}
            \mathrm{T}^{(\cdot)}(x,x)=\frac{1}{2\pi\mathrm{i}}\oint_{\Gamma_{\lambda}}\frac{w}{w-\lambda_x}\mathrm{d}w=\lambda_x.
        \end{equation}
        For $m=x+1$, if $\lambda_x\neq\lambda_{x+1}$, then the residue calculation gives
        \begin{equation}\label{TidentityAdjacentDistinct}
            \mathrm{T}^{(\cdot)}(x,x+1)=-\frac{\lambda_{x}}{2\pi\mathrm{i}}\oint_{\Gamma_{\lambda}}\frac{w}{(w-\lambda_x)(w-\lambda_{x+1})}\mathrm{d}w=-\frac{\lambda_{x+1}\lambda_{x}}{\lambda_{x+1}-\lambda_x}+\frac{\lambda_{x}^2}{\lambda_{x+1}-\lambda_{x}}=-\lambda_x.
        \end{equation}
        If $\lambda_x=\lambda_{x+1}$, then the integrand has a double pole at $w=\lambda_x$, and
        \begin{equation}\label{TidentityAdjacentRepeated}
            \mathrm{T}^{(\cdot)}(x,x+1)=-\frac{\lambda_{x}}{2\pi\mathrm{i}}\oint_{\Gamma_{\lambda}}\frac{w}{(w-\lambda_x)^2}\mathrm{d}w=-\lambda_x.
        \end{equation}
        Thus $\mathrm{T}^{(\cdot)}(x,m)=-\mathsf{L}(x,m)$, where $\mathsf{L}$ is the generator from Definition~\ref{contprocess}, and consequently
        \begin{equation}\label{continuousbackwardT}
            \frac{\mathrm{d}}{\mathrm{d}t}\mathrm{T}^{\psi^{\textnormal{C}}(\cdot,t)}(x,m)=\left(\mathsf{L}\mathrm{T}^{\psi^{\textnormal{C}}(\cdot,t)}\right)(x,m),\qquad \mathrm{T}^{\psi^{\textnormal{C}}(\cdot,0)}(x,m)=\mathbf{1}_{\{m=x\}}.
        \end{equation}
     For fixed $m$ and $S>0$, the identity $\mathrm{T}^{\psi^{\textnormal{C}}(\cdot,t)}(x,m)=0$ for $x>m$ gives
        \begin{equation}\label{continuousTbounded}
            \sup_{0\leq t\leq S}\sup_{x\in\mathbb{Z}_+}\left|\mathrm{T}^{\psi^{\textnormal{C}}(\cdot,t)}(x,m)\right|=\max_{0\leq x\leq m}\sup_{0\leq t\leq S}\left|\mathrm{T}^{\psi^{\textnormal{C}}(\cdot,t)}(x,m)\right|<\infty,
        \end{equation}
        where the final bound follows from the continuity of the contour integral in $t$ on compact time intervals. Therefore, by uniqueness of bounded solutions to the backward equation \eqref{continuousbackwardT} we obtain $\mathrm{T}^{\psi^{\textnormal{C}}(\cdot,t)}(x,m)=\mathbb{P}_{x}[\mathscr{X}(t)=m]$. Moreover, since $\psi^{\textnormal{C}}(\cdot,t)$ is entire for each $t$, the kernel $\mathrm{T}^{\psi^{\textnormal{C}}(\cdot,t)}$ is well-defined for all choices of contour $\Gamma_{\lambda}$ enclosing the spatial poles, and in particular for the nested admissible contours used above.
    
        \textbf{Discrete one-step kernels.}
	        In the case of Definition~\ref{mixprocess}, for each $t\in\mathbb{Z}_+$, we set
        \begin{equation}\label{discreteKernelDefinition}
            \psi^{\textnormal{D}}(w,t)\overset{\textnormal{def}}{=}\frac{(1-\alpha_t w)^{1-\tau_t}}{(1+\beta_t w)^{\tau_t}}\quad\text{and then}\quad\mathrm{T}^{\psi^{\textnormal{D}}(\cdot,t)}(x,m)=-\frac{1}{\lambda_{m}2\pi\mathrm{i}}\oint_{\Gamma_{\lambda}}\frac{\rho_x(w)}{\rho_{m+1}(w)}\frac{(1-\alpha_t w)^{1-\tau_t}}{(1+\beta_t w)^{\tau_t}}\mathrm{d}w.
        \end{equation}
	        The continuous-time function $\psi^{\textnormal{C}}(\cdot,t)$ is entire. The only poles that can arise from the discrete-time functions $\{\psi^{\textnormal{D}}(\cdot,t)\}_{t\in\mathbb{Z}_+}$ occur at $w=-1/\beta_t$ when $\tau_t=1$, and by \textbf{(A1)}, $-1/\beta_t\le -1/\mathfrak{U}$. Hence $\psi^{\textnormal{C}}(\cdot,t),\psi^{\textnormal{D}}(\cdot,t)\in\mathrm{Hol}(\mathbb{H}_{-1/\mathfrak{U}})$ for all relevant $t$, and the corresponding kernels are well-defined. We also use a parameterised version of the discrete kernel. For any $x,m\in\mathbb{Z}_+$ and any $\underline{\lambda}=\left(\lambda_0,\lambda_1,\ldots\right)\in[\mathfrak{L},\mathfrak{B}_1]^{\mathbb{Z}_+}$, denote the kernel $\mathrm{T}^{\psi^{\textnormal{D}}(\cdot,t)}$ with spatial parameters $\underline{\lambda}$ by
        \begin{equation}\label{parameterisedDiscreteKernel}
            \mathrm{T}^{\psi^{\textnormal{D}}(\cdot,t)}(x,m:\underline{\lambda})=-\frac{\prod_{z=x}^{m-1}\lambda_{z}}{2\pi\mathrm{i}}\oint_{\Gamma_{\lambda}}\frac{(1-\alpha_t w)^{1-\tau_t}(1+\beta_t w)^{-\tau_t}}{\prod_{z=x}^{m}(\lambda_z-w)}\mathrm{d}w,\quad x\le m.
        \end{equation}
        Assume first that the finitely many rates $\lambda_x,\ldots,\lambda_m$ are pairwise distinct. By Cauchy's residue theorem, the only contributions to the above contour integral arise from the simple poles at these spatial parameters, and therefore
        \begin{equation}\label{discreteresidueformula}
            \mathrm{T}^{\psi^{\textnormal{D}}(\cdot,t)}(x,m:\underline{\lambda}) =\prod_{z=x}^{m-1}\lambda_{z}\left[\sum_{r=x}^{m}\frac{(1-\alpha_t \lambda_{r})^{1-\tau_t}(1+\beta_t \lambda_{r})^{-\tau_t}}{\prod_{\substack{k=x\\k\ne r}}^{m}\left(\lambda_k-\lambda_{r}\right)}\right].
        \end{equation}
        This can be written in terms of the divided difference of the function $\psi^{\textnormal{D}}(\cdot,t)$, in the standard sense of \cite{deBoor2005},
        \begin{equation}\label{dividdiff}
            \mathrm{T}^{\psi^{\textnormal{D}}(\cdot,t)}(x,m:\underline{\lambda}) =\left(-1\right)^{m-x}\left(\prod_{z=x}^{m-1}\lambda_{z}\right)\psi^{\textnormal{D}}[\lambda_x,\ldots,\lambda_m](t),
        \end{equation}
        where, for distinct nodes, the divided difference $\psi^{\textnormal{D}}[\lambda_x,\ldots,\lambda_m](t)$ is defined recursively by the standard divided-difference recursion \cite{deBoor2005}
        \begin{equation}\label{divdiffrecursion}
            \begin{aligned}
                \psi^{\textnormal{D}}[\lambda_x](t)&=\psi^{\textnormal{D}}(\lambda_x,t),\\
                \psi^{\textnormal{D}}[\lambda_x,\ldots,\lambda_m](t)&=\frac{\psi^{\textnormal{D}}[\lambda_{x+1},\ldots,\lambda_m](t)-\psi^{\textnormal{D}}[\lambda_x,\ldots,\lambda_{m-1}](t)}{\lambda_m-\lambda_x},\quad m>x.
            \end{aligned}
        \end{equation}
        The restriction to distinct rates is only used to make the residue computation \eqref{discreteresidueformula} and the recursive divided-difference formula \eqref{divdiffrecursion} literal. Both sides of \eqref{dividdiff}, after interpreting divided differences by their continuous extension in the nodes, are continuous functions of the finite vector $(\lambda_x,\ldots,\lambda_m)\in[\mathfrak{L},\mathfrak{B}_1]^{m-x+1}$: the contour may be chosen from the admissible family \eqref{admissiblecontours} using only the constants in \textbf{(A1)}, so it remains a positive distance from the compact spatial interval, and divided differences have a standard continuous extension to repeated nodes; see \cite{deBoor2005}. Since finite vectors with pairwise distinct coordinates are dense in $[\mathfrak{L},\mathfrak{B}_1]^{m-x+1}$, \eqref{dividdiff} extends to arbitrary spatial rates satisfying \textbf{(A1)}.
        For $\tau_t=0$, and first for distinct adjacent nodes, we can explicitly compute the divided differences for all $x\in\mathbb{Z}_+$ as
        \begin{equation}\label{divdiff0}
            \begin{aligned}
                \psi^{\textnormal{D}}[\lambda_x](t)&=(1-\alpha_t \lambda_x),\\
                \psi^{\textnormal{D}}[\lambda_x,\lambda_{x+1}](t)&=\frac{(1-\alpha_t \lambda_{x+1})-(1-\alpha_t \lambda_x)}{\lambda_{x+1}-\lambda_x}=-\alpha_t.
            \end{aligned}
        \end{equation}
        The first-order divided difference in \eqref{divdiff0} is independent of the nodes, and hence all divided differences of order at least two vanish. For example, for distinct $\lambda_x$ and $\lambda_{x+2}$,
        \begin{equation*}
            \begin{aligned}
                \psi^{\textnormal{D}}[\lambda_x,\lambda_{x+1},\lambda_{x+2}](t)&=\frac{-\alpha_t+\alpha_t}{\lambda_{x+2}-\lambda_x}=0,\\
                \psi^{\textnormal{D}}[\lambda_x,\ldots,\lambda_{m}](t)&=0,\quad\text{for all $m>x+1$}.
            \end{aligned}
        \end{equation*}
        By the continuous extension just described, the same values hold when some of the spatial rates coincide. Substituting the non-zero values \eqref{divdiff0} into the divided-difference identity \eqref{dividdiff}, we obtain the Bernoulli one-step kernel
        \begin{equation}\label{BernoulliTkernel}
            \mathrm{T}^{\psi^{\textnormal{D}}(\cdot,t)}(x,m:\underline{\lambda}) =\left\{
            \begin{aligned}
                &1-\alpha_t \lambda_x,\quad m=x,\\
                &\alpha_t\lambda_x,\quad m=x+1,\\
                &0,\quad m<x\text{ or }m>x+1,
            \end{aligned}
            \right.
        \end{equation}
        For $\tau_t=1$, and first for distinct nodes, we inductively compute the divided differences for all $x\in\mathbb{Z}_+$ as
        \begin{equation}\label{geometricDivDiff}
            \psi^{\textnormal{D}}[\lambda_x,\ldots,\lambda_{m}](t)=-\frac{1}{\beta_t}\prod_{z=x}^{m}\left(\frac{-\beta_{t}}{1+\beta_t \lambda_z}\right),\quad m\ge x.
        \end{equation}
	        For the initial case $m=x$, this gives
        \begin{equation*}
            \psi^{\textnormal{D}}[\lambda_x](t)=\frac{1}{1+\beta_t \lambda_x},
        \end{equation*}
        then assuming the inductive hypothesis holds for divided differences of order $n-1$ with $n\geq1$, we have
        \begin{equation*}
            \begin{aligned}
                \psi^{\textnormal{D}}[\lambda_x,\ldots,\lambda_{x+n}](t)
                &=\frac{1}{\lambda_{x+n}-\lambda_x}\left[-\frac{1}{\beta_t}\prod_{z=x+1}^{x+n}\left(\frac{-\beta_{t}}{1+\beta_t \lambda_z}\right)+\frac{1}{\beta_t}\prod_{z=x}^{x+n-1}\left(\frac{-\beta_{t}}{1+\beta_t \lambda_z}\right)\right]\\
                &=-\frac{1}{\beta_t}\prod_{z=x}^{x+n}\left(\frac{-\beta_{t}}{1+\beta_t \lambda_z}\right).
            \end{aligned}
        \end{equation*}
        Again, the identity \eqref{geometricDivDiff} extends to repeated spatial rates by continuity of divided differences in their nodes. Similarly, substituting these values into the divided-difference identity \eqref{dividdiff}, we obtain the geometric one-step kernel
        \begin{equation}\label{geometricTkernel}
            \mathrm{T}^{\psi^{\textnormal{D}}(\cdot,t)}(x,m:\underline{\lambda}) =\left\{
            \begin{aligned}
                &\frac{1}{1+\beta_t\lambda_m}\prod_{z=x}^{m-1}\left(\frac{\beta_{t}\lambda_z}{1+\beta_t \lambda_z}\right),\quad m\geq x,\\
                &0,\quad m<x.
            \end{aligned}
            \right.
        \end{equation}
        The Bernoulli kernel \eqref{BernoulliTkernel} and the geometric kernel \eqref{geometricTkernel} are exactly the one-step transition probabilities in Definition~\ref{mixprocess} for the respective choices $\tau_t=0$ and $\tau_t=1$. Thus the one-step transition probabilities coincide with the kernels $\mathrm{T}^{\psi^{\textnormal{D}}(\cdot,t)}$ for the corresponding choices of $\psi^{\textnormal{D}}(\cdot,t)$.
        
        \textbf{Transition kernels.}
	        Applying the multiplicative composition rule \eqref{multiplicativeT} to the continuous time factors gives that, in Definition~\ref{contprocess}, for $s,t\in\mathbb{R}_+$ with $s\leq t$ and any strictly increasing sequence $s<t_1<\cdots<t_n<t$, we have
        \begin{equation}\label{continuouskernelcomposition}
            \mathbb{P}[\mathscr{X}(t)=y(t)|\mathscr{X}(s)=x]
            =\left(\mathrm{T}^{\psi^{\textnormal{C}}(\cdot,t_1-s)}\mathrm{T}^{\psi^{\textnormal{C}}(\cdot,t_2-t_1)}\cdots \mathrm{T}^{\psi^{\textnormal{C}}(\cdot,t-t_n)}\right)(x,y(t))
            =\mathrm{T}^{\psi^{\textnormal{C}}(\cdot,t-s)}(x,y(t)).
        \end{equation}
        Here the last equality in \eqref{continuouskernelcomposition} follows from the pointwise identity $\psi^{\textnormal{C}}(w,t)\psi^{\textnormal{C}}(w,s)=\psi^{\textnormal{C}}(w,t+s)$. In the case of Definition~\ref{mixprocess}, for $s,t\in\mathbb{Z}_+$ with $s\le t$, and with the empty product interpreted as $1$, the same multiplicative composition rule \eqref{multiplicativeT} gives
        \begin{equation}\label{discretekernelcomposition}
            \mathbb{P}[\mathscr{X}(t)=y(t)|\mathscr{X}(s)=x]
            =\left(\mathrm{T}^{\psi^{\textnormal{D}}(\cdot,s)}\mathrm{T}^{\psi^{\textnormal{D}}(\cdot,s+1)}\cdots \mathrm{T}^{\psi^{\textnormal{D}}(\cdot,t-1)}\right)(x,y(t))
            =\mathrm{T}^{\prod_{n=s}^{t-1}\psi^{\textnormal{D}}(\cdot,n)}(x,y(t)).
        \end{equation}
	        Finally, setting the endpoint in \eqref{defofTf} equal to $y(t)$, and collecting all components of the integrand that depend on $t$ and the endpoint $y(t)$ into an exponent, recovers
	        \begin{equation}\label{phaseLogExtraction}
	            \mathfrak{F}_y(w,t)=\log\left(\psi^{\textnormal{C}}(w,t)\right)\mathbf{1}_{\{\mathscr{T}=\mathbb{R}_+\}}+\sum_{n=0}^{t-1}\log\left(\psi^{\textnormal{D}}(w,n)\right)\mathbf{1}_{\{\mathscr{T}=\mathbb{Z}_+\}}-\sum_{k=0}^{y(t)}\log\left(1-\frac{w}{\lambda_k}\right),
        \end{equation}
        coinciding with Definition~\ref{form} when we substitute the appropriate definitions. Equivalently, in the exact contour identity \eqref{cir}, the exponential is the single-valued product
        \begin{equation}\label{phaseproductform}
            \mathrm{e}^{\mathfrak{F}_y(w,t)}=\left(\mathrm{e}^{-tw}\prod_{k=0}^{y(t)}\left(1-\frac{w}{\lambda_k}\right)^{-1}\right)\mathbf{1}_{\{\mathscr{T}=\mathbb{R}_+\}}+\left(\prod_{n=0}^{t-1}\frac{(1-\alpha_n w)^{1-\tau_n}}{(1+\beta_n w)^{\tau_n}}\prod_{k=0}^{y(t)}\left(1-\frac{w}{\lambda_k}\right)^{-1}\right)\mathbf{1}_{\{\mathscr{T}=\mathbb{Z}_+\}},
        \end{equation}
    completing the proof.
    \end{proof}

    \subsection{Analysis of the phase function}
    \noindent We next analyse the deterministic scale $\gamma(t)$ from Definition~\ref{gammadef}, which measures the residual inhomogeneity remaining after centering at the deterministic point $\mathscr{Z}(t)$. The quantity $\gamma(t)$ in some sense captures two distinct sources of scale: deterministic inhomogeneity and fluctuations. 
    
    \begin{definition}[Endpoint window]\label{def:endpoint-window}
        For $T\in\mathscr{T}$,  denote the finite-time endpoint window by
        \begin{equation*}
            \mathscr{C}_{T}
            \overset{\textnormal{def}}{=}
            \left\{m\in\mathbb{Z}_+:\, |m-\mathscr{Z}(T)|\leq \gamma(T)^{1/2}\Sigma(T)^2\right\}.
        \end{equation*}
    \end{definition}
            Since $\gamma(T)^{1/2}$ is larger than $\gamma(T)$ once $\gamma(T)$ is small, the window $\mathscr{C}_{T}$ is deliberately wider than the intrinsic localization scale $\gamma(T)\Sigma(T)^2$ studied later in Section~\ref{MOM}. Its role is to contain the relevant fluctuations with room to spare, so that the analytic estimates can be applied uniformly over all terminal endpoints that arise in the summations.

     We define the open ball centered at $c\in\mathbb{C}$, of radius $r>0$, by
    \begin{equation*}
        \mathcal{B}_{r}(c)\overset{\textnormal{def}}{=}\{z\in\mathbb{C}: |z-c|<r\}.
    \end{equation*}
    \noindent In the following  technical proposition we isolate the deterministic analytic estimates for $\mathfrak{F}_y$, with endpoint-uniform bounds on the window $\mathscr{C}_{t}$. These estimates will be used repeatedly to control stationary-point perturbations and contour expansions.
    \begin{proposition}[Properties of the phase function]\label{specificcasesub}
        Assuming \textbf{(A1)} and \textbf{(A2)}, the phase function $\mathfrak{F}_y$ from Definition~\ref{form} satisfies the following properties:
        \begin{enumerate}[label=(\roman*)]
            \item\label{phaseprop:i} For constants $\mathfrak{L}$ and $\mathfrak{U}$ as defined in \textbf{(A1)}, let $\mathrm{R}_{-}\overset{\textnormal{def}}{=}\min(\mathfrak{L},\mathfrak{U}^{-1})$. Then $\mathfrak{F}_y(\cdot,t)$ is holomorphic in $\mathcal{B}_{\mathrm{R}_{-}}(0)$, for all $y\in\mathscr{P}$ and all $t\in\mathscr{T}$.
            \item\label{phaseprop:ii} There exists $\mathrm{T}_1\in\mathscr{T}$ such that, for all $y\in\mathscr{P}$ and all $t>\mathrm{T}_1$,
			        \begin{equation*}
			            \begin{aligned}
			                \mathrm{K}_{\mathrm{l}}(t)\left(y(t)-\mathscr{Z}(t)\right)-2\gamma(t)\mathscr{Z}(t)-\frac{3\Lambda}{2}
			                &\leq\mathfrak{F}_{y}^{(1)}(0,t),\\
			                \mathfrak{F}_{y}^{(1)}(0,t)
			                &\leq
			                \mathrm{K}_{\mathrm{u}}(t)\left(y(t)-\mathscr{Z}(t)\right)+2\gamma(t)\mathscr{Z}(t)+\frac{3\Lambda}{2}.
			            \end{aligned}
		                \end{equation*}
                where for each $t\in\mathscr{T}$, $\mathrm{K}_{\mathrm{l}}(t)=\mathfrak{B}_1^{-1}$ and $\mathrm{K}_{\mathrm{u}}(t)=\mathfrak{L}^{-1}$ whenever $y(t)\geq\mathscr{Z}(t)$ and conversely $\mathrm{K}_{\mathrm{l}}(t)=\mathfrak{L}^{-1}$ and $\mathrm{K}_{\mathrm{u}}(t)=\mathfrak{B}_1^{-1}$ whenever $y(t)<\mathscr{Z}(t)$. Moreover, with the constant $\mathrm{C}_1$ appearing in \textbf{(iv)}, there exists $\mathrm{T}_{\mathscr{C},1}\in\mathscr{T}$ such that, for every $y\in\mathscr{P}$, whenever $t>\mathrm{T}_{\mathscr{C},1}$ and $y(t)\in\mathscr{C}_{t}$,
                \begin{equation*}
                    \left|\mathfrak{F}_{y}^{(1)}(0,t)\right|\leq \left(\mathfrak{L}^{-1}\mathrm{C}_{1}^{-1}+2\right)\gamma(t)^{1/2}\mathscr{Z}(t)+\frac{3\Lambda}{2}.
                \end{equation*}
            \item\label{phaseprop:iii} Let $\mathrm{R}_2\overset{\textnormal{def}}{=}\mathrm{R}_{-}/2$. There exist constants $\mathrm{M}_1,\mathrm{M}_2>0$ and a threshold $\mathrm{T}_{2,\mathscr{C}}\in\mathscr{T}$ such that
	                \begin{equation*}
	                    \left|\mathfrak{F}_{y}^{(n)}(w,t)\right|\leq (n-1)!\mathrm{M}_1\mathrm{M}_2^{n}\Sigma(t)^2,
	                \end{equation*}
                    for every $t>\mathrm{T}_{2,\mathscr{C}}$, every $n\geq1$, every $w\in\mathcal{B}_{\mathrm{R}_2}(0)$, and every path $y\in\mathscr{P}$ with $y(t)\in\mathscr{C}_{t}$.
            \item\label{phaseprop:iv} There exist constants $\mathrm{C}_1,\mathrm{C}_2>0$ and $\mathrm{T}_3\in\mathscr{T}$ such that
                \begin{equation*}
                    \mathrm{C}_1\Sigma(t)^2\leq\mathscr{Z}(t)\leq \mathrm{C}_2\Sigma(t)^2,\quad\text{for all $t>\mathrm{T}_3$.}
                \end{equation*}
                Furthermore, there exist constants $\mathrm{C}_{\mathscr{C}}>0$ and $\mathrm{T}_{4,\mathscr{C}}\in\mathscr{T}$ such that, for every $y\in\mathscr{P}$, whenever $t>\mathrm{T}_{4,\mathscr{C}}$ and $y(t)\in\mathscr{C}_{t}$,
                \begin{equation*}
                    \left|\mathfrak{F}_{y}^{(2)}(w,t)-\Sigma(t)^2\right|
                    \leq \mathrm{C}_{\mathscr{C}}\left(\gamma(t)^{1/2}+|w|\right)\Sigma(t)^2,\quad\text{for all $w\in\mathcal{B}_{\mathrm{R}_2}(0)$.}
                \end{equation*}
        \end{enumerate}
        \begin{proof}
            We split the proof into five steps. The first step identifies a common holomorphic neighbourhood of the origin. The second step proves the first-derivative estimate by separating the endpoint displacement $y(t)-\mathscr{Z}(t)$ from the deterministic centring error at $\mathscr{Z}(t)$. The third step compares $\Sigma(t)^2$ with $\mathscr{Z}(t)$; this is the scale comparison used in the later derivative bounds. The fourth step proves the endpoint-uniform higher-derivative estimate. The final step proves the strengthened estimates on the finite-time window $\mathscr{C}_{t}$, where the endpoint error is controlled at the wider scale $\gamma(t)^{1/2}\Sigma(t)^2$.
         
            \textbf{Step 1: Holomorphicity.} Recall from Definition~\ref{form} that every singularity of $\mathfrak{F}_y(\cdot,t)$ comes from one of the logarithmic factors. The spatial logarithms have singularities at $w=\lambda_n$. By \textbf{(A1)}, these lie to the right of $\mathfrak{L}$. In the discrete-time regime, the Bernoulli logarithms have singularities at $w=\alpha_n^{-1}$, which lie to the right of $\mathfrak{B}_2$, and the geometric logarithms have singularities at $w=-\beta_n^{-1}$, whose moduli are at least $\mathfrak{U}^{-1}$. The resulting separation from the origin is illustrated schematically in Figure~\ref{phaseholomorphicdiskfig}. Hence every singularity lies outside $\mathcal{B}_{\mathrm{R}_{-}}(0)$, where
            \begin{equation*}
                \mathrm{R}_{-}\overset{\textnormal{def}}{=}\min(\mathfrak{L},\mathfrak{U}^{-1}).
            \end{equation*}
            \begin{figure}[ht!]
                \centering
                \begin{tikzpicture}[x=1cm,y=1cm,line cap=round,line join=round]
                    \definecolor{poleRed}{RGB}{190,30,30}
                    \definecolor{poleBlue}{RGB}{25,55,210}
                    \definecolor{diskFill}{RGB}{230,245,255}
                    \definecolor{diskEdge}{RGB}{25,110,180}
                    \definecolor{axisGray}{RGB}{70,70,70}

                    \draw[->,axisGray] (-3.8,0) -- (7.1,0) node[below right] {\small $\mathfrak{Re}\,w$};
                    \draw[->,axisGray] (0,-2.15) -- (0,2.15) node[above left] {\small $\mathfrak{Im}\,w$};

                    \draw[axisGray] (-0.70,-0.13) -- (-0.70,0.13);
                    \node[axisGray,below] at (-0.70,-0.18) {\small $-\mathfrak{U}^{-1}$};

                    \draw[axisGray] (1,-0.13) -- (1,0.13);
                    \node[axisGray,below] at (1,-0.18) {\small $\mathfrak{L}$};

                    \draw[axisGray] (3.94,-0.13) -- (3.94,0.13);
                    \node[axisGray,below] at (3.94,-0.18) {\small $\mathfrak{B}_1$};

                    \draw[axisGray] (4.70,-0.13) -- (4.70,0.13);
                    \node[axisGray,below] at (4.70,-0.18) {\small $\mathfrak{B}_2$};

                    \draw[poleRed,thick] (-2.95,-0.07) -- (-2.81,0.07);
                    \draw[poleRed,thick] (-2.95,0.07) -- (-2.81,-0.07);
                    \draw[poleRed,thick] (-2.55,-0.07) -- (-2.41,0.07);
                    \draw[poleRed,thick] (-2.55,0.07) -- (-2.41,-0.07);
                    \draw[poleRed,thick] (-2.15,-0.07) -- (-2.01,0.07);
                    \draw[poleRed,thick] (-2.15,0.07) -- (-2.01,-0.07);
                    \draw[poleRed,thin,decorate,decoration={brace,amplitude=3pt,raise=2pt}] (-2.99,0.12) -- (-1.97,0.12);
                    \node[poleRed,anchor=south] at (-2.48,0.42) {\small $-\beta_n^{-1}$};

                    \draw[black,thick] (1.45,-0.07) -- (1.59,0.07);
                    \draw[black,thick] (1.45,0.07) -- (1.59,-0.07);
                    \draw[black,thick] (1.95,-0.07) -- (2.09,0.07);
                    \draw[black,thick] (1.95,0.07) -- (2.09,-0.07);
                    \draw[black,thick] (2.45,-0.07) -- (2.59,0.07);
                    \draw[black,thick] (2.45,0.07) -- (2.59,-0.07);
                    \draw[black,thick] (2.95,-0.07) -- (3.09,0.07);
                    \draw[black,thick] (2.95,0.07) -- (3.09,-0.07);
                    \draw[axisGray,thin,decorate,decoration={brace,amplitude=3pt,raise=2pt}] (1.42,0.12) -- (3.12,0.12);
                    \node[axisGray,anchor=south] at (2.27,0.42) {\small $\lambda_n$};

                    \draw[green!70!black,thick] (5.15,-0.07) -- (5.29,0.07);
                    \draw[green!70!black,thick] (5.15,0.07) -- (5.29,-0.07);
                    \draw[green!70!black,thick] (5.65,-0.07) -- (5.79,0.07);
                    \draw[green!70!black,thick] (5.65,0.07) -- (5.79,-0.07);
                    \draw[green!70!black,thick] (6.15,-0.07) -- (6.29,0.07);
                    \draw[green!70!black,thick] (6.15,0.07) -- (6.29,-0.07);
                    \draw[green!70!black,thick] (6.65,-0.07) -- (6.79,0.07);
                    \draw[green!70!black,thick] (6.65,0.07) -- (6.79,-0.07);
                    \draw[green!70!black,thin,decorate,decoration={brace,amplitude=3pt,raise=2pt}] (5.12,0.12) -- (6.82,0.12);
                    \node[green!70!black,anchor=south] at (5.97,0.42) {\small $\alpha_n^{-1}$};
                    
                    \fill[diskFill,opacity=0.4] (0,0) circle (0.70); 
                    \draw[thick,diskEdge] (0,0) circle (0.70);
                    \node[below=2pt] at (0,0) {\small $0$};
                    \draw[diskEdge,thin] (0.50,0.50) -- (0.98,0.98);
                    \node[diskEdge,anchor=south west] at (0.98,0.98) {\small $\mathcal{B}_{\mathrm{R}_{-}}(0)$};

                \end{tikzpicture}
                \caption{Sketch of relevant singularities on the real axis: the $\lambda_n$ sit in $(\mathfrak{L},\mathfrak{B}_1)$, the Bernoulli points $\alpha_n^{-1}$ are right of $\mathfrak{B}_2$, and geometric points $-\beta_n^{-1}$ have modulus at least $\mathfrak{U}^{-1}$. The shaded disk $\mathcal{B}_{\mathrm{R}_{-}}(0)$ contains no branch points, so $\mathfrak{F}_y(\cdot,t)$ is holomorphic there. Pole locations may depend on $n$ and $t$, but the inequalities are uniform by \textbf{(A1)}.}
                \label{phaseholomorphicdiskfig}
           
            \end{figure}
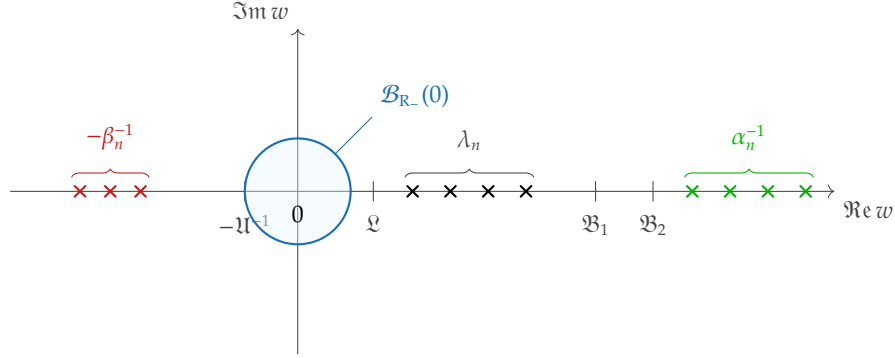
            As an illustration that branch cuts can be chosen disjoint from $\mathcal{B}_{\mathrm{R}_{-}}(0)$, one may take cuts along $[\lambda_n,\infty)$ for each spatial $\log(1-w/\lambda_n)$, along $[\alpha_n^{-1},\infty)$ for each Bernoulli $\log(1-\alpha_n w)$, and along $(-\infty,-\beta_n^{-1}]$ for each geometric $\log(1+\beta_n w)$. Under \textbf{(A1)} none of these rays meets the disk. Fixing values at $w=0$ as in Definition~\ref{form}, each summand is then holomorphic on $\mathcal{B}_{\mathrm{R}_{-}}(0)$, since the disk is simply connected and each logarithmic argument is holomorphic and nowhere vanishing there. For each fixed $t$, the phase is a finite sum of these logarithms and, in continuous time, one additional linear term. Hence $\mathfrak{F}_y(\cdot,t)$ is holomorphic in $\mathcal{B}_{\mathrm{R}_{-}}(0)$, proving \textbf{(i)}.
          
            \textbf{Step 2: The first derivative at the origin.} Differentiating the phase at the origin gives
            \begin{equation}\label{firstderivativeone}
                \mathfrak{F}_{y}^{(1)}(0,t)=\sum_{n=0}^{y(t)}\frac{1}{\lambda_{n}}-\mathbf{1}_{\{\mathscr{T}=\mathbb{R}_+\}}t-\mathbf{1}_{\{\mathscr{T}=\mathbb{Z}_+\}}\sum_{n=0}^{t-1}\left((1-\tau_{n})\alpha_{n}+\tau_{n}\beta_{n}\right).
            \end{equation}
            By Definition~\ref{themeananddeviation}, the temporal contribution in \eqref{firstderivativeone} is precisely $\Lambda\mathrm{z}(t)$. Thus the first derivative can be written in the simpler form
            \begin{equation}\label{expressionforthefirst}
                \mathfrak{F}_{y}^{(1)}(0,t)
                =
                \sum_{n=0}^{y(t)}\frac{1}{\lambda_n}-\Lambda\mathrm{z}(t).
            \end{equation}
            Let $\mathfrak{A}=\Lambda^{-1}\min(1,\mathfrak{U}^{-1},\mathfrak{L})$ and $\mathfrak{B}=\Lambda^{-1}\max(1,\mathfrak{B}_2^{-1},\mathfrak{U})$. Then, by the definition of $\mathrm{z}(t)$ and \textbf{(A1)}, for all $t\in\mathscr{T}$ we have
            \begin{equation}\label{zlinearbounds}
                \mathfrak{A}t\leq \mathrm{z}(t)\leq \mathfrak{B}t.
            \end{equation}
            Since $\mathfrak{A}>0$, \eqref{zlinearbounds} implies $\mathrm{z}(t)\to\infty$ and $\mathscr{Z}(t)\to\infty$ as $t\to\infty$, reflecting the linear time scaling of the deterministic mean path. Choose $\mathrm{T}_0\in\mathscr{T}$ such that, for all $t\geq\mathrm{T}_0$,
            \begin{equation}\label{Zlinearbounds}
                \frac{1}{2}\mathrm{z}(t)\leq\mathscr{Z}(t)\leq\frac{3}{2}\mathrm{z}(t),\qquad
                \mathscr{Z}(t)+1\leq2\mathscr{Z}(t).
            \end{equation}
            Such a choice is possible because $\mathscr{Z}(t)=\lfloor\mathrm{z}(t)\rceil$ differs from $\mathrm{z}(t)$ by at most $1/2$ and $\mathrm{z}(t)\to\infty$. Combining \eqref{zlinearbounds} and \eqref{Zlinearbounds} also gives
            \begin{equation}\label{Zlinearintbounds}
                \frac{\mathfrak{A}}{2}t\leq\mathscr{Z}(t)\leq\frac{3}{2}\mathfrak{B}t,\qquad t\geq\mathrm{T}_0.
            \end{equation}
            We now split \eqref{expressionforthefirst} at the deterministic integer centre $\mathscr{Z}(t)$:
            \begin{equation}\label{centeringequation}
                \mathfrak{F}_{y}^{(1)}(0,t)=\left(\sum_{n=0}^{y(t)}\frac{1}{\lambda_n}-\sum_{n=0}^{\mathscr{Z}(t)}\frac{1}{\lambda_n}\right)+\left(\sum_{n=0}^{\mathscr{Z}(t)}\frac{1}{\lambda_n}-\Lambda\mathrm{z}(t)\right).
            \end{equation}
            By \textbf{(A1)}, we have $\mathfrak{B}_1^{-1}\leq\lambda_n^{-1}\leq\mathfrak{L}^{-1}$ for every $n\in\mathbb{Z}_+$. The expression in the first bracket measures the deviation between the given path $y$ and the deterministic mean path $\mathscr{Z}(t)$, scaled by the spatial parameter, and is directly proportional to their difference. To formalise this, for any $t\in\mathscr{T}$, define the following functions:
            \begin{equation*}
                \mathrm{K}_{\mathrm{l}}(t)=
                \begin{cases}
                    \mathfrak{B}_1^{-1}, & \text{if } y(t)\geq\mathscr{Z}(t), \\
                    \mathfrak{L}^{-1}, & \text{if } y(t)<\mathscr{Z}(t),
                \end{cases}
                \qquad
                \mathrm{K}_{\mathrm{u}}(t)=
                \begin{cases}
                    \mathfrak{L}^{-1}, & \text{if } y(t)\geq\mathscr{Z}(t), \\
                    \mathfrak{B}_1^{-1}, & \text{if } y(t)<\mathscr{Z}(t).
                \end{cases}
            \end{equation*}
            This choice accounts for the direction of the displacement $y(t)-\mathscr{Z}(t)$. When $y(t)\geq\mathscr{Z}(t)$, the bracketed term is a positive sum over $y(t)-\mathscr{Z}(t)$ terms, each inverse-bounded by $\mathfrak{B}_1^{-1}$ from below and by $\mathfrak{L}^{-1}$ from above. Conversely, for $y(t)<\mathscr{Z}(t)$, the bracket becomes negative, corresponding to $-(\mathscr{Z}(t)-y(t))$ terms, and multiplying this by the larger reciprocal-rate bound yields the correct lower estimate. Thus, both cases are captured uniformly by these conventions.
       
            \begin{equation}\label{first-derivative-displacement}
                \mathrm{K}_{\mathrm{l}}(t)\left(y(t)-\mathscr{Z}(t)\right)
                \leq
                \sum_{n=0}^{y(t)}\frac{1}{\lambda_n}
                -
                \sum_{n=0}^{\mathscr{Z}(t)}\frac{1}{\lambda_n}
                \leq
                \mathrm{K}_{\mathrm{u}}(t)\left(y(t)-\mathscr{Z}(t)\right).
            \end{equation}
            The second bracket is deterministic and can be directly estimated via the scale function $\gamma(t)$. Adding and subtracting the empirical average at $\mathscr{Z}(t)$ gives
            \begin{equation*}
                \sum_{n=0}^{\mathscr{Z}(t)}\frac{1}{\lambda_n}-\Lambda\mathrm{z}(t)
                =
                \left(\frac{1}{\mathscr{Z}(t)+1}\sum_{n=0}^{\mathscr{Z}(t)}\frac{1}{\lambda_n}-\Lambda\right)\left(\mathscr{Z}(t)+1\right)
                +
                \Lambda\left(\mathscr{Z}(t)+1-\mathrm{z}(t)\right).
            \end{equation*}
            The first term is bounded by $\gamma(t)(\mathscr{Z}(t)+1)$ by Definition~\ref{gammadef}. The second term is bounded by $3\Lambda/2$: indeed, the nearest-integer convention gives $|\mathscr{Z}(t)-\mathrm{z}(t)|\leq1/2$, and upon using the triangle inequality, we get
            \begin{equation}\label{centering-error-bound}
                \left|\sum_{n=0}^{\mathscr{Z}(t)}\frac{1}{\lambda_n}-\Lambda\mathrm{z}(t)\right|\leq\gamma(t)\left(\mathscr{Z}(t)+1\right)+\frac{3\Lambda}{2}.
            \end{equation}
            By \eqref{Zlinearbounds}, the deterministic centring error is bounded by
            \begin{equation*}
                \left|\sum_{n=0}^{\mathscr{Z}(t)}\frac{1}{\lambda_n}-\Lambda\mathrm{z}(t)\right|\leq2\gamma(t)\mathscr{Z}(t)+\frac{3\Lambda}{2},\quad\text{for all $t\geq \mathrm{T}_0$.}
            \end{equation*}
            Combining the above and the previous estimate \eqref{first-derivative-displacement} with the formula for the first derivative of the phase given in \eqref{centeringequation} proves the first display in \textbf{(ii)}.
         
            \textbf{Step 3: Comparing $\Sigma(t)^2$ and $\mathscr{Z}(t)$.} We now prove the scale comparison in \textbf{(iv)}. The purpose of this step is to show that the deterministic fluctuation scale $\Sigma(t)^2$ is comparable with the deterministic endpoint $\mathscr{Z}(t)$ by constants depending only on the parameters in \textbf{(A1)} and on $\Lambda$. In the continuous-time regime, Definition~\ref{themeananddeviation} gives, for all $t\in\mathscr{T}$,
            \begin{equation*}
                \Sigma(t)^2=\sum_{n=0}^{\mathscr{Z}(t)}\frac{1}{\lambda_n^2}.
            \end{equation*}
            By \textbf{(A1)}\textnormal{(i)}, $\mathfrak{L}<\lambda_n<\mathfrak{B}_1$ for every $n\in\mathbb{Z}_+$. Taking reciprocals reverses the inequalities and gives $\mathfrak{B}_1^{-1}<\lambda_n^{-1}<\mathfrak{L}^{-1}$; hence $\mathfrak{B}_1^{-2}<\lambda_n^{-2}<\mathfrak{L}^{-2}$. Since the sum has $\mathscr{Z}(t)+1$ terms, and since $t>\mathrm{T}_0$ implies $\mathscr{Z}(t)+1\leq2\mathscr{Z}(t)$ by \eqref{Zlinearbounds}, we obtain
            \begin{equation}\label{condefsigm}
                \mathfrak{B}_1^{-2}\mathscr{Z}(t)\leq\Sigma(t)^2\leq2\mathfrak{L}^{-2}\mathscr{Z}(t),\quad\text{for $t>\mathrm{T}_0$}.
            \end{equation}
            In the discrete-time regime, the definition of $\Sigma(t)^2$ contains the same spatial contribution together with the inhomogeneous temporal variance contribution:
            \begin{equation}
                \Sigma(t)^2\overset{\textnormal{def}}{=}\sum_{n=0}^{\mathscr{Z}(t)}\frac{1}{\lambda_n^2}-\sum_{n=0}^{t-1}\big(\alpha_n^2(1-\tau_n)-\beta_n^2\tau_n\big).
            \end{equation}
            For the upper bound, the spatial part is bounded exactly as above by $2\mathfrak{L}^{-2}\mathscr{Z}(t)$ for $t>\mathrm{T}_0$. The negative Bernoulli term can only decrease $\Sigma(t)^2$, but the coarser bound obtained by taking absolute values is sufficient here. Indeed, \textbf{(A1)}\textnormal{(ii)} gives $\alpha_n^2\leq\mathfrak{B}_2^{-2}$ and $\beta_n^2\leq\mathfrak{U}^2$, while $0\leq\tau_n\leq1$. Therefore the total temporal square contribution is bounded above by $(\mathfrak{B}_2^{-2}+\mathfrak{U}^2)t$. Combining this with $t\leq2\mathfrak{A}^{-1}\mathscr{Z}(t)$ from \eqref{Zlinearintbounds} yields
            \begin{equation}\label{numbertwotwotwo}
                \Sigma(t)^2\leq2\mathfrak{L}^{-2}\mathscr{Z}(t)+\left(\mathfrak{B}_2^{-2}+\mathfrak{U}^{2}\right)t\leq\left(2\mathfrak{L}^{-2}+\frac{2(\mathfrak{B}_2^{-2}+\mathfrak{U}^{2})}{\mathfrak{A}}\right)\mathscr{Z}(t),\qquad t>\mathrm{T}_0.
            \end{equation}
            It remains to prove the lower bound in the discrete-time regime. The only possible loss is the Bernoulli part of the temporal variance, because it enters with a negative sign. We compare this negative quadratic term with the corresponding first-order Bernoulli drift, and then recover a positive multiple of the deterministic temporal drift $\Lambda\mathrm{z}(t)$, up to the reciprocal-rate centring error already isolated in \eqref{centering-error-bound}. Set
            \begin{equation*}
                \mathrm{d}_0\overset{\textnormal{def}}{=}\min\left(\mathfrak{B}_1^{-1}-\mathfrak{B}_2^{-1},\mathfrak{B}_1^{-1}\right)>0,
            \end{equation*}
            where positivity follows from $\mathfrak{B}_1<\mathfrak{B}_2$ in \textbf{(A1)}. By Definition~\ref{themeananddeviation}, in the discrete-time regime,
            \begin{equation*}
                \Sigma(t)^2=\sum_{n=0}^{\mathscr{Z}(t)}\frac{1}{\lambda_n^2}-\sum_{n=0}^{t-1}\alpha_n^2(1-\tau_n)+\sum_{n=0}^{t-1}\beta_n^2\tau_n.
            \end{equation*}
            The final sum is non-negative because $\beta_n>0$ and $\tau_n\geq0$, so dropping it gives a valid lower bound. For the spatial part, \textbf{(A1)}\textnormal{(ii)} gives $\lambda_n<\mathfrak{B}_1$, and hence $\lambda_n^{-1}>\mathfrak{B}_1^{-1}$. Multiplying this inequality by the positive number $\lambda_n^{-1}$ gives $\lambda_n^{-2}\geq\mathfrak{B}_1^{-1}\lambda_n^{-1}$. For the Bernoulli part, \textbf{(A1)}\textnormal{(ii)} gives $0<\alpha_n<\mathfrak{B}_2^{-1}$, and multiplying by the positive number $\alpha_n$ gives $\alpha_n^2\leq\mathfrak{B}_2^{-1}\alpha_n$. Since this Bernoulli contribution appears with a negative sign, replacing $\alpha_n^2$ by the larger quantity $\mathfrak{B}_2^{-1}\alpha_n$ preserves the lower bound. Therefore
            \begin{equation}\label{sigma-lower-first-step}
                \Sigma(t)^2\geq\mathfrak{B}_1^{-1}\sum_{n=0}^{\mathscr{Z}(t)}\frac{1}{\lambda_n}-\mathfrak{B}_2^{-1}\sum_{n=0}^{t-1}(1-\tau_n)\alpha_n.
            \end{equation}
            We now insert the deterministic centre $\Lambda\mathrm{z}(t)$ into the spatial sum. By Definition~\ref{themeananddeviation}, the discrete-time deterministic drift is
            \begin{equation*}
                \Lambda\mathrm{z}(t)=\sum_{n=0}^{t-1}(1-\tau_n)\alpha_n+\sum_{n=0}^{t-1}\tau_n\beta_n.
            \end{equation*}
            Hence the spatial sum admits the exact decomposition
            \begin{equation}\label{spatial-sum-splitting}
                \sum_{n=0}^{\mathscr{Z}(t)}\frac{1}{\lambda_n}=\sum_{n=0}^{t-1}(1-\tau_n)\alpha_n+\sum_{n=0}^{t-1}\tau_n\beta_n+\left(\sum_{n=0}^{\mathscr{Z}(t)}\frac{1}{\lambda_n}-\Lambda\mathrm{z}(t)\right).
            \end{equation}
            Substituting this into \eqref{sigma-lower-first-step} separates the positive temporal drift from the centring error:
            \begin{equation}\label{sigma-lower-separated}
                \Sigma(t)^2\geq\left(\mathfrak{B}_1^{-1}-\mathfrak{B}_2^{-1}\right)\sum_{n=0}^{t-1}(1-\tau_n)\alpha_n+\mathfrak{B}_1^{-1}\sum_{n=0}^{t-1}\tau_n\beta_n+\mathfrak{B}_1^{-1}\left(\sum_{n=0}^{\mathscr{Z}(t)}\frac{1}{\lambda_n}-\Lambda\mathrm{z}(t)\right).
            \end{equation}
            The first two sums in \eqref{sigma-lower-separated} are non-negative. By the definition of $\mathrm{d}_0$, their coefficients are both bounded below by $\mathrm{d}_0$. Hence
            \begin{equation*}
                \left(\mathfrak{B}_1^{-1}-\mathfrak{B}_2^{-1}\right)\sum_{n=0}^{t-1}(1-\tau_n)\alpha_n+\mathfrak{B}_1^{-1}\sum_{n=0}^{t-1}\tau_n\beta_n\geq\mathrm{d}_0\Lambda\mathrm{z}(t).
            \end{equation*}
            Finally, the remaining centring error may have either sign, so for a lower bound we estimate it by the negative of its absolute value. Combining this with \eqref{sigma-lower-separated} gives
            \begin{equation}\label{lowerboundforsigma}
                \Sigma(t)^2\geq\mathrm{d}_0\Lambda\mathrm{z}(t)-\mathfrak{B}_1^{-1}\left|\sum_{n=0}^{\mathscr{Z}(t)}\frac{1}{\lambda_n}-\Lambda\mathrm{z}(t)\right|.
            \end{equation}
            The lower estimate \eqref{lowerboundforsigma} separates the positive multiple of the deterministic drift $\Lambda\mathrm{z}(t)$ from the reciprocal-rate centring error controlled by the empirical-average assumption. We now control the absolute-value term in \eqref{lowerboundforsigma}. By \textbf{(A2)}, there exists $\mathrm{T}_{3,1}\in\mathscr{T}$ such that, for all $t>\mathrm{T}_{3,1}$,
            \begin{equation*}
                \left|\frac{1}{\mathscr{Z}(t)+1}\sum_{n=0}^{\mathscr{Z}(t)}\frac{1}{\lambda_n}-\Lambda\right|
                \leq
                \frac{\mathfrak{B}_1\mathrm{d}_0\Lambda}{10}.
            \end{equation*}
            Choose $\mathrm{T}_{3,2}\in\mathscr{T}$ with $\mathrm{T}_{3,2}\geq\max(\mathrm{T}_0,\mathrm{T}_{3,1})$ such that $\mathscr{Z}(t)+1\leq2\mathrm{z}(t)$ and $3\Lambda/2\leq3\mathfrak{B}_1\mathrm{d}_0\Lambda\mathrm{z}(t)/10$ for all $t>\mathrm{T}_{3,2}$. Such a choice is possible because $\mathrm{z}(t)\geq\mathfrak{A}t$ by \eqref{zlinearbounds} and $\mathscr{Z}(t)\leq3\mathrm{z}(t)/2$ by \eqref{Zlinearbounds}. Using the same add-and-subtract decomposition as in \eqref{centering-error-bound}, these choices give
            \begin{equation*}
                \left|\sum_{n=0}^{\mathscr{Z}(t)}\frac{1}{\lambda_n}-\Lambda\mathrm{z}(t)\right|
                \leq
                \frac{\mathfrak{B}_1\mathrm{d}_0\Lambda}{2}\mathrm{z}(t),
                \qquad t>\mathrm{T}_{3,2}.
            \end{equation*}
            Substituting this into \eqref{lowerboundforsigma} gives
            \begin{equation}\label{numberoneoneone}
                \Sigma(t)^2
                \geq
                \frac{\mathrm{d}_0\Lambda}{2}\mathrm{z}(t)
                \geq
                \frac{\mathrm{d}_0\Lambda}{3}\mathscr{Z}(t),
                \qquad t>\mathrm{T}_{3,2},
            \end{equation}
            where the final inequality uses $\mathscr{Z}(t)\leq3\mathrm{z}(t)/2$ from \eqref{Zlinearbounds}. Combining the continuous-time lower estimate \eqref{condefsigm}, the discrete-time lower estimate \eqref{numberoneoneone}, and the upper estimate \eqref{numbertwotwotwo}, we have a lower bound on $\Sigma(t)^2$ with coefficient at least $\min(\mathfrak{B}_1^{-2},\mathrm{d}_0\Lambda/3)$ and an upper bound on $\Sigma(t)^2$ with coefficient at most $2\mathfrak{L}^{-2}+2(\mathfrak{B}_2^{-2}+\mathfrak{U}^2)/\mathfrak{A}$. To state the comparison in the direction required by \textbf{(iv)}, define
            \begin{equation*}
                \mathrm{C}_1\overset{\textnormal{def}}{=} \left(2\mathfrak{L}^{-2}+\frac{2(\mathfrak{B}_2^{-2}+\mathfrak{U}^{2})}{\mathfrak{A}}\right)^{-1},
                \qquad
                \mathrm{C}_2\overset{\textnormal{def}}{=}\max\left(\mathfrak{B}_1^{2},\frac{3}{\mathrm{d}_0\Lambda}\right).
            \end{equation*}
            We may therefore choose the lower time in the scale comparison explicitly as follows: in the continuous-time regime take $\mathrm{T}_3\geq\mathrm{T}_0$, while in the discrete-time regime take $\mathrm{T}_3\geq\max(\mathrm{T}_0,\mathrm{T}_{3,2})$. With this choice,
            \begin{equation}\label{sigmazexplicit}
                \mathrm{C}_1\Sigma(t)^2\leq \mathscr{Z}(t)\leq \mathrm{C}_2\Sigma(t)^2,\qquad t>\mathrm{T}_3.
            \end{equation}
            This proves the first assertion in \textbf{(iv)}.
	         
                \textbf{Step 4: Higher derivatives.} We next prove \textbf{(iii)}. Since the statement requires a bound uniform over $w\in\mathcal{B}_{\mathrm{R}_2}(0)$, we must estimate the derivatives at a general point in this ball, rather than only at the origin. Set $\mathrm{R}_2=\mathrm{R}_{-}/2$. The first task is to turn \textbf{(A1)} into explicit lower bounds on every denominator which can occur after differentiating the logarithmic factors. For the spatial factors, \textbf{(A1)} gives $\lambda_m>\mathfrak{L}\geq\mathrm{R}_{-}$, and therefore, for $w\in\mathcal{B}_{\mathrm{R}_2}(0)$, the reverse triangle inequality gives
	            \begin{equation*}
	                |\lambda_m-w|\geq \lambda_m-|w|\geq \mathrm{R}_{-}-\mathrm{R}_2=\frac{\mathrm{R}_{-}}{2}.
	            \end{equation*}
		            For the Bernoulli factors in the discrete-time regime, \textbf{(A1)}\textnormal{(ii)} gives $0<\alpha_r<\mathfrak{B}_2^{-1}$. Taking reciprocals gives $\alpha_r^{-1}>\mathfrak{B}_2>\mathfrak{B}_1>\mathfrak{L}\geq\mathrm{R}_{-}$, and hence the same reverse-triangle estimate gives
	            \begin{equation*}
	                |\alpha_r^{-1}-w|\geq\alpha_r^{-1}-|w|\geq\mathrm{R}_{-}-\mathrm{R}_2=\frac{\mathrm{R}_{-}}{2}.
	            \end{equation*}
		            For the geometric factors, \textbf{(A1)}\textnormal{(ii)} gives $0<\beta_r<\mathfrak{U}$. Taking reciprocals gives $\beta_r^{-1}>\mathfrak{U}^{-1}\geq\mathrm{R}_{-}$, and hence
	            \begin{equation*}
	                |\beta_r^{-1}+w|\geq\beta_r^{-1}-|w|\geq\mathrm{R}_{-}-\mathrm{R}_2=\frac{\mathrm{R}_{-}}{2}.
	            \end{equation*}
	            Differentiating the phase term by term now gives, for every $n\geq1$ and every $w\in\mathcal{B}_{\mathrm{R}_2}(0)$,
		            \begin{align}\label{equationforderva}
		                \frac{1}{(n-1)!}\mathfrak{F}_{y}^{(n)}(w,t)
		                ={}&-t\mathbf{1}_{\{\mathscr{T}=\mathbb{R}_+,\,n=1\}}+\sum_{m=0}^{y(t)}\frac{1}{(\lambda_m-w)^n}\notag\\
		                &-\mathbf{1}_{\{\mathscr{T}=\mathbb{Z}_+\}}\sum_{r=0}^{t-1}\left((1-\tau_r)\frac{\alpha_r^n}{(1-\alpha_r w)^n}+(-1)^{n-1}\tau_r\frac{\beta_r^n}{(1+\beta_r w)^n}\right).
		            \end{align}
		            The first term in \eqref{equationforderva} comes from the continuous-time linear term $-tw$, and therefore appears only when $n=1$. The spatial contribution follows from differentiating $-\log(1-w/\lambda_m)$. The Bernoulli contribution follows from differentiating $\log(1-\alpha_r w)$, and the geometric contribution follows from differentiating $-\log(1+\beta_r w)$; in the displayed formula the geometric contribution appears inside the parentheses with the factor $(-1)^{n-1}$, together with the common outer minus sign. The denominator bounds above imply, for every admissible $w$,
	                    \begin{equation*}
	                        \begin{aligned}
	                            \left|\frac{1}{(\lambda_m-w)^n}\right|&\leq\left(\frac{2}{\mathrm{R}_{-}}\right)^n,\\
				                \left|\frac{\alpha_r^n}{(1-\alpha_r w)^n}\right|=\left|\frac{1}{(\alpha_r^{-1}-w)^n}\right|&\leq\left(\frac{2}{\mathrm{R}_{-}}\right)^n,\\
				                \left|\frac{\beta_r^n}{(1+\beta_r w)^n}\right|=\left|\frac{1}{(\beta_r^{-1}+w)^n}\right|&\leq\left(\frac{2}{\mathrm{R}_{-}}\right)^n.
	                        \end{aligned}
	                    \end{equation*}
	            Next, define
	            \begin{equation*}
	                \mathrm{M}_2\overset{\textnormal{def}}{=}\max\left(1,\frac{2}{\mathrm{R}_{-}}\right).
	            \end{equation*}
	            It follows from \eqref{equationforderva} that each logarithmic summand contributes at most $(n-1)!\mathrm{M}_2^n$ to the absolute value of the $n$-th derivative. The possible continuous-time linear contribution is also bounded by $(n-1)!\mathrm{M}_2^nt$, since it is equal to $t$ when $n=1$ and is zero when $n\geq2$. There are $y(t)+1$ spatial logarithmic terms; in the discrete-time regime there are exactly $t$ temporal logarithmic terms, and in the continuous-time regime the factor $t$ has already been included to control the linear term. Therefore, in both regimes,
	            \begin{equation*}
	                \left|\mathfrak{F}_{y}^{(n)}(w,t)\right|\leq(n-1)!\mathrm{M}_2^n\left(y(t)+1+t\right),\qquad n\geq1,\quad w\in\mathcal{B}_{\mathrm{R}_2}(0).
	            \end{equation*}
	            It remains to convert the factor $y(t)+1+t$ into the fluctuation scale $\Sigma(t)^2$ uniformly on the endpoint window. By \textbf{(A2)} and the scale comparison just proved, $\gamma(t)\to0$. Choose $\mathrm{T}_{2,\mathscr{C}}\in\mathscr{T}$ satisfying $\mathrm{T}_{2,\mathscr{C}}\geq\max(\mathrm{T}_0,\mathrm{T}_3)$, $\mathscr{Z}(t)\geq1$, and $\gamma(t)^{1/2}\leq\mathrm{C}_1$ for all $t>\mathrm{T}_{2,\mathscr{C}}$. If $y(t)\in\mathscr{C}_{t}$, then the definition of $\mathscr{C}_{t}$ and the lower scale comparison in \eqref{sigmazexplicit}, equivalently $\Sigma(t)^2\leq\mathrm{C}_1^{-1}\mathscr{Z}(t)$, give
            \begin{equation*}
                y(t)+1
                \leq
                \mathscr{Z}(t)+\gamma(t)^{1/2}\Sigma(t)^2+1
                \leq
                \mathscr{Z}(t)+\mathrm{C}_1^{-1}\gamma(t)^{1/2}\mathscr{Z}(t)+1
                \leq
                3\mathscr{Z}(t).
            \end{equation*}
            Combining this deterministic endpoint bound with $t\leq2\mathfrak{A}^{-1}\mathscr{Z}(t)$ from \eqref{Zlinearintbounds} gives $y(t)+1+t\leq(3+2\mathfrak{A}^{-1})\mathscr{Z}(t)$. Define
            \begin{equation*}
                \mathrm{M}_1\overset{\textnormal{def}}{=}\mathrm{C}_2\left(3+\frac{2}{\mathfrak{A}}\right).
            \end{equation*}
            Applying the upper scale comparison $\mathscr{Z}(t)\leq\mathrm{C}_2\Sigma(t)^2$ from \eqref{sigmazexplicit} to the bound above proves \textbf{(iii)}, uniformly over all endpoints in $\mathscr{C}_{t}$.
            
            \textbf{Step 5: The strengthened estimates on $\mathscr{C}_{t}$.} Assume now that $y(t)\in\mathscr{C}_{t}$. Then, by definition of the endpoint window,
            \begin{equation*}
                |y(t)-\mathscr{Z}(t)|\leq\gamma(t)^{1/2}\Sigma(t)^2.
            \end{equation*}
            Choose $\mathrm{T}_{\mathscr{C},1}\geq\max(\mathrm{T}_1,\mathrm{T}_3)$ satisfying $\gamma(t)\leq\gamma(t)^{1/2}$ for all $t>\mathrm{T}_{\mathscr{C},1}$. Taking absolute values in the first display of \textbf{(ii)} gives a coefficient $\mathfrak{L}^{-1}$ in front of $|y(t)-\mathscr{Z}(t)|$, because $\mathfrak{L}^{-1}$ is the larger of the two reciprocal-rate bounds $\mathfrak{B}_1^{-1}$ and $\mathfrak{L}^{-1}$. Using the endpoint-window bound $|y(t)-\mathscr{Z}(t)|\leq\gamma(t)^{1/2}\Sigma(t)^2$, the lower scale comparison $\Sigma(t)^2\leq \mathrm{C}_1^{-1}\mathscr{Z}(t)$ from \eqref{sigmazexplicit}, and $\gamma(t)\leq\gamma(t)^{1/2}$, we obtain, for $t>\mathrm{T}_{\mathscr{C},1}$,
            \begin{equation*}
                \left|\mathfrak{F}_{y}^{(1)}(0,t)\right|\leq\mathfrak{L}^{-1}\gamma(t)^{1/2}\Sigma(t)^2+2\gamma(t)\mathscr{Z}(t)+\frac{3\Lambda}{2}\leq\left(\mathfrak{L}^{-1}\mathrm{C}_1^{-1}+2\right)\gamma(t)^{1/2}\mathscr{Z}(t)+\frac{3\Lambda}{2}.
            \end{equation*}
            This proves the strengthened first-derivative estimate in \textbf{(ii)} with constants written explicitly; in particular, the bound is uniform over the endpoint window $\mathscr{C}_{t}$.
            It remains to compare $\mathfrak{F}_{y}^{(2)}(w,t)$ with $\Sigma(t)^2$. At the origin, the temporal terms in $\mathfrak{F}_{y}^{(2)}(0,t)$ and $\Sigma(t)^2$ are identical, so only the endpoint of the spatial sum changes:
            \begin{equation}\label{formulaforsecond}
                \mathfrak{F}_{y}^{(2)}(0,t)-\Sigma(t)^2
                =
                \sum_{n=0}^{y(t)}\frac{1}{\lambda_n^2}
                -
                \sum_{n=0}^{\mathscr{Z}(t)}\frac{1}{\lambda_n^2}.
            \end{equation}
            This difference has the same sign structure as the displacement term in the first-derivative estimate. If $y(t)\geq\mathscr{Z}(t)$, it is the sum of $y(t)-\mathscr{Z}(t)$ non-negative terms; if $y(t)<\mathscr{Z}(t)$, it is minus the sum of $\mathscr{Z}(t)-y(t)$ non-negative terms. Since \textbf{(A1)} gives $\lambda_n^{-2}\leq\mathfrak{L}^{-2}$ for every $n$, both cases give
            \begin{equation}\label{sumerrorterm}
                \left|\mathfrak{F}_{y}^{(2)}(0,t)-\Sigma(t)^2\right|
                \leq
                \mathfrak{L}^{-2}\left|y(t)-\mathscr{Z}(t)\right|.
            \end{equation}
            Hence the endpoint control from $y(t)\in\mathscr{C}_{t}$ yields
            \begin{equation*}
                \left|\mathfrak{F}_{y}^{(2)}(0,t)-\Sigma(t)^2\right|
                \leq
                \mathfrak{L}^{-2}\gamma(t)^{1/2}\Sigma(t)^2.
            \end{equation*}
            For a general point $w\in\mathcal{B}_{\mathrm{R}_2}(0)$, the line segment $\{sw:0\leq s\leq1\}$ stays inside $\mathcal{B}_{\mathrm{R}_2}(0)$. Therefore the fundamental theorem of calculus gives
            \begin{equation*}
                \mathfrak{F}_{y}^{(2)}(w,t)-\mathfrak{F}_{y}^{(2)}(0,t)
                =
                \int_0^1\frac{\mathrm{d}}{\mathrm{d}s}\mathfrak{F}_{y}^{(2)}(sw,t)\,\mathrm{d}s
                =
                w\int_0^1\mathfrak{F}_{y}^{(3)}(sw,t)\,\mathrm{d}s.
            \end{equation*}
            Applying the uniform endpoint version of \textbf{(iii)} with $n=3$ along this segment gives, for $t>\mathrm{T}_{2,\mathscr{C}}$,
            \begin{equation*}
                \left|\mathfrak{F}_{y}^{(2)}(w,t)-\mathfrak{F}_{y}^{(2)}(0,t)\right|
                \leq
                2\mathrm{M}_1\mathrm{M}_2^3|w|\Sigma(t)^2.
            \end{equation*}
            Combining the last two displays proves the final assertion in \textbf{(iv)} with
            \begin{equation*}
                \mathrm{C}_{\mathscr{C}}\overset{\textnormal{def}}{=}\mathfrak{L}^{-2}+2\mathrm{M}_1\mathrm{M}_2^3,
                \qquad
                \mathrm{T}_{4,\mathscr{C}}\overset{\textnormal{def}}{=}\max(\mathrm{T}_{\mathscr{C},1},\mathrm{T}_{2,\mathscr{C}}),
            \end{equation*}
            completing the proof.
	        \end{proof}
	    \end{proposition}
        \begin{figure}[ht!]
            \centering
            \includegraphics[width=0.7\textwidth]{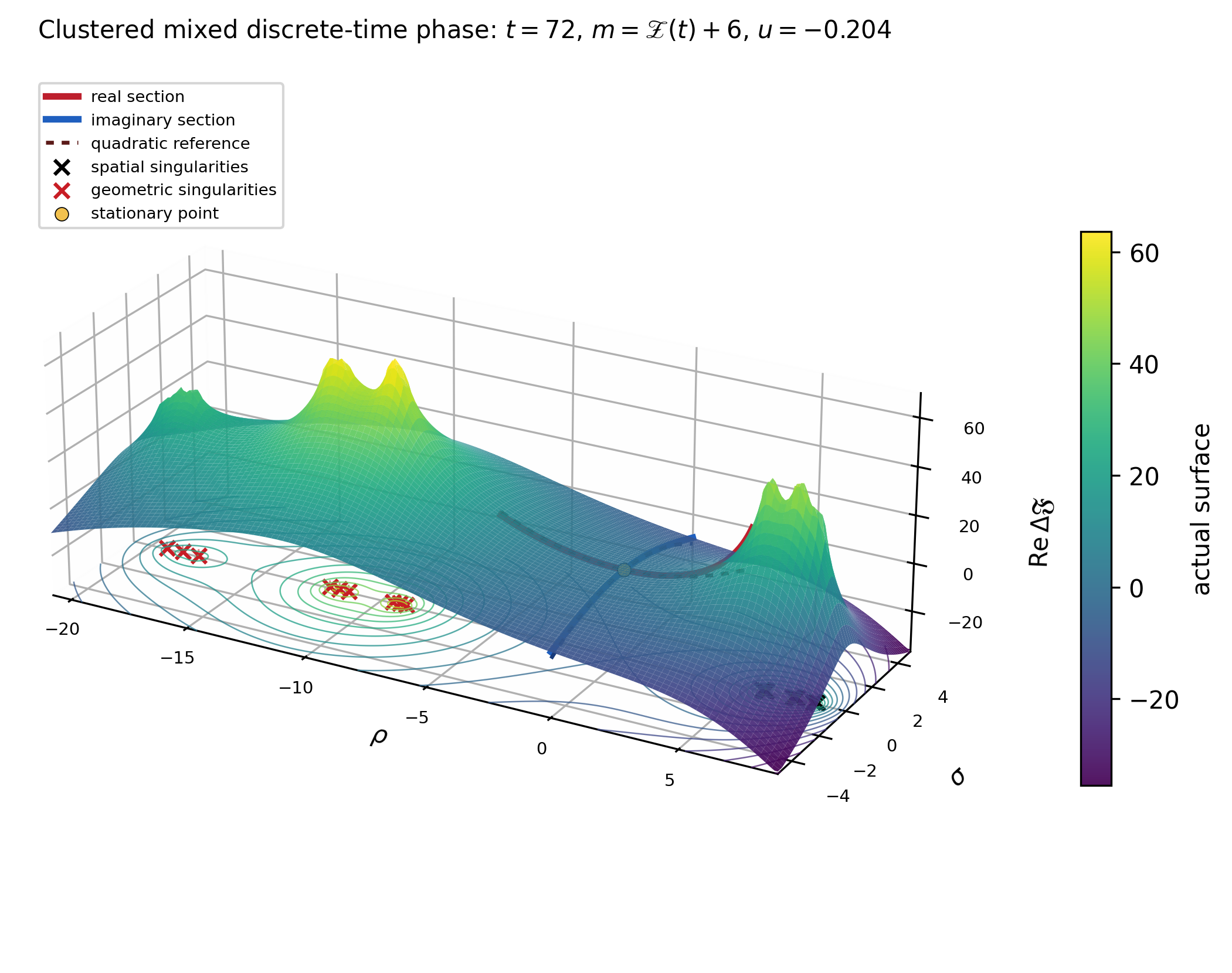}
            \caption{Finite-time surface plot of the mixed discrete-time phase for representative clustered parameters satisfying \textbf{(A1)}. The stationary point lies between the spatial singularities $w=\lambda_n$ and the geometric singularities $w=-1/\beta_r$, and the solid and dotted cross-sections compare the actual phase with its local quadratic approximation.}
            \label{phasequadraticsurface}
        \end{figure}

        \noindent The endpoint window $\mathscr{C}_{t}$ is wide enough to contain the endpoint fluctuations that will be supplied probabilistically in Section~\ref{MOM}, but narrow enough that the constants in Proposition~\ref{specificcasesub} do not depend on which endpoint in the summation is chosen. Section~\ref{ASP} uses the quadratic saddle picture to construct and approximate the stationary point, Section~\ref{CM} builds descent contours compatible with that geometry, and Section~\ref{decouplingsection} uses both inputs to derive the one-particle expansion on this deterministic endpoint window. Section~\ref{MOM} then proves that the random endpoints lie in the same window with high probability.
    \subsection{Approximate stationary point}\label{ASP}
    To perform steepest-descent analysis on the contour-integral representation of the transition probabilities obtained in Proposition~\ref{contourrepresentation}, we need a robust approximation of the stationary point of the phase in the admissible path regime. In the present inhomogeneous setting the phase function may have several stationary points, reflecting the finite collection of spatial and temporal singularities. The next proposition therefore shows that, once the terminal endpoint lies in $\mathscr{C}_{t}$, the stationary-point equation has one relevant real solution near the origin and that this solution is controlled by the linearized equation at the origin with an explicit quadratic error.
         \begin{proposition}\label{axpro}
            Assume \textbf{(A1)} and \textbf{(A2)}. There exist deterministic constants $\mathrm{T}_{\mathrm{A}}\in\mathscr{T}$ and $\mathrm{C}_{\mathrm{A}}>0$ such that, for every $t>\mathrm{T}_{\mathrm{A}}$ and every path $y\in\mathscr{P}$ satisfying $y(t)\in\mathscr{C}_{t}$, the minimizer in Definition~\ref{definitionstationary} is a unique real number, denoted by $\mathrm{u}_y(t)$. More explicitly,
	        \begin{equation}\label{firstconditiontosat}
			              \mathfrak{F}_y^{(1)}\!\left(\mathrm{u}_y(t),t\right)=0,\quad
			              |\mathrm{u}_y(t)|=\min_{\substack{v\in\mathbb{C}\\ \mathfrak{F}^{(1)}_y(v,t)=0}}|v|,\quad\text{and}\quad\mathscr{U}_{y}(t)=-\Sigma(t)\,\mathrm{u}_y(t).
	       \end{equation}
	       Moreover, for every such $t$ and $y$,
	       \begin{equation}\label{secondconditiontosat}
	           \left|\mathrm{u}_y(t)+\frac{\mathfrak{F}_y^{(1)}(0,t)}{\mathfrak{F}_y^{(2)}(0,t)}\right|
	           \leq
	           \mathrm{C}_{\mathrm{A}}\!\left|\frac{\mathfrak{F}_y^{(1)}(0,t)}{\mathfrak{F}_y^{(2)}(0,t)}\right|^{2}.
	       \end{equation}
        \end{proposition}
    \begin{proof}
    The proof is in three parts. \textbf{Step~1} rescales the stationary-point equation by $\mathfrak{F}_y^{(2)}(0,t)$, proves that this denominator is bounded away from zero on $\mathscr{C}_t$, and records the resulting linearized problem at the origin together with a smallness estimate for its root. \textbf{Step~2} treats the remaining nonlinearity via a Taylor expansion and a sign-bracketing argument to produce a real stationary point and the quadratic approximation error. \textbf{Step~3} upgrades this to uniqueness on a fixed neighbourhood and identifies the point selected by the minimisation rule in Definition~\ref{definitionstationary}.
   
    \textbf{Step 1: Normalization at the origin.} By \textbf{(iv)} of Proposition~\ref{specificcasesub}, there exist constants $\mathrm{C}_{\mathscr{C}}>0$ and $\mathrm{T}_{4,\mathscr{C}}\in\mathscr{T}$ such that, whenever $t>\mathrm{T}_{4,\mathscr{C}}$, $y(t)\in\mathscr{C}_{t}$ and $w\in\mathcal{B}_{\mathrm{R}_2}(0)$,
    \begin{equation*}
        \left|\mathfrak{F}_{y}^{(2)}(w,t)-\Sigma(t)^2\right|\leq\mathrm{C}_{\mathscr{C}}\left(\gamma(t)^{1/2}+|w|\right)\Sigma(t)^2.
    \end{equation*}
    By \textbf{(A2)} and \textbf{(iv)} of Proposition~\ref{specificcasesub}, $\gamma(t)\to0$. We may therefore choose a deterministic threshold $\mathrm{T}_{\mathrm{A},0}\in\mathscr{T}$ such that $\mathrm{T}_{\mathrm{A},0}>\mathrm{T}_{4,\mathscr{C}}$ and $\mathrm{C}_{\mathscr{C}}\gamma(t)^{1/2}\leq1/2$ for all $t>\mathrm{T}_{\mathrm{A},0}$. Applying the preceding display at $w=0$ then gives the two-sided denominator bound
    \begin{equation}\label{axpro-denominator}
        \frac{1}{2}\Sigma(t)^2\leq\mathfrak{F}_y^{(2)}(0,t)\leq\frac{3}{2}\Sigma(t)^2,
        \qquad t>\mathrm{T}_{\mathrm{A},0},\quad y(t)\in\mathscr{C}_{t}.
    \end{equation}
    In particular, $\mathfrak{F}_y^{(2)}(0,t)>0$ on the full endpoint window, so the following normalization is well defined there:
    \begin{equation}\label{axpro-J-def}
        \mathfrak{J}_{y}(w,t)\overset{\textnormal{def}}{=}\frac{\mathfrak{F}_y^{(1)}(w,t)}{\mathfrak{F}_y^{(2)}(0,t)},\qquad a_{0,y}(t)\overset{\textnormal{def}}{=}\mathfrak{J}_{y}(0,t)=\frac{\mathfrak{F}_y^{(1)}(0,t)}{\mathfrak{F}_y^{(2)}(0,t)}.
    \end{equation}
    By \textbf{(i)} of Proposition~\ref{specificcasesub}, the phase is holomorphic in a fixed neighbourhood of the origin; hence $\mathfrak{J}_y(\cdot,t)$ is holomorphic there as well. Moreover,
    \begin{equation}\label{axpro-J-linear-data}
        \mathfrak{J}_{y}^{(1)}(0,t)=\frac{\mathfrak{F}_y^{(2)}(0,t)}{\mathfrak{F}_y^{(2)}(0,t)}=1.
    \end{equation}
    Thus the linear part of $\mathfrak{J}_y(w,t)$ at the origin is exactly $a_{0,y}(t)+w$. We next record the endpoint-uniform size of the zero of this linear approximation. Let
    \begin{equation}\label{axpro-BA-def}
        \mathrm{A}_{\mathrm{A}}\overset{\textnormal{def}}{=}\mathfrak{L}^{-1}\mathrm{C}_{1}^{-1}+2,
        \qquad
        \mathrm{B}_{\mathrm{A}}\overset{\textnormal{def}}{=}2\mathrm{A}_{\mathrm{A}}\mathrm{C}_{2}+3\Lambda,
    \end{equation}
    where $\mathrm{C}_1,\mathrm{C}_2$ are the scale constants from \textbf{(iv)} of Proposition~\ref{specificcasesub}. The endpoint first-derivative estimate in \textbf{(ii)} of that proposition gives, for $t>\mathrm{T}_{\mathscr{C},1}$ and $y(t)\in\mathscr{C}_{t}$,
    \begin{equation*}
        \left|\mathfrak{F}_{y}^{(1)}(0,t)\right|\leq \mathrm{A}_{\mathrm{A}}\gamma(t)^{1/2}\mathscr{Z}(t)+\frac{3\Lambda}{2}.
    \end{equation*}
    Dividing this inequality by the lower denominator bound in \eqref{axpro-denominator}, and then using $\mathscr{Z}(t)\leq\mathrm{C}_2\Sigma(t)^2$ from \textbf{(iv)} of Proposition~\ref{specificcasesub}, gives
    \begin{equation}\label{axpro-a0-bound}
        |a_{0,y}(t)|\leq2\mathrm{A}_{\mathrm{A}}\mathrm{C}_{2}\gamma(t)^{1/2}+3\Lambda\Sigma(t)^{-2}.
    \end{equation}
    Finally, the definition of $\gamma(t)$ gives $\gamma(t)\geq\Sigma(t)^{-1/2}$, and \textbf{(iv)} of Proposition~\ref{specificcasesub} implies $\Sigma(t)\to\infty$. Increasing $\mathrm{T}_{\mathrm{A},0}$ if necessary so that $\Sigma(t)\geq1$ for all $t>\mathrm{T}_{\mathrm{A},0}$, we have $\Sigma(t)^{-2}\leq\Sigma(t)^{-1/4}\leq\gamma(t)^{1/2}$ on this range. Therefore \eqref{axpro-a0-bound} yields the uniform Newton-scale estimate
    \begin{equation}\label{axpro-a0-small}
        |a_{0,y}(t)|\leq\mathrm{B}_{\mathrm{A}}\gamma(t)^{1/2},
        \qquad t>\max(\mathrm{T}_{\mathrm{A},0},\mathrm{T}_{\mathscr{C},1}),\quad y(t)\in\mathscr{C}_{t}.
    \end{equation}
   
    \textbf{Step 2: Bracketing a real zero.} We now use the second derivative of the normalized function to control the nonlinear error. Let $\mathrm{M}_1,\mathrm{M}_2$ and $\mathrm{R}_2$ be the constants from \textbf{(iii)} of Proposition~\ref{specificcasesub}. Increasing the lower time so that $t$ exceeds both $\mathrm{T}_{2,\mathscr{C}}$ and $\mathrm{T}_{\mathrm{A},0}$, the uniform endpoint derivative bound in \textbf{(iii)} with $n=3$, together with \eqref{axpro-denominator}, gives, for every $w\in\mathcal{B}_{\mathrm{R}_2}(0)$,
    \begin{equation}\label{axpro-Jsecond}
        \left|\mathfrak{J}_{y}^{(2)}(w,t)\right|=\left|\frac{\mathfrak{F}_{y}^{(3)}(w,t)}{\mathfrak{F}_y^{(2)}(0,t)}\right|\leq\frac{2\mathrm{M}_1\mathrm{M}_2^3\Sigma(t)^2}{\frac{1}{2}\Sigma(t)^2}=4\mathrm{M}_1\mathrm{M}_2^3\overset{\textnormal{def}}{=}\mathrm{D}_{\mathrm{J}}.
    \end{equation}
    Define the deterministic constants
    \begin{equation}\label{axpro-CA-def}
        \mathrm{r}_{\mathrm{A}}\overset{\textnormal{def}}{=}\frac{\mathrm{R}_2}{2},\qquad \mathrm{C}_{\mathrm{A}}\overset{\textnormal{def}}{=}4\mathrm{D}_{\mathrm{J}}.
    \end{equation}
    Since $\gamma(t)\to0$, the estimate \eqref{axpro-a0-small} allows us to increase the lower time to a deterministic threshold $\mathrm{T}_{\mathrm{A},1}$ such that, for every $t>\mathrm{T}_{\mathrm{A},1}$ and every $y(t)\in\mathscr{C}_{t}$,
    \begin{equation}\label{axpro-small-a}
        |a_{0,y}(t)|\leq\frac{\mathrm{r}_{\mathrm{A}}}{4},\qquad \mathrm{C}_{\mathrm{A}}|a_{0,y}(t)|\leq1.
    \end{equation}
    Fix $t>\mathrm{T}_{\mathrm{A},1}$ and a path $y$ with $y(t)\in\mathscr{C}_{t}$. If $w$ belongs to the real interval
    \begin{equation*}
        \bigl[-a_{0,y}(t)-\mathrm{C}_{\mathrm{A}}|a_{0,y}(t)|^2,\,-a_{0,y}(t)+\mathrm{C}_{\mathrm{A}}|a_{0,y}(t)|^2\bigr],
    \end{equation*}
    then \eqref{axpro-small-a} gives
    \begin{equation}\label{axpro-w-small}
        |w|\leq|a_{0,y}(t)|+\mathrm{C}_{\mathrm{A}}|a_{0,y}(t)|^2\leq2|a_{0,y}(t)|\leq\frac{\mathrm{r}_{\mathrm{A}}}{2}<\mathrm{R}_2.
    \end{equation}
    Hence the whole line segment from $0$ to $w$ lies inside $\mathcal{B}_{\mathrm{R}_2}(0)$. Taylor's theorem with integral remainder, together with \eqref{axpro-J-linear-data} and \eqref{axpro-Jsecond}, gives
    \begin{equation}\label{axpro-remainder}
        \mathfrak{J}_y(w,t)=a_{0,y}(t)+w+w^2\int_0^1(1-s)\mathfrak{J}_y^{(2)}(sw,t)\,\mathrm{d}s,
    \end{equation}
    and therefore
    \begin{equation}\label{axpro-remainder-bound}
        \left|\mathfrak{J}_{y}(w,t)-a_{0,y}(t)-w\right|\leq\frac{\mathrm{D}_{\mathrm{J}}}{2}|w|^2\leq2\mathrm{D}_{\mathrm{J}}|a_{0,y}(t)|^2=\frac{\mathrm{C}_{\mathrm{A}}}{2}|a_{0,y}(t)|^2.
    \end{equation}
    Evaluating \eqref{axpro-remainder-bound} at the endpoints of the interval gives
    \begin{equation}\label{axpro-sign-bracket}
        \mathfrak{J}_{y}\bigl(-a_{0,y}(t)-\mathrm{C}_{\mathrm{A}}|a_{0,y}(t)|^2,t\bigr)\leq-\frac{\mathrm{C}_{\mathrm{A}}}{2}|a_{0,y}(t)|^2,
        \qquad
        \mathfrak{J}_{y}\bigl(-a_{0,y}(t)+\mathrm{C}_{\mathrm{A}}|a_{0,y}(t)|^2,t\bigr)\geq\frac{\mathrm{C}_{\mathrm{A}}}{2}|a_{0,y}(t)|^2.
    \end{equation}
    The function $\mathfrak{J}_y(\cdot,t)$ is real-valued on this real interval because all parameters in the phase are real and the interval lies inside the holomorphic ball from \textbf{(i)} of Proposition~\ref{specificcasesub}. If $a_{0,y}(t)=0$, then $\mathfrak{J}_y(0,t)=0$ and we set $\mathrm{u}_y(t)=0$. If $a_{0,y}(t)\neq0$, the two inequalities in \eqref{axpro-sign-bracket} and the intermediate value theorem give a real zero $\mathrm{u}_y(t)$ in this interval.
    \begin{figure}[ht!]
        \centering
        \begin{tikzpicture}[x=6.1cm,y=3.0cm,>=Stealth,line cap=round,line join=round]
            \definecolor{rootBlue}{RGB}{30,80,190}
            \definecolor{rootRed}{RGB}{190,45,35}
            \definecolor{softGray}{RGB}{110,110,110}
            \draw[->,softGray] (-0.92,0) -- (0.28,0) node[right] {$w$};
            \draw[->,softGray] (0,-0.42) -- (0,0.72) node[above] {$\mathfrak{J}_y(w,t)$};
            \draw[dashed,softGray,domain=-0.82:0.18,samples=2] plot (\x,{0.42+\x});
            \draw[thick,black,domain=-0.82:0.18,samples=160,smooth] plot (\x,{0.42+\x+0.18*(\x)^2});
            \draw[densely dotted,rootRed] (-0.52,-0.18) -- (-0.52,0.10);
            \draw[densely dotted,softGray] (-0.42,-0.18) -- (-0.42,0.09);
            \draw[densely dotted,rootBlue] (-0.32,-0.18) -- (-0.32,0.18);
            \fill[rootRed] (-0.52,0) circle (0.7pt);
            \fill[black] (-0.452,0) circle (0.7pt);
            \fill[rootBlue] (-0.32,0) circle (0.7pt);
            \node[font=\scriptsize,rootRed,anchor=east] at (-0.6,-0.1) {$-a_{0,y}(t)-\mathrm{C}_{\mathrm{A}}|a_{0,y}(t)|^2$};
            \node[font=\scriptsize,softGray,anchor=north] at (-0.42,-0.20) {$-a_{0,y}(t)$};
            \node[font=\scriptsize,rootBlue,anchor=west] at (-0.30,-0.15) {$-a_{0,y}(t)+\mathrm{C}_{\mathrm{A}}|a_{0,y}(t)|^2$};
            \node[font=\scriptsize,anchor=south east] at (-0.55,0.37) {$a_{0,y}(t)+w$};
            \node[font=\scriptsize,anchor=south west] at (-0.18,0.33) {$\mathfrak{J}_y(w,t)$};
            \draw[decorate,decoration={brace,amplitude=4pt,mirror},softGray] (-0.52,-0.34) -- (-0.32,-0.34);
            \node[font=\scriptsize,softGray,anchor=north] at (-0.42,-0.39) {bracketing interval};
        \end{tikzpicture}
        \caption{Schematic of the bracketing argument in the proof of Proposition~\ref{axpro}, drawn for $a_{0,y}(t)>0$. The dashed line is the linear approximation $a_{0,y}(t)+w$, whose zero is $-a_{0,y}(t)$. The normalized derivative $\mathfrak{J}_y(w,t)$ differs from this line by at most a quadratic error on the displayed interval, so its values at the two endpoints have opposite signs. The case $a_{0,y}(t)<0$ is the reflected picture.}
        \label{approxstationarypointfig}
    \end{figure}
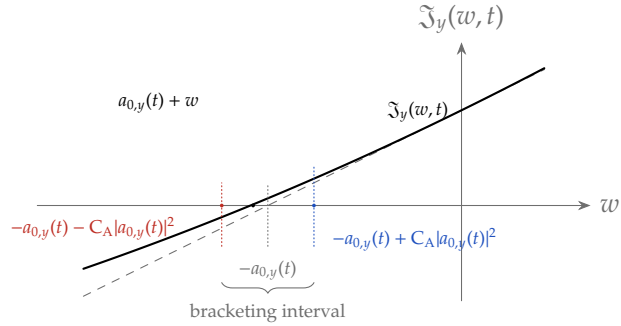
    \noindent Since $\mathfrak{F}_y^{(2)}(0,t)>0$ on this range, the zeros of $\mathfrak{J}_y(\cdot,t)$ are exactly the zeros of $\mathfrak{F}_y^{(1)}(\cdot,t)$. Thus the zero just constructed satisfies $\mathfrak{F}_y^{(1)}(\mathrm{u}_y(t),t)=0$. The interval localization \eqref{axpro-remainder-bound} evaluated at $w=\mathrm{u}_y(t)$ also gives
    \begin{equation}\label{axpro-u-localized}
        \left|\mathrm{u}_y(t)+a_{0,y}(t)\right|\leq\mathrm{C}_{\mathrm{A}}|a_{0,y}(t)|^2.
    \end{equation}
    Substituting the definition of $a_{0,y}(t)$ from \eqref{axpro-J-def} proves \eqref{secondconditiontosat}. Combining \eqref{axpro-u-localized}, \eqref{axpro-small-a} and \eqref{axpro-a0-small} gives
    \begin{equation}\label{axpro-u-small}
        |\mathrm{u}_y(t)|\leq |a_{0,y}(t)|+\mathrm{C}_{\mathrm{A}}|a_{0,y}(t)|^2\leq2|a_{0,y}(t)|\leq2\mathrm{B}_{\mathrm{A}}\gamma(t)^{1/2}.
    \end{equation}
   
    \textbf{Step 3: Local uniqueness.} It remains to show that the minimisation rule in Definition~\ref{definitionstationary} selects exactly the zero constructed above. We first prove injectivity of $\mathfrak{J}_y(\cdot,t)$ in a fixed ball around the origin. For $|z|\leq\mathrm{R}_2$, the line segment from $0$ to $z$ is contained in $\mathcal{B}_{\mathrm{R}_2}(0)$, and the fundamental theorem of calculus gives
    \begin{equation}\label{axpro-Jprime-bound}
        \mathfrak{J}_y^{(1)}(z,t)-1=z\int_0^1\mathfrak{J}_y^{(2)}(sz,t)\,\mathrm{d}s,
        \quad\text{and then via \eqref{axpro-Jsecond},}\quad
        \left|\mathfrak{J}_y^{(1)}(z,t)-1\right|\leq\mathrm{D}_{\mathrm{J}}|z|.
    \end{equation}
    Define
    \begin{equation}\label{axpro-injectivity-radius}
        \mathrm{r}_{\mathrm{J}}\overset{\textnormal{def}}{=}\min\left(\frac{\mathrm{R}_2}{2},\frac{1}{2\mathrm{D}_{\mathrm{J}}}\right).
    \end{equation}
    If $z_1,z_2\in\mathcal{B}_{\mathrm{r}_{\mathrm{J}}}(0)$, the segment joining them is contained in $\mathcal{B}_{\mathrm{R}_2}(0)$. Applying the fundamental theorem of calculus along this segment and subtracting the identity map gives
    \begin{equation}\label{axpro-injectivity-expansion}
        \mathfrak{J}_y(z_1,t)-\mathfrak{J}_y(z_2,t)=(z_1-z_2)+(z_1-z_2)\int_0^1\left(\mathfrak{J}_y^{(1)}(z_2+s(z_1-z_2),t)-1\right)\,\mathrm{d}s.
    \end{equation}
    By \eqref{axpro-Jprime-bound} and \eqref{axpro-injectivity-radius}, the integral in \eqref{axpro-injectivity-expansion} has modulus at most $1/2$. The reverse triangle inequality therefore gives
    \begin{equation}\label{axpro-injective-bound}
        \left|\mathfrak{J}_y(z_1,t)-\mathfrak{J}_y(z_2,t)\right|\geq\frac{1}{2}|z_1-z_2|.
    \end{equation}
    Thus $\mathfrak{J}_y(\cdot,t)$ is injective on $\mathcal{B}_{\mathrm{r}_{\mathrm{J}}}(0)$ throughout the endpoint window. Increasing the final deterministic lower time once more, using \eqref{axpro-u-small} and $\gamma(t)\to0$, we may also assume that $|\mathrm{u}_y(t)|<\mathrm{r}_{\mathrm{J}}$ for every $y(t)\in\mathscr{C}_{t}$. Fix $t$ beyond this final threshold and $y(t)\in\mathscr{C}_{t}$. Since $\mathrm{u}_y(t)$ is a zero of $\mathfrak{J}_y(\cdot,t)$, any zero with modulus at most $|\mathrm{u}_y(t)|$ lies in the compact closed disk $\{|z|\leq |\mathrm{u}_y(t)|\}\subset\mathcal{B}_{\mathrm{r}_{\mathrm{J}}}(0)$. The zero set of the holomorphic function $\mathfrak{J}_y(\cdot,t)$ in this disk is closed and non-empty, so we may choose a zero $v$ of minimal modulus in this disk. Every zero outside this disk has strictly larger modulus than $|\mathrm{u}_y(t)|$, and hence larger modulus than $|v|$, so $v$ is a global minimum-modulus zero. Since $\mathfrak{J}_y(v,t)=0=\mathfrak{J}_y(\mathrm{u}_y(t),t)$ and both points lie in $\mathcal{B}_{\mathrm{r}_{\mathrm{J}}}(0)$, injectivity from \eqref{axpro-injective-bound} forces $v=\mathrm{u}_y(t)$. Therefore the minimizer in Definition~\ref{definitionstationary} is unique and equals $\mathrm{u}_y(t)$, so
    \begin{equation*}
        \mathscr{U}_y(t)=-\Sigma(t)\,\mathrm{u}_y(t).
    \end{equation*}
    Taking $\mathrm{T}_{\mathrm{A}}$ to dominate all deterministic thresholds introduced above completes the proof.
    \end{proof}

\subsection{Contour manipulation}\label{CM}
In this subsection we construct explicit integration contours for the contour integral representation of the transition kernel from Proposition~\ref{contourrepresentation}, written in exponential form with phase $\mathfrak{F}_y(\cdot,t)$ from Definition~\ref{form}. The deterministic phase estimates in Proposition~\ref{specificcasesub} supply holomorphicity, stable second-order structure near the origin, and uniform derivative bounds on the endpoint window $\mathscr{C}_{t}$, while Proposition~\ref{axpro} (together with the minimisation rule in Definition~\ref{definitionstationary}) identifies the real stationary point $\mathrm{u}_y(t)$ at which the integrand should be expanded. The geometric problem is therefore to choose a closed contour passing through (or near) $\mathrm{u}_y(t)$ along which $\mathfrak{Re}[\mathfrak{F}_y]$ decays strongly enough away from the stationary point to make the Cauchy tail exponentially small, while keeping the deformation inside the domain where the preceding estimates are valid.

 We do not seek a globally exact classical steepest-descent contour. Instead we work with the weaker notion of a descent curve. Concretely, we first build explicit global circular arcs on which $\mathfrak{Re}[\mathfrak{F}_y]$ decreases outwards from $\mathrm{u}_y(t)$, then splice these to a local segment that approximates the true steepest-descent direction through $\mathrm{u}_y(t)$ closely enough for the fine Taylor expansion carried out later. The concatenated contour may be deformed from the original Cauchy contour in \eqref{cir} and splits naturally into a small neighbourhood of $\mathrm{u}_y(t)$ and a complementary remainder on which the exponential factor is strictly smaller. All constructions are uniform over terminal endpoints $y(t)\in\mathscr{C}_{t}$; probabilistic input guaranteeing that random endpoints lie in $\mathscr{C}_{t}$ with asymptotically full probability is supplied later in Section~\ref{MOM}.
\medskip

\begin{definition}[Descent curve] \label{descentcurvedefn}
Let $\Omega\subset\mathbb C$, let $\mathfrak{F}:\Omega\to\mathbb C$, and let $I\subset\mathbb R$ be an interval. Let $\Phi:I\to\mathbb C$ be a continuous injective map. We say that $\Phi$ is a descent curve for $\mathfrak{F}$ in $\Omega$, from $x\in\Omega$, if there exists $\phi\in \Phi^{-1}(\Omega)$ such that $\Phi(\phi)=x$, and, writing $\mathcal{K}_{\Omega,x}$ for the connected component of $\Phi^{-1}(\Omega)$ containing $\phi$, the following conditions hold. First, the function
\[
\theta \mapsto \mathfrak{Re}\bigl[\mathfrak{F}(\Phi(\theta))\bigr],\qquad \theta\in \mathcal{K}_{\Omega,x},
\]
is locally absolutely continuous on $\mathcal{K}_{\Omega,x}$. Second, there exists a countable, possibly empty, set $\{p_n\}\subset \mathcal{K}_{\Omega,x}$ such that $\theta \mapsto \mathfrak{Re}\bigl[\mathfrak{F}(\Phi(\theta))\bigr]$ is differentiable at every point of $\mathcal{K}_{\Omega,x}^\circ\setminus \{p_n\}$, and
\[
\frac{\mathrm{d}}{\mathrm{d}\theta}\mathfrak{Re}\bigl[\mathfrak{F}(\Phi(\theta))\bigr]>0,\qquad
\theta\in \left(\mathcal{K}_{\Omega,x}^\circ\cap(-\infty,\phi)\right)\setminus \{p_n\},
\]
while
\[
\frac{\mathrm{d}}{\mathrm{d}\theta}\mathfrak{Re}\bigl[\mathfrak{F}(\Phi(\theta))\bigr]<0,\qquad
\theta\in \left(\mathcal{K}_{\Omega,x}^\circ\cap(\phi,\infty)\right)\setminus \{p_n\}.
\]
The restricted curve $\Phi|_{\mathcal{K}_{\Omega,x}}$ is called the descent segment of $\Phi$ from $x$.
\end{definition}

 Notice that for an injective $\Phi:I\to\mathbb{C}$ to be a steepest descent contour for $\mathfrak{F}$ holomorphic in $\Omega$, from $x\in\Omega$, it must be a descent curve and also satisfy the level-set condition
\begin{equation*}
\mathfrak{Im}\left[\mathfrak{F}(\Phi(\theta))-\mathfrak{F}(x)\right]=\textnormal{constant}.
\end{equation*}
As we will see in this subsection, descent curves require less structure to construct than exact steepest-descent contours and still satisfy the following key property.

\begin{proposition}\label{prop:descent-endpoint-max}
Suppose that $\Phi:I\to\mathbb C$ is a descent curve for $\mathfrak{F}:\Omega\to\mathbb C$ in $\Omega\subset\mathbb C$, from $x=\Phi(\phi)\in\Omega$, in the sense of Definition \ref{descentcurvedefn}. Let $\mathcal{K}_{\Omega,x}$ be the parameter interval of the descent segment from $x$. Let $J$ be a non-empty compact subset of $\mathcal{K}_{\Omega,x}$, and define the one-sided nearest points
\[
\theta_-\overset{\mathrm{def}}{=}\max\bigl(J\cap(-\infty,\phi]\bigr),\qquad
\theta_+\overset{\mathrm{def}}{=}\min\bigl(J\cap[\phi,\infty)\bigr),
\]
whenever the corresponding sets are non-empty. Then
\[
\max_{\theta\in J}\mathfrak{Re}\bigl[\mathfrak{F}(\Phi(\theta))\bigr]
=
\max\left\{
\mathfrak{Re}\bigl[\mathfrak{F}(\Phi(\theta_-))\bigr],
\mathfrak{Re}\bigl[\mathfrak{F}(\Phi(\theta_+))\bigr]
\right\},
\]
where an empty one-sided candidate is omitted from the maximum.
\end{proposition}

\begin{proof}
To ease notation let us write
\[
g(\theta)=\Re\bigl[\mathfrak{F}(\Phi(\theta))\bigr],\qquad \theta\in \mathcal{K}_{\Omega,x}.
\]
By Definition \ref{descentcurvedefn}, $g$ is locally absolutely continuous on $\mathcal{K}_{\Omega,x}$, and hence continuous on $\mathcal{K}_{\Omega,x}$. Therefore $g$ attains its maximum on the compact set $J$. The intersections $J\cap[\phi,\infty)$ and $J\cap(-\infty,\phi]$ are compact whenever they are non-empty, so the one-sided nearest points $\theta_+$ and $\theta_-$ exist in the cases appearing in the statement.

Assume first that $J\cap[\phi,\infty)\neq\varnothing$. Let $\theta\in J\cap[\phi,\infty)$. If $\theta=\theta_+$, there is nothing to prove. If $\theta>\theta_+$, then the interval $[\theta_+,\theta]$ is contained in $\mathcal{K}_{\Omega,x}$, since $\mathcal{K}_{\Omega,x}$ is an interval. By local absolute continuity and the fundamental theorem of calculus for absolutely continuous functions,
\[
g(\theta)-g(\theta_+)=\int_{\theta_+}^{\theta} g'(s)\,\mathrm{d}s.
\]
The derivative satisfies $g'(s)<0$ for all $s\in(\theta_+,\theta)\setminus \{p_n\}$, where $\{p_n\}$ is countable. Since $\{p_n\}$ has Lebesgue measure zero, it follows that $g'(s)<0$ for almost every $s\in(\theta_+,\theta)$. Hence
\[
\int_{\theta_+}^{\theta} g'(s)\,\mathrm{d}s<0,
\]
and therefore $g(\theta)<g(\theta_+)$. Thus the maximum of $g$ over $J\cap[\phi,\infty)$ is attained at $\theta_+$. The left side is completely analogous. Since
\[
J=\bigl(J\cap(-\infty,\phi]\bigr)\cup\bigl(J\cap[\phi,\infty)\bigr),
\]
the maximum over $J$ is the larger of the two one-sided maxima, with any empty side omitted. This proves the claim.
\end{proof}

 In the following proposition we construct a simple but effective descent curve for the phase function $\mathfrak{F}_y(\cdot,t)$ from Definition~\ref{form}.
\begin{proposition}\label{contour1}
Let $\mathrm{R}_{\mathrm{C}}\overset{\textnormal{def}}{=}(\mathfrak{B}_1+\mathfrak{B}_2)/2$.  Then for all $y\in\mathscr{P}$, all $\mathrm{u}\in\mathbb{R}$ with $|\mathrm{u}|<\min\{\mathfrak{L}/2,\mathfrak{U}^{-1}/2,(\mathfrak{B}_2-\mathfrak{B}_1)/4\}$, and all $t\in\mathscr{T}$, the map $\Phi_{u}^{\mathrm{C}}:[-\pi,\pi)\to\mathbb{C}$ defined by
\begin{equation*}
\Phi_{u}^{\mathrm{C}}(\theta)\overset{\textnormal{def}}{=}\mathrm{u}+\mathrm{R}_{\mathrm{C}}\bigl((1-\cos\theta)+\mathrm{i}\sin\theta\bigr)
\end{equation*}
is injective, traces the circle of radius $\mathrm{R}_{\mathrm{C}}$ intersecting $\mathrm{u}$ and $\mathrm{u}+2\mathrm{R}_{\mathrm{C}}$ on the real axis, centred at $\mathrm{u}+\mathrm{R}_{\mathrm{C}}$, and is a descent curve for $\mathfrak{F}_{y}(\cdot,t)$ in the domain
\begin{equation*}
\Omega\overset{\textnormal{def}}{=}
\begin{cases}
\mathbb{C}\setminus\{\lambda_n:n\in\mathbb{Z}_+\},&\textnormal{if }\mathscr{T}=\mathbb{R}_+,\\
\mathbb{C}\setminus\Big(\{\lambda_n:n\in\mathbb{Z}_+\}\cup\{\alpha_n^{-1}:\tau_n=0\}\cup\{-\beta_n^{-1}:\tau_n=1\}\Big),&\textnormal{if }\mathscr{T}=\mathbb{Z}_+,
\end{cases}
\end{equation*}
from $x=\mathrm{u}$ with $\phi=0$, in the sense of Definition~\ref{descentcurvedefn}.
\begin{proof}
	Fix $t\in\mathscr{T}$ and $y\in\mathscr{P}$. Write $\Phi(\theta)\overset{\textnormal{def}}{=}\Phi_{u}^{\mathrm{C}}(\theta)=\phi(\theta)+\mathrm{i}\psi(\theta)$, where
	\begin{equation*}
	\phi(\theta)=\mathrm{u}+\mathrm{R}_{\mathrm{C}}(1-\cos\theta)\quad\textnormal{and}\quad\psi(\theta)=\mathrm{R}_{\mathrm{C}}\sin\theta.
	\end{equation*}

 The map $\theta\mapsto(1-\cos\theta,\sin\theta)$ is injective on $[-\pi,\pi)$, since the only pair of distinct angles in $[-\pi,\pi]$ with the same image is $(-\pi,\pi)$, and $\pi$ is excluded. Scaling by $\mathrm{R}_{\mathrm{C}}>0$ and translating by $\mathrm{u}$ preserves injectivity, so $\Phi$ is injective on $[-\pi,\pi)$.

     Then $\phi^{(1)}(\theta)=\mathrm{R}_{\mathrm{C}}\sin\theta$ and $\psi^{(1)}(\theta)=\mathrm{R}_{\mathrm{C}}\cos\theta$. Let $\{p_n\}_{n\in\mathbb{Z}_+}\subset[-\pi,\pi)$ denote the (possibly empty) set of parameters for which $\Phi(p_n)\notin\Omega$; equivalently, $\Phi(p_n)$ is one of the spatial poles, or, in discrete time, one of the temporal branch points $\alpha_n^{-1}$ or $-\beta_n^{-1}$. On $[-\pi,\pi)\setminus\{p_n\}$ we have $\Phi(\theta)\in\Omega$ and $\Phi(\theta)$ avoids every pole of the logarithmic terms in the definition of $\mathfrak{F}_y(\cdot,t)$, so $\theta\mapsto \mathfrak{Re}[\mathfrak{F}_y(\Phi(\theta),t)]$ is differentiable there. Differentiating term-by-term using $\frac{\mathrm{d}}{\mathrm{d}\theta}\log|\Phi(\theta)-a|=\mathfrak{Re}\left(\frac{\Phi^{(1)}(\theta)}{\Phi(\theta)-a}\right)$ yields
\begin{equation}\label{dReF-along-curve}
\begin{aligned}
    \frac{\mathrm{d}}{\mathrm{d}\theta}\mathfrak{Re}[\mathfrak{F}_y(\Phi(\theta),t)]=&-\mathbf{1}_{\{\mathscr{T}=\mathbb{R}_+\}}t\,\phi^{(1)}(\theta)+\sum_{k=0}^{y(t)}\mathfrak{Re}\left[\frac{\Phi^{(1)}(\theta)}{\lambda_k-\Phi(\theta)}\right]\\&-\mathbf{1}_{\{\mathscr{T}=\mathbb{Z}_+\}}\sum_{n=0}^{t-1}(1-\tau_n)\mathfrak{Re}\left[\frac{\Phi^{(1)}(\theta)}{\alpha_n^{-1}-\Phi(\theta)}\right]\\&-\mathbf{1}_{\{\mathscr{T}=\mathbb{Z}_+\}}\sum_{n=0}^{t-1}\tau_n\,\mathfrak{Re}\left[\frac{\Phi^{(1)}(\theta)}{\Phi(\theta)+\beta_n^{-1}}\right].
\end{aligned}
\end{equation}
In continuous time ($\mathscr{T}=\mathbb{R}_+$) the first bracketed contribution is $t\ge0$; in discrete time ($\mathscr{T}=\mathbb{Z}_+$) it vanishes and only the spatial and temporal sums remain.
For $a\in\mathbb{R}$, a direct computation shows
\begin{equation*}
    \begin{aligned}
        &\mathfrak{Re}\left[\frac{\Phi^{(1)}(\theta)}{a-\Phi(\theta)}\right]=\frac{(a-\phi(\theta))\phi^{(1)}(\theta)-\psi(\theta)\psi^{(1)}(\theta)}{|a-\Phi(\theta)|^2},\\&\mathfrak{Re}\left[\frac{\Phi^{(1)}(\theta)}{\Phi(\theta)+a}\right]=\frac{(\phi(\theta)+a)\phi^{(1)}(\theta)+\psi(\theta)\psi^{(1)}(\theta)}{|\Phi(\theta)+a|^2}.
    \end{aligned}
\end{equation*}
	Using $\phi^{(1)}(\theta)=\mathrm{R}_{\mathrm{C}}\sin\theta$ and $\psi(\theta)\psi^{(1)}(\theta)=\mathrm{R}_{\mathrm{C}}^2\sin\theta\cos\theta$ and the identity $\phi(\theta)+\mathrm{R}_{\mathrm{C}}\cos\theta=\mathrm{u}+\mathrm{R}_{\mathrm{C}}$, we obtain
\begin{equation*}
    \begin{aligned}
	        &\mathfrak{Re}\left[\frac{\Phi^{(1)}(\theta)}{\lambda_k-\Phi(\theta)}\right]=-\mathrm{R}_{\mathrm{C}}\sin\theta\,\frac{\mathrm{u}+\mathrm{R}_{\mathrm{C}}-\lambda_k}{|\lambda_k-\Phi(\theta)|^2},\\&\mathfrak{Re}\left[\frac{\Phi^{(1)}(\theta)}{\alpha_n^{-1}-\Phi(\theta)}\right]=\mathrm{R}_{\mathrm{C}}\sin\theta\,\frac{\alpha_n^{-1}-\mathrm{u}-\mathrm{R}_{\mathrm{C}}}{|\alpha_n^{-1}-\Phi(\theta)|^2},\\& \mathfrak{Re}\left[\frac{\Phi^{(1)}(\theta)}{\Phi(\theta)+\beta_n^{-1}}\right]=\mathrm{R}_{\mathrm{C}}\sin\theta\,\frac{\mathrm{u}+\mathrm{R}_{\mathrm{C}}+\beta_n^{-1}}{|\Phi(\theta)+\beta_n^{-1}|^2}.
    \end{aligned}
\end{equation*}
Substituting into \eqref{dReF-along-curve} gives, for all $\theta\in[-\pi,\pi)\setminus\{p_n\}$,
\begin{equation}\label{factorised-derivative}
    \begin{aligned}
	        \frac{\mathrm{d}}{\mathrm{d}\theta}\mathfrak{Re}[\mathfrak{F}_y(\Phi(\theta),t)]=&-\mathrm{R}_{\mathrm{C}}\sin\theta\Bigg[\mathbf{1}_{\{\mathscr{T}=\mathbb{R}_+\}}t+\sum_{k=0}^{y(t)}\frac{\mathrm{u}+\mathrm{R}_{\mathrm{C}}-\lambda_k}{|\lambda_k-\Phi(\theta)|^2}\\&+\mathbf{1}_{\{\mathscr{T}=\mathbb{Z}_+\}}\sum_{n=0}^{t-1}\left((1-\tau_n)\frac{\alpha_n^{-1}-\mathrm{u}-\mathrm{R}_{\mathrm{C}}}{|\alpha_n^{-1}-\Phi(\theta)|^2}+\tau_n\frac{\mathrm{u}+\mathrm{R}_{\mathrm{C}}+\beta_n^{-1}}{|\Phi(\theta)+\beta_n^{-1}|^2}\right)\Bigg].
    \end{aligned}
\end{equation}
By \textbf{(A1)} we have $\mathfrak{B}_1<\mathfrak{B}_2$. Combining this gap with $|\mathrm{u}|<(\mathfrak{B}_2-\mathfrak{B}_1)/4$ gives
\begin{equation*}
	\mathrm{u}+\mathrm{R}_{\mathrm{C}}-\lambda_k\ge \mathrm{u}+\frac{\mathfrak{B}_1+\mathfrak{B}_2}{2}-\mathfrak{B}_1>\frac{\mathfrak{B}_2-\mathfrak{B}_1}{4},
\end{equation*}
for all $k$, using $\lambda_k\le\mathfrak{B}_1$. Moreover, on $\{\tau_n=0\}$ we have $\alpha_n\le\mathfrak{B}_2^{-1}$ and hence $\alpha_n^{-1}\ge\mathfrak{B}_2$, so
\begin{equation*}
	\alpha_n^{-1}-\mathrm{u}-\mathrm{R}_{\mathrm{C}}\ge \mathfrak{B}_2-\mathrm{u}-\frac{\mathfrak{B}_1+\mathfrak{B}_2}{2}>\frac{\mathfrak{B}_2-\mathfrak{B}_1}{4},
\end{equation*}
and on $\{\tau_n=1\}$ we have $\beta_n\le\mathfrak{U}$ and hence $\beta_n^{-1}\ge\mathfrak{U}^{-1}>0$. Also, using $|\mathrm{u}|<(\mathfrak{B}_2-\mathfrak{B}_1)/4$ and $\mathrm{R}_{\mathrm{C}}=(\mathfrak{B}_1+\mathfrak{B}_2)/2$,
\begin{equation*}
	\mathrm{u}+\mathrm{R}_{\mathrm{C}}>-\frac{\mathfrak{B}_2-\mathfrak{B}_1}{4}+\frac{\mathfrak{B}_1+\mathfrak{B}_2}{2}
	=\frac{3\mathfrak{B}_1+\mathfrak{B}_2}{4}>0,
\end{equation*}
so $\mathrm{u}+\mathrm{R}_{\mathrm{C}}+\beta_n^{-1}>0$ for every $n$. Consequently, the bracketed term in \eqref{factorised-derivative} is strictly positive for all $\theta\in[-\pi,\pi)\setminus\{p_n\}$, and therefore the sign of $\frac{\mathrm{d}}{\mathrm{d}\theta}\mathfrak{Re}[\mathfrak{F}_y(\Phi(\theta),t)]$ is determined by $-\sin\theta$. In particular,
\begin{equation*}
    \begin{aligned}
        &\frac{\mathrm{d}}{\mathrm{d}\theta}\mathfrak{Re}[\mathfrak{F}_y(\Phi(\theta),t)]<0\quad\textnormal{for all }\theta\in(0,\pi)\setminus\{p_n\},\\& \frac{\mathrm{d}}{\mathrm{d}\theta}\mathfrak{Re}[\mathfrak{F}_y(\Phi(\theta),t)]>0\quad\textnormal{for all }\theta\in(-\pi,0)\setminus\{p_n\}.
    \end{aligned}
\end{equation*}
Since $\Phi(0)=\mathrm{u}$, these inequalities are exactly the descent-curve inequalities from Definition~\ref{descentcurvedefn} with $\phi=0$ and $x=\mathrm{u}$. Therefore $\Phi_{u}^{\mathrm{C}}$ is a descent curve for $\mathfrak{F}_y(\cdot,t)$ in $\Omega$, from $\mathrm{u}$.
\end{proof}
\end{proposition}

 We will use a segment of $\Phi_{u}^{\mathrm{C}}$ as the integration contour when evaluating \eqref{cir} at the fluctuation scale; it is spliced to the local near-steepest-descent segment from Proposition~\ref{descentcurve} in Corollary~\ref{steepestdescentpath}. However, the subsequent analysis requires a two-term expansion near the stationary point. This means that the local contour must match the steepest-descent contour beyond its tangent direction. The curve $\Phi_{u}^{\mathrm{C}}$ matches the first-order behaviour in the precise sense that, after centering at the stationary point, its linear term is the vertical displacement $\mathrm{i}\theta$. To capture the next correction to the imaginary level set of the phase, we must examine the derivatives $\mathfrak{F}_{y}^{(n)}(\mathrm{u}_{y}(t),t)$ for $n\geq 3$. The endpoint estimates control these derivatives in magnitude, but they do not impose any uniform cancellation pattern on the odd-order derivatives. Therefore, for $y \in \mathscr{P}$ we define
    \begin{equation}\label{defineofN}
        \mathfrak{N}_{y}(t)\overset{\mathrm{def}}{=}\inf\left\{n\in\mathbb{Z}_+:\mathfrak{F}_{y}^{(2n+1)}(\mathrm{u}_{y}(t),t)\ne 0\right\},
    \end{equation}
    with the convention that the infimum is $\infty$ if the set is empty. This quantity identifies the lowest-order non-vanishing odd derivative of $\mathfrak{F}_{y}(\cdot,t)$ at the stationary point. It allows us to track the leading contribution of the odd part of the expansion and, in particular, to quantify the deviation from the exact steepest-descent contour.
    \begin{proposition}\label{descentcurve}
       There exist constants $\epsilon_{\mathrm{SD}}>0$ and $\mathrm{T}_{\mathrm{SD}}\in\mathscr{T}$ such that, for every $t>\mathrm{T}_{\mathrm{SD}}$ and every auxiliary path $y\in\mathscr{P}$ satisfying $y(t)\in\mathscr{C}_{t}$, the following curve, defined for $|\theta|<\epsilon_{\mathrm{SD}}$, is a descent curve for $\mathfrak{F}_{y}(\cdot,t)$ in the ball $\mathcal{B}_{\epsilon_{\mathrm{SD}}}(\mathrm{u}_{y}(t))$, from the stationary point $\mathrm{u}_y(t)$ constructed in Proposition~\ref{axpro}:
       \begin{equation*}
           \Phi_{y}^{\mathrm{SD}}(\theta,t)=
           \begin{cases}
               \mathrm{u}_{y}(t)+\mathrm{i}\theta+\dfrac{\mathfrak{a}_y(t)}{(2\mathfrak{N}_y(t))!}\theta^{2\mathfrak{N}_{y}(t)},&\text{if }\mathfrak{N}_y(t)<\infty,\\
               \mathrm{u}_{y}(t)+\mathrm{i}\theta,&\text{if }\mathfrak{N}_y(t)=\infty.
           \end{cases}
       \end{equation*}
       Here the coefficient $\mathfrak{a}_y(t)$ is defined by
            \begin{equation*}
                \mathfrak{a}_y(t)\overset{\textnormal{def}}{=}
                \begin{cases}
                    \dfrac{(-1)^{\mathfrak{N}_{y}(t)+1}}{(2\mathfrak{N}_{y}(t)+1)}\dfrac{\mathfrak{F}_{y}^{(2\mathfrak{N}_{y}(t)+1)}(\mathrm{u}_{y}(t),t)}{\mathfrak{F}_{y}^{(2)}(\mathrm{u}_{y}(t),t)},&\text{if }\mathfrak{N}_{y}(t)<\infty,\\
	                    0,&\text{if }\mathfrak{N}_{y}(t)=\infty.
	                \end{cases}
	            \end{equation*}
	            Moreover, writing
	            \begin{equation*}
	                \Phi_y^{\mathrm{SD}}(\theta,t)-\mathrm{u}_y(t)=\mathrm{i}\theta+\chi_y(\theta,t),
	            \end{equation*}
	            where $\chi_y(\theta,t)=\mathfrak{a}_y(t)\theta^{2\mathfrak{N}_y(t)}/(2\mathfrak{N}_y(t))!$ if $\mathfrak{N}_y(t)<\infty$ and $\chi_y(\theta,t)=0$ if $\mathfrak{N}_y(t)=\infty$, and where $\chi_y^{(1)}$ denotes differentiation with respect to $\theta$, the same constants may be chosen so that, for all $|\theta|\leq\mathrm{M}_{2}^{-1}$,
	            \begin{equation}\label{eq:SD-chi-bound}
	                |\chi_y(\theta,t)|\leq \mathrm{K}_{\mathrm{SD}}|\theta|^2,\qquad |\chi_y^{(1)}(\theta,t)|\leq \mathrm{K}_{\mathrm{SD}}|\theta|,\qquad \mathrm{K}_{\mathrm{SD}}\overset{\textnormal{def}}{=}2\mathrm{M}_{1}\mathrm{M}_{2}^{3},
	            \end{equation}
	            where $\mathrm{M}_{1}$ and $\mathrm{M}_{2}$ are the constants from \textbf{(iii)} of Proposition~\ref{specificcasesub}.
	            \begin{proof}
	                We prove the result with constants independent of the choice of path $y\in\mathscr{P}$ satisfying $y(t)\in\mathscr{C}_{t}$. The relevant uniformity is endpoint uniformity: the phase $\mathfrak{F}_{y}(\cdot,t)$ depends on the auxiliary path $y$ only through the terminal value $y(t)$, and therefore all estimates below are estimates on the finite-time window $\mathscr{C}_{t}$. By Proposition~\ref{axpro}, after increasing the lower time to dominate $\mathrm{T}_{\mathrm{A}}$, the stationary point $\mathrm{u}_{y}(t)$ exists uniformly over this window. To see that it lies in a fixed small neighbourhood of the origin, we combine the endpoint estimates from Proposition~\ref{specificcasesub} with the quadratic estimate \eqref{secondconditiontosat}. The endpoint first-derivative estimate in \textbf{(ii)} gives, uniformly for all such paths,
	                \begin{equation*}
	                    |\mathfrak{F}_{y}^{(1)}(0,t)|\leq \left(\mathfrak{L}^{-1}\mathrm{C}_{1}^{-1}+2\right)\gamma(t)^{1/2}\mathscr{Z}(t)+\frac{3\Lambda}{2}.
	                \end{equation*}
                Since $\gamma(t)\geq\Sigma(t)^{-1/2}$ by Definition~\ref{gammadef} and $\Sigma(t)\to\infty$ by \textbf{(iv)} of Proposition~\ref{specificcasesub}, we increase the lower time so that $\Sigma(t)\geq1$ and hence $\Sigma(t)^{-2}\leq\gamma(t)^{1/2}$. Combining this with $\mathscr{Z}(t)\leq\mathrm{C}_{2}\Sigma(t)^2$ and the denominator comparison $\mathfrak{F}_{y}^{(2)}(0,t)\geq\Sigma(t)^2/2$ from \textbf{(iv)} on the endpoint window gives a deterministic constant $\mathrm{D}_{\mathrm{SD}}>0$ such that
                \begin{equation*}
                    \left|\frac{\mathfrak{F}_{y}^{(1)}(0,t)}{\mathfrak{F}_{y}^{(2)}(0,t)}\right|\leq \mathrm{D}_{\mathrm{SD}}\gamma(t)^{1/2}.
                \end{equation*}
                Increase the lower time once more so that
                \begin{equation*}
                    \mathrm{C}_{\mathrm{A}}\mathrm{D}_{\mathrm{SD}}\gamma(t)^{1/2}\leq1.
                \end{equation*}
                Then the quadratic estimate \eqref{secondconditiontosat} implies
                \begin{equation*}
                    |\mathrm{u}_{y}(t)|\leq2\mathrm{D}_{\mathrm{SD}}\gamma(t)^{1/2}.
                \end{equation*}
                Increasing the lower time again, using $\gamma(t)\to0$, we obtain a threshold $\mathrm{T}_{\mathrm{SD},0}\in\mathscr{T}$ such that, whenever $t>\mathrm{T}_{\mathrm{SD},0}$ and $y(t)\in\mathscr{C}_{t}$, the stationary point lies in $\mathcal{B}_{\mathrm{R}_{2}/8}(0)$ and satisfies
                \begin{equation}\label{eq:SD-second-lower}
                    \mathfrak{F}_{y}^{(2)}(\mathrm{u}_{y}(t),t)\geq\frac{1}{2}\Sigma(t)^2.
                \end{equation}
	                Indeed, the comparison in \textbf{(iv)} of Proposition~\ref{specificcasesub} applied at $w=\mathrm{u}_{y}(t)$ gives \eqref{eq:SD-second-lower}, because $\mathrm{T}_{\mathrm{SD},0}$ has been chosen so that $\mathrm{C}_{\mathscr{C}}(\gamma(t)^{1/2}+|\mathrm{u}_{y}(t)|)\leq1/2$ on the endpoint window. We next bound the correction term in $\Phi_{y}^{\mathrm{SD}}$. Let $\mathrm{M}_1,\mathrm{M}_2$ be the endpoint constants from \textbf{(iii)} of Proposition~\ref{specificcasesub}. We choose $\mathrm{T}_{\mathrm{SD},0}\geq\mathrm{T}_{2,\mathscr{C}}$, so that the endpoint derivative estimate in \textbf{(iii)} is available on the same time range. Together with the lower bound \eqref{eq:SD-second-lower}, this implies that, for $t>\mathrm{T}_{\mathrm{SD},0}$, $y(t)\in\mathscr{C}_{t}$, $|\theta|\leq\mathrm{M}_{2}^{-1}$ and $\mathfrak{N}_{y}(t)<\infty$,
                \begin{equation}\label{estimateforaa}
                    \left|
                    \frac{\mathfrak{a}_y(t)}{(2\mathfrak{N}_y(t))!}\theta^{2\mathfrak{N}_y(t)}
                    \right|
                    \leq
                    2\mathrm{M}_{1}\mathrm{M}_{2}\left(\mathrm{M}_{2}|\theta|\right)^{2\mathfrak{N}_y(t)}
                    \leq
                    2\mathrm{M}_{1}\mathrm{M}_{2}^{2}|\theta|.
                \end{equation}
                The first inequality is just the definition of $\mathfrak{a}_y(t)$, the lower bound \eqref{eq:SD-second-lower}, and the derivative estimate
                \begin{equation*}
                    \left|\mathfrak{F}_{y}^{(2\mathfrak{N}_{y}(t)+1)}(\mathrm{u}_{y}(t),t)\right|
                    \leq
                    (2\mathfrak{N}_{y}(t))!\mathrm{M}_{1}\mathrm{M}_{2}^{2\mathfrak{N}_{y}(t)+1}\Sigma(t)^2.
	                \end{equation*}
	                The second inequality in \eqref{estimateforaa} uses $\mathfrak{N}_{y}(t)\geq1$ and $\mathrm{M}_{2}|\theta|\leq1$. If $\mathfrak{N}_{y}(t)=\infty$, then the correction term is zero by definition, so \eqref{estimateforaa} holds with left-hand side equal to zero. We next prove the sharper deviation estimate \eqref{eq:SD-chi-bound}, because this is the form needed later when the local contour is expanded in Proposition~\ref{Thm2}. When $\mathfrak{N}_{y}(t)=\infty$, both inequalities in \eqref{eq:SD-chi-bound} are immediate from $\chi_y\equiv0$. Assume therefore that $\mathfrak{N}_{y}(t)<\infty$. The first inequality in \eqref{estimateforaa} gives
	                \begin{equation*}
	                    |\chi_y(\theta,t)|\leq2\mathrm{M}_1\mathrm{M}_2(\mathrm{M}_2|\theta|)^{2\mathfrak{N}_y(t)}.
	                \end{equation*}
	                Since $\mathfrak{N}_{y}(t)\geq1$ and $\mathrm{M}_{2}|\theta|\leq1$, the elementary inequality $r^{2\mathfrak{N}_{y}(t)}\leq r^2$ for $0\leq r\leq1$ gives
	                \begin{equation*}
	                    |\chi_y(\theta,t)|\leq2\mathrm{M}_{1}\mathrm{M}_{2}^{3}|\theta|^2.
	                \end{equation*}
	                For the derivative estimate, differentiating the explicit expression for $\chi_y$ gives
	                \begin{equation*}
	                    \left|\chi_y^{(1)}(\theta,t)\right|=\left|\frac{2\mathfrak{N}_y(t)\mathfrak{a}_y(t)}{(2\mathfrak{N}_y(t))!}\theta^{2\mathfrak{N}_y(t)-1}\right|.
	                \end{equation*}
		                Substituting the definition of $\mathfrak{a}_{y}(t)$ gives a factor $2\mathfrak{N}_{y}(t)/(2\mathfrak{N}_{y}(t)+1)$, which is at most one. Applying the same derivative estimate used in \eqref{estimateforaa}, and then applying the lower curvature bound \eqref{eq:SD-second-lower}, yields
		                \begin{equation*}
			                    \left|\chi_y^{(1)}(\theta,t)\right|\leq\frac{(2\mathfrak{N}_{y}(t))!\mathrm{M}_{1}\mathrm{M}_{2}^{2\mathfrak{N}_{y}(t)+1}\Sigma(t)^2}{(2\mathfrak{N}_{y}(t))!\frac{1}{2}\Sigma(t)^2}|\theta|^{2\mathfrak{N}_y(t)-1}=2\mathrm{M}_{1}\mathrm{M}_{2}^{2}(\mathrm{M}_{2}|\theta|)^{2\mathfrak{N}_y(t)-1}.
		                \end{equation*}
	                Again $\mathfrak{N}_{y}(t)\geq1$ and $\mathrm{M}_{2}|\theta|\leq1$, so $r^{2\mathfrak{N}_{y}(t)-1}\leq r$ for $0\leq r\leq1$, and therefore
	                \begin{equation*}
	                    \left|\chi_y^{(1)}(\theta,t)\right|\leq2\mathrm{M}_{1}\mathrm{M}_{2}^{3}|\theta|.
	                \end{equation*}
		                This proves \eqref{eq:SD-chi-bound}. Since $\mathrm{M}_{2}\geq1$ by the definition of $\mathrm{M}_{2}$ in Proposition~\ref{specificcasesub}, the condition $|\theta|\leq\mathrm{M}_{2}^{-1}$ implies $|\theta|\leq1$. Thus, with
	                \begin{equation*}
	                    \mathrm{A}_{\mathrm{SD}}\overset{\textnormal{def}}{=}1+\mathrm{K}_{\mathrm{SD}},
	                \end{equation*}
	                we have
	                \begin{equation}\label{eq:SD-path-stays-close}
	                    \left|\Phi_{y}^{\mathrm{SD}}(\theta,t)-\mathrm{u}_y(t)\right|
                    \leq\mathrm{A}_{\mathrm{SD}}|\theta|,
                    \qquad
                    t>\mathrm{T}_{\mathrm{SD},0},\quad y(t)\in\mathscr{C}_{t},\quad |\theta|\leq\mathrm{M}_{2}^{-1}.
                \end{equation}
It remains to choose a radius on which the real part decreases. Since $\mathfrak{F}_{y}$ is holomorphic on $\mathcal{B}_{\mathrm{R}_{-}}(0)$ by \textbf{(i)} of Proposition~\ref{specificcasesub} and $\mathrm{u}_{y}(t)\in\mathcal{B}_{\mathrm{R}_{2}/8}(0)$, the path segment stays inside $\mathcal{B}_{\mathrm{R}_{2}}(0)$ whenever $\mathrm{A}_{\mathrm{SD}}|\theta|\leq3\mathrm{R}_{2}/8$. On this range the derivative bounds in \textbf{(iii)} control the Taylor remainder for
                \begin{equation*}
                    \theta\longmapsto\mathfrak{Re}\left[\mathfrak{F}_{y}\left(\Phi_{y}^{\mathrm{SD}}(\theta,t),t\right)\right]
                \end{equation*}
	                around $\theta=0$. To see the leading term without hiding the estimates, write
	                \begin{equation*}
	                    \Phi_{y}^{\mathrm{SD}}(\theta,t)-\mathrm{u}_y(t)=\mathrm{i}\theta+\chi_y(\theta,t),
	                \end{equation*}
	                where $\chi_y(\theta,t)$ is the correction term from \eqref{eq:SD-chi-bound}. Since the Taylor analysis is restricted to $|\theta|\leq1$, \eqref{eq:SD-chi-bound} implies both $|\Phi_{y}^{\mathrm{SD}}(\theta,t)-\mathrm{u}_y(t)|\leq\mathrm{A}_{\mathrm{SD}}|\theta|$ and $|\frac{\mathrm{d}}{\mathrm{d}\theta}\Phi_{y}^{\mathrm{SD}}(\theta,t)|\leq\mathrm{A}_{\mathrm{SD}}$. Since $\mathfrak{F}_{y}^{(1)}(\mathrm{u}_y(t),t)=0$, the first contribution to the derivative of the real part is the quadratic term
                \begin{equation*}
                    \mathfrak{Re}\left[
                    \mathfrak{F}_{y}^{(2)}(\mathrm{u}_y(t),t)
                    \left(\Phi_{y}^{\mathrm{SD}}(\theta,t)-\mathrm{u}_y(t)\right)
                    \frac{\mathrm{d}}{\mathrm{d}\theta}\Phi_{y}^{\mathrm{SD}}(\theta,t)
                    \right]
                    =
                    -\mathfrak{F}_{y}^{(2)}(\mathrm{u}_{y}(t),t)\theta+\mathrm{E}_{\mathrm{quad}}(\theta,t),
                \end{equation*}
	                where the correction is $\mathrm{E}_{\mathrm{quad}}(\theta,t)=\mathfrak{F}_{y}^{(2)}(\mathrm{u}_{y}(t),t)\chi_y(\theta,t)\chi_y^{(1)}(\theta,t)$ and therefore satisfies $|\mathrm{E}_{\mathrm{quad}}(\theta,t)|\leq\mathrm{K}_{\mathrm{SD}}^{2}\mathfrak{F}_{y}^{(2)}(\mathrm{u}_{y}(t),t)|\theta|^2$ on the subinterval $|\theta|\leq1$, by \eqref{eq:SD-chi-bound}. We also include in $\mathrm{T}_{\mathrm{SD},0}$ the condition giving $\mathfrak{F}_{y}^{(2)}(\mathrm{u}_{y}(t),t)\leq2\Sigma(t)^2$ from the endpoint second-derivative comparison in \textbf{(iv)}. The remaining Taylor terms have degree at least three in $\Phi_{y}^{\mathrm{SD}}(\theta,t)-\mathrm{u}_y(t)$; applying Taylor's theorem with integral remainder, \eqref{eq:SD-path-stays-close}, the bound on $\frac{\mathrm{d}}{\mathrm{d}\theta}\Phi_{y}^{\mathrm{SD}}$, and the uniform derivative estimate in \textbf{(iii)}, gives a deterministic constant
                \begin{equation*}
	                    \mathrm{B}_{\mathrm{SD}}\overset{\textnormal{def}}{=}8\mathrm{M}_{1}\mathrm{M}_{2}^{3}\mathrm{A}_{\mathrm{SD}}^{3}+2\mathrm{K}_{\mathrm{SD}}^{2}
                \end{equation*}
                such that, for $t>\mathrm{T}_{\mathrm{SD},0}$, $y(t)\in\mathscr{C}_{t}$ and $|\theta|\leq\min(1,\mathrm{M}_{2}^{-1},3\mathrm{R}_{2}/(8\mathrm{A}_{\mathrm{SD}}))$,
                \begin{equation}\label{taylorbyrem}
                    \frac{\mathrm{d}}{\mathrm{d}\theta}\mathfrak{Re}\left[\mathfrak{F}_{y}\left(\Phi_{y}^{\mathrm{SD}}(\theta,t),t\right)\right]
                    =-\mathfrak{F}_{y}^{(2)}(\mathrm{u}_{y}(t),t)\theta+\mathrm{R}_{\mathrm{SD}}(\theta,t),
                    \qquad
                    |\mathrm{R}_{\mathrm{SD}}(\theta,t)|\leq\mathrm{B}_{\mathrm{SD}}\Sigma(t)^2\theta^2.
                \end{equation}
                The coefficient $\mathfrak{a}_{y}(t)$ is the correction that cancels the first non-zero odd contribution to the imaginary level of the phase. For the descent property, estimate \eqref{taylorbyrem} is the decisive input: after differentiating the real part along the curve, the leading term is the negative quadratic contribution and every other contribution is bounded by a constant times $\Sigma(t)^2\theta^2$, uniformly over $y(t)\in\mathscr{C}_{t}$. We now fix the radius once and for all by setting
                \begin{equation*}
                    \epsilon_{\mathrm{SD}}\overset{\textnormal{def}}{=}
                    \min\left(1,\frac{1}{\mathrm{M}_{2}},\frac{3\mathrm{R}_{2}}{8\mathrm{A}_{\mathrm{SD}}},\frac{1}{4\mathrm{B}_{\mathrm{SD}}}\right).
                \end{equation*}
                If $0<\theta\leq\epsilon_{\mathrm{SD}}$, then \eqref{eq:SD-second-lower} and \eqref{taylorbyrem} give
                \begin{equation*}
                    \frac{\mathrm{d}}{\mathrm{d}\theta}\mathfrak{Re}\left[\mathfrak{F}_{y}\left(\Phi_{y}^{\mathrm{SD}}(\theta,t),t\right)\right]
                    \leq
                    -\frac{1}{2}\Sigma(t)^2\theta+\frac{1}{4}\Sigma(t)^2\theta
                    <0.
                \end{equation*}
                If $-\epsilon_{\mathrm{SD}}\leq\theta<0$, then \eqref{eq:SD-second-lower} and \eqref{taylorbyrem} instead give
                \begin{equation*}
                    \frac{\mathrm{d}}{\mathrm{d}\theta}\mathfrak{Re}\left[\mathfrak{F}_{y}\left(\Phi_{y}^{\mathrm{SD}}(\theta,t),t\right)\right]
                    \geq
                    \frac{1}{2}\Sigma(t)^2|\theta|-\frac{1}{4}\Sigma(t)^2|\theta|
                    >0,
                \end{equation*}
                because now $-\mathfrak{F}_{y}^{(2)}(\mathrm{u}_{y}(t),t)\theta$ is positive and the remainder is bounded by $\frac{1}{4}\Sigma(t)^2|\theta|$. Therefore
                \begin{equation*}
                     \frac{\mathrm{d}}{\mathrm{d}\theta}\mathfrak{Re}\left[\mathfrak{F}_{y}\left(\Phi_{y}^{\mathrm{SD}}(\theta,t),t\right)\right]<0,\quad\text{for all }\theta\in(0,\epsilon_{\mathrm{SD}}),
                \end{equation*}
                and
                \begin{equation*}
                     \frac{\mathrm{d}}{\mathrm{d}\theta}\mathfrak{Re}\left[\mathfrak{F}_{y}\left(\Phi_{y}^{\mathrm{SD}}(\theta,t),t\right)\right]>0,\quad\text{for all }\theta\in(-\epsilon_{\mathrm{SD}},0).
                \end{equation*}
	                These two inequalities give the descent property in $\mathcal{B}_{\epsilon_{\mathrm{SD}}}(\mathrm{u}_{y}(t))$. Taking $\mathrm{T}_{\mathrm{SD}}\overset{\textnormal{def}}{=}\mathrm{T}_{\mathrm{SD},0}$ completes the proof.
	            \end{proof}
    \end{proposition}
    
     Since $\Phi_{y}^{\mathrm{SD}}$ is known to be a descent curve only in a neighbourhood of the stationary point, we splice it with $\Phi_{u}^{\mathrm{C}}$ to obtain the global contour used in the steepest-descent analysis. The statement is formulated on the endpoint window $\mathscr{C}_{t}$, because this is the only range used in Proposition~\ref{Thm2} and it is the range on which the constants are uniform.
    \begin{cor}\label{steepestdescentpath}
        There exist constants $\theta_{\mathrm{C}}>0$ and $\mathrm{T}_{\mathrm{C}}\in\mathscr{T}$ such that, for every $t>\mathrm{T}_{\mathrm{C}}$, every auxiliary path $y\in\mathscr{P}$ satisfying $y(t)\in\mathscr{C}_{t}$, and every $\theta\in(0,\theta_{\mathrm{C}})$, the following concatenation of $\Phi_{y}^{\mathrm{SD}}(\cdot,t)$ and $\Phi_{\mathrm{c}(\theta,t)}^{\mathrm{C}}$ is an injective descent curve for $\mathfrak{F}_{y}(\cdot,t)$ in the logarithmic domain $\Omega$ from Proposition~\ref{contour1}, from the stationary point $\mathrm{u}_{y}(t)$. Set
        \begin{equation*}
            \phi_{\theta}\overset{\textnormal{def}}{=}\mathrm{R}_{\mathrm{C}}\sin\theta,
            \qquad
            U_{\theta}\overset{\textnormal{def}}{=}\frac{\mathrm{R}_{\mathrm{C}}\sin\theta}{\theta}\pi,
        \end{equation*}
        where $\mathrm{R}_{\mathrm{C}}=(\mathfrak{B}_{1}+\mathfrak{B}_{2})/2$ is the circular radius from Proposition~\ref{contour1}, and define $\Phi_y(\cdot,t,\theta):[-U_{\theta},U_{\theta})\to\mathbb{C}$ by
        \begin{equation*}
            \Phi_y(u,t,\theta)\overset{\textnormal{def}}{=}
            \begin{cases}
                \Phi_{\mathrm{c}(\theta,t)}^{\mathrm{C}}\!\left(\dfrac{\theta u}{\mathrm{R}_{\mathrm{C}}\sin\theta}\right),& -U_{\theta}\leq u<-\phi_{\theta},\\
                \Phi_y^{\mathrm{SD}}(u,t),& -\phi_{\theta}\leq u\leq\phi_{\theta},\\
                \Phi_{\mathrm{c}(\theta,t)}^{\mathrm{C}}\!\left(\dfrac{\theta u}{\mathrm{R}_{\mathrm{C}}\sin\theta}\right),& \phi_{\theta}<u<U_{\theta}.
            \end{cases}
        \end{equation*}
        The corresponding closed contour is used in contour integrals with positive orientation, opposite to the increasing-$u$ parametrisation above. When a compact parametrisation is needed, we use the endpoint convention $\Phi_y(U_{\theta},t,\theta)=\Phi_y(-U_{\theta},t,\theta)$. Set
        \begin{equation*}
            \mathrm{c}(\theta,t)\overset{\textnormal{def}}{=}
            \begin{cases}
                \mathrm{u}_y(t)+\dfrac{\mathfrak{a}_y(t)}{(2\mathfrak{N}_y(t))!}\bigl(\mathrm{R}_{\mathrm{C}}\sin\theta\bigr)^{2\mathfrak{N}_y(t)}-\mathrm{R}_{\mathrm{C}}(1-\cos(\theta)),&\text{if }\mathfrak{N}_y(t)<\infty,\\
                \mathrm{u}_y(t)-\mathrm{R}_{\mathrm{C}}(1-\cos(\theta)),&\text{if }\mathfrak{N}_y(t)=\infty.
            \end{cases}
        \end{equation*}
        Moreover, with this positive orientation, the interior of $\Phi_y(\cdot,t,\theta)$ encloses the sequence of spatial rates $(\lambda_{n})_{n\in\mathbb{Z}_+}$.
        \begin{proof}
            The proof has two parts. First we check that the local and circular pieces meet exactly and that the local piece remains inside the neighbourhood from Proposition~\ref{descentcurve}. Then we choose $\theta_{\mathrm{C}}$ and $\mathrm{T}_{\mathrm{C}}$ so that the centre of the circular piece satisfies the hypothesis of Proposition~\ref{contour1}. By the definition of $\mathrm{c}(\theta,t)$, for every $\theta\in(0,\pi)$,
            \begin{equation*}
                \begin{aligned}
                    \Phi_{\mathrm{c}(\theta,t)}^{\mathrm{C}}(\theta)
                    &=
                    \mathrm{c}(\theta,t)+\mathrm{R}_{\mathrm{C}}\left(1-\cos(\theta)+\mathrm{i}\sin(\theta)\right)\\
                    &=
                    \begin{cases}
                        \mathrm{u}_y(t)+\mathrm{i}\mathrm{R}_{\mathrm{C}}\sin(\theta)+\dfrac{\mathfrak{a}_y(t)}{(2\mathfrak{N}_y(t))!}\left(\mathrm{R}_{\mathrm{C}}\sin(\theta)\right)^{2\mathfrak{N}_y(t)},&\text{if }\mathfrak{N}_y(t)<\infty,\\
                        \mathrm{u}_y(t)+\mathrm{i}\mathrm{R}_{\mathrm{C}}\sin(\theta),&\text{if }\mathfrak{N}_y(t)=\infty,
                    \end{cases}\\
                    &=
                    \Phi_{y}^{\mathrm{SD}}(\mathrm{R}_{\mathrm{C}}\sin\theta,t).
                \end{aligned}
            \end{equation*}
            Thus the two pieces join at the point corresponding to $u=\mathrm{R}_{\mathrm{C}}\sin\theta$; the same computation with $-\theta$ gives the join at $u=-\mathrm{R}_{\mathrm{C}}\sin\theta$. Since $\sin(\theta)\leq\theta$ for $\theta>0$, \eqref{eq:SD-path-stays-close} gives, for $0<u\leq\mathrm{R}_{\mathrm{C}}\sin\theta$ and $t>\mathrm{T}_{\mathrm{SD}}$,
            \begin{equation*}
                \left|\Phi_{y}^{\mathrm{SD}}(u,t)-\mathrm{u}_y(t)\right|
                \leq
                \mathrm{A}_{\mathrm{SD}}\mathrm{R}_{\mathrm{C}}\theta.
            \end{equation*}
            We therefore require $\theta\leq\epsilon_{\mathrm{SD}}/(\mathrm{A}_{\mathrm{SD}}\mathrm{R}_{\mathrm{C}})$, which ensures that the whole local piece lies in $\mathcal{B}_{\epsilon_{\mathrm{SD}}}(\mathrm{u}_{y}(t))$ and is a descent segment by Proposition~\ref{descentcurve}. It remains to control the centre $\mathrm{c}(\theta,t)$. By \eqref{estimateforaa}, $\mathrm{R}_{\mathrm{C}}\sin\theta\leq\mathrm{R}_{\mathrm{C}}\theta$, and $1-\cos(\theta)\leq\theta^2/2$, whenever $\theta\leq(\mathrm{M}_{2}\mathrm{R}_{\mathrm{C}})^{-1}$ we have
            \begin{equation}\label{eq:c-centre-bound}
                |\mathrm{c}(\theta,t)|
                \leq
                |\mathrm{u}_{y}(t)|+
                \left(2\mathrm{M}_{1}\mathrm{M}_{2}^{3}\mathrm{R}_{\mathrm{C}}^{2}+\frac{\mathrm{R}_{\mathrm{C}}}{2}\right)\theta^{2}.
            \end{equation}
            We now make the separation from the real singularities explicit. Set
            \begin{equation*}
                \mathrm{d}_{\mathrm{gap}}\overset{\textnormal{def}}{=}\frac{1}{8}\min\left(\mathfrak{L},\mathfrak{U}^{-1},\mathfrak{B}_{2}-\mathfrak{B}_{1}\right).
            \end{equation*}
            By \eqref{axpro-u-small}, the point $\mathrm{u}_{y}(t)$ is bounded by $2\mathrm{B}_{\mathrm{A}}\gamma(t)^{1/2}$ uniformly on $y(t)\in\mathscr{C}_{t}$. Since $\gamma(t)$ can be made arbitrarily small by increasing the deterministic lower time, the endpoint stationary-point construction used at the start of Proposition~\ref{descentcurve} gives a threshold $\mathrm{T}_{\mathrm{C},1}\in\mathscr{T}$ such that
            \begin{equation*}
                |\mathrm{u}_{y}(t)|\leq\mathrm{d}_{\mathrm{gap}},\qquad t>\mathrm{T}_{\mathrm{C},1},\quad y(t)\in\mathscr{C}_{t}.
            \end{equation*}
            Define
            \begin{equation*}
                \theta_{\mathrm{C}}\overset{\textnormal{def}}{=}
                \min\left(
                    \frac{\pi}{2},
                    \frac{1}{\mathrm{M}_{2}\mathrm{R}_{\mathrm{C}}},
                    \frac{\epsilon_{\mathrm{SD}}}{\mathrm{A}_{\mathrm{SD}}\mathrm{R}_{\mathrm{C}}},
                    \frac{\mathrm{d}_{\mathrm{gap}}}{\mathrm{A}_{\mathrm{SD}}\mathrm{R}_{\mathrm{C}}},
                    \left(\frac{\mathrm{d}_{\mathrm{gap}}}{2\mathrm{M}_{1}\mathrm{M}_{2}^{3}\mathrm{R}_{\mathrm{C}}^{2}+\mathrm{R}_{\mathrm{C}}/2}\right)^{1/2},
                    \left(\frac{2\mathrm{d}_{\mathrm{gap}}}{\mathrm{R}_{\mathrm{C}}}\right)^{1/2}
                \right).
            \end{equation*}
            Finally, set $\mathrm{T}_{\mathrm{C}}\overset{\textnormal{def}}{=}\max(\mathrm{T}_{\mathrm{SD}},\mathrm{T}_{\mathrm{C},1})$. Then \eqref{eq:c-centre-bound} implies
            \begin{equation*}
                |\mathrm{c}(\theta,t)|<2\mathrm{d}_{\mathrm{gap}}\leq\frac{\mathfrak{B}_{2}-\mathfrak{B}_{1}}{4},
                \qquad
                t>\mathrm{T}_{\mathrm{C}},\quad y(t)\in\mathscr{C}_{t},\quad 0<\theta<\theta_{\mathrm{C}}.
            \end{equation*}
            This verifies the gap condition required in Proposition~\ref{contour1}. Moreover, whenever
            \begin{equation*}
                \mathrm{R}_{\mathrm{C}}\sin\theta<|u|<\left(\frac{\mathrm{R}_{\mathrm{C}}\sin\theta}{\theta}\right)\pi,
            \end{equation*}
            the angle $\widetilde{\theta}\overset{\textnormal{def}}{=}\theta u/(\mathrm{R}_{\mathrm{C}}\sin\theta)$ lies in $(-\pi,-\theta)\cup(\theta,\pi)\subset[-\pi,\pi)$. Hence the circular segment is exactly the map $\Phi_{\mathrm{c}(\theta,t)}^{\mathrm{C}}(\widetilde{\theta})$ from Proposition~\ref{contour1} with centre $\mathrm{c}(\theta,t)$, and is therefore a descent contour on that segment. Since $\theta/(\mathrm{R}_{\mathrm{C}}\sin\theta)>0$ for $0<\theta<\pi$, this affine change of parameter preserves the sign of the derivative of the real part.
          
             We next record injectivity of the splice. The local piece is injective by Proposition~\ref{descentcurve}, and the two retained circular pieces are injective by Proposition~\ref{contour1}. Their only common points are the two splice endpoints. Indeed, the local piece has imaginary part in $[-\mathrm{R}_{\mathrm{C}}\sin\theta,\mathrm{R}_{\mathrm{C}}\sin\theta]$. On the retained circular arc, if $\theta<|\widetilde{\theta}|<\pi-\theta$, then $|\sin\widetilde{\theta}|>\sin\theta$, so the imaginary part is outside this interval. If instead $\pi-\theta\leq|\widetilde{\theta}|<\pi$, then, since $\theta<\theta_{\mathrm{C}}\leq\pi/2$,
            \begin{equation*}
                \mathfrak{Re}\left[\Phi_{\mathrm{c}(\theta,t)}^{\mathrm{C}}(\widetilde{\theta})\right]
                \geq \mathrm{c}(\theta,t)+\mathrm{R}_{\mathrm{C}}(1+\cos\theta)
                >2\mathrm{d}_{\mathrm{gap}},
            \end{equation*}
            where the last inequality uses $|\mathrm{c}(\theta,t)|<2\mathrm{d}_{\mathrm{gap}}$ and $4\mathrm{d}_{\mathrm{gap}}<\mathrm{R}_{\mathrm{C}}$. By contrast, the local real-part bound in \eqref{eq:SD-path-stays-close} and the choice of $\theta_{\mathrm{C}}$ give $|\mathfrak{Re}[\Phi_y^{\mathrm{SD}}(u,t)]|\leq2\mathrm{d}_{\mathrm{gap}}$ on the local piece. Thus the local piece does not meet the retained circular arc except at the splice endpoints. Combining this injectivity with the descent property of the local segment from Proposition~\ref{descentcurve} and the circular descent property above proves that $\Phi_y(\cdot,t,\theta)$ is a descent contour from $\mathrm{u}_{y}(t)$.

             It remains to justify the enclosure statement for the closed concatenation. This is the topology shown in Figure~\ref{contourconstructionfig}: the red circular part is a large contour around the spatial rates, and the blue segment replaces only a small left-hand arc near the stationary point. To make this precise, first consider the complete circle $\Phi_{\mathrm{c}(\theta,t)}^{\mathrm{C}}([-\pi,\pi))$, viewed as a closed curve with its positive boundary orientation when discussing winding numbers. Its left and right real-axis intersection points are $\mathrm{c}(\theta,t)$ and $\mathrm{c}(\theta,t)+2\mathrm{R}_{\mathrm{C}}$, so the real points lying inside the circle are exactly those between these two values. Since $|\mathrm{c}(\theta,t)|<2\mathrm{d}_{\mathrm{gap}}$, we have $-\mathfrak{U}^{-1}<\mathrm{c}(\theta,t)<\mathfrak{L}$, so the left intersection point is strictly to the right of every geometric pole $-\beta_n^{-1}<-\mathfrak{U}^{-1}$ and strictly to the left of every spatial rate $\lambda_n>\mathfrak{L}$. Moreover, $\mathrm{c}(\theta,t)+2\mathrm{R}_{\mathrm{C}}>\mathfrak{B}_{1}+\mathfrak{B}_{2}-2\mathrm{d}_{\mathrm{gap}}>\mathfrak{B}_{1}$, because $\mathrm{R}_{\mathrm{C}}=(\mathfrak{B}_{1}+\mathfrak{B}_{2})/2$ and $2\mathrm{d}_{\mathrm{gap}}<\mathfrak{B}_{2}$. Hence the complete circle encloses every spatial rate and no geometric pole; equivalently, with positive orientation, its winding number is one around each spatial rate and zero around each geometric pole.
           
             The spliced contour $\Phi_y(\cdot,t,\theta)$ is obtained from this complete circle by removing the small circular arc with $|\widetilde{\theta}|\leq\theta$ and inserting the local steepest-descent segment $\Phi_y^{\mathrm{SD}}$. Both curves have the same endpoints by the matching computation above. Moreover, both lie in the singularity-free strip $-\mathfrak{U}^{-1}<\mathfrak{Re}(w)<\mathfrak{L}$. Indeed, on the omitted circular arc,
            \begin{equation*}
                \left|\mathfrak{Re}\left[\Phi_{\mathrm{c}(\theta,t)}^{\mathrm{C}}(\widetilde{\theta})\right]\right|\leq|\mathrm{c}(\theta,t)|+\mathrm{R}_{\mathrm{C}}(1-\cos\widetilde{\theta})<2\mathrm{d}_{\mathrm{gap}}+\frac{\mathrm{R}_{\mathrm{C}}\theta^2}{2}\leq3\mathrm{d}_{\mathrm{gap}},
            \end{equation*}
            while on the blue local segment, \eqref{eq:SD-path-stays-close} and the choice of $\theta_{\mathrm{C}}$ give
            \begin{equation*}
                \left|\mathfrak{Re}\left[\Phi_y^{\mathrm{SD}}(u,t)\right]\right|\leq|\mathrm{u}_{y}(t)|+\left|\Phi_y^{\mathrm{SD}}(u,t)-\mathrm{u}_{y}(t)\right|\leq2\mathrm{d}_{\mathrm{gap}}.
            \end{equation*}
            Since $3\mathrm{d}_{\mathrm{gap}}<\min(\mathfrak{L},\mathfrak{U}^{-1})$, both bounds place the relevant replacement region inside the strip. This strip contains no spatial rates and no geometric poles. Because the strip is convex, the omitted circular arc can be deformed to the inserted local segment, with endpoints fixed, without crossing any of these real singularities. The replacement therefore preserves the winding numbers computed for the complete circle. Consequently, the interior of the concatenated contour encloses the spatial rates $(\lambda_n)_{n\in\mathbb{Z}_+}$ and excludes the geometric poles.
        \end{proof}
    \end{cor}
    
     Figure~\ref{contourconstructionfig} illustrates schematically how the local segment $\Phi_{y}^{\mathrm{SD}}$ from Proposition~\ref{descentcurve} meets the translated circular arc $\Phi_{\mathrm{c}(\theta,t)}^{\mathrm{C}}$ from Corollary~\ref{steepestdescentpath}, relative to the stationary point $\mathrm{u}_y(t)$ and the axis markers for the geometric poles $-\beta_n^{-1}$, the spatial rates $\lambda_n$, and the Bernoulli zeros $\alpha_n^{-1}$. The three panels correspond to the sign regimes $\mathfrak{a}_y(t)<0$, $\mathfrak{a}_y(t)=0$, and $\mathfrak{a}_y(t)>0$, which determine whether the local correction bends the near-vertical branch to the left, leaves it vertical, or bends it to the right.
    \begin{figure}[ht!]
        \centering
        \begin{subfigure}[t]{0.31\textwidth}
            \centering
            \begin{tikzpicture}[line cap=round,line join=round,>=Stealth]
                \definecolor{contourBlue}{RGB}{35,70,210}
                \definecolor{poleRed}{RGB}{170,30,30}
                \definecolor{zeroBlue}{RGB}{40,110,220}
                \definecolor{auxGray}{RGB}{100,100,100}
                \draw[->,auxGray] (-1.85,0) -- (1.95,0);
                \draw[->,auxGray] (0,-1.15) -- (0,1.15);
                \foreach \x in {-1.52,-1.34,-1.16} {
                    \node[poleRed] at (\x,0) {\tiny$\times$};
                }
                \foreach \x in {0.66,0.86,1.06} {
                    \node[black] at (\x-0.5,0) {\tiny$\times$};
                }
                \foreach \x in {1.42,1.60,1.78} {
                    \node[zeroBlue] at (\x-0.6,0) {\tiny$\times$};
                }
                \draw[contourBlue,thick,domain=-0.8:0.8,samples=200,smooth,variable=\y]
                    plot ({-0.32*(\y*\y)/0.64},{\y});
                \draw[poleRed,thick,domain=0.9273:5.3559,samples=240,smooth,variable=\s]
                    plot ({-0.72 + 1 - cos(\s r)},{sin(\s r)});
                \fill[black] (0,0) circle (1pt);
                \fill[contourBlue] (-0.32,0.8) circle (0.9pt);
                \fill[contourBlue] (-0.32,-0.8) circle (0.9pt);
                \node[font=\scriptsize,below, red] at (-1.34,0) {$-\beta_n^{-1}$};
                \node[font=\scriptsize,below] at (0.86-0.5,0) {$\lambda_n$};
                \node[font=\scriptsize,below, blue] at (1.60-0.6,0) {$\alpha_n^{-1}$};
                \node[font=\scriptsize,below left] at (0,0) {$\mathrm{u}_y(t)$};
                \node[font=\scriptsize] at (0,-1.34) {$\mathfrak{a}_y(t)<0$};
            \end{tikzpicture}
        \end{subfigure}
        \hfill
        \begin{subfigure}[t]{0.31\textwidth}
            \centering
            \begin{tikzpicture}[line cap=round,line join=round,>=Stealth]
                \definecolor{contourBlue}{RGB}{35,70,210}
                \definecolor{poleRed}{RGB}{170,30,30}
                \definecolor{zeroBlue}{RGB}{40,110,220}
                \definecolor{auxGray}{RGB}{100,100,100}
                \draw[->,auxGray] (-1.85,0) -- (1.95,0);
                \draw[->,auxGray] (0,-1.15) -- (0,1.15);
                \foreach \x in {-1.52,-1.34,-1.16} {
                    \node[poleRed] at (\x,0) {\tiny$\times$};
                }
                \foreach \x in {0.66,0.86,1.06} {
                    \node[black] at (\x-0.4,0) {\tiny$\times$};
                }
                \foreach \x in {1.42,1.60,1.78} {
                    \node[zeroBlue] at (\x-0.5,0) {\tiny$\times$};
                }
                \draw[contourBlue,thick] (0,-0.8) -- (0,0.8);
                \draw[poleRed,thick,domain=0.9273:5.3559,samples=240,smooth,variable=\s]
                    plot ({-0.40 + 1 - cos(\s r)},{sin(\s r)});
                \fill[black] (0,0) circle (1pt);
                \fill[contourBlue] (0,0.8) circle (0.9pt);
                \fill[contourBlue] (0,-0.8) circle (0.9pt);
                \node[font=\scriptsize,below, red] at (-1.34,0) {$-\beta_n^{-1}$};
                \node[font=\scriptsize,below] at (0.86-0.4,0) {$\lambda_n$};
                \node[font=\scriptsize,below, blue] at (1.60-0.5,0) {$\alpha_n^{-1}$};
                \node[font=\scriptsize,below left] at (0,0) {$\mathrm{u}_y(t)$};
                \node[font=\scriptsize] at (0,-1.34) {$\mathfrak{a}_y(t)=0$};
            \end{tikzpicture}
        \end{subfigure}
        \hfill
        \begin{subfigure}[t]{0.31\textwidth}
            \centering
            \begin{tikzpicture}[line cap=round,line join=round,>=Stealth]
                \definecolor{contourBlue}{RGB}{35,70,210}
                \definecolor{poleRed}{RGB}{170,30,30}
                \definecolor{zeroBlue}{RGB}{40,110,220}
                \definecolor{auxGray}{RGB}{100,100,100}
                \draw[->,auxGray] (-1.85,0) -- (1.95,0);
                \draw[->,auxGray] (0,-1.15) -- (0,1.15);
                \foreach \x in {-1.52,-1.34,-1.16} {
                    \node[poleRed] at (\x,0) {\tiny$\times$};
                }
                \foreach \x in {0.66,0.86,1.06} {
                    \node[black] at (\x-0.4,0) {\tiny$\times$};
                }
                \foreach \x in {1.42,1.60,1.78} {
                    \node[zeroBlue] at (\x-0.5,0) {\tiny$\times$};
                }
                \draw[contourBlue,thick,domain=-0.8:0.8,samples=200,smooth,variable=\y]
                    plot ({0.28*(\y*\y)/0.64},{\y});
                \draw[poleRed,thick,domain=0.9273:5.3559,samples=240,smooth,variable=\s]
                    plot ({-0.12 + 1 - cos(\s r)},{sin(\s r)});
                \fill[black] (0,0) circle (1pt);
                \fill[contourBlue] (0.28,0.8) circle (0.9pt);
                \fill[contourBlue] (0.28,-0.8) circle (0.9pt);
                \node[font=\scriptsize,below, red] at (-1.34,0) {$-\beta_n^{-1}$};
                \node[font=\scriptsize,below] at (0.86-0.4,0) {$\lambda_n$};
                \node[font=\scriptsize,below, blue] at (1.60-0.5,0) {$\alpha_n^{-1}$};
                \node[font=\scriptsize,below left] at (0,0) {$\mathrm{u}_y(t)$};
                \node[font=\scriptsize] at (0,-1.34) {$\mathfrak{a}_y(t)>0$};
            \end{tikzpicture}
        \end{subfigure}
        \caption{Schematic contour geometry for Corollary~\ref{steepestdescentpath}. In each panel, the blue curve is the local segment $\Phi_{y}^{\mathrm{SD}}$ and the red arc is the circular piece $\Phi_{\mathrm{c}(\theta,t)}^{\mathrm{C}}$ used to close the contour. The black dot marks the stationary point $\mathrm{u}_y(t)$, and the real-axis markers indicate the relative locations of the geometric poles $-\beta_n^{-1}$, the spatial rates $\lambda_n$, and the Bernoulli zeros $\alpha_n^{-1}$. From left to right, the panels show the three sign regimes for $\mathfrak{a}_y(t)$.}
        \label{contourconstructionfig}
    \end{figure}
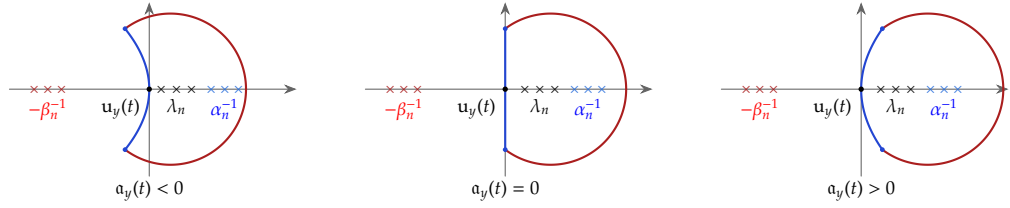

    \begin{lemma}[Endpoint reduction on the complementary contour]\label{lem:Gamma2-Re-max}
    Fix $t\in\mathscr{T}$, a path $y\in\mathscr{P}$ with $y(t)\in\mathscr{C}_{t}$, angles $\theta_{\star}\in(0,\theta_{\mathrm{C}})$ and a cutoff $\vartheta\in(0,\mathrm{R}_{\mathrm{C}}\sin\theta_{\star})$ such that Corollary~\ref{steepestdescentpath} applies at $(t,\theta_{\star})$. Let $\Gamma(t)$ denote the resulting closed contour from that corollary, with the endpoint convention from Corollary~\ref{steepestdescentpath}, parametrized as
    \begin{equation*}
        \Gamma(t)=\left\{\Phi_y(u,t,\theta_{\star}):-\frac{\mathrm{R}_{\mathrm{C}}\sin\theta_{\star}}{\theta_{\star}}\pi\leq u\leq\frac{\mathrm{R}_{\mathrm{C}}\sin\theta_{\star}}{\theta_{\star}}\pi\right\}.
    \end{equation*}
    Define the $\vartheta$-short and $\vartheta$-long pieces
    \begin{equation*}
        \Gamma_{1}^{(\vartheta)}(t)=\left\{\Phi_y(u,t,\theta_{\star}):|u|\leq\vartheta\right\},
        \qquad
        \Gamma_{2}^{(\vartheta)}(t)=\left\{\Phi_y(u,t,\theta_{\star}):\vartheta\leq |u|\leq\frac{\mathrm{R}_{\mathrm{C}}\sin\theta_{\star}}{\theta_{\star}}\pi\right\}.
    \end{equation*}
    Then, we have
    \begin{equation}\label{eq:Gamma2-Re-max}
        \max_{z\in\Gamma_{2}^{(\vartheta)}(t)}\mathfrak{Re}\bigl[\mathfrak{F}_{y}(z,t)\bigr]
        \leq
        \max_{\sigma\in\{-1,1\}}\mathfrak{Re}\bigl[\mathfrak{F}_{y}\bigl(\Phi_y(\sigma\vartheta,t,\theta_{\star}),t\bigr)\bigr].
    \end{equation}
    \end{lemma}
    \begin{proof}
    Write $U\overset{\textnormal{def}}{=}\mathrm{R}_{\mathrm{C}}\sin(\theta_{\star})/\theta_{\star}\cdot\pi$, so $\Gamma(t)=\{\Phi_y(u,t,\theta_{\star}):-U\leq u\leq U\}$ and $\Gamma_{2}^{(\vartheta)}(t)=\{\Phi_y(u,t,\theta_{\star}):\vartheta\leq |u|\leq U\}$, with the two endpoint values $u=-U$ and $u=U$ identified.
    On each local outward interval where $|u|\leq\mathrm{R}_{\mathrm{C}}\sin\theta_{\star}$ and $\Phi_y(u,t,\theta_{\star})=\Phi_{y}^{\mathrm{SD}}(u,t)$, the map $u\mapsto\Phi_y(u,t,\theta_{\star})$ is a descent curve for $\mathfrak{F}_y(\cdot,t)$ from $\mathrm{u}_y(t)$ in the sense of Definition~\ref{descentcurvedefn}, with interior stationary parameter $u=0$.
    Hence Proposition~\ref{prop:descent-endpoint-max} applies on the compact intervals
    \begin{equation*}
        J_+\overset{\textnormal{def}}{=}[\vartheta,\mathrm{R}_{\mathrm{C}}\sin\theta_{\star}]
        \quad\text{and}\quad
        J_-\overset{\textnormal{def}}{=}[-\mathrm{R}_{\mathrm{C}}\sin\theta_{\star},-\vartheta],
    \end{equation*}
    showing that the maximum of $\mathfrak{Re}[\mathfrak{F}_y(\Phi_y(\cdot,t,\theta_{\star}),t)]$ over each local outward branch is attained at $u=\pm\vartheta$. On each circular subarc of $\Gamma(t)$, the trace agrees with a translated circular template from Proposition~\ref{contour1} with centre $\mathrm{c}(\theta_{\star},t)$ from Corollary~\ref{steepestdescentpath}, hence is a descent curve in its native angle variable. The final sentence in the proof of Corollary~\ref{steepestdescentpath} shows that the affine reparametrisation in $u$ preserves the sign of the derivative of $\mathfrak{Re}[\mathfrak{F}_y]$ along that circular piece. Therefore Proposition~\ref{prop:descent-endpoint-max} applies on each compact circular subinterval contained in $\Gamma_{2}^{(\vartheta)}(t)$, and continuity at the identified endpoint $u=\pm U$ covers the closed circular piece. Hence the maximum of $\mathfrak{Re}[\mathfrak{F}_y]$ over the circular part of $\Gamma_{2}^{(\vartheta)}(t)$ is attained at the splice parameters $u=\pm\mathrm{R}_{\mathrm{C}}\sin\theta_{\star}$. Finally, since $\vartheta<\mathrm{R}_{\mathrm{C}}\sin\theta_{\star}$, the same descent monotonicity along the local outward branches implies
    \begin{equation*}
        \mathfrak{Re}\bigl[\mathfrak{F}_y(\Phi_y(\pm\mathrm{R}_{\mathrm{C}}\sin\theta_{\star},t,\theta_{\star}),t)\bigr]
        \leq
        \mathfrak{Re}\bigl[\mathfrak{F}_y(\Phi_y(\pm\vartheta,t,\theta_{\star}),t)\bigr],
    \end{equation*}
    because moving from $|u|=\vartheta$ to $|u|=\mathrm{R}_{\mathrm{C}}\sin\theta_{\star}$ along $\Phi_{y}^{\mathrm{SD}}(\cdot,t)$ moves away from $\mathrm{u}_y(t)$ along a descent segment. Combining the local and circular cases yields \eqref{eq:Gamma2-Re-max}.
    \end{proof}
    \medskip

    \subsection{Decoupling initial data}\label{decouplingsection}
    This section completes the deterministic one-particle analysis of the Markov kernel. Throughout the section $y\in\mathscr{P}$ is a path, and the endpoint condition is imposed by requiring $y(t)\in\mathscr{C}_{t}$. On this endpoint window, the uniform versions of \textbf{(ii)}--\textbf{(iv)} in Proposition~\ref{specificcasesub} and the stationary-point construction in Proposition~\ref{axpro} identify and approximate $\mathrm{u}_{y}(t)$ with constants independent of the particular path satisfying $y(t)\in\mathscr{C}_{t}$. The key outcome of this subsection is a kind of decoupling principle. The next proposition gives the relevant expansion uniformly over the deterministic endpoint window. 
    
    \begin{proposition}\label{Thm2}
    Assume \textbf{(A1)} and \textbf{(A2)}. There exist constants $\mathrm{T}_{\mathrm{dec}}\in\mathscr{T}$ and $\mathrm{c}_{\varphi}>0$ such that, for every $t>\mathrm{T}_{\mathrm{dec}}$ and every path $y\in\mathscr{P}$ satisfying $y(t)\in\mathscr{C}_{t}$, the stationary point $\mathrm{u}_{y}(t)$ is well defined and belongs to $\mathcal{B}_{\mathrm{R}_{2}}(0)$. Moreover, for every $x\in\mathbb{Z}_+$, with the explicit constant $\mathrm{C}_{x}$ defined in \eqref{eq:Cx-explicit}, independently of $t$ and $y$,
        \begin{align}\label{ointc}
            \left|\frac{\mathbb{P}_{x}[\mathscr{X}(t)=y(t)]}{\mathbb{P}_{0}[\mathscr{X}(t)=y(t)]}-\sum_{n=0}^{x}\frac{\rho_{x}^{(n)}\left(\mathrm{u}_y(t)\right)}{n!}\eta_{n,y}(t)\right|\leq \mathrm{C}_{x}\mathrm{e}^{-\mathrm{c}_{\varphi}\Sigma(t)^{\frac{2}{3}}}.
        \end{align}
        The coefficients $\eta_{n,y}(t)$ are defined by \eqref{eq:eta-def} in the proof below.
        Finally, for each fixed $n\in\mathbb{Z}_+$ there is a lower time $\mathrm{T}_{\mathrm{dec},n}\geq\mathrm{T}_{\mathrm{dec}}$ such that, with the convention $(-1)!!=1$ and with the explicit constant $\mathrm{C}_{n}$ defined in \eqref{eq:Cn-explicit}, independently of $t$ and $y$,
        \begin{equation}\label{leadingordereta}
            \begin{aligned}
                &\left|\eta_{n,y}(t)-\frac{(-1)^{\frac{n}{2}}(n-1)!!}{\Sigma(t)^{n}}\right|\leq \mathrm{C}_{n}\gamma(t)^{1/2}\Sigma(t)^{-n},\quad\text{if $n$ is even,}\\
                &\left|\eta_{n,y}(t)\right|\leq \mathrm{C}_{n}\Sigma(t)^{-n-1},\quad\text{if $n$ is odd.}
            \end{aligned}
        \end{equation}
        The bound \eqref{ointc} holds for every $t>\mathrm{T}_{\mathrm{dec}}$ and every path $y\in\mathscr{P}$ satisfying $y(t)\in\mathscr{C}_{t}$. For each fixed $n$, the coefficient bounds in \eqref{leadingordereta} hold for every $t>\mathrm{T}_{\mathrm{dec},n}$ and every such path.
    \end{proposition}
    \begin{proof}
        \medskip

        \noindent\textbf{Outline.} The proof is long and technical so we give an outline. We fix a deterministic lower time $\mathrm{T}_{\mathrm{dec}}$ by finitely many increases, each justified by uniform bounds on the endpoint window $\{y(t)\in\mathscr{C}_{t}\}$ from Proposition~\ref{specificcasesub} and Proposition~\ref{axpro}.
        The argument is organised in seven labelled segments below (Step~4 is subdivided into Steps~4a--4c). Step~1 fixes the deterministic threshold and deforms \eqref{cir} to the descent contour from Corollary~\ref{steepestdescentpath}. Step~2 splits that contour into the central arc $\Gamma_{1,y}(t)$ and the complementary arc $\Gamma_{2,y}(t)$, with Lemma~\ref{lem:Gamma2-Re-max} in Section~\ref{CM} providing the endpoint reduction on $\Gamma_{2,y}(t)$. Step~3 expands the polynomial prefactor $\rho_x$ on $\Gamma_{1,y}(t)$ and introduces the normalised coefficients $\eta_{n,y}(t)$. Steps~4a--4c estimate the local integrals $\mathcal{I}_{n,y}(t)$ by passing to the steepest-descent parameter, expanding the phase, and evaluating the resulting Gaussian main term with explicit remainders. Step~5 proves exponential smallness of the complementary integral $I_{2,x,y}(t)$. Step~6 proves the denominator lower bound, transfers the local integral estimates to $\eta_{n,y}(t)$, and replaces the exact curvature by $\Sigma(t)^2$, yielding \eqref{leadingordereta}. Step~7 combines the exact decomposition with the complementary-contour estimate to obtain \eqref{ointc}.
        \begin{center}
        \begin{tikzpicture}[x=1cm,y=1cm,>=Stealth,font=\scriptsize,
            stepbox/.style={draw=black!55,fill=black!3,rounded corners=2pt,align=center,inner xsep=4pt,inner ysep=3pt,minimum width=2.18cm,minimum height=0.80cm},
            arrow/.style={->,draw=black!65,line width=0.45pt},
            tag/.style={font=\tiny,text=black!65,inner sep=1pt}]
            \node[stepbox,fill=blue!5] (s1) at (0,0) {\textbf{Step 1}\\threshold and\\contour deformation};
            \node[stepbox,fill=blue!5] (s2) at (2.65,0) {\textbf{Step 2}\\central/complement\\split};
            \node[stepbox,fill=blue!5] (s3) at (5.30,0) {\textbf{Step 3}\\prefactor expansion\\and $\eta_{n,y}$};
            \node[stepbox,fill=green!7] (s4) at (7.95,0) {\textbf{Steps 4a--c}\\local integrals\\and Gaussian term};
            \node[stepbox,fill=red!5] (s5) at (7.95,-1.35) {\textbf{Step 5}\\complementary\\tail $I_{2,x,y}$};
            \node[stepbox,fill=orange!8] (s6) at (5.30,-1.35) {\textbf{Step 6}\\denominator and\\coefficients};
            \node[stepbox,fill=orange!8] (s7) at (2.65,-1.35) {\textbf{Step 7}\\final ratio\\estimate};
            \draw[arrow] (s1) -- (s2);
            \draw[arrow] (s2) -- (s3);
            \draw[arrow] (s3) -- (s4);
            \draw[arrow] (s4) -- (s5);
            \draw[arrow] (s5) -- (s6);
            \draw[arrow] (s6) -- (s7);
        \end{tikzpicture}
        \end{center}
        
        \textbf{Step 1: threshold and contour deformation.}
        We construct the threshold $\mathrm{T}_{\mathrm{dec}}$ appearing in the statement by starting from a deterministic lower time that dominates the uniform endpoint-window estimates in Proposition~\ref{specificcasesub} and Proposition~\ref{axpro}, and then increasing it finitely many further times as the argument demands; each increase is justified by a uniform inequality on $\{y(t)\in\mathscr{C}_{t}\}$, hence is independent of the particular path $y$. By Proposition~\ref{axpro}, once the lower time dominates $\mathrm{T}_{\mathrm{A}}$, the point $\mathrm{u}_{y}(t)$ is the unique real stationary point identified by Definition~\ref{definitionstationary} for every path satisfying $y(t)\in\mathscr{C}_{t}$. The localization estimate \eqref{axpro-u-small} and the fact that $\gamma(t)\to0$ allow us to increase $\mathrm{T}_{\mathrm{dec}}$ so that $\mathrm{u}_{y}(t)\in\mathcal{B}_{\mathrm{R}_{2}}(0)$ throughout the endpoint window. We shall take
            \begin{equation*}
                \mathrm{c}_{\varphi}\overset{\textnormal{def}}{=}\frac{1}{8\,2^{1/3}}.
            \end{equation*}
            We now fix $t>\mathrm{T}_{\mathrm{dec}}$ and a path $y\in\mathscr{P}$ satisfying $y(t)\in\mathscr{C}_{t}$, and fix the splicing angle
            \begin{equation*}
                \theta_{\star}\overset{\textnormal{def}}{=}\frac{1}{2}\theta_{\mathrm{C}},
                \quad\text{and set}\quad
                \phi_{\star}\overset{\textnormal{def}}{=}\mathrm{R}_{\mathrm{C}}\sin(\theta_{\star}),
            \end{equation*}
            where $\theta_{\mathrm{C}}$ is the constant from Corollary~\ref{steepestdescentpath} and $\mathrm{R}_{\mathrm{C}}=(\mathfrak{B}_{1}+\mathfrak{B}_{2})/2$ is the circular radius from Proposition~\ref{contour1}. After increasing $\mathrm{T}_{\mathrm{dec}}$ so that $\mathrm{T}_{\mathrm{dec}}\geq\mathrm{T}_{\mathrm{C}}$, that corollary gives the descent contour
            \begin{equation*}
                \Gamma_y(t)\overset{\textnormal{def}}{=}
                \left\{\Phi_y(u,t,\theta_{\star}): -\frac{\phi_{\star}}{\theta_{\star}}\pi\leq u<\frac{\phi_{\star}}{\theta_{\star}}\pi\right\}.
            \end{equation*}
	            This contour is the positively oriented contour determined by Corollary~\ref{steepestdescentpath}; when $u$ increases, the parametrisation gives the descent direction and therefore runs in the opposite orientation. The contour constructed in Corollary~\ref{steepestdescentpath} is homologous to the original Cauchy contour from \eqref{cir} in the domain obtained by deleting the spatial poles $\{\lambda_0,\ldots,\lambda_{y(t)}\}$ and the geometric poles. The Bernoulli factors are zeros in the branch-free product form of $\mathrm{e}^{\mathfrak{F}_y}$, and hence do not obstruct this deformation. More explicitly, both contours have winding number one around each spatial pole $\lambda_0,\ldots,\lambda_{y(t)}$ and winding number zero around every geometric pole; the connecting homotopy stays inside the holomorphic domain of the full integrand away from those spatial poles. Consequently, Cauchy's theorem applied to the meromorphic integrand $\rho_x(w)\mathrm{e}^{\mathfrak{F}_y(w,t)}$ gives the same residue sum and permits the deformation in \eqref{cir}:
            \begin{equation*}
                \mathbb{P}_{x}\left[\mathscr{X}(t)=y(t)\right]=-\frac{1}{\lambda_{y(t)}2\pi \mathrm{i}}\oint_{\Gamma_y(t)}\rho_{x}(w)\mathrm{e}^{\mathfrak{F}_{y}(w,t)}\mathrm{d}w.
            \end{equation*}
            We now choose the local cutoff
            \begin{equation*}
                \theta_y(t)\overset{\textnormal{def}}{=}\left|\mathfrak{F}_{y}^{(2)}(\mathrm{u}_y(t),t)\right|^{-1/3}.
            \end{equation*}
            By the lower comparison in \eqref{eq:SD-second-lower}, $\mathfrak{F}_{y}^{(2)}(\mathrm{u}_y(t),t)\geq\Sigma(t)^2/2$ on the endpoint window, and $\Sigma(t)\to\infty$ by \textbf{(iv)} of Proposition~\ref{specificcasesub}. Hence $\theta_y(t)\leq 2^{1/3}\Sigma(t)^{-2/3}$ and $\theta_y(t)\to0$ uniformly on the endpoint window. We choose $\mathrm{T}_{\mathrm{dec}}$ to dominate the lower time for which $2^{1/3}\Sigma(t)^{-2/3}<\phi_{\star}$; then $\theta_y(t)<\phi_{\star}$ for every $t>\mathrm{T}_{\mathrm{dec}}$ and every path $y\in\mathscr{P}$ satisfying $y(t)\in\mathscr{C}_{t}$.
           
            \textbf{Step 2: central arc versus complementary arc.}
            We split the contour integral into two parts: for $i\in\{1,2\}$ and $x\in\mathbb{Z}_+$,
            \begin{equation*}
                I_{i,x,y}(t)\overset{\textnormal{def}}{=}-\frac{1}{\lambda_{y(t)}2\pi \mathrm{i}}\int_{\Gamma_{i,y}(t)}\rho_{x}(w)\mathrm{e}^{\mathfrak{F}_{y}(w,t)}\mathrm{d}w,
            \end{equation*}
            where we define the two segments of $\Gamma_y(t)$, $\Gamma_{1,y}(t)$ and $\Gamma_{2,y}(t)$ as follows:
            \begin{equation*}
                \Gamma_{1,y}(t)=\left\{\Phi_y(u,t,\theta_{\star}):u\in\left[-\theta_y(t),\theta_y(t)\right] \right\}\quad\text{and}\quad\Gamma_{2,y}(t)=\Gamma_y(t)\setminus\Gamma_{1,y}(t).
            \end{equation*}
           
            \textbf{Step 3: Taylor expansion of the prefactor and coefficient normalisation.}
            Since $\rho_x$ is a polynomial of degree $x$, its Taylor expansion around $\mathrm{u}_{y}(t)$ is finite. Expanding $\rho_x$ inside the integral defining $I_{1,x,y}(t)$ therefore gives
            \begin{equation}\label{firstpart}
            I_{1,x,y}(t)=\sum_{n=0}^{x}\frac{\rho_x^{(n)}(\mathrm{u}_{y}(t))}{n!}\mathcal{I}_{n,y}(t),
            \end{equation}
            where for each $n\in\mathbb{Z}_+$ we define
            \begin{equation}\label{defofmathcalI}
                \mathcal{I}_{n,y}(t)\overset{\textnormal{def}}{=}-\frac{1}{\lambda_{y(t)}2\pi\mathrm{i}}\int_{\Gamma_{1,y}(t)}\left(w-\mathrm{u}_{y}(t)\right)^{n}\mathrm{e}^{\mathfrak{F}_{y}(w,t)}\mathrm{d}w.
            \end{equation}
            Equivalently,
            \begin{equation*}
                I_{1,x,y}(t)=\left(\sum_{n=0}^{x}\frac{\rho_x^{(n)}(\mathrm{u}_{y}(t))}{n!}\eta_{n,y}(t)\right)\mathbb{P}_0\left[\mathscr{X}(t)=y(t)\right],
            \end{equation*}
            where $\eta_{n,y}(t)$ is defined by
            \begin{equation}\label{eq:eta-def}
                \mathcal{I}_{n,y}(t)=\eta_{n,y}(t)\mathbb{P}_{0}\left[\mathscr{X}(t)=y(t)\right].
            \end{equation}
            This quotient is legitimate on the time range ultimately used in the proposition: Step~6 proves the explicit lower bound \eqref{eq:P0-lower}, and hence $\mathbb{P}_{0}[\mathscr{X}(t)=y(t)]>0$, uniformly on the endpoint window after the stated lower-time increase.
            Subsequently, setting $x=0$, we can use this decomposition to obtain
            \begin{equation*}
                \mathbb{P}_{0}\left[\mathscr{X}(t)=y(t)\right]=\mathcal{I}_{0,y}(t)+I_{2,0,y}(t),
            \end{equation*}
            and therefore, whenever $|I_{2,0,y}(t)|\leq |\mathcal{I}_{0,y}(t)|/2$, the coefficient $\eta_{n,y}(t)$ differs from the local ratio by the explicit bound
            \begin{equation}\label{etaratio}
            \left|\eta_{n,y}(t)-\frac{\mathcal{I}_{n,y}(t)}{\mathcal{I}_{0,y}(t)}\right|
            \leq
            2\left|\frac{\mathcal{I}_{n,y}(t)}{\mathcal{I}_{0,y}(t)}\right|
            \left|\frac{I_{2,0,y}(t)}{\mathcal{I}_{0,y}(t)}\right|.
            \end{equation}
            \medskip

            \noindent The regime $|I_{2,0,y}(t)|\leq|\mathcal{I}_{0,y}(t)|/2$ used here is recorded in \eqref{eq:I20-dominated-by-local} at the start of Step~6.
	            Indeed, \eqref{etaratio} is obtained by writing $\eta_{n,y}(t)=\mathcal{I}_{n,y}(t)/(\mathcal{I}_{0,y}(t)+I_{2,0,y}(t))$ and applying the elementary inequality $|(1+z)^{-1}-1|\leq2|z|$ for $|z|\leq1/2$.

            \textbf{Step 4a: Local parametrisation and Taylor coefficients of the phase.}
            It remains to estimate the local integrals $\mathcal{I}_{n,y}(t)$ and the complementary contour integral $I_{2,x,y}(t)$ explicitly. For the local part, we write $\mathcal{I}_{n,y}(t)$ in parametric form:
            \begin{equation}\label{parametriccontourint}
	                \begin{aligned}
                        \mathcal{I}_{n,y}(t)
                        =
                        \frac{1}{\lambda_{y(t)}2\pi \mathrm{i}}
                        \int_{-\theta_y(t)}^{\theta_y(t)}
                        \left(\Phi_y(u,t,\theta_{\star})-\mathrm{u}_{y}(t)\right)^n
                        \frac{\mathrm{d}}{\mathrm{d}u}\Phi_y(u,t,\theta_{\star})
                        \mathrm{e}^{\mathfrak{F}_{y}(\Phi_y(u,t,\theta_{\star}),t)}
                        \mathrm{d}u.
                    \end{aligned}
            \end{equation}
            In \eqref{defofmathcalI}, the integral over $\Gamma_{1,y}(t)$ uses the positive orientation inherited from $\Gamma_y(t)$. On the local segment, the increasing-$u$ parametrisation $u\mapsto\Phi_y(u,t,\theta_{\star})$ runs through the same geometric arc in the reverse direction from the contour orientation on $\Gamma_{1,y}(t)$; therefore the sign from reversing the parametrisation cancels the minus sign in \eqref{defofmathcalI}, giving \eqref{parametriccontourint}. By construction of $\Phi$, for $t\in\mathscr{T}$ such that $\theta_y(t)<\phi_{\star}$, we are integrating over the segment of the contour corresponding to $\Phi_{y}^{\mathrm{SD}}$.
         
            \emph{The case $\mathfrak{N}_y(t)=\infty$.}
            We treat this case in one self-contained paragraph, since all subsequent formulas simplify drastically. The definition of $\mathfrak{N}_y(t)$ gives $\mathfrak{F}_{y}^{(2m+1)}(\mathrm{u}_y(t),t)=0$ for every $m\geq1$, $\Phi_{y}^{\mathrm{SD}}(u,t)=\mathrm{u}_y(t)+\mathrm{i}u$, $\mathfrak{a}_y(t)=0$, and $\chi_y(\theta,t)=0$. Consequently $\Phi_y(-u,t,\theta_{\star})-\mathrm{u}_y(t)=-(\Phi_y(u,t,\theta_{\star})-\mathrm{u}_y(t))$, all odd Taylor coefficients $G_{i,y}(t)$ of $H_y$ vanish, and $H_y(-u,t)=H_y(u,t)$ on the local segment. Hence $\exp(H_y(u,t))$ is even in $u$ while $(\Phi_y(u,t,\theta_{\star})-\mathrm{u}_y(t))^n=\mathrm{i}^n u^n$ has parity $n$, so the integrand in \eqref{parametriccontourint} is odd whenever $n$ is odd and $\mathcal{I}_{n,y}(t)=0$ for every odd $n$. The remaining estimates for even $n$ then reduce to the special case of the formulas below with $\mathfrak{a}_y(t)=0$ and the single term $k=0$ retained; explicitly, the parenthesis in the rescaled integral \eqref{parametricthing} collapses to the single vertical term $\mathrm{i}^{n+1}u^n/\mathfrak{F}_y^{(2)}(\mathrm{u}_y(t),t)^{(n+1)/2}$. Thus the local bounds \eqref{eq:local-remainder-explicit}, \eqref{eq:local-even-In}, \eqref{eq:local-odd-In}, and the normalising bound \eqref{eq:I0-local-bound} below also hold in the case $\mathfrak{N}_y(t)=\infty$, with the odd local integrals equal to zero and with the correction terms deleted.
            
            \emph{The case $\mathfrak{N}_y(t)<\infty$.}
            We now assume $\mathfrak{N}_y(t)<\infty$ throughout the remainder of Step~4. (A few formulas below still display both branches of $\mathfrak{N}_y(t)<\infty$/$\mathfrak{N}_y(t)=\infty$ for compactness; under the standing assumption of this case only the $\mathfrak{N}_y(t)<\infty$ branch is in active use, the other being a place-holder that vanishes by the simplification of the preceding paragraph.) Substituting the formula for $\Phi_{y}^{\mathrm{SD}}$ into $\left(\Phi_y(u,t,\theta_{\star})-\mathrm{u}_{y}(t)\right)^n\frac{\mathrm{d}}{\mathrm{d}u}\Phi_y(u,t,\theta_{\star})$ yields the expansion
			            \begin{equation}\label{integrandexpansion}
	                    \begin{aligned}
	                    \left(\Phi_y(u,t,\theta_{\star})-\mathrm{u}_{y}(t)\right)^n\frac{\mathrm{d}}{\mathrm{d}u}\Phi_y(u,t,\theta_{\star})
	                    =
			            \sum_{k=0}^{n+1}\mathrm{i}^{n-k+1}\left(\frac{\mathfrak{a}_y(t)}{(2\mathfrak{N}_y(t))!}\right)^{k}
			            \mathfrak{c}_{n,k,y}(t)u^{2\mathfrak{N}_y(t)k+n-k},
	                    \end{aligned}
			            \end{equation}
            where to reduce the complexity of later sums, we define
            \begin{equation*}
                \mathfrak{c}_{n,k,y}(t)\overset{\textnormal{def}}{=}\binom{n}{k}\mathbf{1}_{\{0\leq k\leq n\}}+2\mathfrak{N}_y(t)\binom{n}{k-1}\mathbf{1}_{\{1\leq k\leq n+1\}}.
            \end{equation*}
            To organise the local expansion, set
            \begin{equation*}
                H_y(u,t)\overset{\textnormal{def}}{=}\mathfrak{F}_y(\Phi_y(u,t,\theta_{\star}),t)-\mathfrak{F}_y(\mathrm{u}_y(t),t).
            \end{equation*}
            By \textbf{(i)} of Proposition~\ref{specificcasesub}, $\mathfrak{F}_{y}(\cdot,t)$ is holomorphic in $\mathcal{B}_{\mathrm{R}_{-}}(0)$. Since $\theta_y(t)<\phi_{\star}$ and the local segment is exactly $\Phi_{y}^{\mathrm{SD}}$, \eqref{eq:SD-path-stays-close} shows that the whole segment $\Gamma_{1,y}(t)$ lies in $\mathcal{B}_{\mathrm{R}_{2}}(0)$ once $\mathrm{A}_{\mathrm{SD}}\theta_y(t)\leq3\mathrm{R}_{2}/8$. We include this deterministic condition in the definition of $\mathrm{T}_{\mathrm{dec}}$, again possible uniformly over paths $y\in\mathscr{P}$ satisfying $y(t)\in\mathscr{C}_{t}$ because $\theta_y(t)\leq2^{1/3}\Sigma(t)^{-2/3}$. Therefore $H_y(\cdot,t)$ is holomorphic on a fixed neighbourhood of the origin, and Taylor's theorem gives
	            \begin{equation}\label{globalapproximate}
		                \mathfrak{F}_{y}(\Phi_y(u,t,\theta_{\star}),t)=\mathfrak{F}_{y}(\mathrm{u}_y(t),t)+\sum_{n=1}^{\infty}G_{n,y}(t)\frac{u^n}{n!},
	            \end{equation}
            where the Taylor coefficients are defined by
            \begin{equation*}
	                G_{n,y}(t)\overset{\textnormal{def}}{=}\frac{\mathrm{d}^n}{\mathrm{d}u^n}\mathfrak{F}_y(\Phi_y(u,t,\theta_{\star}),t)\Bigg|_{u=0}.
            \end{equation*}
            For $n\geq 1$, these are also the Taylor coefficients of $H_y(\cdot,t)$ at the origin. By Fa\`a di Bruno's formula via exponential Bell polynomials,
            \begin{equation*}
                G_{n,y}(t)=\sum_{k=1}^{n}\mathfrak{F}_y^{(k)}(\mathrm{u}_y(t),t)B_{n,k,y}(\Phi,t),\qquad\text{for all $n\geq 1$,}
            \end{equation*}
            where the path-dependent exponential Bell polynomials are defined, for each $n\geq1$ and $k\geq1$, by
            \begin{equation}\label{bellsum}
                \begin{aligned}
                    &B_{n,k,y}(\Phi,t)=n!\sum_{\underline{j}\in\mathcal{R}_{n,k}}\prod_{i=1}^{n-k+1}\frac{\Phi_y^{(i)}(0,t)^{j_i}}{(i!)^{j_i}j_i!},\quad&\text{for $1\leq k\leq n$}\\&B_{n,k,y}(\Phi,t)=0&\text{ otherwise.}
                \end{aligned}
            \end{equation}
            and the combinatorial index set is
            \begin{equation*}
                \mathcal{R}_{n,k}\overset{\textnormal{def}}{=}\left\{\underline{j}=(j_1,\ldots,j_{n-k+1})\in\mathbb{Z}_+^{n-k+1}:\sum_{i=1}^{n-k+1}j_i=k\text{ and }\sum_{i=1}^{n-k+1}ij_i=n\right\}.
            \end{equation*}
	            In the finite-$\mathfrak{N}_y(t)$ case, the definition of $\Phi_{y}^{\mathrm{SD}}$ gives, for all $n\geq1$,
	            \begin{equation}\label{Phi}
	                \Phi_y^{(n)}(0,t)=\mathrm{i}\,\mathbf{1}_{\{n=1\}}+\mathfrak{a}_y(t)\,\mathbf{1}_{\{n=2\mathfrak{N}_y(t)\}}.
	            \end{equation}
            Hence all derivatives of $\Phi_y$ at the origin vanish except the first and the $2\mathfrak{N}_y(t)$-th.
          
            \emph{Combinatorial specialisation:}
            We now specialise the Bell-polynomial sum in \eqref{bellsum} using the sparse derivative structure in \eqref{Phi}. Since $\Phi_y^{(i)}(0,t)=0$ unless $i=1$ or $i=2\mathfrak{N}_y(t)$, a nonzero term in \eqref{bellsum} must satisfy $j_i=0$ for every $i\notin\{1,2\mathfrak{N}_y(t)\}$. Thus the two partition constraints defining $\mathcal{R}_{n,k}$ reduce to
            \begin{equation*}
                j_1+j_{2\mathfrak{N}_y(t)}=k,\qquad j_1+2\mathfrak{N}_y(t)j_{2\mathfrak{N}_y(t)}=n.
            \end{equation*}
            Subtracting the first identity from the second shows that any contribution using the $2\mathfrak{N}_y(t)$-th derivative must have
            \begin{equation*}
                j_{2\mathfrak{N}_y(t)}=\frac{n-k}{2\mathfrak{N}_y(t)-1}\overset{\textnormal{def}}{=}r_{n,k,y}(t).
            \end{equation*}
            This quotient must belong to $\mathbb{Z}_+\setminus\{0\}$ for the correction term; if it does not, no partition with a higher-derivative factor exists. Substituting \eqref{Phi} into \eqref{bellsum} therefore gives
            \begin{equation*}
                B_{n,k,y}(\Phi,t)=\mathrm{i}^{n}\mathbf{1}_{\{n=k\}}+R_{n,k,y}(t).
            \end{equation*}
            where $R_{n,k,y}(t)$ is given by
            \begin{equation*}
                R_{n,k,y}(t)\overset{\textnormal{def}}{=}\begin{cases}\dfrac{n!}{k!}\binom{k}{r_{n,k,y}(t)}\left(\dfrac{\mathfrak{a}_y(t)}{(2\mathfrak{N}_y(t))!}\right)^{r_{n,k,y}(t)}\mathrm{i}^{k-r_{n,k,y}(t)},&\text{if }r_{n,k,y}(t)\in\mathbb{Z}_+\setminus\{0\},\\0,&\text{otherwise.}\end{cases}
            \end{equation*}
            The term $\mathrm{i}^{n}\mathbf{1}_{\{n=k\}}$ corresponds to the pure first-derivative partition $j_1=k=n$, which is the case $r_{n,k,y}(t)=0$ and is deliberately separated from $R_{n,k,y}(t)$. The remainder $R_{n,k,y}(t)$ collects exactly the partitions with at least one block of size $2\mathfrak{N}_y(t)$: when $r_{n,k,y}(t)\in\mathbb{Z}_+\setminus\{0\}$ there are $r_{n,k,y}(t)$ such blocks and $k-r_{n,k,y}(t)$ first-derivative blocks, producing the displayed multinomial factor; otherwise the Bell constraints have no solution contributing to the correction term, and the remainder is zero.
        
            \textbf{Step 4b: Bell polynomials, tail bounds on $B_{n,y}$, and the Gaussian prefactor.}
            Therefore, using the stationary-point identity $\mathfrak{F}_y^{(1)}(\mathrm{u}_y(t),t)=0$, we obtain
            \begin{equation}\label{eq:Snclean}
			                G_{n,y}(t)=\mathrm{i}^{n}\mathfrak{F}_y^{(n)}(\mathrm{u}_y(t),t)+\sum_{k=2}^{n}\mathfrak{F}_y^{(k)}(\mathrm{u}_y(t),t)R_{n,k,y}(t).
            \end{equation}
            We now read off the low-order coefficients from \eqref{eq:Snclean}. If $n<2\mathfrak{N}_y(t)$ and $k\geq2$, then $n-k<2\mathfrak{N}_y(t)-1$, so the quotient $r_{n,k,y}(t)=(n-k)/(2\mathfrak{N}_y(t)-1)$ cannot be a positive integer. Hence $R_{n,k,y}(t)=0$ for every $k\geq2$, and
            \begin{equation*}
	                G_{n,y}(t)=\mathrm{i}^{n}\mathfrak{F}_y^{(n)}(\mathrm{u}_y(t),t),\qquad\text{for all $n<2\mathfrak{N}_y(t)$.}
            \end{equation*}
            When $n=2\mathfrak{N}_y(t)$, the only possible extra term is the one with $k=1$, and that term again vanishes because $\mathfrak{F}_y^{(1)}(\mathrm{u}_y(t),t)=0$. Hence
            \begin{equation*}
	                G_{2\mathfrak{N}_y(t),y}(t)=\mathrm{i}^{2\mathfrak{N}_y(t)}\mathfrak{F}_y^{(2\mathfrak{N}_y(t))}(\mathrm{u}_y(t),t).
            \end{equation*}
            Finally, when $n=2\mathfrak{N}_y(t)+1$, the only non-zero contribution in the sum in \eqref{eq:Snclean} comes from $k=2$, so
            \begin{equation*}
	                G_{2\mathfrak{N}_y(t)+1,y}(t)=\mathrm{i}^{2\mathfrak{N}_y(t)+1}\mathfrak{F}_y^{(2\mathfrak{N}_y(t)+1)}(\mathrm{u}_y(t),t)+\mathrm{i}(2\mathfrak{N}_y(t)+1)\mathfrak{a}_y(t)\mathfrak{F}_y^{(2)}(\mathrm{u}_y(t),t)=0,
            \end{equation*}
            where the final equality is exactly the definition of $\mathfrak{a}_y(t)$. Moreover, by the definition of $\mathfrak{N}_y(t)$,
            \begin{equation*}
                \mathfrak{F}_{y}^{(2m+1)}(\mathrm{u}_y(t),t)=0,\qquad\text{for all $m<\mathfrak{N}_y(t)$.}
            \end{equation*}
            Combining the previous identities, we obtain
            \begin{equation*}
	                G_{n,y}(t)=0,\qquad\text{for all odd $n\leq 2\mathfrak{N}_y(t)+1$ in the finite-$\mathfrak{N}_y(t)$ case.}
            \end{equation*}
            We now isolate the Gaussian quadratic term. Since
            \begin{equation*}
	                G_{2,y}(t)=\mathrm{i}^{2}\mathfrak{F}_y^{(2)}(\mathrm{u}_y(t),t)=-\mathfrak{F}_y^{(2)}(\mathrm{u}_y(t),t),
            \end{equation*}
            we now define the truncated Taylor coefficients as
            \begin{equation*}
	                \widetilde{G}_{n,y}(t)\overset{\textnormal{def}}{=}G_{n,y}(t)\mathbf{1}_{\{n\geq 3\}}.
            \end{equation*}
            Then the truncated Bell expansion is given by
            \begin{equation}\label{bellexpansiontwo}
	                H_y(u,t)=-\mathfrak{F}_y^{(2)}(\mathrm{u}_y(t),t)\frac{u^2}{2}+\sum_{n=1}^{\infty}\widetilde{G}_{n,y}(t)\frac{u^n}{n!},
            \end{equation}
            and applying the generating function for complete exponential Bell polynomials gives
            \begin{equation}\label{generatingexpansion}
	                \exp\left(\sum_{n=1}^{\infty}\widetilde{G}_{n,y}(t)\frac{u^n}{n!}\right)=\sum_{n=0}^{\infty}B_{n,y}(t)\frac{u^n}{n!},
            \end{equation}
            where, in this context, the Bell coefficients are defined by
            \begin{equation*}
	                B_{n,y}(t)\overset{\textnormal{def}}{=}n!\sum_{\underline{j}\in\mathcal{R}_n}\prod_{i=1}^{n}\frac{\widetilde{G}_{i,y}(t)^{j_i}}{(i!)^{j_i}j_i!},\quad
            \end{equation*}
            with the combinatorial index set $\mathcal{R}_n$, given by
            \begin{equation*}
                \mathcal{R}_n\overset{\textnormal{def}}{=}\left\{\underline{j}=(j_1,\ldots,j_n)\in\mathbb{Z}_+^n:\sum_{i=1}^{n}ij_i=n\right\}.
            \end{equation*}
            Therefore, the exponent in \eqref{parametriccontourint} excluding the $\exp(\mathfrak{F}_y(\mathrm{u}_y(t),t))$ term, can be written as
            \begin{equation*}
                \exp\left(H_y(u,t)\right)=\left(\sum_{n=0}^{\infty}B_{n,y}(t)\frac{u^n}{n!}\right)\exp\left(-\mathfrak{F}_y^{(2)}(\mathrm{u}_y(t),t)\frac{u^2}{2}\right).
            \end{equation*}
	            Since $\widetilde{G}_{i,y}(t)=0$ for $i=1,2$ and for all odd $i\leq 2\mathfrak{N}_y(t)+1$ (recall we have assumed $\mathfrak{N}_y(t)<\infty$ since Step~4a), if $n$ is odd and $n\leq 2\mathfrak{N}_y(t)+1$, then every non-zero term in the Bell expansion of $B_{n,y}(t)$ would have to involve only even indices, which is impossible. Therefore
            \begin{equation*}
                B_{n,y}(t)=0,\qquad\text{for all odd $n$ with $n\leq 2\mathfrak{N}_y(t)+1$.}
            \end{equation*}
	            To control the higher coefficients, we use the fixed holomorphic neighbourhood just obtained and Cauchy's estimate. Choose
	            \begin{equation*}
	                \mathrm{r}_{\mathrm{H}}\overset{\textnormal{def}}{=}\frac{1}{2}\min\left(\phi_{\star},\mathrm{M}_{2}^{-1},\frac{3\mathrm{R}_{2}}{8\mathrm{A}_{\mathrm{SD}}}\right).
	            \end{equation*}
	            Increasing $\mathrm{T}_{\mathrm{dec}}$ so that $\theta_y(t)\leq\mathrm{r}_{\mathrm{H}}$, the disk $\mathcal{B}_{2\mathrm{r}_{\mathrm{H}}}(0)$ in the $u$-plane is mapped into $\mathcal{B}_{\mathrm{R}_{2}}(0)$ by the local contour parametrisation. Indeed, on this disk we use the analytic continuation of the local polynomial branch $\Phi_y^{\mathrm{SD}}$, not the globally spliced map, and the proof of \eqref{eq:SD-path-stays-close} applies with $|u|$ in place of the real variable $|\theta|$, because the correction term is a monomial in $u$ and $2\mathrm{r}_{\mathrm{H}}\leq\mathrm{M}_{2}^{-1}$. Explicitly, for every complex $u$ with $|u|\leq2\mathrm{r}_{\mathrm{H}}$,
                \begin{equation*}
                    |\Phi_y^{\mathrm{SD}}(u,t)-\mathrm{u}_y(t)|\leq \mathrm{A}_{\mathrm{SD}}|u|\leq 2\mathrm{A}_{\mathrm{SD}}\mathrm{r}_{\mathrm{H}}\leq \frac{3\mathrm{R}_{2}}{8}.
                \end{equation*}
                Together with $\mathrm{u}_y(t)\in\mathcal{B}_{\mathrm{R}_{2}/8}(0)$, this gives $\Phi_y^{\mathrm{SD}}(u,t)\in\mathcal{B}_{\mathrm{R}_{2}}(0)$ for every complex $u$ with $|u|\leq2\mathrm{r}_{\mathrm{H}}$. Hence all points used in Cauchy's estimate lie in the region where \textbf{(iii)} of Proposition~\ref{specificcasesub} applies. In particular, for the corresponding analytic continuation of $H_y$, the segment joining $\mathrm{u}_y(t)$ to $\Phi_y^{\mathrm{SD}}(u,t)$ lies inside $\mathcal{B}_{\mathrm{R}_{2}}(0)$ whenever $|u|\leq2\mathrm{r}_{\mathrm{H}}$, and so the first-derivative bound in \textbf{(iii)} gives
	            \begin{equation*}
	                |H_y(u,t)|\leq\frac{9}{4}\mathrm{M}_{1}\mathrm{M}_{2}\mathrm{R}_{2}\left|\mathfrak{F}_y^{(2)}(\mathrm{u}_y(t),t)\right|.
	            \end{equation*}
	            Here the factor $9\mathrm{R}_{2}/8$ bounds the distance from $\mathrm{u}_y(t)$ to any point of $\mathcal{B}_{\mathrm{R}_{2}}(0)$, while the extra factor $2$ converts $\Sigma(t)^2$ to $\mathfrak{F}_y^{(2)}(\mathrm{u}_y(t),t)$ using \eqref{eq:SD-second-lower}. Thus Cauchy's estimate on the circle $|u|=\mathrm{r}_{\mathrm{H}}$ gives the following fully explicit coefficient bound. Define
	            \begin{equation*}
	                \mathrm{K}_{\mathrm{H}}\overset{\textnormal{def}}{=}\frac{9}{4}\mathrm{M}_{1}\mathrm{M}_{2}\mathrm{R}_{2},\qquad \mathrm{M}_{\mathrm{H}}\overset{\textnormal{def}}{=}\max\left(1,\mathrm{r}_{\mathrm{H}}^{-1},\mathrm{K}_{\mathrm{H}}\right)^2.
	            \end{equation*}
	            Then, since $\mathrm{K}_{\mathrm{H}}\mathrm{r}_{\mathrm{H}}^{-n}\leq\mathrm{M}_{\mathrm{H}}^n$ for every $n\geq1$, we have
	            \begin{equation}\label{eq:Gn-bound}
		                |G_{n,y}(t)|\leq n!\mathrm{M}_{\mathrm{H}}^n\left|\mathfrak{F}_y^{(2)}(\mathrm{u}_y(t),t)\right|,
	                \qquad\text{for all $n\geq3$, $t>\mathrm{T}_{\mathrm{dec}}$ and $y(t)\in\mathscr{C}_{t}$.}
	            \end{equation}
	            We now insert this coefficient bound into the complete Bell polynomial defining $B_{n,y}(t)$. First we record which indices can actually contribute. By definition, $\widetilde{G}_{i,y}(t)=0$ for $i=1,2$. In addition, the parity cancellation proved above gives $G_{i,y}(t)=0$ for every odd $i$ in the infinite-$\mathfrak{N}_y(t)$ case, and for every odd $i\leq2\mathfrak{N}_y(t)+1$ in the finite case. Therefore a term in the Bell expansion for $B_{n,y}(t)$ can be nonzero only if every index with $j_i>0$ belongs to
	            \begin{equation*}
	                A_{n,y}(t)=
                    \begin{cases}
                    \{i\in\{1,\ldots,n\}:\text{$i\geq4$ and $i$ is even}\},&\text{if $\mathfrak{N}_y(t)=\infty$,}\\
                    \{i\in\{1,\ldots,n\}:\text{$i\geq4$ and either $i$ is even or $i\geq2\mathfrak{N}_y(t)+3$ is odd}\},&\text{if $\mathfrak{N}_y(t)<\infty$.}
                    \end{cases}
	            \end{equation*}
	            Thus the factors with $i\notin A_{n,y}(t)$ either vanish, if $j_i>0$, or contribute $1$, if $j_i=0$. For the remaining active indices we apply the bound for $G_{i,y}(t)$ with $i\geq3$. The factorial $i!$ in this bound cancels the denominator $(i!)^{j_i}$ in the Bell polynomial, and the constraint $\sum_{i=1}^{n}ij_i=n$ in the definition of $\mathcal{R}_n$ gives
	            \begin{equation*}
	                \prod_{i\in A_{n,y}(t)}\mathrm{M}_{\mathrm{H}}^{ij_i}=\mathrm{M}_{\mathrm{H}}^{\sum_{i\in A_{n,y}(t)}ij_i}\leq \mathrm{M}_{\mathrm{H}}^{n}.
	            \end{equation*}
	            Consequently, for every $t>\mathrm{T}_{\mathrm{dec}}$ and $y(t)\in\mathscr{C}_{t}$,
	            \begin{equation*}
	                |B_{n,y}(t)|\leq n!\mathrm{M}_{\mathrm{H}}^n\sum_{\underline{j}\in\mathcal{R}_n}\left(\prod_{i=1}^{n}\frac{1}{j_i!}\right)\left|\mathfrak{F}_y^{(2)}(\mathrm{u}_y(t),t)\right|^{\sum_{i\in A_{n,y}(t)}j_i}.
	            \end{equation*}
		            It remains only to estimate the number of active factors appearing in each nonzero Bell term. We increase $\mathrm{T}_{\mathrm{dec}}$, using the lower bound \eqref{eq:SD-second-lower} and $\Sigma(t)\to\infty$ from \textbf{(iv)} of Proposition~\ref{specificcasesub}, so that
	            \begin{equation*}
	                \mathfrak{F}_y^{(2)}(\mathrm{u}_y(t),t)\geq1,\qquad t>\mathrm{T}_{\mathrm{dec}},\quad y(t)\in\mathscr{C}_{t}.
	            \end{equation*}
	            This allows us to replace a smaller power of $\mathfrak{F}_y^{(2)}(\mathrm{u}_y(t),t)$ by a larger one. If $n$ is even, then every active index is at least $4$, so, for each $\underline{j}\in\mathcal{R}_n$,
	            \begin{equation*}
	                4\sum_{i\in A_{n,y}(t)}j_i\leq \sum_{i=1}^{n}ij_i=n,\quad\text{which implies}\quad\sum_{i\in A_{n,y}(t)}j_i\leq \frac{n}{4}.
	            \end{equation*}
	            If $\mathfrak{N}_y(t)<\infty$, $n$ is odd, and $n\geq 2\mathfrak{N}_y(t)+3$, then any nonzero term must contain at least one active odd index; otherwise the sum $\sum_{i=1}^{n}ij_i$ would be even. The smallest possible active odd index is $2\mathfrak{N}_y(t)+3$. After selecting one such odd index, every remaining active index is at least $4$, and hence
	            \begin{equation*}
	                n=\sum_{i=1}^{n}ij_i\geq (2\mathfrak{N}_y(t)+3)+4\left(\sum_{i\in A_{n,y}(t)}j_i-1\right)=2\mathfrak{N}_y(t)-1+4\sum_{i\in A_{n,y}(t)}j_i,
	            \end{equation*}
            and therefore
            \begin{equation*}
                \sum_{i\in A_{n,y}(t)}j_i\leq \frac{n-2\mathfrak{N}_y(t)+1}{4}.
            \end{equation*}
	            We now make the remaining combinatorial constant explicit. Since $\mathcal{R}_n$ is the full partition index set for the complete Bell polynomial, the sum over $\mathcal{R}_n$ is the coefficient of $z^n$ in the exponential generating product
	            \begin{equation*}
		                \sum_{\underline{j}\in\mathcal{R}_n}\prod_{i=1}^{n}\frac{1}{j_i!}=[z^n]\prod_{i=1}^{\infty}\sum_{j_i=0}^{\infty}\frac{z^{ij_i}}{j_i!}=[z^n]\exp\left(\sum_{i=1}^{\infty}z^i\right)=[z^n]\exp\left(\frac{z}{1-z}\right).
	            \end{equation*}
	            Expanding the last expression once more gives the closed form
	            \begin{equation*}
	                \sum_{\underline{j}\in\mathcal{R}_n}\prod_{i=1}^{n}\frac{1}{j_i!}=\sum_{k=1}^{n}\frac{1}{k!}\binom{n-1}{k-1},\qquad n\geq1.
	            \end{equation*}
	            We therefore define the explicit constant
	            \begin{equation*}
	                \mathrm{B}_{n}\overset{\textnormal{def}}{=}n!\mathrm{M}_{\mathrm{H}}^n\sum_{k=1}^{n}\frac{1}{k!}\binom{n-1}{k-1},\qquad n\geq1.
	            \end{equation*}
	            This constant is independent of $t$ and $y$. Moreover, using $\binom{n-1}{k-1}\leq n^{k-1}/(k-1)!$ and then setting $m=k-1$, we have
	            \begin{equation*}
	                \sum_{k=1}^{n}\frac{1}{k!}\binom{n-1}{k-1}\leq\sum_{m=0}^{\infty}\frac{n^m}{(m!)^2}\leq\sum_{m=0}^{\infty}\frac{(2\sqrt{n})^{2m}}{(2m)!}\leq \mathrm{e}^{2\sqrt{n}},
	            \end{equation*}
	            and hence $\mathrm{B}_{n}\leq n!\mathrm{M}_{\mathrm{H}}^n\mathrm{e}^{2\sqrt{n}}$. Thus the dependence of the Bell coefficient bound on $n$ is fully explicit. Combining this with the two preceding active-index estimates gives
	            \begin{equation}\label{absB}
	                \begin{aligned}
	                        &|B_{n,y}(t)|\leq \mathrm{B}_{n}\left|\mathfrak{F}_y^{(2)}(\mathrm{u}_y(t),t)\right|^{\frac{n}{4}},\qquad&\text{for all even $n\geq1$,}\\
                        &|B_{n,y}(t)|\leq \mathrm{B}_{n}\left|\mathfrak{F}_y^{(2)}(\mathrm{u}_y(t),t)\right|^{\frac{n-2\mathfrak{N}_y(t)+1}{4}},\qquad&\text{if $\mathfrak{N}_y(t)<\infty$ and $n>2\mathfrak{N}_y(t)+1$ is odd.}
                    \end{aligned}
            \end{equation}
           
            \textbf{Step 4c: Rescaling, main terms, Gaussian moments, and local error bounds.}
	            Continuing first in the finite-$\mathfrak{N}_y(t)$ case, set
                \begin{equation*}
                    \mathrm{A}_y(t)\overset{\textnormal{def}}{=}\theta_y(t)\sqrt{\mathfrak{F}_y^{(2)}(\mathrm{u}_y(t),t)}
                    =
                    \mathfrak{F}_y^{(2)}(\mathrm{u}_y(t),t)^{1/6}.
                \end{equation*}
                We now rescale by the quadratic scale. Making the change of variables
	            \begin{equation*}
	                u=\mathfrak{F}_y^{(2)}(\mathrm{u}_y(t),t)^{-1/2}z
	            \end{equation*}
	            in \eqref{parametriccontourint}, and then renaming $z$ back to $u$, gives
	            \begin{equation}\label{parametricthing}
	                \begin{aligned}
	                    \mathcal{I}_{n,y}(t)
	                    &=
	                    \frac{\exp\left(\mathfrak{F}_{y}(\mathrm{u}_{y}(t),t)\right)}{\lambda_{y(t)}2\pi\mathrm{i}}
		                    \int_{-\mathrm{A}_y(t)}^{\mathrm{A}_y(t)}
	                    \Bigg(
	                    \sum_{k=0}^{n+1}
		                    \mathrm{i}^{n-k+1}\left(\frac{\mathfrak{a}_y(t)}{(2\mathfrak{N}_y(t))!}\right)^k\mathfrak{c}_{n,k,y}(t)
	                    \frac{u^{2\mathfrak{N}_y(t)k+n-k}}{\mathfrak{F}_y^{(2)}(\mathrm{u}_y(t),t)^{\frac{2\mathfrak{N}_y(t)k+n-k+1}{2}}}
	                    \Bigg)\\
	                    &\quad\times
	                    \Bigg(
	                    \sum_{m=0}^{\infty}B_{m,y}(t)\frac{u^m}{m!\,\mathfrak{F}_y^{(2)}(\mathrm{u}_y(t),t)^{m/2}}
	                    \Bigg)
	                    \exp\left(-\frac{u^2}{2}\right)\mathrm{d}u.
	                \end{aligned}
	            \end{equation}
                (The case $\mathfrak{N}_y(t)=\infty$ has already been handled in Step~4a.) The Taylor series is evaluated at the original local parameter $\theta=\mathfrak{F}_y^{(2)}(\mathrm{u}_y(t),t)^{-1/2}u$. Since $|u|\leq\mathrm{A}_y(t)$ is equivalent to $|\theta|\leq\theta_y(t)$, the expansion remains inside the fixed holomorphic neighbourhood used in Step~4b. By the lower bound \eqref{eq:SD-second-lower}, $\mathrm{A}_y(t)=\mathfrak{F}_y^{(2)}(\mathrm{u}_y(t),t)^{1/6}\geq2^{-1/6}\Sigma(t)^{1/3}$. Thus, for every fixed $\mathrm{L}\geq1$, choosing $\mathrm{T}_{\mathrm{dec}}$ so that $2^{-1/6}\Sigma(t)^{1/3}\geq\mathrm{L}$ ensures that the rescaled integration interval contains $[-\mathrm{L},\mathrm{L}]$ for every $t>\mathrm{T}_{\mathrm{dec}}$ and every path $y\in\mathscr{P}$ satisfying $y(t)\in\mathscr{C}_{t}$.
               
                 The estimates below repeatedly use the two facts established above: the definition $\mathrm{A}_y(t)=\mathfrak{F}_y^{(2)}(\mathrm{u}_y(t),t)^{1/6}$ and the lower curvature bound $\mathfrak{F}_y^{(2)}(\mathrm{u}_y(t),t)\geq\Sigma(t)^2/2$ from \eqref{eq:SD-second-lower}. First, define the explicit Gaussian moment constants, for $r\in\mathbb{Z}_+$
            \begin{equation*}
                \mathrm{G}_{r}\overset{\textnormal{def}}{=}\int_{-\infty}^{\infty}|u|^{r}\exp\left(-\frac{u^2}{2}\right)\mathrm{d}u=2^{\frac{r+1}{2}}\Gamma\left(\frac{r+1}{2}\right),\quad\text{and set}\quad\mathrm{R}_{\mathrm{B}}\overset{\textnormal{def}}{=}2\mathrm{e}^{3}\mathrm{M}_{\mathrm{H}}^{2}.
            \end{equation*}
            We claim that, after increasing \(\mathrm{T}_{\mathrm{dec}}\), the Bell factor in \eqref{parametricthing} satisfies
            \begin{equation}\label{eq:Bell-tail-explicit}
                \left|\sum_{m=0}^{\infty}B_{m,y}(t)\frac{u^m}{m!\,\mathfrak{F}_y^{(2)}(\mathrm{u}_y(t),t)^{m/2}}-1\right|\leq \mathrm{R}_{\mathrm{B}}u^2\mathfrak{F}_y^{(2)}(\mathrm{u}_y(t),t)^{-1/2},\qquad |u|\leq\mathrm{A}_y(t),
            \end{equation}
            and, on the same range,
            \begin{equation}\label{eq:Bell-factor-two}
                \left|\sum_{m=0}^{\infty}B_{m,y}(t)\frac{u^m}{m!\,\mathfrak{F}_y^{(2)}(\mathrm{u}_y(t),t)^{m/2}}\right|\leq2.
            \end{equation}
            To prove \eqref{eq:Bell-tail-explicit}, first note that $B_{1,y}(t)=B_{2,y}(t)=0$ because $\widetilde{G}_{1,y}(t)=\widetilde{G}_{2,y}(t)=0$. For $m\geq3$, combine the vanishing of the inadmissible odd Bell coefficients with \eqref{absB}. If $m$ is even, \eqref{absB} gives the power $m/4$. If $m$ is odd and $\mathfrak{N}_y(t)<\infty$ with $m>2\mathfrak{N}_y(t)+1$, then the exponent in \eqref{absB} is at most $m/4$; in all other odd cases, $B_{m,y}(t)=0$. Thus, once $\mathfrak{F}_y^{(2)}(\mathrm{u}_y(t),t)\geq1$, the single bound $|B_{m,y}(t)|\leq\mathrm{B}_{m}\mathfrak{F}_y^{(2)}(\mathrm{u}_y(t),t)^{m/4}$ holds for every $m\geq3$. The explicit formula for $\mathrm{B}_{m}$ gives $\mathrm{B}_{m}/m!\leq \mathrm{e}(\mathrm{e}\mathrm{M}_{\mathrm{H}})^m$, because $\mathrm{e}^{2\sqrt{m}}\leq\mathrm{e}^{m+1}$ for $m\geq1$. Hence, for $|u|\leq\mathrm{A}_y(t)=\mathfrak{F}_y^{(2)}(\mathrm{u}_y(t),t)^{1/6}$,
            \begin{equation*}
                \begin{aligned}
                    \left|\sum_{m=1}^{\infty}B_{m,y}(t)\frac{u^m}{m!\,\mathfrak{F}_y^{(2)}(\mathrm{u}_y(t),t)^{m/2}}\right|
                    &\leq\mathrm{e}\sum_{m=3}^{\infty}\left(\mathrm{e}\mathrm{M}_{\mathrm{H}}|u|\mathfrak{F}_y^{(2)}(\mathrm{u}_y(t),t)^{-1/4}\right)^m.
                \end{aligned}
            \end{equation*}
            Since $|u|\mathfrak{F}_y^{(2)}(\mathrm{u}_y(t),t)^{-1/4}\leq\mathfrak{F}_y^{(2)}(\mathrm{u}_y(t),t)^{-1/12}$, we increase $\mathrm{T}_{\mathrm{dec}}$ so that $\mathrm{e}\mathrm{M}_{\mathrm{H}}\mathfrak{F}_y^{(2)}(\mathrm{u}_y(t),t)^{-1/12}\leq1/2$. Enlarging the resulting constant, the same geometric-series estimate gives
            \begin{equation*}
                \mathrm{e}\sum_{m=2}^{\infty}\left(\mathrm{e}\mathrm{M}_{\mathrm{H}}|u|\mathfrak{F}_y^{(2)}(\mathrm{u}_y(t),t)^{-1/4}\right)^m\leq2\mathrm{e}^{3}\mathrm{M}_{\mathrm{H}}^2u^2\mathfrak{F}_y^{(2)}(\mathrm{u}_y(t),t)^{-1/2},
            \end{equation*}
            which proves \eqref{eq:Bell-tail-explicit}. Increasing \(\mathrm{T}_{\mathrm{dec}}\) again so that \(\mathrm{R}_{\mathrm{B}}\mathfrak{F}_y^{(2)}(\mathrm{u}_y(t),t)^{-1/6}\leq1\), and using \(|u|\leq\mathfrak{F}_y^{(2)}(\mathrm{u}_y(t),t)^{1/6}\), gives \eqref{eq:Bell-factor-two}.
	           
	           We also need an explicit bound for the deviation of the local contour from the vertical tangent. Write
	            \begin{equation*}
	                \Phi_y^{\mathrm{SD}}(\theta,t)-\mathrm{u}_y(t)=\mathrm{i}\theta+\chi_y(\theta,t),
	            \end{equation*}
	            where \(\chi_y(\theta,t)=\mathfrak{a}_y(t)\theta^{2\mathfrak{N}_y(t)}/(2\mathfrak{N}_y(t))!\) if \(\mathfrak{N}_y(t)<\infty\), and \(\chi_y(\theta,t)=0\) if \(\mathfrak{N}_y(t)=\infty\). Since $\mathrm{T}_{\mathrm{dec}}\geq\mathrm{T}_{\mathrm{C}}\geq\mathrm{T}_{\mathrm{SD}}$, the contour-deviation estimate \eqref{eq:SD-chi-bound} from Proposition~\ref{descentcurve} gives, whenever \(|\theta|\leq\mathrm{M}_{2}^{-1}\),
	            \begin{equation*}
	                |\chi_y(\theta,t)|\leq \mathrm{K}_{\mathrm{SD}}|\theta|^2,
	                \qquad
	                |\chi_y^{(1)}(\theta,t)|\leq \mathrm{K}_{\mathrm{SD}}|\theta|.
	            \end{equation*}
	         We increase the common threshold so that \(\theta_y(t)\leq\min(1,\mathrm{M}_{2}^{-1})\). Then, for each fixed \(n\in\mathbb{Z}_+\), the product coming from the parametrisation satisfies
            \begin{equation}\label{eq:poly-correction-bound}
                \left|\left(\Phi_y^{\mathrm{SD}}(\theta,t)-\mathrm{u}_y(t)\right)^n\frac{\mathrm{d}}{\mathrm{d}\theta}\Phi_y^{\mathrm{SD}}(\theta,t)-\mathrm{i}^{n+1}\theta^n\right|\leq \mathrm{K}_{\mathrm{poly},n}|\theta|^{n+1},
            \end{equation}
            where
            \begin{equation*}
                \mathrm{K}_{\mathrm{poly},n}\overset{\textnormal{def}}{=}\mathrm{K}_{\mathrm{SD}}\left(1+n(1+\mathrm{K}_{\mathrm{SD}})(2+\mathrm{K}_{\mathrm{SD}})^n\right).
            \end{equation*}
            Indeed, subtract and add \((\mathrm{i}\theta)^n(\mathrm{i}+\chi_y^{(1)}(\theta,t))\). The term containing \(\chi_y^{(1)}\) is bounded by \(\mathrm{K}_{\mathrm{SD}}|\theta|^{n+1}\). For the remaining term, the elementary inequality \(|a^n-b^n|\leq n|a-b|(|a|+|b|)^{n-1}\), applied with \(a=\mathrm{i}\theta+\chi_y(\theta,t)\) and \(b=\mathrm{i}\theta\), together with \(|\mathrm{i}\theta+\chi_y(\theta,t)|\leq(1+\mathrm{K}_{\mathrm{SD}})|\theta|\), gives the displayed bound.
         
             We next identify the lowest-order term for even \(n\). The contribution obtained by keeping the vertical factor \(\mathrm{i}^{n+1}\theta^n\) and the zeroth Bell coefficient is
            \begin{equation*}
                \mathcal{I}^{\mathrm{main,even}}_{n,y}(t)\overset{\textnormal{def}}{=}\frac{\mathrm{i}^{n}\exp\left(\mathfrak{F}_{y}(\mathrm{u}_{y}(t),t)\right)}{\lambda_{y(t)}2\pi\mathfrak{F}_y^{(2)}(\mathrm{u}_y(t),t)^{\frac{n+1}{2}}}\int_{-\mathrm{A}_y(t)}^{\mathrm{A}_y(t)}u^n\exp\left(-\frac{u^2}{2}\right)\mathrm{d}u.
            \end{equation*}
            If \(n\) is odd, this term is identically zero because the integration interval is symmetric and \(u^n\exp(-u^2/2)\) is odd. This observation avoids any need for a uniform bound on \(\mathfrak{N}_y(t)\): all odd contributions are treated as errors forced either by the Bell tail \eqref{eq:Bell-tail-explicit} or by the contour correction \eqref{eq:poly-correction-bound}. Combining \eqref{eq:Bell-tail-explicit}, \eqref{eq:Bell-factor-two}, and \eqref{eq:poly-correction-bound}, and then using the Gaussian moments \(\mathrm{G}_{n+1}\) and \(\mathrm{G}_{n+2}\), gives, for every \(n\in\mathbb{Z}_+\),
            \begin{equation}\label{eq:local-remainder-explicit}
                \left|\mathcal{I}_{n,y}(t)-\mathcal{I}^{\mathrm{main,even}}_{n,y}(t)\right|\leq\frac{\mathrm{D}_{\mathrm{rem},n}}{\lambda_{y(t)}\mathfrak{F}_y^{(2)}(\mathrm{u}_y(t),t)^{\frac{n+2}{2}}}\exp\left(\mathfrak{F}_{y}(\mathrm{u}_{y}(t),t)\right),
            \end{equation}
            where the explicit constant is
            \begin{equation*}
                \mathrm{D}_{\mathrm{rem},n}\overset{\textnormal{def}}{=}\frac{1}{2\pi}\left(\mathrm{R}_{\mathrm{B}}\mathrm{G}_{n+2}+2\mathrm{K}_{\mathrm{poly},n}\mathrm{G}_{n+1}\right).
            \end{equation*}
            To see \eqref{eq:local-remainder-explicit}, decompose the difference from the leading vertical Gaussian term as
            \begin{equation*}
                \textnormal{vertical factor}\times\textnormal{Bell tail}
                \quad+\quad
                \textnormal{contour correction}\times\textnormal{full Bell factor}.
            \end{equation*}
            The first product is controlled by \eqref{eq:Bell-tail-explicit}: the vertical part multiplied by the Bell tail contributes at most
            \begin{equation*}
                \frac{\mathrm{R}_{\mathrm{B}}\mathrm{G}_{n+2}}{2\pi\lambda_{y(t)}\mathfrak{F}_y^{(2)}(\mathrm{u}_y(t),t)^{\frac{n+2}{2}}}\exp\left(\mathfrak{F}_y(\mathrm{u}_y(t),t)\right).
            \end{equation*}
            The contour-correction part is bounded by \eqref{eq:poly-correction-bound}; after the change of variables $\theta=\mathfrak{F}_y^{(2)}(\mathrm{u}_y(t),t)^{-1/2}u$, the extra power $|\theta|^{n+1}\mathrm{d}\theta$ gives $\mathfrak{F}_y^{(2)}(\mathrm{u}_y(t),t)^{-(n+2)/2}|u|^{n+1}\mathrm{d}u$, and \eqref{eq:Bell-factor-two} supplies the factor $2$. We can now evaluate the even leading term using Gaussian moments. After increasing \(\mathrm{T}_{\mathrm{dec}}\), we have \(\mathrm{A}_y(t)\geq1\) for every \(t>\mathrm{T}_{\mathrm{dec}}\) and \(y(t)\in\mathscr{C}_{t}\). For \(p\in\mathbb{Z}_+\), set
            \begin{equation*}
                \mathrm{C}_{\mathrm{tail},p}\overset{\textnormal{def}}{=}2\sum_{\ell=0}^{p}\binom{p}{\ell}\ell!.
            \end{equation*}
            If \(\mathrm{A}\geq1\), then, for \(u=\mathrm{A}+v\) with \(v\geq0\), using \((\mathrm{A}+v)^p\leq\sum_{\ell=0}^{p}\binom{p}{\ell}\mathrm{A}^{p-\ell}v^{\ell}\leq\mathrm{A}^{p}\sum_{\ell=0}^{p}\binom{p}{\ell}v^{\ell}\), and using \(\exp(-(\mathrm{A}+v)^2/2)\leq\exp(-\mathrm{A}^2/2)\exp(-v)\), we get
            \begin{equation*}
                2\int_{\mathrm{A}}^{\infty}u^p\exp\left(-\frac{u^2}{2}\right)\mathrm{d}u\leq \mathrm{C}_{\mathrm{tail},p}\mathrm{A}^{p}\exp\left(-\frac{\mathrm{A}^2}{2}\right).
            \end{equation*}
            For each fixed even \(n\), define a further lower time \(\mathrm{T}_{\mathrm{dec},n}^{\mathrm{tail}}\geq\mathrm{T}_{\mathrm{dec}}\) so that
            \begin{equation*}
                \mathrm{C}_{\mathrm{tail},n}\mathfrak{F}_y^{(2)}(\mathrm{u}_y(t),t)^{n/6}\exp\left(-\frac{1}{2}\mathfrak{F}_y^{(2)}(\mathrm{u}_y(t),t)^{1/3}\right)\leq \mathfrak{F}_y^{(2)}(\mathrm{u}_y(t),t)^{-1/2},
                \qquad t>\mathrm{T}_{\mathrm{dec},n}^{\mathrm{tail}},\quad y(t)\in\mathscr{C}_{t}.
            \end{equation*}
            This is possible because $\mathfrak{F}_y^{(2)}(\mathrm{u}_y(t),t)\geq\Sigma(t)^2/2$ and $\Sigma(t)\to\infty$. Thus, for even $n$, the truncation error in replacing the finite Gaussian moment by the full Gaussian moment contributes at most
            \begin{equation*}
                \frac{1}{2\pi\lambda_{y(t)}\mathfrak{F}_y^{(2)}(\mathrm{u}_y(t),t)^{\frac{n+2}{2}}}\exp\left(\mathfrak{F}_y(\mathrm{u}_y(t),t)\right).
            \end{equation*}
            Combining this truncation estimate with \eqref{eq:local-remainder-explicit}, and defining
            \begin{equation*}
                \mathrm{D}_{\mathrm{loc},n}\overset{\textnormal{def}}{=}\mathrm{D}_{\mathrm{rem},n}+\frac{1}{2\pi},
            \end{equation*}
	            gives, for even \(n\),
		            \begin{equation}\label{eq:local-even-In}
		                \begin{aligned}
		                    &\left|\mathcal{I}_{n,y}(t)-\frac{\mathrm{i}^{n}(n-1)!!}{\lambda_{y(t)}\sqrt{2\pi}\,\mathfrak{F}_y^{(2)}(\mathrm{u}_y(t),t)^{\frac{n+1}{2}}}\exp\left(\mathfrak{F}_{y}(\mathrm{u}_{y}(t),t)\right)\right|\\
		                    &\qquad\leq\frac{\mathrm{D}_{\mathrm{loc},n}}{\lambda_{y(t)}\,\mathfrak{F}_y^{(2)}(\mathrm{u}_y(t),t)^{\frac{n+2}{2}}}\exp\left(\mathfrak{F}_{y}(\mathrm{u}_{y}(t),t)\right).
		                \end{aligned}
		            \end{equation}
            For the special case \(n=0\), the additional truncation threshold is absorbed into the common $\mathrm{T}_{\mathrm{dec}}$, because the resulting estimate is used below to lower-bound the normalising probability.
            If \(n\) is odd, then \(\mathcal{I}^{\mathrm{main,even}}_{n,y}(t)=0\), and \eqref{eq:local-remainder-explicit} gives directly
	            \begin{equation}\label{eq:local-odd-In}
	                \left|\mathcal{I}_{n,y}(t)\right|\leq\frac{\mathrm{D}_{\mathrm{rem},n}}{\lambda_{y(t)}\,\mathfrak{F}_y^{(2)}(\mathrm{u}_y(t),t)^{\frac{n+2}{2}}}\exp\left(\mathfrak{F}_{y}(\mathrm{u}_{y}(t),t)\right).
	            \end{equation}
	            In particular, taking \(n=0\) in \eqref{eq:local-even-In}, we obtain
	            \begin{equation}\label{eq:I0-local-bound}
	                \left|\mathcal{I}_{0,y}(t)-\frac{1}{\lambda_{y(t)}\sqrt{2\pi\,\mathfrak{F}_y^{(2)}(\mathrm{u}_y(t),t)}}\exp\left(\mathfrak{F}_{y}(\mathrm{u}_{y}(t),t)\right)\right|\leq\frac{\mathrm{D}_{\mathrm{loc},0}}{\lambda_{y(t)}\,\mathfrak{F}_y^{(2)}(\mathrm{u}_y(t),t)}\exp\left(\mathfrak{F}_{y}(\mathrm{u}_{y}(t),t)\right).
	            \end{equation}
          
            \textbf{Step 5: exponential smallness on the complementary arc.}
	            We next estimate the complementary contour integral $I_{2,x,y}(t)$. We apply Lemma~\ref{lem:Gamma2-Re-max} with $\vartheta=\theta_y(t)$. Its hypotheses are all in force: $y(t)\in\mathscr{C}_{t}$ by assumption; $\theta_{\star}\in(0,\theta_{\mathrm C})$ by construction; Corollary~\ref{steepestdescentpath} applies because $\mathrm{T}_{\mathrm{dec}}\geq\mathrm{T}_{\mathrm C}$; and Step~1 imposed $\theta_y(t)<\phi_{\star}=\mathrm{R}_{\mathrm{C}}\sin\theta_{\star}$. Hence the phase bound \eqref{eq:Gamma2-Re-max} holds. Moreover $\Gamma_{2}^{(\theta_y(t))}(t)$ agrees with the complementary contour $\Gamma_{2,y}(t)$ from Step~2 up to the two cutoff endpoints, which do not affect the contour integral; the endpoint-reduction estimate is applied on the corresponding closure.
		            Since $\theta_y(t)<\phi_{\star}$, both splice points $\Phi_y(\pm\theta_y(t),t,\theta_{\star})$ lie on the local segment where \eqref{globalapproximate}--\eqref{bellexpansiontwo} holds. Fix $\sigma\in\{-1,1\}$. The coefficient bound \eqref{eq:Gn-bound}, together with the threshold condition $\mathrm{M}_{\mathrm{H}}\theta_y(t)\leq1/2$, gives the traceable tail bound
		            \begin{equation*}
			                \left|\sum_{m=3}^{\infty}G_{m,y}(t)\frac{(\sigma\theta_y(t))^m}{m!}\right|
		                \leq
		                2\mathrm{M}_{\mathrm{H}}^3\,\mathfrak{F}_y^{(2)}(\mathrm{u}_y(t),t)\theta_y(t)^3.
		            \end{equation*}
		            Since $\theta_y(t)=\mathfrak{F}_y^{(2)}(\mathrm{u}_y(t),t)^{-1/3}$, we can rewrite the right-hand side as
	                \begin{equation*}
	                    2\mathrm{M}_{\mathrm{H}}^3\mathfrak{F}_y^{(2)}(\mathrm{u}_y(t),t)^{-1/3}
	                    \cdot
	                    \mathfrak{F}_y^{(2)}(\mathrm{u}_y(t),t)\theta_y(t)^2.
	                \end{equation*}
	                Increasing $\mathrm{T}_{\mathrm{dec}}$ so that
	                \begin{equation*}
	                    2\mathrm{M}_{\mathrm{H}}^3\mathfrak{F}_y^{(2)}(\mathrm{u}_y(t),t)^{-1/3}\leq\frac{1}{4},
	                    \qquad
	                    t>\mathrm{T}_{\mathrm{dec}},\quad y(t)\in\mathscr{C}_{t},
	                \end{equation*}
                we obtain, for this same fixed $\sigma$,
	            \begin{equation*}
	                \mathfrak{Re}\left[H_y(\sigma\theta_y(t),t)\right]
	                \leq
		                -\frac{1}{2}\mathfrak{F}_y^{(2)}(\mathrm{u}_y(t),t)\theta_y(t)^2
		                +
		                2\mathrm{M}_{\mathrm{H}}^3\mathfrak{F}_y^{(2)}(\mathrm{u}_y(t),t)\theta_y(t)^3
		                \leq
		                -\frac{1}{4}\mathfrak{F}_y^{(2)}(\mathrm{u}_y(t),t)\theta_y(t)^2.
		            \end{equation*}
	            Therefore
	            \begin{equation*}
	                \mathfrak{Re}\left[\mathfrak{F}_y(\Phi_y(\sigma\theta_y(t),t,\theta_{\star}),t)\right]
	                \leq
	                \mathfrak{F}_{y}(\mathrm{u}_{y}(t),t)-\frac{1}{4}\mathfrak{F}_y^{(2)}(\mathrm{u}_y(t),t)^{1/3}.
	            \end{equation*}
		            Since the right-hand side is independent of $\sigma$, the endpoint maximum satisfies the same estimate:
                    \begin{equation*}
                        \max_{\sigma\in\{-1,1\}}\mathfrak{Re}\left[\mathfrak{F}_y(\Phi_y(\sigma\theta_y(t),t,\theta_{\star}),t)\right]\leq\mathfrak{F}_{y}(\mathrm{u}_{y}(t),t)-\frac{1}{4}\mathfrak{F}_y^{(2)}(\mathrm{u}_y(t),t)^{1/3}.
                    \end{equation*}
	                    We now make the contour and polynomial prefactor constants explicit. The parameter interval for $\Gamma_y(t)$ has length
	                    \begin{equation*}
	                        2\mathrm{U}_{\Gamma},\qquad \mathrm{U}_{\Gamma}\overset{\textnormal{def}}{=}\frac{\phi_{\star}}{\theta_{\star}}\pi.
	                    \end{equation*}
	                    On the local part the derivative of the parametrisation is bounded by $\mathrm{A}_{\mathrm{SD}}$, while on the circular part the change of variables in Corollary~\ref{steepestdescentpath} gives the bound $\theta_{\star}/\sin(\theta_{\star})=\mathrm{R}_{\mathrm{C}}\theta_{\star}/\phi_{\star}$. Thus we set
	                    \begin{equation*}
	                        \mathrm{A}_{\Gamma}\overset{\textnormal{def}}{=}\max\left(\mathrm{A}_{\mathrm{SD}},\frac{\mathrm{R}_{\mathrm{C}}\theta_{\star}}{\phi_{\star}}\right),\qquad \mathrm{L}_{\Gamma}\overset{\textnormal{def}}{=}2\mathrm{U}_{\Gamma}\mathrm{A}_{\Gamma},\qquad \mathrm{R}_{\Gamma}\overset{\textnormal{def}}{=}\frac{\mathrm{R}_{2}}{8}+\mathrm{U}_{\Gamma}\mathrm{A}_{\Gamma}.
	                    \end{equation*}
	                    Since $|\mathrm{u}_{y}(t)|\leq\mathrm{R}_{2}/8$, every point of $\Gamma_y(t)$ has modulus at most $\mathrm{R}_{\Gamma}$, and the length of $\Gamma_y(t)$ is at most $\mathrm{L}_{\Gamma}$. Also, by \textbf{(A1)}, $\lambda_k\geq\mathfrak{L}$ for every $k$, and hence
	                    \begin{equation*}
	                        |\rho_x(w)|=\prod_{k=0}^{x-1}\left|1-\frac{w}{\lambda_k}\right|\leq\left(1+\frac{\mathrm{R}_{\Gamma}}{\mathfrak{L}}\right)^x,\qquad |w|\leq\mathrm{R}_{\Gamma}.
	                    \end{equation*}
	                    Therefore, with
	                    \begin{equation*}
	                        \mathrm{C}_{x,\Gamma}\overset{\textnormal{def}}{=}\frac{\mathrm{L}_{\Gamma}}{2\pi\mathfrak{L}}\left(1+\frac{\mathrm{R}_{\Gamma}}{\mathfrak{L}}\right)^x,\qquad x\in\mathbb{Z}_+,
	                    \end{equation*}
	                    the prefactor and contour length in the definition of $I_{2,x,y}(t)$ are bounded explicitly, uniformly over $t$ and $y(t)\in\mathscr{C}_{t}$:
		            \begin{equation}\label{boundforI2}
		                \left|I_{2,x,y}(t)\right|
		                \leq
		                \mathrm{C}_{x,\Gamma}\exp\left(\max_{\sigma\in\{-1,1\}}\mathfrak{Re}\left[\mathfrak{F}_y(\Phi_y(\sigma\theta_y(t),t,\theta_{\star}),t)\right]\right).
	            \end{equation}
	            Combining the previous displays yields the pointwise estimate
	            \begin{equation*}
	                \left|I_{2,x,y}(t)\right|
	                \leq
	                \mathrm{C}_{x,\Gamma}\exp\left(\mathfrak{F}_{y}(\mathrm{u}_{y}(t),t)-\frac{1}{4}\mathfrak{F}_y^{(2)}(\mathrm{u}_y(t),t)^{1/3}\right).
	            \end{equation*}
				            In anticipation of normalising this pointwise bound by $\mathbb{P}_{0}[\mathscr{X}(t)=y(t)]$, define the statement constant
				            \begin{equation}\label{eq:Cx-explicit}
				                \mathrm{C}_{x}\overset{\textnormal{def}}{=}2\sqrt{2\pi}\,\mathfrak{B}_{1}\mathrm{C}_{x,\Gamma},\qquad x\in\mathbb{Z}_+.
				            \end{equation}
			            Increasing $\mathrm{T}_{\mathrm{dec}}$ so that
			            \begin{equation*}
			                \mathfrak{F}_y^{(2)}(\mathrm{u}_y(t),t)^{1/2}\exp\left(-\frac{1}{4}\mathfrak{F}_y^{(2)}(\mathrm{u}_y(t),t)^{1/3}\right)\leq\exp\left(-\frac{1}{8}\mathfrak{F}_y^{(2)}(\mathrm{u}_y(t),t)^{1/3}\right),
			            \end{equation*}
			            for every $t>\mathrm{T}_{\mathrm{dec}}$ and $y(t)\in\mathscr{C}_{t}$, and then using the lower curvature estimate \eqref{eq:SD-second-lower}, we get
			            \begin{equation*}
				                \exp\left(-\frac{1}{8}\mathfrak{F}_y^{(2)}(\mathrm{u}_y(t),t)^{1/3}\right)\leq\mathrm{e}^{-\mathrm{c}_{\varphi}\Sigma(t)^{2/3}},
				            \end{equation*}
				            because $\mathrm{c}_{\varphi}=1/(8\,2^{1/3})$. The only missing ingredient for the normalised form is a lower bound on $\mathbb{P}_{0}[\mathscr{X}(t)=y(t)]$, which we prove next from the $n=0$ local integral.
	            
                \textbf{Step 6: from local integrals to $\eta_{n,y}(t)$ and $\Sigma(t)$.}
	            The denominator in the definition of $\eta_{n,y}(t)$ is the full transition probability from the origin, not only the local Gaussian integral. We therefore first prove an explicit lower bound for this denominator. Let
	            \begin{equation*}
	                \mathcal{L}_{0,y}(t)\overset{\textnormal{def}}{=}\left(\lambda_{y(t)}\sqrt{2\pi\,\mathfrak{F}_y^{(2)}(\mathrm{u}_y(t),t)}\right)^{-1}\exp\left(\mathfrak{F}_y(\mathrm{u}_y(t),t)\right)
	            \end{equation*}
	            denote the leading term in \eqref{eq:I0-local-bound}. Dividing the error term in \eqref{eq:I0-local-bound} by $\mathcal{L}_{0,y}(t)$ gives
	            \begin{equation*}
	                \left|\mathcal{I}_{0,y}(t)-\mathcal{L}_{0,y}(t)\right|\leq\sqrt{2\pi}\,\mathrm{D}_{\mathrm{loc},0}\mathfrak{F}_y^{(2)}(\mathrm{u}_y(t),t)^{-1/2}\mathcal{L}_{0,y}(t).
	            \end{equation*}
	            Similarly, the pointwise estimate from Step~5 and the upper rate bound $\lambda_{y(t)}\leq\mathfrak{B}_1$ from \textbf{(A1)} give
	            \begin{equation*}
	                |I_{2,0,y}(t)|\leq\sqrt{2\pi}\,\mathfrak{B}_{1}\mathrm{C}_{0,\Gamma}\mathfrak{F}_y^{(2)}(\mathrm{u}_y(t),t)^{1/2}\exp\left(-\frac{1}{4}\mathfrak{F}_y^{(2)}(\mathrm{u}_y(t),t)^{1/3}\right)\mathcal{L}_{0,y}(t).
	            \end{equation*}
	            Increasing $\mathrm{T}_{\mathrm{dec}}$, using \eqref{eq:SD-second-lower} and $\Sigma(t)\to\infty$, we may assume that both coefficients on the right-hand sides of the preceding two displays are at most $1/4$, uniformly for $y(t)\in\mathscr{C}_{t}$. Hence
	            \begin{equation}\label{eq:P0-lower}
	                |\mathcal{I}_{0,y}(t)|\geq\frac{3}{4}\mathcal{L}_{0,y}(t),\qquad \mathbb{P}_{0}[\mathscr{X}(t)=y(t)]\geq\frac{1}{2}\mathcal{L}_{0,y}(t),
	            \end{equation}
	            and also
	            \begin{equation}\label{eq:I20-dominated-by-local}
	                \left|I_{2,0,y}(t)\right|\leq \frac{1}{2}\left|\mathcal{I}_{0,y}(t)\right|,\qquad t>\mathrm{T}_{\mathrm{dec}},\quad y(t)\in\mathscr{C}_{t}.
	            \end{equation}
	            The second inequality in \eqref{eq:P0-lower} follows because $\mathcal{L}_{0,y}(t)>0$ is real, $\mathbb{P}_{0}[\mathscr{X}(t)=y(t)]$ is a probability, and the identity $\mathbb{P}_{0}[\mathscr{X}(t)=y(t)]=\mathcal{I}_{0,y}(t)+I_{2,0,y}(t)$ together with the two one-quarter bounds gives $\left|\mathbb{P}_{0}[\mathscr{X}(t)=y(t)]-\mathcal{L}_{0,y}(t)\right|\leq\mathcal{L}_{0,y}(t)/2$. The domination \eqref{eq:I20-dominated-by-local} is stronger than what is needed for \eqref{etaratio}, but it records explicitly that the complementary contour cannot change the leading local mass. Combining \eqref{eq:P0-lower} with the pointwise estimate for $I_{2,x,y}(t)$ in Step~5 gives, for every fixed $x\in\mathbb{Z}_+$,
	            \begin{equation*}
	                \frac{|I_{2,x,y}(t)|}{\mathbb{P}_{0}[\mathscr{X}(t)=y(t)]}\leq2\sqrt{2\pi}\,\mathfrak{B}_{1}\mathrm{C}_{x,\Gamma}\mathfrak{F}_y^{(2)}(\mathrm{u}_y(t),t)^{1/2}\exp\left(-\frac{1}{4}\mathfrak{F}_y^{(2)}(\mathrm{u}_y(t),t)^{1/3}\right).
	            \end{equation*}
	            The final threshold condition imposed at the end of Step~5 then turns this into the normalised complementary-contour estimate
	            \begin{equation}\label{eq:I2-small}
	                \left|I_{2,x,y}(t)\right|\leq\mathrm{C}_{x}\mathrm{e}^{-\mathrm{c}_{\varphi}\Sigma(t)^{2/3}}\mathbb{P}_{0}\left[\mathscr{X}(t)=y(t)\right],
	            \end{equation}
	            with $\mathrm{C}_{x}$ exactly as defined in \eqref{eq:Cx-explicit}. We now pass from local integrals to the coefficients $\eta_{n,y}(t)$. For even $n$, set $\mathrm{Q}_{n}\overset{\textnormal{def}}{=}(n-1)!!$, with the convention $\mathrm{Q}_{0}=(-1)!!=1$. The leading term in \eqref{eq:local-even-In} can then be written as
	            \begin{equation*}
	                (-1)^{\frac{n}{2}}\mathrm{Q}_{n}\mathfrak{F}_{y}^{(2)}(\mathrm{u}_{y}(t),t)^{-\frac{n}{2}}\mathcal{L}_{0,y}(t).
	            \end{equation*}
	            Dividing the local error estimate \eqref{eq:local-even-In} by $\mathcal{L}_{0,y}(t)$ gives the explicit relative error
	            \begin{equation*}
	                \left|\mathcal{I}_{n,y}(t)-(-1)^{\frac{n}{2}}\mathrm{Q}_{n}\mathfrak{F}_{y}^{(2)}(\mathrm{u}_{y}(t),t)^{-\frac{n}{2}}\mathcal{L}_{0,y}(t)\right|\leq\sqrt{2\pi}\,\mathrm{D}_{\mathrm{loc},n}\mathfrak{F}_{y}^{(2)}(\mathrm{u}_{y}(t),t)^{-\frac{n+1}{2}}\mathcal{L}_{0,y}(t).
	            \end{equation*}
            Since $\eta_{n,y}(t)=\mathcal{I}_{n,y}(t)/\mathbb{P}_{0}[\mathscr{X}(t)=y(t)]$, we separate the numerator error from the denominator normalisation. With
                \begin{equation*}
                    \mathrm{M}_{n,y}(t)\overset{\textnormal{def}}{=}(-1)^{\frac{n}{2}}\mathrm{Q}_{n}\mathfrak{F}_{y}^{(2)}(\mathrm{u}_{y}(t),t)^{-\frac{n}{2}},
                \end{equation*}
                the triangle inequality gives
                \begin{equation*}
                    \left|\eta_{n,y}(t)-\mathrm{M}_{n,y}(t)\right|
                    \leq
                    \frac{|\mathcal{I}_{n,y}(t)-\mathrm{M}_{n,y}(t)\mathcal{L}_{0,y}(t)|}{\mathbb{P}_{0}[\mathscr{X}(t)=y(t)]}
                    +|\mathrm{M}_{n,y}(t)|\left|\frac{\mathcal{L}_{0,y}(t)}{\mathbb{P}_{0}[\mathscr{X}(t)=y(t)]}-1\right|.
                \end{equation*}
                The first term is bounded by the local estimate above and \eqref{eq:P0-lower}. The second term is bounded by writing $\mathbb{P}_{0}[\mathscr{X}(t)=y(t)]=\mathcal{L}_{0,y}(t)+(\mathcal{I}_{0,y}(t)-\mathcal{L}_{0,y}(t))+I_{2,0,y}(t)$ and using the two denominator errors displayed before \eqref{eq:P0-lower}. Therefore
	            \begin{equation}\label{eq:eta-even-F2}
	                \left|\eta_{n,y}(t)-(-1)^{\frac{n}{2}}\mathrm{Q}_{n}\mathfrak{F}_{y}^{(2)}(\mathrm{u}_{y}(t),t)^{-\frac{n}{2}}\right|\leq \mathrm{D}_{\eta,n}^{(0)}\mathfrak{F}_{y}^{(2)}(\mathrm{u}_{y}(t),t)^{-\frac{n+1}{2}}+\mathrm{D}_{\eta,n}^{(1)}\mathrm{e}^{-\mathrm{c}_{\varphi}\Sigma(t)^{2/3}}\mathfrak{F}_{y}^{(2)}(\mathrm{u}_{y}(t),t)^{-\frac{n}{2}},
	            \end{equation}
	            where the explicit constants are
	            \begin{equation*}
	                \mathrm{D}_{\eta,n}^{(0)}\overset{\textnormal{def}}{=}2\sqrt{2\pi}\left(\mathrm{D}_{\mathrm{loc},n}+\mathrm{Q}_{n}\mathrm{D}_{\mathrm{loc},0}\right),\qquad \mathrm{D}_{\eta,n}^{(1)}\overset{\textnormal{def}}{=}2\sqrt{2\pi}\,\mathfrak{B}_{1}\mathrm{Q}_{n}\mathrm{C}_{0,\Gamma}.
	            \end{equation*}
	            Indeed, the first term in $\mathrm{D}_{\eta,n}^{(0)}$ comes from the numerator error in \eqref{eq:local-even-In}; the second comes from replacing $\mathbb{P}_{0}[\mathscr{X}(t)=y(t)]$ by $\mathcal{L}_{0,y}(t)$ in the denominator; and $\mathrm{D}_{\eta,n}^{(1)}$ is the contribution of the complementary contour $I_{2,0,y}(t)$, after using the exponential absorption from Step~5.
                \medskip

	            \noindent It remains to replace the exact curvature $\mathfrak{F}_{y}^{(2)}(\mathrm{u}_y(t),t)$ by the deterministic scale $\Sigma(t)^2$. By \textbf{(iv)} from Proposition~\ref{specificcasesub},
	            \begin{equation*}
	                \left|\mathfrak{F}_{y}^{(2)}(\mathrm{u}_{y}(t),t)-\Sigma(t)^2\right|\leq \mathrm{C}_{\mathscr{C}}\left(\gamma(t)^{1/2}+\left|\mathrm{u}_y(t)\right|\right)\Sigma(t)^2.
	            \end{equation*}
	            The stationary-point smallness estimate \eqref{axpro-u-small} gives $|\mathrm{u}_y(t)|\leq2\mathrm{B}_{\mathrm{A}}\gamma(t)^{1/2}$ after increasing $\mathrm{T}_{\mathrm{dec}}$ to dominate the lower time in Proposition~\ref{axpro}. Hence
	            \begin{equation*}
	                \left|\frac{\mathfrak{F}_{y}^{(2)}(\mathrm{u}_y(t),t)}{\Sigma(t)^2}-1\right|\leq\mathrm{C}_{\mathscr{C}}\left(1+2\mathrm{B}_{\mathrm{A}}\right)\gamma(t)^{1/2}.
	            \end{equation*}
	            Increasing $\mathrm{T}_{\mathrm{dec}}$ once more, we may assume that the right-hand side is at most $1/2$. On this range the ratio $\mathfrak{F}_{y}^{(2)}(\mathrm{u}_y(t),t)/\Sigma(t)^2$ lies in $[1/2,3/2]$, and the mean-value theorem applied to $r\mapsto r^{-n/2}$ gives
	            \begin{equation*}
	                \left|\mathfrak{F}_{y}^{(2)}(\mathrm{u}_y(t),t)^{-\frac{n}{2}}-\Sigma(t)^{-n}\right|\leq n2^{\frac{n}{2}}\mathrm{C}_{\mathscr{C}}\left(1+2\mathrm{B}_{\mathrm{A}}\right)\gamma(t)^{1/2}\Sigma(t)^{-n}.
	            \end{equation*}
                Multiplying this curvature-replacement bound by $\mathrm{Q}_{n}=(n-1)!!$ gives the first contribution to the even constant in \eqref{eq:Cn-explicit}.
	            Also, \eqref{eq:SD-second-lower} gives $\mathfrak{F}_{y}^{(2)}(\mathrm{u}_y(t),t)^{-(n+1)/2}\leq2^{(n+1)/2}\Sigma(t)^{-n-1}$ and $\mathfrak{F}_{y}^{(2)}(\mathrm{u}_y(t),t)^{-n/2}\leq2^{n/2}\Sigma(t)^{-n}$. Since $\gamma(t)\geq\Sigma(t)^{-1/2}$ by Definition~\ref{gammadef}, we have $\gamma(t)^{1/2}\geq\Sigma(t)^{-1/4}$; after increasing the $n$-dependent lower time so that $\Sigma(t)\geq1$, this implies $\Sigma(t)^{-1}\leq\gamma(t)^{1/2}$. Increasing this lower time once more so that $\mathrm{e}^{-\mathrm{c}_{\varphi}\Sigma(t)^{2/3}}\leq\gamma(t)^{1/2}$, the even estimate \eqref{eq:eta-even-F2} becomes the even part of \eqref{leadingordereta}. We denote by $\mathrm{T}_{\mathrm{dec},n}$ the maximum of $\mathrm{T}_{\mathrm{dec}}$, $\mathrm{T}_{\mathrm{dec},n}^{\mathrm{tail}}$ when $n$ is even, and the finitely many additional $n$-dependent lower times used in this paragraph. If $n$ is odd, the leading Gaussian contribution is absent because the symmetric integral of $u^n\exp(-u^2/2)$ is zero. The local estimate \eqref{eq:local-odd-In}, divided by the denominator lower bound in \eqref{eq:P0-lower}, gives directly
	            \begin{equation*}
	                \left|\eta_{n,y}(t)\right|\leq\mathrm{D}_{\eta,n}^{\mathrm{odd}}\mathfrak{F}_{y}^{(2)}(\mathrm{u}_{y}(t),t)^{-\frac{n+1}{2}},\qquad \mathrm{D}_{\eta,n}^{\mathrm{odd}}\overset{\textnormal{def}}{=}2\sqrt{2\pi}\,\mathrm{D}_{\mathrm{rem},n}.
	            \end{equation*}
	            Applying \eqref{eq:SD-second-lower} yields $\left|\eta_{n,y}(t)\right|\leq2^{(n+1)/2}\mathrm{D}_{\eta,n}^{\mathrm{odd}}\Sigma(t)^{-n-1}$, which is the odd part of \eqref{leadingordereta}; for odd $n$ we take $\mathrm{T}_{\mathrm{dec},n}$ to be the maximum of $\mathrm{T}_{\mathrm{dec}}$ and the lower times used in the denominator and local-remainder estimates above. Collecting the constants used in the even and odd cases, the constant in the statement may be taken to be
	            \begin{equation}\label{eq:Cn-explicit}
	                \mathrm{C}_{n}\overset{\textnormal{def}}{=}
	                \begin{cases}
	                    n2^{\frac{n}{2}}(n-1)!!\,\mathrm{C}_{\mathscr{C}}\left(1+2\mathrm{B}_{\mathrm{A}}\right)+2^{\frac{n+1}{2}}\mathrm{D}_{\eta,n}^{(0)}+2^{\frac{n}{2}}\mathrm{D}_{\eta,n}^{(1)},&\text{if $n$ is even,}\\
	                    2^{\frac{n+1}{2}}\mathrm{D}_{\eta,n}^{\mathrm{odd}},&\text{if $n$ is odd.}
	                \end{cases}
	            \end{equation}
        
            \textbf{Step 7: assembling \eqref{ointc}.}
	            Finally, the exact decomposition above gives
	            \begin{equation*}
	                \mathbb{P}_{x}\left[\mathscr{X}(t)=y(t)\right]=\left(\sum_{n=0}^{x}\frac{\rho_x^{(n)}(\mathrm{u}_{y}(t))}{n!}\eta_{n,y}(t)\right)\mathbb{P}_0\left[\mathscr{X}(t)=y(t)\right]+I_{2,x,y}(t).
	            \end{equation*}
		            Dividing by $\mathbb{P}_0[\mathscr{X}(t)=y(t)]$ and applying \eqref{eq:I2-small} gives \eqref{ointc} with the explicit constant $\mathrm{C}_{x}$ from \eqref{eq:Cx-explicit}.
      
        \textbf{On the thresholds.}
        The proof enlarges a starting lower time finitely many times, always using bounds uniform on $\{y(t)\in\mathscr{C}_{t}\}$.
        The $x$-dependent expansion \eqref{ointc} uses the common threshold $\mathrm{T}_{\mathrm{dec}}$, whose cumulative requirements are:
        \textup{(i)} applicability of Corollary~\ref{steepestdescentpath} and the endpoint-window inputs from Propositions~\ref{specificcasesub} and~\ref{axpro}, including $\mathrm{u}_y(t)\in\mathcal{B}_{\mathrm{R}_2}(0)$ and $\theta_y(t)<\phi_{\star}$;
        \textup{(ii)} containment of the local contour segment in $\mathcal{B}_{\mathrm{R}_2}(0)$ and the polydisk/Cauchy regime for $H_y(\cdot,t)$ from Step~4;
	        \textup{(iii)} the $n=0$ local bound needed for the denominator lower estimate;
        \textup{(iv)} the complementary-contour regime for Lemma~\ref{lem:Gamma2-Re-max} at cutoff $\vartheta=\theta_y(t)$, the endpoint reduction \eqref{eq:Gamma2-Re-max}, and the successive comparisons leading to \eqref{eq:I2-small};
        \textup{(v)} the domination \eqref{eq:I20-dominated-by-local} used to normalise the complementary-contour error.
        Any deterministic time satisfying \textup{(i)}--\textup{(v)} may be taken as $\mathrm{T}_{\mathrm{dec}}$. For the coefficient estimate \eqref{leadingordereta}, a fixed $n$ may require the additional finite-time conditions from the Gaussian truncation and curvature-replacement steps; their maximum with $\mathrm{T}_{\mathrm{dec}}$ is the threshold $\mathrm{T}_{\mathrm{dec},n}$ appearing in the statement.
        \end{proof}

     To see \eqref{ointc} concretely, it is useful to examine the homogeneous unit-rate process, where the transition probabilities are explicit. In this case the transition ratio can be expanded directly in the same Taylor basis as the decoupling sum. Example \ref{ex:homogeneous-unit-rates} shows how the saddle coordinate, the Taylor coefficients of $\rho_x$, and the even--odd hierarchy in \eqref{leadingordereta} fit together in the simplest possible model.

    \begin{example}[Homogeneous unit rates: the decoupling expansion in Poisson form]\label{ex:homogeneous-unit-rates}
        Work in continuous time with $\lambda_n\equiv 1$, fix $x\in\mathbb{Z}_+$, and restrict to integer times $t\geq1$. Then $\mathscr{Z}(t)=t$ and $\Sigma(t)^2=t+1$. Let $y$ be a path with $y(t)\in\mathscr{C}_t$. After increasing the lower time in a way depending only on the fixed value of $x$, the endpoint window gives $y(t)\geq x$. The stationary equation gives
        \begin{equation*}
            \mathrm{u}_y(t)=1-\frac{y(t)+1}{t},\qquad
            \mathscr{U}_y(t)=\frac{\sqrt{t+1}}{t}(y(t)+1-t),\qquad
            \rho_x(w)=(1-w)^x.
        \end{equation*}
        Since $\mathscr{X}(t)-x$ has distribution $\mathrm{Poisson}(t)$ under $\mathbb{P}_x$, the Poisson masses cancel and
        \begin{equation}\label{eq:homogeneous-ratio}
            \frac{\mathbb{P}_x[\mathscr{X}(t)=y(t)]}{\mathbb{P}_0[\mathscr{X}(t)=y(t)]}
            =\frac{y(t)!}{(y(t)-x)!t^x}
            =\prod_{j=0}^{x-1}\frac{y(t)-j}{t}.
        \end{equation}
        We now write this exact ratio in the Taylor basis appearing in \eqref{ointc}. Since
        \begin{equation}\label{eq:homogeneous-rho-derivatives}
            \frac{\rho_x^{(n)}(\mathrm{u}_{y}(t))}{n!}
            =(-1)^n\binom{x}{n}\left(\frac{y(t)+1}{t}\right)^{x-n},
        \end{equation}
        coefficient extraction gives
        \begin{equation}\label{eq:homogeneous-coefficients}
            \frac{\mathbb{P}_x[\mathscr{X}(t)=y(t)]}{\mathbb{P}_0[\mathscr{X}(t)=y(t)]}
            =
            \sum_{n=0}^{x}(-1)^n\binom{x}{n}\left(\frac{y(t)+1}{t}\right)^{x-n}\!\!\!\!
            \tilde{\eta}_{n,y}(t),
        \end{equation}
        where
        \begin{equation}\label{eq:homogeneous-first-etas}
            \tilde{\eta}_{n,y}(t)
            \overset{\textnormal{def}}{=}
            n![\zeta^n]\left\{
                \exp\left(\frac{y(t)+1}{t}\zeta\right)
                \left(1-\frac{\zeta}{t}\right)^{y(t)}
            \right\}
            =
            \frac{n!}{t^n}\sum_{\ell=0}^{n}\binom{y(t)}{\ell}
            \frac{(-1)^\ell(y(t)+1)^{n-\ell}}{(n-\ell)!}.
        \end{equation}
        To verify \eqref{eq:homogeneous-coefficients}, multiply its right-hand side by $s^x/x!$ and sum over $x$. The Taylor factors in \eqref{eq:homogeneous-rho-derivatives} give
        \begin{equation*}
            \exp\left(\frac{y(t)+1}{t}s\right)
            \sum_{n\geq0}\tilde{\eta}_{n,y}(t)\frac{(-s)^n}{n!}
            =
            \exp\left(\frac{y(t)+1}{t}s\right)
            \exp\left(-\frac{y(t)+1}{t}s\right)
            \left(1+\frac{s}{t}\right)^{y(t)}
            =
            \left(1+\frac{s}{t}\right)^{y(t)},
        \end{equation*}
        which is the exponential generating function of the falling-factorial ratio in \eqref{eq:homogeneous-ratio}. The first coefficients already display the mechanism behind \eqref{leadingordereta}:
        \begin{equation}\label{eq:homogeneous-low-coefficients}
            \tilde{\eta}_{0,y}(t)=1,\qquad
            \tilde{\eta}_{1,y}(t)=\frac{1}{t},\qquad
            \tilde{\eta}_{2,y}(t)=\frac{1-y(t)}{t^2}
        \end{equation}
        and
        \begin{equation*}
            \tilde{\eta}_{3,y}(t)=\frac{1-5y(t)}{t^3},\qquad
            \tilde{\eta}_{4,y}(t)=\frac{3y(t)^2-20y(t)+1}{t^4}.
        \end{equation*}
        In the central case $y(t)=t$, these identities compare directly with the Gaussian powers appearing in \eqref{leadingordereta}:
        \begin{equation*}
            \left|\tilde{\eta}_{2,y}(t)+\Sigma(t)^{-2}\right|=\frac{1}{t^2(t+1)},\qquad \left|\tilde{\eta}_{4,y}(t)-3\Sigma(t)^{-4}\right|\leq\frac{69}{t^3},\qquad t\geq1.
        \end{equation*}
        The odd coefficients displayed above are smaller than the corresponding even Gaussian scale in this central example. For comparison with the local coefficients in \eqref{defofmathcalI}, the homogeneous local contour can also be written explicitly. On the local steepest-descent circle, set
        \begin{equation*}
            w(\phi)=1-\frac{y(t)+1}{t}\mathrm{e}^{\mathrm{i}\phi},
            \qquad
            \phi_t=\arcsin\left(\frac{t}{y(t)+1}\mathfrak{F}_y^{(2)}(\mathrm{u}_y(t),t)^{-1/3}\right),
        \end{equation*}
        after the lower time has been chosen so that the argument of the arcsine is at most one. The change of variables in \eqref{defofmathcalI}, restricted to the symmetric arc $|\phi|\leq\phi_t$, gives
        \begin{equation*}
            \eta_{n,y}(t)
            =
            \frac{y(t)!\,(y(t)+1)^{n-y(t)}}{2\pi t^n}
            \sum_{k=0}^{n}(-1)^k\binom{n}{k}
            J_{y(t)-k}\big(y(t)+1,\phi_t\big),
        \end{equation*}
        where
        \begin{equation*}
            J_r(A,\phi)\overset{\textnormal{def}}{=}
            \int_{-\phi}^{\phi}\exp\!\left(A\mathrm{e}^{\mathrm{i}\psi}\right)\mathrm{e}^{-\mathrm{i}r\psi}\,\mathrm{d}\psi
            =
            \sum_{q=0}^{\infty}\frac{A^q}{q!}\,\mathrm{S}_{q-r}(\phi),
            \qquad
            \mathrm{S}_{\ell}(\phi)=
            \begin{cases}
                2\sin(\ell\phi)/\ell,&\ell\neq0,\\
                2\phi,&\ell=0.
            \end{cases}
        \end{equation*}
        Thus the only distinction between the local coefficients $\eta_{n,y}(t)$ and the full-ratio coefficients $\tilde{\eta}_{n,y}(t)$ is the angular range of integration. Indeed, since $J_r(A,\pi)=2\pi A^r/r!$ for $r\geq0$, for each fixed $n$ there is a lower time after which $n\leq y(t)$, and replacing $\phi_t$ by $\pi$ in the preceding display gives exactly \eqref{eq:homogeneous-first-etas}. For instance, at $t=1000$ and $y(t)=1000$, the rounded values are
        \begin{equation*}
        \begin{array}{c|cc}
        n & \eta_{n,y}(t) & \tilde{\eta}_{n,y}(t)\\
        \hline
        0 & 0.99846 & 1\\
        1 & 1.00975\times10^{-3} & 1.00000\times10^{-3}\\
        2 & -9.80590\times10^{-4} & -9.99000\times10^{-4}\\
        3 & -4.90769\times10^{-6} & -4.99900\times10^{-6}\\
        4 & 2.75659\times10^{-6} & 2.98000\times10^{-6}
        \end{array}
        \end{equation*}
        The exact coefficients $\tilde{\eta}_{n,y}(t)$ in \eqref{eq:homogeneous-first-etas} come from the full Poisson transition ratio, whereas $\eta_{n,y}(t)$ is defined through the local contour integral in Step~3 of Proposition~\ref{Thm2}. 
    \end{example}

    \noindent The same local estimates also give a kind of local limit theorem which, although not used directly in the proof, could be of independent interest.
    
    \begin{cor}[Stationary-point local limit]\label{stationarypointlocallimit}
    Assume \textbf{(A1)} and \textbf{(A2)}. For every $x\in\mathbb{Z}_+$, there exist $\mathrm{T}_{\mathrm{loc},x}\in\mathscr{T}$ and $\mathrm{K}_{x}>0$ such that, for every $t>\mathrm{T}_{\mathrm{loc},x}$ and every path $y\in\mathscr{P}$ satisfying $y(t)\in\mathscr{C}_{t}$, the minimizer in Definition~\ref{definitionstationary} is unique and equals $\mathrm{u}_y(t)$, with
        \begin{equation}\label{stationarypointlocallimitbound}
            \left|\mathbb{P}_{x}\left[\mathscr{X}(t)=y(t)\right]-\frac{\rho_{x}(\mathrm{u}_{y}(t))}{\lambda_{y(t)}\sqrt{2\pi}\,\Sigma(t)}\exp\left(-\frac{1}{2}\mathscr{U}_{y}(t)^2\right)\right|\leq\frac{\mathrm{K}_{x}\gamma(t)^{1/2}}{\Sigma(t)}.
        \end{equation}
        Here $\mathrm{K}_{x}$ depends only on $x$ and the fixed structural constants from Proposition~\ref{axpro} and Proposition~\ref{specificcasesub}, and in particular is independent of the choice of $y$ once $y(t)\in\mathscr{C}_{t}$.
    \end{cor}
    \begin{proof}
        Fix $x\in\mathbb{Z}_+$. We construct $\mathrm{T}_{\mathrm{loc},x}$ by starting from the threshold $\mathrm{T}_{\mathrm{dec}}$ in Proposition~\ref{Thm2} and increasing it finitely many times; these increases may depend on this fixed $x$, but remain independent of $y$ because all endpoint estimates invoked below are uniform on $\mathscr{C}_{t}$. Let $t>\mathrm{T}_{\mathrm{loc},x}$, and let $y\in\mathscr{P}$ satisfy $y(t)\in\mathscr{C}_{t}$. The proof of Proposition~\ref{Thm2} gives the exact decomposition
        \begin{equation*}
            \mathbb{P}_{x}\left[\mathscr{X}(t)=y(t)\right]=\sum_{n=0}^{x}\frac{\rho_x^{(n)}(\mathrm{u}_{y}(t))}{n!}\mathcal{I}_{n,y}(t)+I_{2,x,y}(t).
        \end{equation*}
        We first isolate the exact saddle approximation. Since \eqref{axpro-u-small} and $\gamma(t)\to0$ put $\mathrm{u}_{y}(t)$ in a fixed compact subset of $\mathcal{B}_{\mathrm{R}_{2}}(0)$ after increasing the lower time, the finitely many derivatives $\rho_x^{(n)}(\mathrm{u}_{y}(t))$, $0\leq n\leq x$, are bounded by a constant depending only on $x$. The $n=0$ estimate \eqref{eq:I0-local-bound}, together with the lower rate bound $\lambda_{y(t)}\geq\mathfrak{L}$ from \textbf{(A1)} and the curvature lower bound \eqref{eq:SD-second-lower}, gives
        \begin{equation*}
            \left|\rho_x(\mathrm{u}_{y}(t))\mathcal{I}_{0,y}(t)-\frac{\rho_{x}(\mathrm{u}_{y}(t))}{\lambda_{y(t)}\sqrt{2\pi\,\mathfrak{F}_{y}^{(2)}(\mathrm{u}_{y}(t),t)}}\exp\left(\mathfrak{F}_{y}(\mathrm{u}_{y}(t),t)\right)\right|\leq\frac{\mathrm{K}_{x,0}}{\Sigma(t)^2}.
        \end{equation*}
        Here we have also used $\exp(\mathfrak{F}_{y}(\mathrm{u}_{y}(t),t))\leq1$. To justify this last inequality, note that $\mathfrak{F}_{y}(0,t)=0$, $\mathfrak{F}_{y}^{(1)}(\mathrm{u}_{y}(t),t)=0$ by \eqref{firstconditiontosat}, and \eqref{eq:SD-second-lower} together with \textbf{(iv)} of Proposition~\ref{specificcasesub} makes $\mathfrak{F}_{y}^{(2)}(w,t)$ non-negative on the real segment between $0$ and $\mathrm{u}_{y}(t)$ after increasing $\mathrm{T}_{\mathrm{loc},x}$. This segment lies inside $\mathcal{B}_{\mathrm{R}_{2}}(0)$, so the logarithmic branches in Definition~\ref{form} are real-valued along it. Thus the real function $w\mapsto\mathfrak{F}_{y}(w,t)$ is convex on that segment and has a stationary point at $\mathrm{u}_{y}(t)$, so $\mathfrak{F}_{y}(\mathrm{u}_{y}(t),t)\leq\mathfrak{F}_{y}(0,t)=0$. It remains to control the terms with $n\geq1$ and the complementary contour. For each fixed $1\leq n\leq x$, the local estimates \eqref{eq:local-even-In} and \eqref{eq:local-odd-In}, together with the coefficient bounds preceding them, imply
        \begin{equation*}
            \left|\mathcal{I}_{n,y}(t)\right|\leq\frac{\mathrm{K}_{x,n}}{\Sigma(t)^2},
        \end{equation*}
        after increasing the lower time, uniformly over $y(t)\in\mathscr{C}_{t}$. Indeed, every such term contains at least one additional inverse power of $\mathfrak{F}_{y}^{(2)}(\mathrm{u}_{y}(t),t)^{1/2}$ compared with $\mathcal{I}_{0,y}(t)$, and \eqref{eq:SD-second-lower} converts this into at least one additional inverse power of $\Sigma(t)$. Also, \eqref{eq:I2-small}, the trivial bound $\mathbb{P}_{0}[\mathscr{X}(t)=y(t)]\leq1$, and the fact that $\mathrm{e}^{-\mathrm{c}_{\varphi}\Sigma(t)^{2/3}}\leq\Sigma(t)^{-2}$ after increasing $\mathrm{T}_{\mathrm{loc},x}$ show that $|I_{2,x,y}(t)|\leq\mathrm{K}_{x,\mathrm{tail}}\Sigma(t)^{-2}$. Summing the finitely many bounds for $0\leq n\leq x$ yields the exact-saddle estimate
        \begin{equation}\label{stationarypointlocallimitexact}
            \left|\mathbb{P}_{x}\left[\mathscr{X}(t)=y(t)\right]-\frac{\rho_{x}(\mathrm{u}_{y}(t))}{\lambda_{y(t)}\sqrt{2\pi\,\mathfrak{F}_{y}^{(2)}(\mathrm{u}_{y}(t),t)}}\exp\left(\mathfrak{F}_{y}(\mathrm{u}_{y}(t),t)\right)\right|\leq\frac{\mathrm{K}_{x,\mathrm{ex}}}{\Sigma(t)^2}.
        \end{equation}
        \noindent We now replace the exact saddle height and curvature by the stationary-point coordinate. On the present time range, \eqref{firstconditiontosat} gives $\mathscr{U}_{y}(t)=-\Sigma(t)\mathrm{u}_{y}(t)$. Since $\mathfrak{F}_{y}(0,t)=0$ and $\mathfrak{F}_{y}^{(1)}(\mathrm{u}_{y}(t),t)=0$, the fundamental theorem of calculus gives the exact identity
        \begin{equation*}
            \mathfrak{F}_{y}(\mathrm{u}_{y}(t),t)+\frac{1}{2}\Sigma(t)^2\mathrm{u}_{y}(t)^2=-\mathrm{u}_{y}(t)^2\int_{0}^{1}r\left(\mathfrak{F}_{y}^{(2)}(r\mathrm{u}_{y}(t),t)-\Sigma(t)^2\right)\mathrm{d}r.
        \end{equation*}
        The real segment from $0$ to $\mathrm{u}_{y}(t)$ lies inside $\mathcal{B}_{\mathrm{R}_{2}}(0)$ after the same lower-time increase used above. Applying \textbf{(iv)} of Proposition~\ref{specificcasesub} on this segment and using \eqref{axpro-u-small} gives a constant $\mathrm{K}_{\mathrm{ph}}>0$, independent of $t$ and $y$, such that
        \begin{equation*}
            \left|\mathfrak{F}_{y}(\mathrm{u}_{y}(t),t)+\frac{1}{2}\mathscr{U}_{y}(t)^2\right|\leq\mathrm{K}_{\mathrm{ph}}\gamma(t)^{1/2}\mathscr{U}_{y}(t)^2.
        \end{equation*}
        Similarly, \textbf{(iv)} and \eqref{axpro-u-small} imply
        \begin{equation*}
            \left|\frac{\mathfrak{F}_{y}^{(2)}(\mathrm{u}_{y}(t),t)}{\Sigma(t)^2}-1\right|\leq\mathrm{K}_{\mathrm{curv}}\gamma(t)^{1/2},
        \end{equation*}
        for a constant $\mathrm{K}_{\mathrm{curv}}>0$ independent of $t$ and $y$. Increasing $\mathrm{T}_{\mathrm{loc},x}$ so that the right-hand side is at most $1/2$, the mean-value theorem applied to $r\mapsto r^{-1/2}$ on $[1/2,3/2]$ gives
        \begin{equation*}
            \left|\frac{\Sigma(t)}{\sqrt{\mathfrak{F}_{y}^{(2)}(\mathrm{u}_{y}(t),t)}}-1\right|\leq\mathrm{K}_{\mathrm{curv}}^{\prime}\gamma(t)^{1/2}.
        \end{equation*}
        The phase-height estimate above still contains the factor $\mathscr{U}_{y}(t)^2$, which can be absorbed into a uniform upper bound after the Gaussian decay is retained. Writing $\mathfrak{F}_{y}(\mathrm{u}_{y}(t),t)=-\mathscr{U}_{y}(t)^2/2+\mathcal{E}_{y}(t)$, we have $|\mathcal{E}_{y}(t)|\leq\mathrm{K}_{\mathrm{ph}}\gamma(t)^{1/2}\mathscr{U}_{y}(t)^2$. After increasing $\mathrm{T}_{\mathrm{loc},x}$ so that $\mathrm{K}_{\mathrm{ph}}\gamma(t)^{1/2}\leq1/4$, the elementary estimate $|\mathrm{e}^{z}-1|\leq |z|\mathrm{e}^{|z|}$ yields
        \begin{equation*}
            \left|\exp\left(\mathfrak{F}_{y}(\mathrm{u}_{y}(t),t)\right)-\exp\left(-\frac{1}{2}\mathscr{U}_{y}(t)^2\right)\right|
            \leq
            \mathrm{K}_{\mathrm{ph}}\gamma(t)^{1/2}\mathscr{U}_{y}(t)^2\exp\left(-\frac{1}{4}\mathscr{U}_{y}(t)^2\right).
        \end{equation*}
        Since $\sup_{r\in\mathbb{R}}r^2\mathrm{e}^{-r^2/4}<\infty$, the prefactor $\mathscr{U}_{y}(t)^2\exp(-\mathscr{U}_{y}(t)^2/4)$ is bounded uniformly in $y(t)\in\mathscr{C}_{t}$ and $t>\mathrm{T}_{\mathrm{loc},x}$, and therefore there exists $\mathrm{K}_{\mathrm{ph}}^{\prime}>0$, independent of $t$ and $y$, such that
        \begin{equation*}
            \left|\exp\left(\mathfrak{F}_{y}(\mathrm{u}_{y}(t),t)\right)-\exp\left(-\frac{1}{2}\mathscr{U}_{y}(t)^2\right)\right|\leq\mathrm{K}_{\mathrm{ph}}^{\prime}\gamma(t)^{1/2}.
        \end{equation*}
        The polynomial $\rho_x$ is bounded on the fixed compact set containing all $\mathrm{u}_{y}(t)$, and \textbf{(A1)} gives $\lambda_{y(t)}\geq\mathfrak{L}$. Combining the preceding display with the curvature comparison gives
        \begin{align*}
            &\left|\frac{\rho_{x}(\mathrm{u}_{y}(t))}{\lambda_{y(t)}\sqrt{2\pi\,\mathfrak{F}_{y}^{(2)}(\mathrm{u}_{y}(t),t)}}\exp\left(\mathfrak{F}_{y}(\mathrm{u}_{y}(t),t)\right)-\frac{\rho_{x}(\mathrm{u}_{y}(t))}{\lambda_{y(t)}\sqrt{2\pi}\,\Sigma(t)}\exp\left(-\frac{1}{2}\mathscr{U}_{y}(t)^2\right)\right|\\
            &\qquad\leq\frac{\mathrm{K}_{x,\mathrm{G}}\gamma(t)^{1/2}}{\Sigma(t)}.
        \end{align*}
        Finally, the triangle inequality combines \eqref{stationarypointlocallimitexact}, which bounds the difference from the exact saddle approximation by $\mathrm{K}_{x,\mathrm{ex}}\Sigma(t)^{-2}$, with the Gaussian-coordinate comparison above, which bounds the replacement of that saddle approximation by $\mathrm{K}_{x,\mathrm{G}}\gamma(t)^{1/2}/\Sigma(t)$. By Definition~\ref{gammadef}, $\gamma(t)\geq\Sigma(t)^{-1/2}$, so $\Sigma(t)^{-2}\leq\gamma(t)^{1/2}\Sigma(t)^{-1}$ whenever $\Sigma(t)\geq1$. Increasing $\mathrm{T}_{\mathrm{loc},x}$ once more, we may assume $\mathrm{K}_{x,\mathrm{ex}}\Sigma(t)^{-2}\leq\tfrac{1}{2}\mathrm{K}_{x,\mathrm{G}}\gamma(t)^{1/2}/\Sigma(t)$, because $\Sigma(t)\to\infty$ by \textbf{(iv)} of Proposition~\ref{specificcasesub}. Summing the two displays with the enlarged constant $\mathrm{K}_{x}\overset{\textnormal{def}}{=}\mathrm{K}_{x,\mathrm{G}}+\mathrm{K}_{x,\mathrm{ex}}$ proves \eqref{stationarypointlocallimitbound}.
    \end{proof}
    \section{The process associated with the stationary point}\label{stationarypointprocess}
    \subsection{Large Deviations}\label{MOM}
    This section supplies the probabilistic localization input needed to apply the deterministic endpoint-window expansion from Section~\ref{decouplingsection} to random terminal states. The phase estimates from Section~\ref{glo}, the stationary-point construction from Section~\ref{ASP}, the contour geometry from Section~\ref{CM}, and the one-particle decoupling theorem from Section~\ref{decouplingsection} are all conditional on endpoint control. The purpose here is to show that, under \textbf{(A1)} and \textbf{(A2)}, with probability one there is a finite random time after which the process remains in the finite-time endpoint window $\mathscr{C}_{t}$ from Definition \ref{def:endpoint-window}. Consequently, the deterministic expansion may be applied to random endpoints with the complementary event contributing only a vanishing tail.
    
     The argument has two levels. Proposition~\ref{eigenfunction} provides the exponential observables from which the one-time estimate is derived. The quantitative one-time bound is then obtained in Proposition~\ref{localizationtail} while Proposition~\ref{moments} then upgrades this one-time input to a path tail localization statement.
    \begin{proposition}\label{eigenfunction}
    For all $x\in\mathbb{Z}_+$ and $t\in\mathscr{T}$, the following identity holds for every $w\in\mathbb{C}$ in the continuous-time regime and for every $w\in\mathcal{B}_{\mathfrak{L}+\mathfrak{U}^{-1}}(\mathfrak{L})$ in the discrete-time regime:
    \begin{equation}\label{charac}
        \mathbb{E}_x\left[\rho_{\mathscr{X}(t)}(w)\right]
        =
        \begin{cases}
            \rho_x(w)e^{-tw},&\textnormal{if }\mathscr{T}=\mathbb{R}_+,\\
            \rho_x(w)\displaystyle\prod_{n=0}^{t-1}\dfrac{(1-\alpha_nw)^{1-\tau_n}}{(1+\beta_nw)^{\tau_n}},&\textnormal{if }\mathscr{T}=\mathbb{Z}_+.
        \end{cases}
    \end{equation}
    In the discrete-time formula the powers are ordinary integer powers, since $\tau_n\in\{0,1\}$, and the empty product at $t=0$ is equal to one.
    \begin{proof}
        We prove the identity directly from the transition mechanisms in Definition~\ref{contprocess} and Definition~\ref{mixprocess}. The only algebraic input needed throughout the proof is the one-step recursion following immediately from the definition of the characteristic polynomials in \eqref{characteristic},
        \begin{equation}\label{rhoeigenrecursion}
            \rho_{y+1}(w)=\rho_y(w)\left(1-\frac{w}{\lambda_y}\right),\qquad y\in\mathbb{Z}_+.
        \end{equation}
        We first consider the continuous-time regime. The generator $\mathsf{L}$ of the pure-birth chain in Definition~\ref{contprocess} acts on a test function $f$ by $\mathsf{L}f(y)=\lambda_y(f(y+1)-f(y))$. Substituting $f(y)=\rho_y(w)$ and using \eqref{rhoeigenrecursion} gives the exact eigenfunction relation
        \begin{equation*}
            (\mathsf{L}\rho_{\cdot}(w))(y)=\lambda_y\left(\rho_{y+1}(w)-\rho_y(w)\right)=-w\rho_y(w).
        \end{equation*}
        Since the function $y\mapsto\rho_y(w)$ is not bounded for arbitrary complex $w$, we justify the preceding generator computation by the standard stopping argument. Let $N(t)$ be the number of jumps made by the process up to time $t$. By the upper rate bound in \textbf{(A1)}, $N(t)$ is stochastically dominated by a Poisson random variable with mean $\mathfrak{B}_1t$. On the event $\{N(t)=k\}$ the process is at most at $x+k$, and the lower rate bound in \textbf{(A1)} together with \eqref{rhoeigenrecursion} gives
        \begin{equation*}
            |\rho_{x+k}(w)|\leq|\rho_x(w)|\left(1+\frac{|w|}{\mathfrak{L}}\right)^k.
        \end{equation*}
        The right-hand side has finite expectation under the dominating Poisson law for every fixed $t$ and $w$. Therefore the stopped Dynkin formula for the process killed when it leaves $\{x,\ldots,x+m\}$ may be passed to the limit $m\rightarrow\infty$ by dominated convergence. Thus the function $h(t)\overset{\textnormal{def}}{=}\mathbb{E}_x[\rho_{\mathscr{X}(t)}(w)]$ is finite and satisfies the integral form of the backward equation, equivalently
        \begin{equation*}
            h'(t)=-wh(t),\qquad h(0)=\rho_x(w).
        \end{equation*}
        Solving this scalar differential equation gives $h(t)=\rho_x(w)e^{-tw}$, proving \eqref{charac} in the continuous-time regime. We now turn to the discrete-time regime. In the Bernoulli case $\tau_t=0$, conditioning on $\mathscr{X}(t)=y$ and using \eqref{rhoeigenrecursion} gives
        \begin{equation*}
            \mathbb{E}\left[\rho_{\mathscr{X}(t+1)}(w)\mid\mathscr{X}(t)=y\right]=(1-\lambda_y\alpha_t)\rho_y(w)+\lambda_y\alpha_t\rho_{y+1}(w)=(1-\alpha_tw)\rho_y(w).
        \end{equation*}
        This proves the one-step identity in the Bernoulli case. In the geometric case $\tau_t=1$, again condition on $\mathscr{X}(t)=y$. For $k\geq0$, set
        \begin{equation*}
            A_k(w)\overset{\textnormal{def}}{=}\prod_{i=y}^{y+k-1}\frac{\beta_t(\lambda_i-w)}{1+\beta_t\lambda_i},
        \end{equation*}
        with the empty product $A_0(w)=1$. We first check that the tail of this product vanishes. Since $w\in\mathcal{B}_{\mathfrak{L}+\mathfrak{U}^{-1}}(\mathfrak{L})$, for every $i\in\mathbb{Z}_+$ the triangle inequality and the lower bound $\lambda_i\geq\mathfrak{L}$ give
        \begin{equation*}
            |\lambda_i-w|\leq|\lambda_i-\mathfrak{L}|+|w-\mathfrak{L}|<\lambda_i+\mathfrak{U}^{-1}.
        \end{equation*}
        By the upper bound $\beta_t<\mathfrak{U}$ from \textbf{(A1)}, this implies $|\lambda_i-w|<\lambda_i+\beta_t^{-1}$, or equivalently
        \begin{equation*}
            \left|\frac{\beta_t(\lambda_i-w)}{1+\beta_t\lambda_i}\right|<1.
        \end{equation*}
        The map $\lambda\mapsto|\beta_t(\lambda-w)|(1+\beta_t\lambda)^{-1}$ is continuous on the compact interval $[\mathfrak{L},\mathfrak{B}_1]$ and is strictly smaller than one at every point of that interval. Hence there exists $\mathrm{q}(w,t)<1$ such that
        \begin{equation*}
            \left|\frac{\beta_t(\lambda_i-w)}{1+\beta_t\lambda_i}\right|\leq\mathrm{q}(w,t),\qquad i\in\mathbb{Z}_+.
        \end{equation*}
        Consequently, $|A_k(w)|\leq \mathrm{q}(w,t)^k$ and therefore $\lim_{k\to\infty}A_k(w)=0$. The reason for introducing $A_k(w)$ is that it makes the geometric sum telescope. Indeed, for every $k\geq0$,
        \begin{equation*}
            A_k(w)-A_{k+1}(w)=A_k(w)\left(1-\frac{\beta_t(\lambda_{y+k}-w)}{1+\beta_t\lambda_{y+k}}\right)=A_k(w)\frac{1+\beta_tw}{1+\beta_t\lambda_{y+k}}.
        \end{equation*}
        Combining the geometric transition probabilities from Definition~\ref{mixprocess} with \eqref{rhoeigenrecursion}, we obtain, for every $m\geq0$,
        \begin{align*}
            \sum_{k=0}^{m}\mathbb{P}\left[\mathscr{X}(t+1)=y+k\mid\mathscr{X}(t)=y\right]\rho_{y+k}(w)
            &=\rho_y(w)\sum_{k=0}^{m}\frac{A_k(w)}{1+\beta_t\lambda_{y+k}}\\
            &=\frac{\rho_y(w)}{1+\beta_tw}\left(1-A_{m+1}(w)\right).
        \end{align*}
        Passing to the limit as $m$ tends to infinity is justified by the bound $|A_k(w)|\leq \mathrm{q}(w,t)^k$, and gives
        \begin{equation*}
            \mathbb{E}\left[\rho_{\mathscr{X}(t+1)}(w)\mid\mathscr{X}(t)=y\right]=\frac{\rho_y(w)}{1+\beta_tw}.
        \end{equation*}
        Thus, in both discrete one-step regimes,
        \begin{equation*}
            \mathbb{E}\left[\rho_{\mathscr{X}(t+1)}(w)\mid\mathscr{X}(t)\right]=\rho_{\mathscr{X}(t)}(w)\frac{(1-\alpha_tw)^{1-\tau_t}}{(1+\beta_tw)^{\tau_t}}.
        \end{equation*}
        Iterating this identity by the tower property, starting from $\mathscr{X}(0)=x$, gives the discrete-time part of \eqref{charac}. This completes the proof.
    \end{proof}
    \end{proposition}
	\noindent For vector endpoints we use the coordinatewise endpoint window defined as follows,
	\begin{equation*}
	    \mathscr{C}_{t}^{N}\overset{\textnormal{def}}{=}\left\{\mathbf{m}=(m_1,\ldots,m_N)\in\mathbb{Z}_+^N:m_i\in\mathscr{C}_{t}\textnormal{ for every }i\in\llbracket N\rrbracket\right\}.
	\end{equation*}
\begin{proposition}[One-time tail localization]\label{localizationtail}
Assume \textbf{(A1)} and \textbf{(A2)}. Then there exist constants $\kappa,c>0$ and $\mathrm{T}_{\mathrm{LD}}\in\mathscr{T}$ such that, for every $x\in\mathbb{Z}_+$, there exists $\mathrm{M}_x>0$ such that
\begin{equation}\label{localizationtailbound}
    \mathbb{P}_x\left[\left|\mathscr{X}(t)-\mathscr{Z}(t)\right|\geq \kappa\gamma(t)\Sigma(t)^2\right]\leq \mathrm{M}_x\mathrm{e}^{-c\Sigma(t)},\qquad t\geq\mathrm{T}_{\mathrm{LD}}.
\end{equation}
In particular,
\begin{equation}\label{endpointwindowtailbound}
    \mathbb{P}_x\left[\mathscr{X}(t)\notin\mathscr{C}_{t}\right]\leq \mathrm{M}_x\mathrm{e}^{-c\Sigma(t)},\qquad t\geq\mathrm{T}_{\mathrm{LD}}.
\end{equation}
\end{proposition}
\begin{proof}
We prove the estimate by applying Markov's inequality to the eigenfunction from Proposition~\ref{eigenfunction}. The proof has four steps. First we choose the two deterministic boundary paths at distance $\gamma(t)\Sigma(t)^2$ from $\mathscr{Z}(t)$ and choose a deterministic lower time $\mathrm{T}_{\mathrm{LD}}$ so that, for every $t\geq\mathrm{T}_{\mathrm{LD}}$, the corresponding endpoints lie in the wider window $\mathscr{C}_{t}$. Second we turn the upper and lower tails into bounds involving the phase $\mathfrak{F}_y$. Third we use the new first-derivative estimate in \textbf{(ii)} of Proposition~\ref{specificcasesub} to prove that the two boundary phases have stationary points on the correct sides of the origin. Finally we evaluate the phase at those stationary points and obtain an exponentially small bound. Throughout the proof, $\mathrm{T}_{\mathrm{LD}}$ is increased finitely many times; the final value is the threshold appearing in the statement.

\textbf{Step 1: The intrinsic boundary paths.}
Let $\mathrm{T}_1$ be the threshold from \textbf{(ii)} of Proposition~\ref{specificcasesub}, and let $\mathrm{C}_1,\mathrm{C}_2,\mathrm{T}_3$ be the constants and threshold from \textbf{(iv)} of the same proposition. The sign-sensitive coefficient in the boundary estimates below is $\mathfrak{B}_1^{-1}$, because the lower bound at the upper endpoint and the upper bound at the lower endpoint both use the $\mathfrak{B}_1^{-1}$ branch of the first display in \textbf{(ii)}. Define
\begin{equation*}
    \mathrm{A}_{\star}\overset{\textnormal{def}}{=}2\mathrm{C}_2+1,
    \qquad
    \kappa\overset{\textnormal{def}}{=}4\mathfrak{B}_1\mathrm{A}_{\star}.
\end{equation*}
The value $\kappa=4\mathfrak{B}_1\mathrm{A}_{\star}$ is chosen so that the deterministic displacement of the endpoint dominates the centring error allowed in \textbf{(ii)} of Proposition~\ref{specificcasesub}; this is used explicitly in \eqref{boundary-displacement-bounds}--\eqref{boundary-F-bound}.
\medskip

\noindent By Definition~\ref{gammadef}, $\gamma(t)\geq\Sigma(t)^{-1/2}$ and $\gamma(t)\downarrow0$. Also, \textbf{(iv)} of Proposition~\ref{specificcasesub} gives
\begin{equation}\label{localization-scale-comparison}
    \mathrm{C}_1\Sigma(t)^2\leq\mathscr{Z}(t)\leq\mathrm{C}_2\Sigma(t)^2,
    \qquad t>\mathrm{T}_3.
\end{equation}
Consequently, since $\gamma(t)$ can be made arbitrarily small and $\Sigma(t)$ arbitrarily large after increasing the deterministic lower time, we may assume throughout the proof that
\begin{equation}\label{window-preconditions}
    \kappa\gamma(t)\Sigma(t)^2\leq\frac{1}{2}\mathscr{Z}(t),
    \qquad
    \gamma(t)\Sigma(t)^2\geq1,
    \qquad
    \mathfrak{B}_1^{-1}+\frac{3\Lambda}{2}\leq\mathrm{A}_{\star}\gamma(t)\Sigma(t)^2.
\end{equation}
The first inequality ensures that the lower boundary remains in $\mathbb{Z}_+$, while the last two inequalities allow the integer-rounding errors and the additive constant $3\Lambda/2$ in \textbf{(ii)} of Proposition~\ref{specificcasesub} to be absorbed into the intrinsic window.
\medskip

\noindent For $t\in\mathscr{T}$, set
\begin{equation*}
    \mathrm{Z}_{+}(t)\overset{\textnormal{def}}{=}\mathscr{Z}(t)+\kappa\gamma(t)\Sigma(t)^2,
    \qquad
    \mathrm{Z}_{-}(t)\overset{\textnormal{def}}{=}\mathscr{Z}(t)-\kappa\gamma(t)\Sigma(t)^2,
\end{equation*}
and define the integer boundary paths
\begin{equation*}
    \mathscr{Z}_{+}(t)\overset{\textnormal{def}}{=}\left\lceil\mathrm{Z}_{+}(t)\right\rceil,
    \qquad
    \mathscr{Z}_{-}(t)\overset{\textnormal{def}}{=}\left\lfloor\mathrm{Z}_{-}(t)\right\rfloor.
\end{equation*}
For the estimates below we use these displayed values only after the final lower time; before that time we choose arbitrary non-negative integer values, so that $\mathscr{Z}_{\pm}\in\mathscr{P}$.
After increasing $\mathrm{T}_{\mathrm{LD}}$ so that \eqref{window-preconditions} holds for every $t\geq\mathrm{T}_{\mathrm{LD}}$, the first inequality in \eqref{window-preconditions} gives $\mathscr{Z}_{-}(t)\geq0$. Moreover,
\begin{equation}\label{boundary-displacement-bounds}
    \mathscr{Z}_{+}(t)-\mathscr{Z}(t)\geq\kappa\gamma(t)\Sigma(t)^2,
    \qquad
    \mathscr{Z}(t)-\mathscr{Z}_{-}(t)\geq\kappa\gamma(t)\Sigma(t)^2,
\end{equation}
and, using $\gamma(t)\Sigma(t)^2\geq1$,
\begin{equation}\label{boundary-F-bound}
    |\mathscr{Z}_{\pm}(t)-\mathscr{Z}(t)|\leq(\kappa+1)\gamma(t)\Sigma(t)^2.
\end{equation}
We include in the final choice of $\mathrm{T}_{\mathrm{LD}}$ the deterministic condition $(\kappa+1)\gamma(t)^{1/2}\leq1$ for every $t\geq\mathrm{T}_{\mathrm{LD}}$, which is possible because $\gamma(t)\to0$. Then \eqref{boundary-F-bound} gives
\begin{equation*}
    |\mathscr{Z}_{\pm}(t)-\mathscr{Z}(t)|
    \leq
    \gamma(t)^{1/2}\Sigma(t)^2,
\end{equation*}
so both endpoints $\mathscr{Z}_{+}(t)$ and $\mathscr{Z}_{-}(t)$ belong to $\mathscr{C}_{t}$. Thus the two boundary endpoints satisfy exactly the finite-time hypothesis of Proposition~\ref{axpro}. No separate qualitative path argument is needed at this point: the stationary-point construction in that proposition is now uniform over all endpoints in $\mathscr{C}_{t}$.

\textbf{Step 2: Markov bounds from the eigenfunction.}
Since $\mathscr{X}(t)$ and $\mathscr{Z}(t)$ are integer-valued, the event in the statement is contained in the union of the two one-sided events
\begin{equation}\label{partitionProb}
    \left\{\left|\mathscr{X}(t)-\mathscr{Z}(t)\right|\geq \kappa\gamma(t)\Sigma(t)^2\right\}
    \subseteq
    \left\{\mathscr{X}(t)\geq\mathscr{Z}_{+}(t)\right\}
    \cup
    \left\{\mathscr{X}(t)\leq\mathscr{Z}_{-}(t)\right\}.
\end{equation}
This inclusion is all that is needed. We now bound each event by choosing the sign of the spectral parameter so that $m\mapsto\rho_m(w)$ is monotone in the correct direction. For the spectral parameters below set
\begin{equation*}
    \mathcal{D}_{-}\overset{\textnormal{def}}{=}
    \begin{cases}
        (-\infty,0),&\textnormal{in continuous time},\\
        (-\mathfrak{U}^{-1},0),&\textnormal{in discrete time},
    \end{cases}
    \qquad
    \mathcal{D}_{+}\overset{\textnormal{def}}{=}(0,\mathfrak{L}).
\end{equation*}
Thus both choices lie in the spectral domain of Proposition~\ref{eigenfunction}.
\medskip

\noindent Recall from \eqref{characteristic} that
\begin{equation*}
    \rho_{m+1}(w)=\rho_m(w)\left(1-\frac{w}{\lambda_m}\right).
\end{equation*}
If $w\in\mathcal{D}_{-}$, then every factor $1-w/\lambda_m$ is at least one, so $m\mapsto\rho_m(w)$ is increasing and positive. Therefore
\begin{equation*}
    \left\{\mathscr{X}(t)\geq\mathscr{Z}_{+}(t)\right\}
    =
    \left\{\rho_{\mathscr{X}(t)}(w)\geq\rho_{\mathscr{Z}_{+}(t)}(w)\right\},
    \qquad w\in\mathcal{D}_{-}.
\end{equation*}
If $w\in\mathcal{D}_{+}$, then \textbf{(A1)} gives $0<1-w/\lambda_m\leq1$ for every $m$, so $m\mapsto\rho_m(w)$ is decreasing and positive. Hence
\begin{equation*}
    \left\{\mathscr{X}(t)\leq\mathscr{Z}_{-}(t)\right\}
    =
    \left\{\rho_{\mathscr{X}(t)}(w)\geq\rho_{\mathscr{Z}_{-}(t)}(w)\right\},
    \qquad w\in\mathcal{D}_{+}.
\end{equation*}
Markov's inequality then gives
\begin{align}
    \mathbb{P}_x\left[\mathscr{X}(t)\geq\mathscr{Z}_{+}(t)\right]
    &\leq
    \rho_{\mathscr{Z}_{+}(t)}(w)^{-1}
    \mathbb{E}_x\left[\rho_{\mathscr{X}(t)}(w)\right],
    && w\in\mathcal{D}_{-},\label{uppermarkov}\\
    \mathbb{P}_x\left[\mathscr{X}(t)\leq\mathscr{Z}_{-}(t)\right]
    &\leq
    \rho_{\mathscr{Z}_{-}(t)}(w)^{-1}
    \mathbb{E}_x\left[\rho_{\mathscr{X}(t)}(w)\right],
    && w\in\mathcal{D}_{+}.\label{lowermarkov}
\end{align}
Using Proposition~\ref{eigenfunction} and the branch-free product form of $\exp(\mathfrak{F}_y)$ from \eqref{phaseproductform}, the right-hand side can be rewritten exactly. The only small point is that $\mathfrak{F}_{y}(w,t)$ contains the spatial product up to $y(t)$, whereas $\rho_{y(t)}(w)$ contains the product only up to $y(t)-1$. This accounts for the remaining factor:
\begin{equation}\label{tailphaseidentity}
    \rho_{y(t)}(w)^{-1}\mathbb{E}_x\left[\rho_{\mathscr{X}(t)}(w)\right]
    =
    \left(1-\frac{w}{\lambda_{y(t)}}\right)\rho_x(w)
    \exp\left(\mathfrak{F}_y(w,t)\right).
\end{equation}
Taking $y=\mathscr{Z}_{+}$ in \eqref{tailphaseidentity} for the upper tail and $y=\mathscr{Z}_{-}$ for the lower tail, \eqref{uppermarkov}--\eqref{lowermarkov} become
\begin{align}
    \mathbb{P}_x\left[\mathscr{X}(t)\geq\mathscr{Z}_{+}(t)\right]
    &\leq
    \left(1-\frac{w}{\lambda_{\mathscr{Z}_{+}(t)}}\right)\rho_x(w)
    \exp\left(\mathfrak{F}_{\mathscr{Z}_{+}}(w,t)\right),
    && w\in\mathcal{D}_{-},\label{uppertailphase}\\
    \mathbb{P}_x\left[\mathscr{X}(t)\leq\mathscr{Z}_{-}(t)\right]
    &\leq
    \left(1-\frac{w}{\lambda_{\mathscr{Z}_{-}(t)}}\right)\rho_x(w)
    \exp\left(\mathfrak{F}_{\mathscr{Z}_{-}}(w,t)\right),
    && w\in\mathcal{D}_{+}.\label{lowertailphase}
\end{align}

\textbf{Step 3: Sign and size of the first derivative at the two boundaries.}
The upper boundary satisfies $\mathscr{Z}_{+}(t)\geq\mathscr{Z}(t)$, so \textbf{(ii)} of Proposition~\ref{specificcasesub} uses $\mathrm{K}_{\mathrm l}(t)=\mathfrak{B}_1^{-1}$ in the lower bound. Combining that lower bound with \eqref{boundary-displacement-bounds}, \eqref{localization-scale-comparison}, and \eqref{window-preconditions}, we get
\begin{align}
    \mathfrak{F}_{\mathscr{Z}_{+}}^{(1)}(0,t)
    &\geq
    \mathfrak{B}_1^{-1}\left(\mathscr{Z}_{+}(t)-\mathscr{Z}(t)\right)
    -2\gamma(t)\mathscr{Z}(t)-\frac{3\Lambda}{2}\notag\\
    &\geq
    \left(\mathfrak{B}_1^{-1}\kappa-2\mathrm{C}_2\right)\gamma(t)\Sigma(t)^2-\frac{3\Lambda}{2}\notag\\
    &\geq
    \mathrm{A}_{\star}\gamma(t)\Sigma(t)^2.\label{upperfirstderivative}
\end{align}
For the lower boundary, $\mathscr{Z}_{-}(t)<\mathscr{Z}(t)$, and the upper bound in \textbf{(ii)} uses $\mathrm{K}_{\mathrm u}(t)=\mathfrak{B}_1^{-1}$. Combining this upper bound with the lower-boundary inequality in \eqref{boundary-displacement-bounds}, the upper scale comparison in \eqref{localization-scale-comparison}, and the last inequality in \eqref{window-preconditions} gives
\begin{align}
    \mathfrak{F}_{\mathscr{Z}_{-}}^{(1)}(0,t)
    &\leq
    \mathfrak{B}_1^{-1}\left(\mathscr{Z}_{-}(t)-\mathscr{Z}(t)\right)
    +2\gamma(t)\mathscr{Z}(t)+\frac{3\Lambda}{2}\notag\\
    &\leq
    -\left(\mathfrak{B}_1^{-1}\kappa-2\mathrm{C}_2\right)\gamma(t)\Sigma(t)^2+\frac{3\Lambda}{2}\notag\\
    &\leq
    -\mathrm{A}_{\star}\gamma(t)\Sigma(t)^2.\label{lowerfirstderivative}
\end{align}
This is where the sign information in \textbf{(ii)} enters the proof: the phase at the upper endpoint has positive first derivative at the origin, while the phase at the lower endpoint has negative first derivative at the origin.

Because $\mathscr{Z}_{+}(t),\mathscr{Z}_{-}(t)\in\mathscr{C}_{t}$ for every $t\geq\mathrm{T}_{\mathrm{LD}}$, the uniform endpoint version of the strengthened second-derivative estimate in \textbf{(iv)} of Proposition~\ref{specificcasesub} applies to both phases. We add the condition $\mathrm{C}_{\mathscr{C}}\gamma(t)^{1/2}\leq1/2$ to the final choice of $\mathrm{T}_{\mathrm{LD}}$. Then, for all $t\geq\mathrm{T}_{\mathrm{LD}}$,
\begin{equation}\label{secondderivativepm}
    \frac{1}{2}\Sigma(t)^2
    \leq
    \mathfrak{F}_{\mathscr{Z}_{\pm}}^{(2)}(0,t)
    \leq
    \frac{3}{2}\Sigma(t)^2.
\end{equation}
Indeed, applying \textbf{(iv)} at $w=0$ gives
\begin{equation*}
    \left|\mathfrak{F}_{\mathscr{Z}_{\pm}}^{(2)}(0,t)-\Sigma(t)^2\right|
    \leq
    \mathrm{C}_{\mathscr{C}}\gamma(t)^{1/2}\Sigma(t)^2.
\end{equation*}
The constants are uniform over the two boundary endpoints because both endpoints lie in $\mathscr{C}_{t}$.

\textbf{Step 4: Evaluation at the stationary points.}
For readability, write
\begin{equation*}
    \mathrm{u}_{+}(t)\overset{\textnormal{def}}{=}\mathrm{u}_{\mathscr{Z}_{+}}(t),
    \qquad
    \mathrm{u}_{-}(t)\overset{\textnormal{def}}{=}\mathrm{u}_{\mathscr{Z}_{-}}(t).
\end{equation*}
Since the two boundary endpoints belong to $\mathscr{C}_{t}$, Proposition~\ref{axpro} supplies real stationary points $\mathrm{u}_{\pm}(t)$ after the deterministic lower time $\mathrm{T}_{\mathrm{A}}$. Increasing $\mathrm{T}_{\mathrm{LD}}$ to dominate $\mathrm{T}_{\mathrm{A}}$, these stationary points are defined for every $t\geq\mathrm{T}_{\mathrm{LD}}$ and obey
\begin{equation*}
    \left|\mathrm{u}_{\pm}(t)+\frac{\mathfrak{F}_{\mathscr{Z}_{\pm}}^{(1)}(0,t)}{\mathfrak{F}_{\mathscr{Z}_{\pm}}^{(2)}(0,t)}\right|
    \leq
    \mathrm{C}_{\mathrm{A}}\!\left(\frac{\mathfrak{F}_{\mathscr{Z}_{\pm}}^{(1)}(0,t)}{\mathfrak{F}_{\mathscr{Z}_{\pm}}^{(2)}(0,t)}\right)^{\!2}
\end{equation*}
by \eqref{secondconditiontosat}. Write
\begin{equation*}
    a_{\pm}(t)\overset{\textnormal{def}}{=}
    \frac{\mathfrak{F}_{\mathscr{Z}_{\pm}}^{(1)}(0,t)}{\mathfrak{F}_{\mathscr{Z}_{\pm}}^{(2)}(0,t)}.
\end{equation*}
By \eqref{upperfirstderivative}, \eqref{lowerfirstderivative}, and \eqref{secondderivativepm}, $a_{+}(t)>0$ and $a_{-}(t)<0$ for every $t\geq\mathrm{T}_{\mathrm{LD}}$. The complementary inequalities in \textbf{(ii)} of Proposition~\ref{specificcasesub}, together with \eqref{boundary-F-bound} and \eqref{localization-scale-comparison}, also give a deterministic constant $\mathrm{D}_{\partial}>0$ such that
\begin{equation*}
    \left|\mathfrak{F}_{\mathscr{Z}_{\pm}}^{(1)}(0,t)\right|
    \leq
    \mathrm{D}_{\partial}\left(\gamma(t)\Sigma(t)^2+1\right)
\end{equation*}
for both signs. Hence \eqref{secondderivativepm} gives $|a_{\pm}(t)|\leq \mathrm{D}_{a}(\gamma(t)+\Sigma(t)^{-2})$ for a deterministic constant $\mathrm{D}_{a}>0$, and so $a_{\pm}(t)\to0$. We add to the final choice of $\mathrm{T}_{\mathrm{LD}}$ the two conditions $2\mathrm{C}_{\mathrm{A}}|a_{\pm}(t)|<1$ and $|a_{\pm}(t)|<\mathrm{R}_{-}/3$ for both signs. Then the explicit bound gives $|\mathrm{u}_{\pm}(t)+a_{\pm}(t)|\leq \mathrm{C}_{\mathrm{A}}a_{\pm}(t)^2<|a_{\pm}(t)|/2$, so $\mathrm{u}_{\pm}(t)$ and $-a_{\pm}(t)$ have the same sign. The same inequality also gives $|\mathrm{u}_{\pm}(t)|\leq |a_{\pm}(t)|+\mathrm{C}_{\mathrm{A}}a_{\pm}(t)^2<\mathrm{R}_{-}/2$. Therefore
\begin{equation}\label{stationary-signs}
    \mathrm{u}_{+}(t)<0<\mathrm{u}_{-}(t),
    \qquad
    |\mathrm{u}_{+}(t)|,|\mathrm{u}_{-}(t)|<\frac{1}{2}\mathrm{R}_{-},
\end{equation}
where $\mathrm{R}_{-}$ is the holomorphic radius from \textbf{(i)} of Proposition~\ref{specificcasesub}. Since $\mathrm{R}_{-}\leq\mathfrak{U}^{-1}$ and $\mathrm{R}_{-}\leq\mathfrak{L}$, \eqref{stationary-signs} puts $\mathrm{u}_{+}(t)$ in $\mathcal{D}_{-}$ and $\mathrm{u}_{-}(t)$ in $\mathcal{D}_{+}$. Substituting these choices into \eqref{uppertailphase}--\eqref{lowertailphase} and using \eqref{partitionProb} gives
\begin{equation}\label{pointwisetailviastationarypoint}
\begin{aligned}
    &\mathbb{P}_x\left[\left|\mathscr{X}(t)-\mathscr{Z}(t)\right|\geq\kappa\gamma(t)\Sigma(t)^2\right]\\
    &\quad\leq
    \left(1-\frac{\mathrm{u}_{+}(t)}{\lambda_{\mathscr{Z}_{+}(t)}}\right)
    \rho_x(\mathrm{u}_{+}(t))
    \exp\left(\mathfrak{F}_{\mathscr{Z}_{+}}(\mathrm{u}_{+}(t),t)\right)\\
    &\qquad+
    \left(1-\frac{\mathrm{u}_{-}(t)}{\lambda_{\mathscr{Z}_{-}(t)}}\right)
    \rho_x(\mathrm{u}_{-}(t))
    \exp\left(\mathfrak{F}_{\mathscr{Z}_{-}}(\mathrm{u}_{-}(t),t)\right).
\end{aligned}
\end{equation}
We now estimate the two exponential terms. Since $\mathfrak{F}_{\mathscr{Z}_{\pm}}(0,t)=0$ and $\mathfrak{F}_{\mathscr{Z}_{\pm}}^{(1)}(\mathrm{u}_{\pm}(t),t)=0$, Taylor's theorem on the real segment between $0$ and $\mathrm{u}_{\pm}(t)$ gives
\begin{equation}\label{stationary-taylor}
    0
    =
    \mathfrak{F}_{\mathscr{Z}_{\pm}}(\mathrm{u}_{\pm}(t),t)
    +
    \frac{1}{2}
    \mathfrak{F}_{\mathscr{Z}_{\pm}}^{(2)}(\mathrm{v}_{\pm}(t),t)
    \mathrm{u}_{\pm}(t)^2,
\end{equation}
for some real $\mathrm{v}_{\pm}(t)$ between $0$ and $\mathrm{u}_{\pm}(t)$. By \eqref{stationary-signs}, the whole segment lies in $\mathcal{B}_{\mathrm{R}_{-}/2}(0)$. Applying the strengthened second-derivative estimate in \textbf{(iv)} of Proposition~\ref{specificcasesub} along this segment gives
\begin{equation*}
    \left|\mathfrak{F}_{\mathscr{Z}_{\pm}}^{(2)}(w,t)-\Sigma(t)^2\right|
    \leq
    \mathrm{C}_{\mathscr{C}}\left(\gamma(t)^{1/2}+|w|\right)\Sigma(t)^2
\end{equation*}
for every real $w$ between $0$ and $\mathrm{u}_{\pm}(t)$. We add the deterministic condition $\mathrm{C}_{\mathscr{C}}(\gamma(t)^{1/2}+|\mathrm{u}_{\pm}(t)|)\leq1/2$ for both signs to the final choice of $\mathrm{T}_{\mathrm{LD}}$, which is possible because $\gamma(t)\to0$ and $\mathrm{u}_{\pm}(t)\to0$. Therefore, for every $t\geq\mathrm{T}_{\mathrm{LD}}$,
\begin{equation}\label{segment-second-derivative}
    \frac{1}{2}\Sigma(t)^2
    \leq
    \mathfrak{F}_{\mathscr{Z}_{\pm}}^{(2)}(w,t)
    \leq
    2\Sigma(t)^2
\end{equation}
for every real $w$ between $0$ and $\mathrm{u}_{\pm}(t)$.

 The same segment estimate applied to the derivative gives a lower bound on the distance of the stationary point from the origin. Since
\begin{equation*}
    \mathfrak{F}_{\mathscr{Z}_{\pm}}^{(1)}(\mathrm{u}_{\pm}(t),t)
    -
    \mathfrak{F}_{\mathscr{Z}_{\pm}}^{(1)}(0,t)
    =
    \mathfrak{F}_{\mathscr{Z}_{\pm}}^{(2)}(\widetilde{\mathrm{v}}_{\pm}(t),t)\mathrm{u}_{\pm}(t)
\end{equation*}
for some $\widetilde{\mathrm{v}}_{\pm}(t)$ between $0$ and $\mathrm{u}_{\pm}(t)$, \eqref{upperfirstderivative}, \eqref{lowerfirstderivative}, and \eqref{segment-second-derivative} imply
\begin{equation}\label{stationary-lower-bound}
    |\mathrm{u}_{\pm}(t)|
    \geq
    \frac{\mathrm{A}_{\star}}{2}\gamma(t).
\end{equation}
Substituting \eqref{segment-second-derivative} and \eqref{stationary-lower-bound} into \eqref{stationary-taylor}, we obtain
\begin{equation}\label{phaseupperbound}
    \mathfrak{F}_{\mathscr{Z}_{\pm}}(\mathrm{u}_{\pm}(t),t)
    \leq
    -\frac{\mathrm{A}_{\star}^{2}}{16}\gamma(t)^2\Sigma(t)^2
    \leq
    -\frac{\mathrm{A}_{\star}^{2}}{16}\Sigma(t),
\end{equation}
where the final inequality uses $\gamma(t)\geq\Sigma(t)^{-1/2}$, again from Definition~\ref{gammadef}.

 It remains only to control the non-exponential factors in \eqref{pointwisetailviastationarypoint}. The bound \eqref{stationary-signs} places $\mathrm{u}_{+}(t)$ and $\mathrm{u}_{-}(t)$ in a fixed compact subset of $\mathcal{B}_{\mathrm{R}_{-}}(0)$. Since $x$ is fixed, $\rho_x$ is a fixed polynomial, and \textbf{(A1)} gives $\lambda_n\geq\mathfrak{L}$, there exists $\mathrm{M}_x>0$ such that
\begin{equation}\label{rhoxbounded}
    \left|1-\frac{\mathrm{u}_{+}(t)}{\lambda_{\mathscr{Z}_{+}(t)}}\right|
    |\rho_x(\mathrm{u}_{+}(t))|
    +
    \left|1-\frac{\mathrm{u}_{-}(t)}{\lambda_{\mathscr{Z}_{-}(t)}}\right|
    |\rho_x(\mathrm{u}_{-}(t))|
    \leq
    \mathrm{M}_x
\end{equation}
for every $t\geq\mathrm{T}_{\mathrm{LD}}$. Combining \eqref{pointwisetailviastationarypoint}, \eqref{phaseupperbound}, and \eqref{rhoxbounded} proves \eqref{localizationtailbound}, with $c=\mathrm{A}_{\star}^{2}/16$ and with $\mathrm{T}_{\mathrm{LD}}$ chosen to dominate all thresholds introduced above. We now increase $\mathrm{T}_{\mathrm{LD}}$ once more, if necessary, so that $\kappa\gamma(t)\leq\gamma(t)^{1/2}$ for every $t\geq\mathrm{T}_{\mathrm{LD}}$, which is possible because $\gamma(t)\to0$. If $\mathscr{X}(t)\notin\mathscr{C}_{t}$, then by the definition of the endpoint window $\left|\mathscr{X}(t)-\mathscr{Z}(t)\right|>\gamma(t)^{1/2}\Sigma(t)^2$, and hence $\left|\mathscr{X}(t)-\mathscr{Z}(t)\right|\geq\kappa\gamma(t)\Sigma(t)^2$. Thus \eqref{endpointwindowtailbound} follows directly from \eqref{localizationtailbound}, which completes the proof.
\end{proof}

\begin{proposition}[Almost-sure eventual endpoint-window localization]\label{moments}
Assume \textbf{(A1)} and \textbf{(A2)}. For every $x\in\mathbb{Z}_+$,
\begin{equation*}
    \mathbb{P}_{x}\left[\exists\,T_0\in\mathscr{T}\textnormal{ such that }\mathscr{X}(s)\in\mathscr{C}_{s}\textnormal{ for every }s\geq T_0\right]=1.
\end{equation*}
\end{proposition}
\begin{proof}
Proposition~\ref{localizationtail} gives an exponentially small estimate at one fixed time, whereas the event of interest concerns the whole future tail of the path. In discrete time the upgrade is a direct summation over the future time indices. In continuous time there is one extra issue: a path could, in principle, cross the endpoint window between two deterministic sampling times. We therefore introduce a fine deterministic mesh and separately bound the probability that the birth process makes too many jumps inside one mesh interval.

\textbf{Step 1: A fixed intrinsic window.}
Let $\kappa,c>0$ and $\mathrm{T}_{\mathrm{tail}}\in\mathscr{T}$ be the constants and lower time supplied by Proposition~\ref{localizationtail}. Fix
\begin{equation*}
    \mathrm{K}\overset{\textnormal{def}}{=}8\kappa.
\end{equation*}
For $u\in\mathscr{T}$, introduce the two one-time events
\begin{equation*}
    \mathfrak{I}_{\mathrm{K}}(u)
    \overset{\textnormal{def}}{=}
    \left\{\left|\mathscr{X}(u)-\mathscr{Z}(u)\right|\geq \mathrm{K}\gamma(u)\Sigma(u)^2\right\},
    \qquad
    \mathfrak{I}_{\kappa}(u)
    \overset{\textnormal{def}}{=}
    \left\{\left|\mathscr{X}(u)-\mathscr{Z}(u)\right|\geq \kappa\gamma(u)\Sigma(u)^2\right\}.
\end{equation*}
Since $\mathrm{K}>\kappa$, the inclusion $\mathfrak{I}_{\mathrm{K}}(u)\subseteq\mathfrak{I}_{\kappa}(u)$ holds for every $u$. Hence, for each fixed starting point $x\in\mathbb{Z}_+$, Proposition~\ref{localizationtail} gives a constant $\mathrm{M}_x>0$ such that
\begin{equation}\label{pathwise-one-time-tail}
    \mathbb{P}_x\left[\mathfrak{I}_{\mathrm{K}}(u)\right]
    \leq
    \mathbb{P}_x\left[\mathfrak{I}_{\kappa}(u)\right]
    \leq
    \mathrm{M}_x\mathrm{e}^{-c\Sigma(u)},
    \qquad u\geq \mathrm{T}_{\mathrm{tail}}.
\end{equation}
Since $\gamma(u)\to0$, we may increase $\mathrm{T}_{\mathrm{tail}}$ so that $\mathrm{K}\gamma(u)^{1/2}\leq1$ for all $u\geq\mathrm{T}_{\mathrm{tail}}$. For such $u$, the inclusion
\begin{equation}\label{endpoint-window-inclusion}
    \left\{\mathscr{X}(u)\notin\mathscr{C}_{u}\right\}
    \subseteq
    \mathfrak{I}_{\mathrm{K}}(u)
\end{equation}
holds, because leaving $\mathscr{C}_{u}$ means that $|\mathscr{X}(u)-\mathscr{Z}(u)|>\gamma(u)^{1/2}\Sigma(u)^2$, and the latter threshold is at least $\mathrm{K}\gamma(u)\Sigma(u)^2$.
We shall repeatedly use a deterministic lower growth bound for $\Sigma$. Set
\begin{equation*}
    \mathrm{a}_0\overset{\textnormal{def}}{=}\frac{\min(1,\mathfrak{U}^{-1},\mathfrak{L})}{\Lambda}.
\end{equation*}
In continuous time, Definition~\ref{themeananddeviation} gives $\mathrm{z}(u)=u/\Lambda$. In discrete time, the same definition gives $\mathrm{z}(u)=\Lambda^{-1}\sum_{r=0}^{u-1}(\alpha_r(1-\tau_r)+\beta_r\tau_r)$, and the lower bounds on $\alpha_r$ and $\beta_r$ in \textbf{(A1)} imply $\mathrm{z}(u)\geq \min(\mathfrak{U}^{-1},\mathfrak{L})u/\Lambda$. Therefore, in both regimes, $\mathrm{z}(u)\geq \mathrm{a}_0u$. Since $\mathscr{Z}(u)=\lfloor \mathrm{z}(u)\rceil$ differs from $\mathrm{z}(u)$ by at most $1/2$, there exists $\mathrm{T}_{\mathrm{z}}\in\mathscr{T}$ such that
\begin{equation}\label{Z-linear-lower-bound}
    \mathscr{Z}(u)\geq \frac{\mathrm{a}_0}{2}u,\qquad u>\mathrm{T}_{\mathrm{z}}.
\end{equation}
Combining this with the upper comparison $\mathscr{Z}(u)\leq \mathrm{C}_2\Sigma(u)^2$ from \textbf{(iv)} of Proposition~\ref{specificcasesub}, we obtain
\begin{equation}\label{Sigma-lower-growth}
    \Sigma(u)\geq \mathrm{c}_0u^{1/2},\qquad u>\max(\mathrm{T}_{\mathrm{z}},\mathrm{T}_3),
\end{equation}
where $\mathrm{c}_0\overset{\textnormal{def}}{=}(\mathrm{a}_0/(2\mathrm{C}_2))^{1/2}$ and $\mathrm{T}_3$ is the lower time in \textbf{(iv)} of Proposition~\ref{specificcasesub}.

\textbf{Step 2: The discrete-time upgrade.}
Assume first that $\mathscr{T}=\mathbb{Z}_+$. For $t\in\mathbb{Z}_+$, let
\begin{equation*}
    \mathcal{A}_{\mathrm{K}}(t)\overset{\textnormal{def}}{=}\bigcup_{s=t}^{\infty}\mathfrak{I}_{\mathrm{K}}(s).
\end{equation*}
The event $\mathcal{A}_{\mathrm{K}}(t)$ says that the path exits the $\mathrm{K}\gamma\Sigma^2$ window at least once after time $t$. By \eqref{endpoint-window-inclusion}, the event that $\mathscr{X}(s)$ leaves $\mathscr{C}_{s}$ for some $s\geq t$ is contained in $\mathcal{A}_{\mathrm{K}}(t)$. By Boole's inequality, \eqref{pathwise-one-time-tail}, and \eqref{Sigma-lower-growth}, for every $t>\max(\mathrm{T}_{\mathrm{tail}},\mathrm{T}_{\mathrm{z}},\mathrm{T}_3)$,
\begin{equation*}
    \mathbb{P}_x\left[\mathcal{A}_{\mathrm{K}}(t)\right]
    \leq
    \sum_{s=t}^{\infty}\mathbb{P}_x\left[\mathfrak{I}_{\kappa}(s)\right]
    \leq
    \mathrm{M}_x\sum_{s=t}^{\infty}\mathrm{e}^{-c\Sigma(s)}
    \leq
    \mathrm{M}_x\sum_{s=t}^{\infty}\mathrm{e}^{-c\mathrm{c}_0s^{1/2}}.
\end{equation*}
The series $\sum_{s=1}^{\infty}\mathrm{e}^{-c\mathrm{c}_0s^{1/2}}$ converges by comparison with an integral of the same stretched-exponential form, and hence the tail displayed above tends to zero as $t\to\infty$. Since the events $\mathcal{A}_{\mathrm{K}}(t)$ decrease as $t$ increases, continuity of probability from above gives
\begin{equation*}
    \mathbb{P}_x\left[\bigcap_{t=1}^{\infty}\mathcal{A}_{\mathrm{K}}(t)\right]
    =
    \lim_{t\to\infty}\mathbb{P}_x\left[\mathcal{A}_{\mathrm{K}}(t)\right]
    =
    0.
\end{equation*}
Thus the probability that the path leaves $\mathscr{C}_{s}$ at some discrete time $s\geq t$ tends to zero as $t\to\infty$. The almost-sure statement of the proposition then follows in the discrete-time regime from continuity of probability from above applied to the decreasing events $\{\exists\,s\geq t:\mathscr{X}(s)\notin\mathscr{C}_{s}\}$.

\textbf{Step 3: The continuous-time mesh.}
It remains to prove the same statement when $\mathscr{T}=\mathbb{R}_+$. For $t\in\mathbb{Z}_+$ define
\begin{equation*}
    \mathcal{A}_{\mathrm{K}}(t)\overset{\textnormal{def}}{=}\left\{\exists\,s\geq t\text{ such that }\mathfrak{I}_{\mathrm{K}}(s)\text{ occurs}\right\}.
\end{equation*}
Fix $\sigma\in(0,1)$, set $r_0\overset{\textnormal{def}}{=}0$, and define
\begin{equation*}
    r_{n+1}-r_n\overset{\textnormal{def}}{=}\frac{1}{\mathfrak{B}_1(n+1)^{\sigma}},\qquad n\geq0.
\end{equation*}
Because $\sigma<1$, the sequence $r_n$ tends to infinity, and the intervals
\begin{equation*}
    \mathcal{J}_n(t)\overset{\textnormal{def}}{=}[t+r_n,t+r_{n+1}],
    \qquad
    q_n(t)\overset{\textnormal{def}}{=}t+r_{n+1},
\end{equation*}
cover the whole half-line $[t,\infty)$. We also set
\begin{equation*}
    \mathrm{L}(t)\overset{\textnormal{def}}{=}\lceil\Sigma(t)\rceil,
    \qquad
    \mathcal{N}_n(t)\overset{\textnormal{def}}{=}
    \left\{\sup_{s_1,s_2\in\mathcal{J}_n(t)}
    \left|\mathscr{X}(s_2)-\mathscr{X}(s_1)\right|<\mathrm{L}(t)\right\}.
\end{equation*}
The event $\mathcal{N}_n(t)$ is the good-oscillation event on the interval $\mathcal{J}_n(t)$. We next choose a deterministic lower time $\mathrm{T}_{\mathrm{mesh}}$ such that, for every $t>\mathrm{T}_{\mathrm{mesh}}$, a violation of the $\mathrm{K}$-window somewhere inside $\mathcal{J}_n(t)$ either forces a violation of the $\kappa$-window at the right endpoint $q_n(t)$ or else the good-oscillation event fails.
\medskip

\noindent The deterministic centre moves by a bounded amount over every mesh interval. Indeed, in continuous time $\mathrm{z}(u)=u/\Lambda$, and therefore, for any $s\in\mathcal{J}_n(t)$,
\begin{equation}\label{center-mesh-motion}
    \left|\mathscr{Z}(q_n(t))-\mathscr{Z}(s)\right|
    \leq
    \left|\mathrm{z}(q_n(t))-\mathrm{z}(s)\right|+1
    \leq
    \frac{1}{\Lambda\mathfrak{B}_1}+1
    \overset{\textnormal{def}}{=}\mathrm{D}_{\mathrm{Z}}.
\end{equation}
Here the additional $1$ accounts for the two nearest-integer rounding errors in $\mathscr{Z}=\lfloor\mathrm{z}\rceil$. Since, in continuous time, $\Sigma(u)^2=\sum_{r=0}^{\mathscr{Z}(u)}\lambda_r^{-2}$, \eqref{center-mesh-motion} and the lower spatial-rate bound in \textbf{(A1)} imply
\begin{equation}\label{Sigma-mesh-motion}
    0\leq \Sigma(q_n(t))^2-\Sigma(s)^2\leq \left(\lceil\mathrm{D}_{\mathrm{Z}}\rceil+1\right)\mathfrak{L}^{-2}\overset{\textnormal{def}}{=}\mathrm{D}_{\Sigma}.
\end{equation}
The left inequality uses the monotonicity of $\mathscr{Z}$, and the right inequality uses the fact that at most $\lceil\mathrm{D}_{\mathrm{Z}}\rceil+1$ additional summands can be inserted between $\mathscr{Z}(s)$ and $\mathscr{Z}(q_n(t))$. Since $\Sigma(t)\to\infty$ by \eqref{Sigma-lower-growth}, there exists $\mathrm{T}_{\Sigma}\in\mathbb{Z}_+$ such that, whenever $t>\mathrm{T}_{\Sigma}$, $n\geq0$, and $s\in\mathcal{J}_n(t)$,
\begin{equation}\label{Sigma-comparable-on-mesh}
    \Sigma(s)^2\geq\frac{1}{2}\Sigma(q_n(t))^2.
\end{equation}
Moreover, $\gamma$ is non-increasing by Definition~\ref{gammadef}, because it is the maximum of two tail suprema over $[t,\infty)$. Hence, for $s\leq q_n(t)$, \eqref{Sigma-comparable-on-mesh} gives
\begin{equation}\label{window-comparison-on-mesh}
    \mathrm{K}\gamma(s)\Sigma(s)^2
    \geq
    4\kappa\gamma(q_n(t))\Sigma(q_n(t))^2.
\end{equation}
We shall also use that the deterministic and path oscillation errors are smaller than one copy of the $\kappa$-window. By Definition~\ref{gammadef}, $\gamma(u)\geq \Sigma(u)^{-1/2}$; since $\Sigma(q_n(t))\geq\Sigma(t)$ in continuous time, increasing $\mathrm{T}_{\mathrm{mesh}}$ if necessary gives, uniformly in $n$ and for every $t>\mathrm{T}_{\mathrm{mesh}}$,
\begin{equation}\label{mesh-error-absorbed}
    \mathrm{L}(t)+\mathrm{D}_{\mathrm{Z}}
    \leq
    \kappa\gamma(q_n(t))\Sigma(q_n(t))^2.
\end{equation}
Now suppose that $\mathfrak{I}_{\mathrm{K}}(s)$ occurs for some $s\in\mathcal{J}_n(t)$ and that $\mathcal{N}_n(t)$ also occurs. Combining the triangle inequality with \eqref{center-mesh-motion}, \eqref{window-comparison-on-mesh}, and \eqref{mesh-error-absorbed}, we obtain
\begin{align*}
    \left|\mathscr{X}(q_n(t))-\mathscr{Z}(q_n(t))\right|
    &\geq
    \left|\mathscr{X}(s)-\mathscr{Z}(s)\right|
    -
    \left|\mathscr{X}(q_n(t))-\mathscr{X}(s)\right|
    -
    \left|\mathscr{Z}(q_n(t))-\mathscr{Z}(s)\right|\\
    &\geq
    \mathrm{K}\gamma(s)\Sigma(s)^2-\mathrm{L}(t)-\mathrm{D}_{\mathrm{Z}}\\
    &\geq
    3\kappa\gamma(q_n(t))\Sigma(q_n(t))^2.
\end{align*}
In particular, $\mathfrak{I}_{\kappa}(q_n(t))$ occurs. Thus, for every $t>\mathrm{T}_{\mathrm{mesh}}$,
\begin{equation}\label{continuous-event-inclusion}
    \mathcal{A}_{\mathrm{K}}(t)
    \subseteq
    \bigcup_{n=0}^{\infty}\mathfrak{I}_{\kappa}(q_n(t))
    \cup
    \bigcup_{n=0}^{\infty}\mathcal{N}_n(t)^c.
\end{equation}

\textbf{Step 4: Summing the continuous-time bounds.}
We estimate the two terms in \eqref{continuous-event-inclusion}. First, by \eqref{pathwise-one-time-tail},
\begin{equation}\label{continuous-tail-sum}
    \sum_{n=0}^{\infty}\mathbb{P}_x\left[\mathfrak{I}_{\kappa}(q_n(t))\right]
    \leq
    \mathrm{M}_x\sum_{n=0}^{\infty}\mathrm{e}^{-c\Sigma(q_n(t))}
\end{equation}
for every $t\geq\mathrm{T}_{\mathrm{tail}}$, since $q_n(t)\geq t$ for every $n\geq0$. For each fixed $n$, the summand in \eqref{continuous-tail-sum} tends to zero as $t\to\infty$, because $q_n(t)\to\infty$ and $\Sigma(u)\to\infty$ by \eqref{Sigma-lower-growth}. To justify passing the limit through the infinite sum, note that \eqref{Sigma-lower-growth} gives
\begin{equation*}
    \mathrm{e}^{-c\Sigma(q_n(t))}
    \leq
    \mathrm{e}^{-c\mathrm{c}_0r_{n+1}^{1/2}}
\end{equation*}
for every $t$ beyond the lower scale threshold in \eqref{Sigma-lower-growth}, since $q_n(t)=t+r_{n+1}\geq r_{n+1}$. Also,
\begin{equation*}
    r_{n+1}
    =
    \frac{1}{\mathfrak{B}_1}\sum_{j=1}^{n+1}j^{-\sigma}
    \geq
    \frac{(n+2)^{1-\sigma}-1}{\mathfrak{B}_1(1-\sigma)},
\end{equation*}
by comparison of the sum with the integral of $u^{-\sigma}$. Hence $\sum_{n\geq0}\exp\{-c\mathrm{c}_0r_{n+1}^{1/2}\}<\infty$, since the displayed lower bound implies a stretched-exponential bound in $n$ after a finite index depending only on $\sigma$ and $\mathfrak{B}_1$. The dominated convergence theorem for series therefore gives
\begin{equation}\label{continuous-tail-sum-limit}
    \lim_{t\to\infty}
    \sum_{n=0}^{\infty}\mathbb{P}_x\left[\mathfrak{I}_{\kappa}(q_n(t))\right]
    =
    0.
\end{equation}

It remains to bound the bad-oscillation term. On $\mathcal{N}_n(t)^c$, the non-decreasing birth process makes at least $\mathrm{L}(t)$ jumps during $\mathcal{J}_n(t)$. By the upper rate bound $\lambda_r<\mathfrak{B}_1$ in \textbf{(A1)}, the jump count on $\mathcal{J}_n(t)$ is stochastically dominated by a Poisson random variable with rate $\mathfrak{B}_1$. One way to see this is to construct a rate-$\mathfrak{B}_1$ Poisson clock and accept a proposed jump from site $r$ with probability $\lambda_r/\mathfrak{B}_1$; the accepted process has the required birth rates and is bounded above by the proposing Poisson process. Thus the number of jumps in $\mathcal{J}_n(t)$ is dominated by a Poisson random variable with mean
\begin{equation*}
    \mu_n\overset{\textnormal{def}}{=}\mathfrak{B}_1|\mathcal{J}_n(t)|=(n+1)^{-\sigma}.
\end{equation*}
Since $\mu_n\leq1$, the Poisson tail is bounded by
\begin{equation*}
    \mathbb{P}\left[\mathrm{Poisson}(\mu_n)\geq \mathrm{L}(t)\right]
    =
    \mathrm{e}^{-\mu_n}\sum_{\ell=\mathrm{L}(t)}^{\infty}\frac{\mu_n^\ell}{\ell!}
    \leq
    \sum_{j=0}^{\infty}\frac{\mu_n^{\mathrm{L}(t)+j}}{\mathrm{L}(t)!j!}
    =
    \mathrm{e}^{\mu_n}\frac{\mu_n^{\mathrm{L}(t)}}{\mathrm{L}(t)!}
    \leq
    \mathrm{e}\,\frac{(n+1)^{-\sigma\mathrm{L}(t)}}{\mathrm{L}(t)!}.
\end{equation*}
Choose $\mathrm{T}_{\mathrm{osc}}\in\mathscr{T}$ such that $\sigma\mathrm{L}(t)>1$ for every $t>\mathrm{T}_{\mathrm{osc}}$. For such $t$, summing this estimate over $n$ gives
\begin{equation}\label{bad-oscillation-sum}
    \sum_{n=0}^{\infty}\mathbb{P}_x\left[\mathcal{N}_n(t)^c\right]
    \leq
    \frac{\mathrm{e}}{\mathrm{L}(t)!}\sum_{n=0}^{\infty}(n+1)^{-\sigma\mathrm{L}(t)}
    =
    \frac{\mathrm{e}\,\zeta(\sigma\mathrm{L}(t))}{\mathrm{L}(t)!},
\end{equation}
and the right-hand side tends to zero as $t\to\infty$, because $\mathrm{L}(t)=\lceil\Sigma(t)\rceil\to\infty$ while $\zeta(\sigma\mathrm{L}(t))\to1$. Combining \eqref{continuous-event-inclusion}, \eqref{continuous-tail-sum-limit}, and \eqref{bad-oscillation-sum}, we obtain
\begin{equation*}
    \lim_{t\to\infty}\mathbb{P}_x\left[\mathcal{A}_{\mathrm{K}}(t)\right]=0.
\end{equation*}
The events $\mathcal{A}_{\mathrm{K}}(t)$ decrease as $t$ increases through the integers, so continuity from above gives
\begin{equation*}
    \mathbb{P}_x\left[\bigcap_{t=1}^{\infty}\mathcal{A}_{\mathrm{K}}(t)\right]=0.
\end{equation*}
Consequently, with probability one, there is a random finite time after which
\begin{equation*}
    \left|\mathscr{X}(s)-\mathscr{Z}(s)\right|<\mathrm{K}\gamma(s)\Sigma(s)^2
\end{equation*}
for every $s$ in the continuous time set. Since $\gamma(s)\to0$, there exists $\mathrm{T}_{\mathscr{C}}\in\mathscr{T}$ such that $\mathrm{K}\gamma(s)^{1/2}\leq1$ for every $s>\mathrm{T}_{\mathscr{C}}$; hence $\mathrm{K}\gamma(s)\Sigma(s)^2\leq\gamma(s)^{1/2}\Sigma(s)^2$ for every $s>\mathrm{T}_{\mathscr{C}}$. Therefore $\mathscr{X}(s)\in\mathscr{C}_{s}$ after the maximum of the random finite time above and $\mathrm{T}_{\mathscr{C}}$, which proves the proposition in the continuous-time regime.
\end{proof}
\begin{figure}[!ht]
    \centering
    \includegraphics[width=0.78\textwidth]{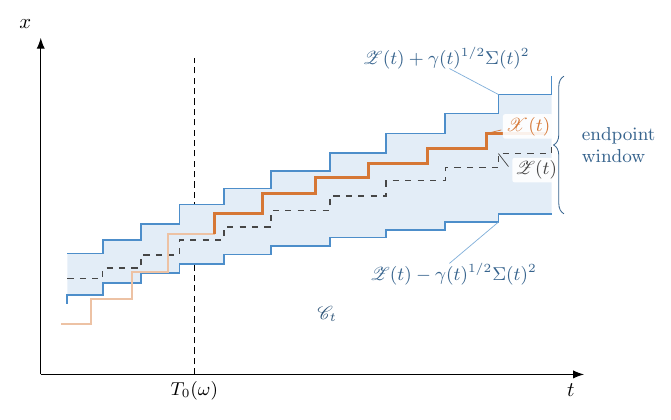}
    \caption{Schematic illustration of Proposition~\ref{moments}. With probability one, there is a finite random time $T_0$ after which the one-particle path remains in the endpoint window $\mathscr{C}_{t}$ around the deterministic centre $\mathscr{Z}(t)$. The blue staircase boundaries reflect the integer-valued endpoint window in the discrete spatial coordinate.}
    \label{fig:prop33-endpoint-cone}
\end{figure}
\subsection{Martingale control}\label{martingalecontrol}
This subsection provides the moment input that is used to lift the one-particle estimates of Section~\ref{decouplingsection} to the multi-particle determinant analysis. The remaining content of Section~\ref{stationarypointprocess} has three components: martingale/moment control for the stationary-point coordinates (Proposition~\ref{boundedmoments}), the Gaussian limit for the stationary-point process (Proposition~\ref{NormalisedProcess}), and the determinant expansion of the Karlin--McGregor semigroup (Proposition~\ref{ThmKM}). Throughout, the endpoint restriction is encoded by the coordinatewise window $\mathscr{C}_{t}^{N}$ introduced before Proposition~\ref{localizationtail}.
    \begin{proposition}\label{boundedmoments}
Assume \textbf{(A1)} and fix $\mathbf{x}\in\mathbb{Z}_{+}^{N}$. For every $p\geq2$ and every $i\in\llbracket N\rrbracket$, there exists $\mathrm{C}_{p,\mathbf{x}}>0$ such that
\begin{equation}\label{eq:first-derivative-moment}
\mathbb{E}_{\mathbf{x}}^N\!\left[\left|\mathfrak{F}^{(1)}_{\mathscr{X}_i}(0,t)\right|^p\right]\leq \mathrm{C}_{p,\mathbf{x}}(1+t)^{p/2},\qquad t\in\mathscr{T}.
\end{equation}
If in addition \textbf{(A2)} holds, then for every $p\geq2$ there exist $\mathrm{T}_{p,N}\in\mathscr{T}$ and $\mathrm{C}_{p,N,\mathbf{x}}>0$ such that
\begin{equation}\label{eq:Up-moment}
\sup_{t\ge \mathrm{T}_{p,N}}\mathbb{E}_{\mathbf{x}}^N\!\left[|\mathbf{U}(t)|^p\mathbf{1}_{\{\mathbf{X}(t)\in\mathscr{C}_{t}^{N}\}}\right]\leq \mathrm{C}_{p,N,\mathbf{x}}.
\end{equation}
		\begin{proof}
			The proof has two parts. First we prove a moment bound for the centered inverse-rate coordinate of one particle. This is the probabilistic input behind \eqref{eq:first-derivative-moment}. Second, under \textbf{(A2)}, we use the endpoint-window stationary-point estimates to transfer that moment control to $\mathbf{U}(t)$ on $\mathscr{C}_{t}^{N}$.
		
        \textbf{Step 1: The inverse-rate martingale.} Set
		\begin{equation*}
		H(m)\overset{\textnormal{def}}{=}\sum_{n=0}^{m-1}\frac{1}{\lambda_n},\qquad m\in\mathbb{Z}_{+}.
		\end{equation*}
		We first prove the required moment estimate for one coordinate process $X$ started from a fixed point $x\in\mathbb{Z}_{+}$. In continuous time, the generator from Definition~\ref{contprocess} is $\mathsf{L}f(m)=\lambda_m(f(m+1)-f(m))$. Since $H$ is unbounded, we localize it: for $K\in\mathbb{N}$ with $K>x$ let $H_K(m)=H(m\wedge K)$ and $\tau_K=\inf\{t\geq0:X(t)\geq K\}$. Dynkin's formula for the bounded test function $H_K$ \cite[Chapter~4]{EthierKurtz1986} gives a martingale, and before $\tau_K$ the generator satisfies
		\begin{equation*}
		(\mathsf{L}H)(m)=\lambda_m\left(\frac{1}{\lambda_m}\right)=1.
		\end{equation*}
		Hence $H(X(t\wedge\tau_K))-H(X(0))-(t\wedge\tau_K)$ is a martingale. The upper rate bound in \textbf{(A1)} gives $\lambda_m<\mathfrak{B}_1$ for every state $m$, so the process can be constructed by thinning a rate-$\mathfrak{B}_1$ Poisson clock: when the dominating clock rings at state $m$, accept the birth with probability $\lambda_m/\mathfrak{B}_1$. This construction gives non-explosion and $\tau_K\uparrow\infty$ almost surely. For each fixed deterministic time, $H(X(t\wedge\tau_K))\leq H(x)+\mathfrak{L}^{-1}\Pi_{\mathfrak{B}_1}(t)$, where $\Pi_{\mathfrak{B}_1}(t)$ is the number of rings of the dominating Poisson clock up to time $t$. Thus the stopped martingale variables are dominated in absolute value by $2H(x)+\mathfrak{L}^{-1}\Pi_{\mathfrak{B}_1}(t)+t$. Hence the stopped martingales are uniformly integrable on each finite time interval, and letting $K\rightarrow\infty$ gives that $H(X(t))-H(x)-t$ is a true martingale. Its jumps have size at most $\mathfrak{L}^{-1}$, and its predictable quadratic variation is
		\begin{equation*}
		\int_{0}^{t}\lambda_{X(s)}\left(\frac{1}{\lambda_{X(s)}}\right)^2\mathrm{d}s=\int_{0}^{t}\frac{1}{\lambda_{X(s)}}\mathrm{d}s\leq\frac{t}{\mathfrak{L}}.
		\end{equation*}
		For the moment estimate we use the Burkholder--Davis--Gundy inequality for c\`adl\`ag martingales in \cite[Chapter~IV, Theorem~48]{Protter2005}. Applied at a fixed time, this theorem controls the martingale by its optional quadratic variation. In the present pure-birth case,
		\begin{equation*}
		[H(X(\cdot))-\cdot,H(X(\cdot))-\cdot]_{t}=\sum_{0<s\leq t}\left(H(X(s))-H(X(s-))\right)^2.
		\end{equation*}
		Each summand is bounded by $\mathfrak{L}^{-2}$, and the number of jumps up to time $t$ is bounded by the same rate-$\mathfrak{B}_1$ Poisson clock used above. Therefore the optional quadratic variation is bounded by $\mathfrak{L}^{-2}\Pi_{\mathfrak{B}_1}(t)$ under this coupling. Since $\mathbb{E}[\Pi_{\mathfrak{B}_1}(t)^{p/2}]\leq \mathrm{C}_{p,\mathfrak{B}_1}(1+t)^{p/2}$ for every $p\geq2$, the c\`adl\`ag BDG inequality applied to $H(X(t))-H(x)-t$ gives the required $p$th moment bound for the martingale part. Absorbing the fixed term $H(x)$ by the triangle inequality, there exists $\mathrm{C}_{p,x}>0$ such that
		\begin{equation}\label{eq:H-moment-continuous}
		\mathbb{E}_{x}\left[\left|H(X(t))-t\right|^p\right]\leq\mathrm{C}_{p,x}(1+t)^{p/2},\qquad t\in\mathbb{R}_{+}.
		\end{equation}
       
        In discrete time, let $X$ be one coordinate process from Definition~\ref{mixprocess}, with natural filtration $(\mathcal{F}_{t})_{t\in\mathbb{Z}_{+}}$. If $\tau_t=0$, then the Bernoulli jump probability is $\alpha_t\lambda_{X(t)}$, and therefore
		\begin{equation*}
		\mathbb{E}\!\left[H(X(t+1))-H(X(t))\mid\mathcal{F}_{t}\right]=\alpha_t.
		\end{equation*}
		If $\tau_t=1$ and $X(t)=m$, then the tail-sum formula for the non-negative increment $H(X(t+1))-H(X(t))$ gives
		\begin{align*}
		\mathbb{E}\!\left[H(X(t+1))-H(X(t))\mid X(t)=m\right]
		&=\sum_{n=m}^{\infty}\frac{1}{\lambda_n}\prod_{r=m}^{n}\frac{\beta_t\lambda_r}{1+\beta_t\lambda_r}\\
		&=\beta_t\sum_{n=m}^{\infty}\left(\prod_{r=m}^{n-1}\frac{\beta_t\lambda_r}{1+\beta_t\lambda_r}-\prod_{r=m}^{n}\frac{\beta_t\lambda_r}{1+\beta_t\lambda_r}\right).
		\end{align*}
		The product at the upper endpoint tends to zero because \textbf{(A1)} gives
		\begin{equation*}
		\frac{\beta_t\lambda_r}{1+\beta_t\lambda_r}\leq\frac{\mathfrak{U}\mathfrak{B}_1}{1+\mathfrak{U}\mathfrak{B}_1}<1.
		\end{equation*}
		Thus the telescoping sum equals $\beta_t$, and the process
		\begin{equation*}
		H(X(t))-\sum_{s=0}^{t-1}\big((1-\tau_s)\alpha_s+\tau_s\beta_s\big)
		\end{equation*}
		is a martingale. Its one-step martingale differences have uniformly bounded conditional $p$th moments. Indeed, a Bernoulli update makes at most one spatial jump, while in a geometric update the probability of at least $r$ jumps is bounded by $q_{\mathrm{geo}}^{r}$ with $q_{\mathrm{geo}}\overset{\textnormal{def}}{=}\mathfrak{U}\mathfrak{B}_1/(1+\mathfrak{U}\mathfrak{B}_1)<1$; multiplying by the factor $\mathfrak{L}^{-1}$ converting spatial jumps into $H$-increments gives a finite bound for the conditional $p$th moment of $H(X(t+1))-H(X(t))$, depending only on $p$ and the constants in \textbf{(A1)}. The deterministic centring term subtracted in the martingale difference is either $\alpha_t$ or $\beta_t$, both bounded above by constants from \textbf{(A1)}, so the martingale differences themselves have conditional $p$th moments bounded by a deterministic constant $\mathrm{D}_{p}$. Let
        \begin{equation*}
            \Delta_s\overset{\textnormal{def}}{=}
            H(X(s+1))-H(X(s))-\big((1-\tau_s)\alpha_s+\tau_s\beta_s\big).
        \end{equation*}
        Burkholder's discrete martingale inequality \cite[Chapter~2, Theorem~2.10]{HallHeyde1980} gives
        \begin{equation*}
            \mathbb{E}\left|\sum_{s=0}^{t-1}\Delta_s\right|^p
            \leq
            \mathrm{C}_p\,\mathbb{E}\left(\sum_{s=0}^{t-1}\Delta_s^2\right)^{p/2}.
        \end{equation*}
        For $t\geq1$, we use $(\sum_{s=0}^{t-1}a_s^2)^{p/2}\leq t^{p/2-1}\sum_{s=0}^{t-1}|a_s|^p$, while $t=0$ is absorbed into the fixed initial value. The uniform conditional moment bound on the $\Delta_s$ therefore gives a bound of order $(1+t)^{p/2}$ for the centred martingale. Absorbing $H(x)$ by the triangle inequality gives a constant $\mathrm{C}_{p,x}>0$ such that
		\begin{equation}\label{eq:H-moment-discrete}
		\mathbb{E}_{x}\!\left[\left|H(X(t))-\sum_{s=0}^{t-1}\big((1-\tau_s)\alpha_s+\tau_s\beta_s\big)\right|^p\right]\leq\mathrm{C}_{p,x}(1+t)^{p/2},\qquad t\in\mathbb{Z}_{+}.
		\end{equation}
		Combining \eqref{eq:H-moment-continuous} and \eqref{eq:H-moment-discrete}, and applying the resulting one-coordinate estimate to $X=\mathscr{X}_i$, we get, in both time regimes,
		\begin{equation}\label{eq:H-centered-moment}
		\mathbb{E}_{\mathbf{x}}^N\!\left[\left|H(\mathscr{X}_i(t))-\mathbf{1}_{\{\mathscr{T}=\mathbb{R}_+\}}t-\mathbf{1}_{\{\mathscr{T}=\mathbb{Z}_+\}}\sum_{s=0}^{t-1}\big((1-\tau_s)\alpha_s+\tau_s\beta_s\big)\right|^p\right]\leq\mathrm{C}_{p,\mathbf{x}}(1+t)^{p/2}.
		\end{equation}
		By the explicit formula for the first derivative of the phase in Definition~\ref{form},
		\begin{equation*}
		\mathfrak{F}^{(1)}_{\mathscr{X}_i}(0,t)=H(\mathscr{X}_i(t))-\mathbf{1}_{\{\mathscr{T}=\mathbb{R}_+\}}t-\mathbf{1}_{\{\mathscr{T}=\mathbb{Z}_+\}}\sum_{s=0}^{t-1}\big((1-\tau_s)\alpha_s+\tau_s\beta_s\big)+\frac{1}{\lambda_{\mathscr{X}_i(t)}}.
		\end{equation*}
		The final term is bounded by $\mathfrak{L}^{-1}$ under \textbf{(A1)}. Therefore \eqref{eq:H-centered-moment} implies \eqref{eq:first-derivative-moment}, after increasing $\mathrm{C}_{p,\mathbf{x}}$.
	
    \textbf{Step 2: Moment bounds for the stationary-point coordinates.}
	Assume now \textbf{(A2)}. On the endpoint event $\{\mathscr{X}_i(t)\in\mathscr{C}_{t}\}$, the stationary-point construction in Proposition~\ref{axpro} is uniform in the terminal value. Let $\mathrm{T}_{\mathrm{U}}\in\mathscr{T}$ be chosen to dominate $\mathrm{T}_{\mathrm{A}}$ and the threshold in \textbf{(iv)} of Proposition~\ref{specificcasesub} at which $\mathrm{C}_{\mathscr{C}}\gamma(t)^{1/2}\leq1/2$. Then the endpoint second-derivative comparison gives $\mathfrak{F}^{(2)}_{\mathscr{X}_i}(0,t)\geq\Sigma(t)^2/2$ whenever $t>\mathrm{T}_{\mathrm{U}}$ and $\mathscr{X}_i(t)\in\mathscr{C}_{t}$. The same comparison also gives $\Sigma(t)^2/\mathfrak{F}^{(2)}_{\mathscr{X}_i}(0,t)\leq1+2\mathrm{C}_{\mathscr{C}}\gamma(t)^{1/2}$, because $(1-r)^{-1}\leq1+2r$ for $0\leq r\leq1/2$. Multiplying \eqref{secondconditiontosat} by $\Sigma(t)$ and using these two consequences of the second-derivative comparison gives
		\begin{equation}\label{eq:U-first-derivative-comparison}
		|\mathscr{U}_i(t)|\leq\frac{\left|\mathfrak{F}^{(1)}_{\mathscr{X}_i}(0,t)\right|}{\Sigma(t)}+\frac{4\mathrm{C}_{\mathrm{A}}\left|\mathfrak{F}^{(1)}_{\mathscr{X}_i}(0,t)\right|^2}{\Sigma(t)^3}+\frac{2\mathrm{C}_{\mathscr{C}}\gamma(t)^{1/2}\left|\mathfrak{F}^{(1)}_{\mathscr{X}_i}(0,t)\right|}{\Sigma(t)}
		\end{equation}
			on this event. This deterministic estimate transfers the martingale bound to the stationary-point coordinate. We next compare $\Sigma(t)^2$ with time. The deterministic centre estimate \eqref{Zlinearintbounds}, proved from Definition~\ref{themeananddeviation} and \textbf{(A1)}, gives a threshold $\mathrm{T}_{\mathrm{Z}}\in\mathscr{T}$ such that $\mathscr{Z}(t)\geq \mathfrak{A}t/2$ for every $t>\mathrm{T}_{\mathrm{Z}}$. Combining this with the upper comparison $\mathscr{Z}(t)\leq \mathrm{C}_2\Sigma(t)^2$ from \textbf{(iv)} of Proposition~\ref{specificcasesub}, there exist constants $\mathrm{c}_{\Sigma}>0$ and $\mathrm{T}_{\Sigma}\in\mathscr{T}$ such that
	\begin{equation*}
	\Sigma(t)^2\geq \mathrm{c}_{\Sigma}(1+t),\qquad t>\mathrm{T}_{\Sigma}.
	\end{equation*}
		Increase the lower time once more so that $\gamma(t)\leq1$. Taking the $p$th power in \eqref{eq:U-first-derivative-comparison}, using $(a+b+c)^p\leq3^{p-1}(a^p+b^p+c^p)$, and applying \eqref{eq:first-derivative-moment} with exponents $p$ and $2p$, gives
		\begin{equation*}
		\sup_{t\ge \mathrm{T}_{p}}\mathbb{E}_{\mathbf{x}}^N\!\left[|\mathscr{U}_i(t)|^p\mathbf{1}_{\{\mathscr{X}_i(t)\in\mathscr{C}_{t}\}}\right]<\infty
		\end{equation*}
			for a deterministic threshold $\mathrm{T}_{p}\in\mathscr{T}$. Indeed, for every $t>\mathrm{T}_{p}$, the first and third terms in \eqref{eq:U-first-derivative-comparison} are controlled by
		\begin{equation*}
		\frac{\mathbb{E}_{\mathbf{x}}^N\!\left[\left|\mathfrak{F}^{(1)}_{\mathscr{X}_i}(0,t)\right|^p\right]}{\Sigma(t)^p}\leq \mathrm{C}_{p,\mathbf{x}}\mathrm{c}_{\Sigma}^{-p/2},
		\end{equation*}
		where we used $\gamma(t)\leq1$ for the third term. The quadratic term is controlled by
		\begin{equation*}
		\frac{\mathbb{E}_{\mathbf{x}}^N\!\left[\left|\mathfrak{F}^{(1)}_{\mathscr{X}_i}(0,t)\right|^{2p}\right]}{\Sigma(t)^{3p}}\leq \mathrm{C}_{2p,\mathbf{x}}\mathrm{c}_{\Sigma}^{-3p/2}(1+t)^{-p/2}\leq \mathrm{C}_{2p,\mathbf{x}}\mathrm{c}_{\Sigma}^{-3p/2}.
		\end{equation*}
		These two estimates use \eqref{eq:first-derivative-moment} and the lower bound $\Sigma(t)^2\geq \mathrm{c}_{\Sigma}(1+t)$, and therefore give a bound independent of $t$. Finally,
	\begin{equation*}
	|\mathbf{U}(t)|^p\le N^{p-1}\sum_{i=1}^{N}|\mathscr{U}_i(t)|^p,
	\end{equation*}
		and the event $\{\mathbf{X}(t)\in\mathscr{C}_{t}^{N}\}$ is contained in each coordinate event $\{\mathscr{X}_i(t)\in\mathscr{C}_{t}\}$. Summing the one-coordinate bounds over $i$ proves \eqref{eq:Up-moment}, after choosing $\mathrm{C}_{p,N,\mathbf{x}}$ larger than the resulting finite sum.
	\end{proof}
    \end{proposition}
    \begin{proposition}[Quantitative Gaussian limit on the endpoint window]\label{NormalisedProcess}
Assume \textbf{(A1)} and \textbf{(A2)}. Fix $\mathbf{x}\in\mathbb{Z}_+^{N}$. Let $H:\mathbb{R}^{N}\to\mathbb{C}$ be continuous and suppose that there are constants $\mathrm{C}_{H}>0$ and $q\geq0$ such that
\begin{equation*}
    |H(\mathbf{v})|\leq \mathrm{C}_{H}(1+|\mathbf{v}|^{q}),
    \qquad
    |H(\mathbf{v})-H(\mathbf{w})|\leq \mathrm{C}_{H}(1+\mathrm{R}^{q})|\mathbf{v}-\mathbf{w}|,
\end{equation*}
whenever $\max(|\mathbf{v}|,|\mathbf{w}|)\leq\mathrm{R}$ and $\mathrm{R}\geq1$; the same constant $\mathrm{C}_{H}$ may be used in both bounds after enlarging it if necessary. There exist $\mathrm{T}_{H,\mathbf{x}}\in\mathscr{T}$ and $\mathrm{C}_{H,\mathbf{x}}>0$ such that, with
\begin{equation}\label{eq:RH-explicit}
    \mathrm{R}_{H,\mathbf{x}}(t)\overset{\textnormal{def}}{=}\mathrm{C}_{H,\mathbf{x}}\left(1+\log(\gamma(t)^{-1})\right)^{-1/4},
\end{equation}
for every $t>\mathrm{T}_{H,\mathbf{x}}$ (in particular, $\gamma(t)<1$ throughout),
\begin{equation}\label{eq:free-U-Gaussian-limit}
    \left|
    \mathbb{E}_{\mathbf{x}}^N\left[
    H\left(\mathbf{U}(t)\right)\mathbf{1}_{\{\mathbf{X}(t)\in\mathscr{C}_{t}^{N}\}}
    \right]
    -
    \frac{1}{(2\pi)^{N/2}}\int_{\mathbb{R}^{N}} H(\mathbf{v})\,\exp\Big(-\tfrac12|\mathbf{v}|^2\Big)\,\mathrm{d}\mathbf{v}
    \right|
    \leq \mathrm{R}_{H,\mathbf{x}}(t).
\end{equation}
\end{proposition}

\begin{proof}
The strategy is to reduce the multivariate test-function bound to a one-particle characteristic-function estimate on a compact frequency window, identify the relevant exponent exactly on the endpoint event $\mathcal{A}_{t}=\{\mathscr{X}(t)\in\mathscr{C}_{t}\}$, control the residual error in $L^{2}$ via the moment input from Proposition~\ref{boundedmoments}, and finally pass from the characteristic-function to the test-function by first cutting off the polynomially growing observable and then applying a compact-frequency smoothing kernel with finite first moment. 

\textbf{Step 1: Reduction to a one-particle characteristic-function estimate.} Fix $\mathrm{R}>0$ and $\mathbf{u}=(u_1,\ldots,u_N)\in[-\mathrm{R},\mathrm{R}]^N$. By independence of the coordinates, and because $\mathscr{U}_{\mathscr{X}_i}(t)$ depends only on $\mathscr{X}_i(t)$, we have
\begin{equation}
\mathbb{E}_{\mathbf{x}}^N\left[\exp\big(\mathrm{i}\langle\mathbf{u},\mathbf{U}(t)\rangle\big)\right]=\prod_{i=1}^{N}\mathbb{E}_{x_i}\left[\exp\big(\mathrm{i}u_i\,\mathscr{U}(t)\big)\right],
\end{equation}
where, in the $i$-th factor, $\mathscr{U}(t)$ denotes the globally defined one-dimensional random variable $\mathscr{U}_{\mathscr{X}}(t)$ for a one-particle process started from $x_i$. For $t$ beyond the deterministic lower time used below, this equals $-\Sigma(t)\,\mathrm{u}_{\mathscr{X}}(t)$ on the endpoint event by Proposition~\ref{axpro}. We prove that, for every fixed starting point $x\in\mathbb{Z}_+$, there are constants $\mathrm{T}_{\mathrm{cf}}(\mathrm{R},x)\in\mathscr{T}$ and $\mathrm{M}_{\mathrm{cf}}(\mathrm{R},x)>0$ such that
\begin{equation}\label{eq:CFgoal}
\sup_{|u|\leq \mathrm{R}}\left|\mathbb{E}_{x}\left[\exp\big(\mathrm{i}u\,\mathscr{U}(t)\big)\right]-\exp\left(-\frac{u^2}{2}\right)\right|
\leq
\mathrm{M}_{\mathrm{cf}}(\mathrm{R},x)\gamma(t)^{1/2},
\qquad t>\mathrm{T}_{\mathrm{cf}}(\mathrm{R},x).
\end{equation}

\textbf{Step 2: Eigenfunction identity at the imaginary spectral value.} Fix $u$ with $|u|\leq\mathrm{R}$ and define $w(u,t)=-\frac{\mathrm{i} u}{\Sigma(t)}$. Since $w(u,t)$ is purely imaginary, $|w(u,t)-\mathfrak{L}|^{2}=\mathfrak{L}^{2}+u^{2}/\Sigma(t)^{2}\leq\mathfrak{L}^{2}+\mathrm{R}^{2}/\Sigma(t)^{2}$, and $\Sigma(t)\to\infty$ by \textbf{(iv)} of Proposition~\ref{specificcasesub}. Choose $\mathrm{T}_{\mathrm{disc}}\in\mathscr{T}$ so that $\mathrm{R}^{2}/\Sigma(t)^{2}<2\mathfrak{L}/\mathfrak{U}+\mathfrak{U}^{-2}$ for every $t>\mathrm{T}_{\mathrm{disc}}$; then $w(u,t)\in\mathcal{B}_{\mathfrak{L}+\mathfrak{U}^{-1}}(\mathfrak{L})$ uniformly in $|u|\leq\mathrm{R}$. Hence Proposition~\ref{eigenfunction} applies in both time regimes and yields, for every $t\in\mathscr{T}$ with $t>\mathrm{T}_{\mathrm{disc}}$,
\begin{multline}\label{eq:eigenID}
\mathbb{E}_x\left[\rho_{\mathscr{X}(t)}(w(u,t))\right]
=
\rho_x(w(u,t))\\
\times
\exp\Bigg(
    -\mathbf{1}_{\{\mathscr{T}=\mathbb{R}_+\}}t\,w(u,t)
    +
    \mathbf{1}_{\{\mathscr{T}=\mathbb{Z}_+\}}
    \sum_{n=0}^{t-1}\Big(
        (1-\tau_n)\log\left(1-\alpha_n w(u,t)\right)
        -
        \tau_n\log\left(1+\beta_n w(u,t)\right)
    \Big)
\Bigg).
\end{multline}

\textbf{Step 3: Exact exponent identity on the endpoint event.} Let
\begin{equation*}
    \mathcal{A}_{t}\overset{\textnormal{def}}{=}\{\mathscr{X}(t)\in\mathscr{C}_{t}\}.
\end{equation*}
By Proposition~\ref{localizationtail}, there exists $\mathrm{T}_{\mathrm{L}}\in\mathscr{T}$ and a constant $\mathrm{M}_x>0$ such that, for all $t>\mathrm{T}_{\mathrm{L}}$,
\begin{equation*}
    \mathbb{P}_{x}(\mathcal{A}_{t}^{\mathrm c})
    \leq
    \mathrm{M}_x\mathrm{e}^{-\mathrm{c}\Sigma(t)}.
\end{equation*}
The contribution of $\mathcal{A}_{t}^{\mathrm c}$ is deferred to Step~5, where it is controlled by Cauchy--Schwarz together with the square-integrability estimate \eqref{eq:Z-L2-bound}; no separate pointwise bound on $\rho_{\mathscr{X}(t)}(w(u,t))$ is needed at this stage. All endpoint-window estimates below are applied only on the event $\mathcal{A}_{t}$, and their lower times are deterministic, because they come from the endpoint versions of Proposition~\ref{specificcasesub} and from the endpoint-uniform construction in Proposition~\ref{axpro}. Using the explicit form $\rho_x(w)=\prod_{k=0}^{x-1}(1-w/\lambda_k)$ from \eqref{characteristic} together with the absolutely convergent series $\log(1+z)=\sum_{m\geq1}(-1)^{m+1}z^m/m$ on $|z|<1$, applied with $z=\mathrm{i} u/(\lambda_k\Sigma(t))$, we obtain
\begin{equation}\label{eq:rhoexp}
\rho_{\mathscr{X}(t)}\left(-\frac{\mathrm{i} u}{\Sigma(t)}\right)=\exp\left(
    \frac{\mathrm{i} u}{\Sigma(t)}\sum_{n=0}^{\mathscr{X}(t)-1}\frac{1}{\lambda_n}+
    \frac{u^2}{2\Sigma(t)^2}\sum_{n=0}^{\mathscr{X}(t)-1}\frac{1}{\lambda_n^2}
    +
    \Xi_1(u,t)
\right),
\end{equation}
where
\begin{equation}
\Xi_1(u,t)
=
-\sum_{m=3}^{\infty}\frac{(-\mathrm{i} u)^m}{m\,\Sigma(t)^m}
\sum_{n=0}^{\mathscr{X}(t)-1}\frac{1}{\lambda_n^m}.
\end{equation}
		By \textbf{(A1)}, $\lambda_n>\mathfrak{L}$ for every $n\in\mathbb{Z}_{+}$. Therefore, whenever $|u|\leq \mathfrak{L}\Sigma(t)/2$, the geometric-series bound applied to the tail in the preceding display gives
\begin{equation}\label{eq:Xi1-bound}
    \left|\Xi_1(u,t)\right|
    \leq
    2\mathfrak{L}^{-3}\frac{|u|^3\mathscr{X}(t)}{\Sigma(t)^3}.
\end{equation}
	Indeed, after factoring out the $m=3$ term, each remaining summand is bounded by the corresponding power of $|u|/(\mathfrak{L}\Sigma(t))\leq1/2$. On $\mathcal{A}_{t}$, the endpoint condition gives $|\mathscr{X}(t)-\mathscr{Z}(t)|\leq\gamma(t)^{1/2}\Sigma(t)^2$. Choose the deterministic lower time so that $\gamma(t)^{1/2}\leq\mathrm{C}_1$ and the lower comparison $\mathscr{Z}(t)\geq \mathrm{C}_1\Sigma(t)^2$ from \textbf{(iv)} of Proposition~\ref{specificcasesub} is valid. Then $\mathscr{X}(t)\leq2\mathscr{Z}(t)$ on $\mathcal{A}_{t}$. The upper comparison $\mathscr{Z}(t)\leq \mathrm{C}_2\Sigma(t)^2$ from the same proposition therefore gives, uniformly for $|u|\leq\mathrm{R}$ and every $t$ beyond this deterministic lower time,
\begin{equation*}
    \left|\Xi_1(u,t)\right|\mathbf{1}_{\mathcal{A}_{t}}
    \leq 4\mathfrak{L}^{-3}\mathrm{C}_2\frac{|u|^3}{\Sigma(t)},
\end{equation*}
which is bounded by a constant depending only on $\mathrm{R}$ times $\Sigma(t)^{-1}$. Next, define
\begin{equation*}
\mathcal{G}_{\mathrm{time}}(u,t)\overset{\textnormal{def}}{=}\exp\left(
\begin{aligned}
&\mathbf{1}_{\{\mathscr{T}=\mathbb{R}_+\}}t\,w(u,t)\\
&\quad-
\mathbf{1}_{\{\mathscr{T}=\mathbb{Z}_+\}}
\sum_{n=0}^{t-1}\left((1-\tau_n)\log\left(1-\alpha_n w(u,t)\right)-\tau_n\log\left(1+\beta_n w(u,t)\right)\right)
\end{aligned}
\right),
\end{equation*}
so that \eqref{eq:eigenID} can be rewritten as
\begin{equation}\label{eq:eigenID-rearranged}
\mathbb{E}_{x}\left[\rho_{\mathscr{X}(t)}\left(w(u,t)\right)\mathcal{G}_{\mathrm{time}}(u,t)\right]=\rho_x\left(w(u,t)\right).
\end{equation}
Expanding $\mathcal{G}_{\mathrm{time}}(u,t)$ at $w=0$ by the same absolutely convergent logarithmic series yields
\begin{equation}\label{eq:TimeExp}
\mathcal{G}_{\mathrm{time}}(u,t)
=
\exp\left(
\begin{aligned}
    &-\mathbf{1}_{\{\mathscr{T}=\mathbb{R}_+\}}\frac{\mathrm{i} u}{\Sigma(t)}\,t\\
    &\quad-
    \mathbf{1}_{\{\mathscr{T}=\mathbb{Z}_+\}}\frac{\mathrm{i} u}{\Sigma(t)}
    \sum_{n=0}^{t-1}\big((1-\tau_n)\alpha_n+\tau_n\beta_n\big)\\
    &\quad-
    \mathbf{1}_{\{\mathscr{T}=\mathbb{Z}_+\}}\frac{u^2}{2\Sigma(t)^2}
    \sum_{n=0}^{t-1}\big((1-\tau_n)\alpha_n^2-\tau_n\beta_n^2\big)
    +
    \Xi_2(u,t)
\end{aligned}
\right),
\end{equation}
where
\begin{equation}
\Xi_2(u,t)
=
\mathbf{1}_{\{\mathscr{T}=\mathbb{Z}_+\}}
\sum_{m=3}^{\infty}\frac{1}{m\,\Sigma(t)^m}
\sum_{n=0}^{t-1}\Big((1-\tau_n)(-\mathrm{i} u\alpha_n)^m-\tau_n(\mathrm{i} u\beta_n)^m\Big).
\end{equation}
By the temporal bounds in \textbf{(A1)}, $\alpha_n,\beta_n\leq\mathfrak{U}$. Hence, whenever $|u|\leq\Sigma(t)/(2\mathfrak{U})$, the same geometric-series estimate gives
\begin{equation}\label{eq:Xi2-bound}
    \left|\Xi_2(u,t)\right|
    \leq
    2\mathfrak{U}^{3}\mathbf{1}_{\{\mathscr{T}=\mathbb{Z}_+\}}\frac{|u|^3t}{\Sigma(t)^3}.
\end{equation}
As in the preceding moment estimate, \textbf{(iv)} of Proposition~\ref{specificcasesub} together with \textbf{(A1)} gives constants $\mathrm{c}_{\Sigma}>0$ and $\mathrm{T}_{\Sigma}\in\mathscr{T}$ such that $\Sigma(t)^2\geq \mathrm{c}_{\Sigma}(1+t)$ for all $t>\mathrm{T}_{\Sigma}$. Choose $\mathrm{T}_{\mathrm{R},2}\in\mathscr{T}$ such that $\mathrm{R}\leq\Sigma(t)/(2\mathfrak{U})$ for all $t>\mathrm{T}_{\mathrm{R},2}$. Therefore, for every $t>\max(\mathrm{T}_{\Sigma},\mathrm{T}_{\mathrm{R},2})$ and every $|u|\leq\mathrm{R}$,
\begin{equation*}
    \left|\Xi_2(u,t)\right|\leq \frac{2\mathfrak{U}^{3}|u|^3}{\mathrm{c}_{\Sigma}\Sigma(t)},
\end{equation*}
which is bounded by a constant depending only on $\mathrm{R}$ times $\Sigma(t)^{-1}$.

\textbf{Step 4: Identification of $\mathscr{U}(t)$ and the explicit error bound.}  Throughout this step, $\mathfrak{F}_{\mathscr{X}}^{(j)}(0,t)$ denotes the $j$th phase derivative $\mathfrak{F}_{y}^{(j)}(0,t)$ of Definition~\ref{form} evaluated along the random terminal value $y=\mathscr{X}$. By Definition~\ref{form} we have
\begin{equation}
\mathfrak{F}_{y}^{(1)}(0,t)=-\mathbf{1}_{\{\mathscr{T}=\mathbb{R}_+\}}\,t+\sum_{n=0}^{y(t)}\frac{1}{\lambda_n}-\mathbf{1}_{\{\mathscr{T}=\mathbb{Z}_+\}}\sum_{n=0}^{t-1}\big((1-\tau_n)\alpha_n+\tau_n\beta_n\big),
\end{equation}
and
\begin{equation}
\mathfrak{F}_{y}^{(2)}(0,t)=\sum_{n=0}^{y(t)}\frac{1}{\lambda_n^2}+\mathbf{1}_{\{\mathscr{T}=\mathbb{Z}_+\}}\sum_{n=0}^{t-1}\big(\tau_n\beta_n^2-(1-\tau_n)\alpha_n^2\big).
\end{equation}
On $\mathcal{A}_{t}$, the terminal endpoint belongs to $\mathscr{C}_{t}$. Therefore the endpoint estimate in \textbf{(iv)} of Proposition~\ref{specificcasesub} gives the concrete comparison
\begin{equation}\label{eq:F2-random-window}
    \left|\mathfrak{F}_{\mathscr{X}}^{(2)}(0,t)-\Sigma(t)^2\right|
    \leq
    \mathrm{C}_{\mathscr{C}}\gamma(t)^{1/2}\Sigma(t)^2
\end{equation}
for every $t>\mathrm{T}_{4,\mathscr{C}}$ on $\mathcal{A}_{t}$. Combining this comparison with the endpoint stationary-point construction from Proposition~\ref{axpro}, we have $\mathscr{U}(t)=-\Sigma(t)\,\mathrm{u}_{\mathscr{X}}(t)$ on $\mathcal{A}_{t}$. Taking $t$ beyond the deterministic maximum of $\mathrm{T}_{\mathrm{A}}$ and $\mathrm{T}_{4,\mathscr{C}}$, \eqref{secondconditiontosat} combined with \textbf{(iv)} of Proposition~\ref{specificcasesub} applied at $w=0$ gives, on $\mathcal{A}_{t}$,
\begin{equation}\label{eq:Uidentify}
    \left|\mathscr{U}(t)-\frac{\mathfrak{F}_{\mathscr{X}}^{(1)}(0,t)}{\Sigma(t)}\right|
    \leq
    \frac{4\mathrm{C}_{\mathrm{A}}\left|\mathfrak{F}_{\mathscr{X}}^{(1)}(0,t)\right|^2}{\Sigma(t)^3}
    +
    \frac{2\mathrm{C}_{\mathscr{C}}\gamma(t)^{1/2}\left|\mathfrak{F}_{\mathscr{X}}^{(1)}(0,t)\right|}{\Sigma(t)}.
\end{equation}
Combining \eqref{eq:rhoexp} and \eqref{eq:TimeExp} gives an exact exponent identity. The only mismatch between the sums in \eqref{eq:rhoexp} and the derivatives in Definition~\ref{form} is the terminal spatial summand. Thus, on $\mathcal{A}_{t}$ and for every $t$ beyond the deterministic lower time fixed in the preceding paragraph,
\begin{equation}\label{eq:explicit-exponent-one}
\rho_{\mathscr{X}(t)}\left(w(u,t)\right)\mathcal{G}_{\mathrm{time}}(u,t)
=
\exp\left(\frac{\mathrm{i} u}{\Sigma(t)}\mathfrak{F}_{\mathscr{X}}^{(1)}(0,t)
+\frac{u^2}{2\Sigma(t)^2}\mathfrak{F}_{\mathscr{X}}^{(2)}(0,t)
+\mathcal{E}_{1,t}(u)\right),
\end{equation}
where the error term is not implicit but is the explicit quantity
\begin{equation}\label{eq:E1-explicit}
\mathcal{E}_{1,t}(u)
\overset{\textnormal{def}}{=}
-\frac{\mathrm{i} u}{\Sigma(t)\lambda_{\mathscr{X}(t)}}
-\frac{u^2}{2\Sigma(t)^2\lambda_{\mathscr{X}(t)}^2}
+\Xi_1(u,t)+\Xi_2(u,t).
\end{equation}
By \textbf{(A1)}, the two terminal summands satisfy
\begin{equation*}
    \frac{|u|}{\Sigma(t)\lambda_{\mathscr{X}(t)}}
    +
    \frac{u^2}{2\Sigma(t)^2\lambda_{\mathscr{X}(t)}^2}
    \leq
    \frac{|u|}{\mathfrak{L}\Sigma(t)}
    +
    \frac{u^2}{2\mathfrak{L}^{2}\Sigma(t)^2}.
\end{equation*}
Together with \eqref{eq:Xi1-bound} and \eqref{eq:Xi2-bound}, this bounds $\mathcal{E}_{1,t}(u)$ explicitly. Replacing $\mathfrak{F}_{\mathscr{X}}^{(1)}(0,t)/\Sigma(t)$ by $\mathscr{U}(t)$ and replacing $\mathfrak{F}_{\mathscr{X}}^{(2)}(0,t)/\Sigma(t)^2$ by $1$ gives
\begin{equation}\label{eq:explicit-exponent-two}
\rho_{\mathscr{X}(t)}\left(w(u,t)\right)\mathcal{G}_{\mathrm{time}}(u,t)
=
\exp\left(\mathrm{i}u\,\mathscr{U}(t)+\frac{u^2}{2}+\mathcal{E}_{2,t}(u)\right),
\end{equation}
where
\begin{equation}\label{eq:E2-explicit}
\mathcal{E}_{2,t}(u)
\overset{\textnormal{def}}{=}
\mathcal{E}_{1,t}(u)
+\mathrm{i}u\left(\frac{\mathfrak{F}_{\mathscr{X}}^{(1)}(0,t)}{\Sigma(t)}-\mathscr{U}(t)\right)
+\frac{u^2}{2}\left(\frac{\mathfrak{F}_{\mathscr{X}}^{(2)}(0,t)}{\Sigma(t)^2}-1\right).
\end{equation}
The preceding estimates give the following completely explicit bound on the endpoint-window event:
\begin{align}\label{eq:E2-bound}
|\mathcal{E}_{2,t}(u)|
&\leq
\frac{|u|}{\mathfrak{L}\Sigma(t)}
+\frac{u^2}{2\mathfrak{L}^{2}\Sigma(t)^2}
+4\mathfrak{L}^{-3}\mathrm{C}_2\frac{|u|^3}{\Sigma(t)}
+\frac{2\mathfrak{U}^{3}|u|^3}{\mathrm{c}_{\Sigma}\Sigma(t)}
+\frac{u^2}{2}\mathrm{C}_{\mathscr{C}}\gamma(t)^{1/2}\notag\\
&\quad+
|u|\left(
\frac{4\mathrm{C}_{\mathrm{A}}\left|\mathfrak{F}_{\mathscr{X}}^{(1)}(0,t)\right|^2}{\Sigma(t)^3}
+\frac{2\mathrm{C}_{\mathscr{C}}\gamma(t)^{1/2}\left|\mathfrak{F}_{\mathscr{X}}^{(1)}(0,t)\right|}{\Sigma(t)}
\right).
\end{align}

\textbf{Step 5: $L^2$ bound on $\mathcal{E}_{2,t}$ and on the endpoint complement.} We now check the probabilistic size of the explicit error uniformly for $|u|\leq\mathrm{R}$. Let $\mathrm{C}_{\mathrm{FD},p}$ denote the constant in \eqref{eq:first-derivative-moment}. Squaring \eqref{eq:E2-bound}, using $(a_1+\cdots+a_m)^2\leq m(a_1^2+\cdots+a_m^2)$, and applying \eqref{eq:first-derivative-moment} with $p=2$ and $p=4$, gives a constant $\mathrm{K}_{\mathrm{E}}(\mathrm{R},x)>0$ and a deterministic lower time $\mathrm{T}_{\mathrm{E}}(\mathrm{R},x)\in\mathscr{T}$ such that, for every $t>\mathrm{T}_{\mathrm{E}}(\mathrm{R},x)$,
    \begin{equation}\label{eq:E2-L2-rate}
        \sup_{|u|\leq\mathrm{R}}
        \mathbb{E}_{x}\left[
        |\mathcal{E}_{2,t}(u)|^2\mathbf{1}_{\mathcal{A}_{t}}
        \right]
        \leq
        \mathrm{K}_{\mathrm{E}}(\mathrm{R},x)\gamma(t).
    \end{equation}
    Indeed, the deterministic terms in \eqref{eq:E2-bound} have squares bounded by constants depending only on $\mathrm{R}$ times $\Sigma(t)^{-2}$ or $\gamma(t)$, and $\Sigma(t)^{-2}\leq\gamma(t)$ after increasing the lower time. The term involving $|\mathfrak{F}_{\mathscr{X}}^{(1)}(0,t)|^2/\Sigma(t)^3$ is controlled by
    \begin{equation*}
        \mathbb{E}_x\left[\frac{|\mathfrak{F}_{\mathscr{X}}^{(1)}(0,t)|^4}{\Sigma(t)^6}\right]
        \leq
        \frac{\mathrm{C}_{\mathrm{FD},4}(1+t)^2}{\Sigma(t)^6}
        \leq
        \frac{\mathrm{C}_{\mathrm{FD},4}}{\mathrm{c}_{\Sigma}^2\Sigma(t)^2}
        \leq
        \frac{\mathrm{C}_{\mathrm{FD},4}}{\mathrm{c}_{\Sigma}^2}\gamma(t),
    \end{equation*}
    while the term involving $\gamma(t)^{1/2}|\mathfrak{F}_{\mathscr{X}}^{(1)}(0,t)|/\Sigma(t)$ is controlled by
    \begin{equation*}
        \gamma(t)\frac{\mathrm{C}_{\mathrm{FD},2}(1+t)}{\Sigma(t)^2}
        \leq
        \frac{\mathrm{C}_{\mathrm{FD},2}}{\mathrm{c}_{\Sigma}}\gamma(t).
    \end{equation*}
    The complement estimate at the beginning of the proof and the inequality $\mathrm{e}^{-\mathrm{c}\Sigma(t)/2}\leq\gamma(t)^{1/2}$, imposed by increasing the deterministic lower time, give constants $\mathrm{K}_{\mathrm{L}}(x)>0$ and $\mathrm{T}_{\mathrm{L},\gamma}(x)\in\mathscr{T}$ such that, for every $t>\mathrm{T}_{\mathrm{L},\gamma}(x)$,
    \begin{equation}\label{eq:Atc-rate}
        \mathbb{P}_{x}(\mathcal{A}_{t}^{\mathrm c})^{1/2}\leq \mathrm{K}_{\mathrm{L}}(x)\gamma(t)^{1/2}.
    \end{equation}

\textbf{Step 6: Uniform $L^2$ bound on $\mathrm{Z}_t(u)$ and the one-particle CF estimate.} We now justify the substitution of \eqref{eq:explicit-exponent-two} inside the expectation without assuming an exponential moment for $\mathcal{E}_{2,t}(u)$. Set
    \begin{equation*}
        \mathrm{Z}_{t}(u)\overset{\textnormal{def}}{=}
        \rho_{\mathscr{X}(t)}\left(-\frac{\mathrm{i}u}{\Sigma(t)}\right)\mathcal{G}_{\mathrm{time}}(u,t).
    \end{equation*}
    We first record a uniform square-integrability bound for this random variable. By \textbf{(A1)} and the inequality $\log(1+v)\geq v/2$ for $0\leq v\leq1$, for every fixed $a>0$ there exists $\mathrm{T}_{a}\in\mathscr{T}$ such that
    \begin{equation*}
        \exp\left(a\frac{\mathscr{X}(t)}{\Sigma(t)^2}\right)
        \leq
        \rho_{\mathscr{X}(t)}\left(-\frac{2\mathfrak{B}_{1}a}{\Sigma(t)^2}\right),
        \qquad t>\mathrm{T}_{a}.
    \end{equation*}
    Indeed, choose $\mathrm{T}_{a}$ so that $2\mathfrak{B}_{1}a/(\mathfrak{L}\Sigma(t)^2)\leq1$ for every $t>\mathrm{T}_{a}$. Then each factor in $\rho_{\mathscr{X}(t)}(-2\mathfrak{B}_{1}a/\Sigma(t)^2)$ is at least $\exp(a/\Sigma(t)^2)$, because $\lambda_n\leq\mathfrak{B}_1$ and $\lambda_n\geq\mathfrak{L}$. In the discrete-time regime we also require $\mathrm{T}_{a}$ to make the real spectral value $-2\mathfrak{B}_{1}a/\Sigma(t)^2$ lie in the ball allowed in Proposition~\ref{eigenfunction}. Applying Proposition~\ref{eigenfunction} at this real value, and using the deterministic comparison
    \begin{equation*}
        t\leq \mathrm{C}_{t,\Sigma}\Sigma(t)^2
    \end{equation*}
    for a constant $\mathrm{C}_{t,\Sigma}>0$ supplied by \textbf{(iv)} of Proposition~\ref{specificcasesub} together with \textbf{(A1)}, gives
    \begin{equation}\label{eq:X-exp-moment}
        \sup_{t>\mathrm{T}_{a}}
        \mathbb{E}_{x}\left[\exp\left(a\frac{\mathscr{X}(t)}{\Sigma(t)^2}\right)\right]
        <\infty.
    \end{equation}
    The last implication is uniform for $t>\mathrm{T}_{a}$: the fixed polynomial $\rho_x$ is bounded on the compact set of spectral values reached after this threshold. In continuous time the temporal factor is bounded by
    \begin{equation*}
        \exp(2\mathfrak{B}_1a\,\mathrm{C}_{t,\Sigma}).
    \end{equation*}
    In discrete time, with $\mathrm{T}_{a}$ also chosen so that $\mathfrak{U}(2\mathfrak{B}_1a)/\Sigma(t)^2\leq1/2$ for every $t>\mathrm{T}_{a}$, the inequalities $\log(1+v)\leq v$ and $-\log(1-v)\leq2v$ give the bound
    \begin{equation*}
        \exp(4\mathfrak{U}\mathfrak{B}_1a\,\mathrm{C}_{t,\Sigma}).
    \end{equation*}
    On the other hand, for the purely imaginary spectral value $-\mathrm{i}u/\Sigma(t)$, the elementary inequality $\log(1+v)\leq v$ for $v\geq0$ gives
    \begin{equation*}
        \left|\rho_{\mathscr{X}(t)}\left(-\frac{\mathrm{i}u}{\Sigma(t)}\right)\right|^2
        \leq
        \exp\left(\frac{u^2}{\mathfrak{L}^2}\frac{\mathscr{X}(t)}{\Sigma(t)^2}\right),
    \end{equation*}
    and the temporal factor satisfies $|\mathcal{G}_{\mathrm{time}}(u,t)|^2\leq \exp(\mathrm{C}_{\mathrm{R}})$ for a finite constant $\mathrm{C}_{\mathrm{R}}$ depending only on $\mathrm{R}$ and the constants in \textbf{(A1)}. In continuous time this modulus is equal to one. In discrete time, the bound follows by applying $\mathfrak{Re}\log(1+\mathrm{i}v)=\frac12\log(1+v^2)\leq v^2/2$ to the factors with $v=\alpha_nu/\Sigma(t)$ and $v=\beta_nu/\Sigma(t)$, followed by $t\leq \mathrm{C}_{t,\Sigma}\Sigma(t)^2$. Combining this deterministic temporal bound with \eqref{eq:X-exp-moment}, for $a=\mathrm{R}^2/\mathfrak{L}^2$, yields
    \begin{equation}\label{eq:Z-L2-bound}
        \sup_{\substack{t>\mathrm{T}_{\mathrm{R}}\\ |u|\leq\mathrm{R}}}\mathbb{E}_{x}\left[|\mathrm{Z}_{t}(u)|^2\right]<\infty
    \end{equation}
    for a threshold $\mathrm{T}_{\mathrm{R}}\in\mathscr{T}$.
    \medskip

    \noindent On $\mathcal{A}_{t}$, equation \eqref{eq:explicit-exponent-two} gives
    \begin{equation*}
        \mathrm{Z}_{t}(u)
        =
        \exp\left(\frac{u^2}{2}\right)
        \exp\big(\mathrm{i}u\,\mathscr{U}(t)\big)
        \exp(\mathcal{E}_{2,t}(u)).
    \end{equation*}
    Let $\mathrm{K}_{\mathrm{Z}}(\mathrm{R},x)$ be a finite constant larger than the square-root of the supremum in \eqref{eq:Z-L2-bound}. On $\mathcal{A}_{t}\cap\{|\mathcal{E}_{2,t}(u)|\leq1\}$, the inequality $|\mathrm{e}^{z}-1|\leq \mathrm{e}|z|$ gives
    \begin{equation*}
        \left|\mathrm{Z}_{t}(u)-\exp\left(\frac{u^2}{2}\right)\exp\big(\mathrm{i}u\,\mathscr{U}(t)\big)\right|
        \leq
        \mathrm{e}^{1+\mathrm{R}^2/2}|\mathcal{E}_{2,t}(u)|.
    \end{equation*}
    On $\mathcal{A}_{t}\cap\{|\mathcal{E}_{2,t}(u)|>1\}$ we use Cauchy's inequality, \eqref{eq:Z-L2-bound}, and \eqref{eq:E2-L2-rate}. Hence \eqref{eq:E2-L2-rate}, \eqref{eq:Atc-rate}, and \eqref{eq:Z-L2-bound} give a constant $\mathrm{K}_{\rho}(\mathrm{R},x)>0$ and a deterministic lower time $\mathrm{T}_{\rho}(\mathrm{R},x)\in\mathscr{T}$ such that, for every $t>\mathrm{T}_{\rho}(\mathrm{R},x)$,
	\begin{equation}\label{eq:CF-rho-rate}
    \sup_{|u|\leq\mathrm{R}}
	\left|
	\rho_x\left(-\frac{\mathrm{i} u}{\Sigma(t)}\right)
	-
	\exp\left(\frac{u^2}{2}\right)\mathbb{E}_{x}\left[\exp\big(\mathrm{i}u\,\mathscr{U}(t)\big)\right]
	\right|
	\leq
    \mathrm{K}_{\rho}(\mathrm{R},x)\gamma(t)^{1/2}.
	\end{equation}
    Since $\rho_x$ is a fixed polynomial with $\rho_x(0)=1$, there are constants $\mathrm{K}_{x,\mathrm{R}}>0$ and $\mathrm{T}_{x,\mathrm{R}}\in\mathscr{T}$ such that, for every $t>\mathrm{T}_{x,\mathrm{R}}$,
    \begin{equation*}
        \sup_{|u|\leq\mathrm{R}}\left|\rho_x\left(-\frac{\mathrm{i} u}{\Sigma(t)}\right)-1\right|
        \leq
        \mathrm{K}_{x,\mathrm{R}}\Sigma(t)^{-1}
        \leq
        \mathrm{K}_{x,\mathrm{R}}\gamma(t)^{1/2}.
    \end{equation*}
    Combining this estimate with \eqref{eq:CF-rho-rate}, and taking $\mathrm{T}_{\mathrm{cf}}(\mathrm{R},x)$ to dominate the lower times in \eqref{eq:CF-rho-rate} and in the last displayed inequality, proves \eqref{eq:CFgoal}.

\textbf{Step 7: From characteristic functions to general test functions.} Consequently, there exist constants $\mathrm{T}_{\mathrm{cf}}(\mathrm{R},\mathbf{x})\in\mathscr{T}$ and $\mathrm{M}_{\mathrm{cf}}(\mathrm{R},\mathbf{x})>0$ such that, for every $t>\mathrm{T}_{\mathrm{cf}}(\mathrm{R},\mathbf{x})$ and every $\mathbf{u}\in[-\mathrm{R},\mathrm{R}]^N$,
\begin{equation}\label{eq:multi-CF-rate}
\left|
\mathbb{E}_{\mathbf{x}}^N\left[\exp\big(\mathrm{i}\langle\mathbf{u},\mathbf{U}(t)\rangle\big)\right]
-
\exp\left(-\frac12|\mathbf{u}|^2\right)
\right|
\leq
\mathrm{M}_{\mathrm{cf}}(\mathrm{R},\mathbf{x})\gamma(t)^{1/2}.
\end{equation}
This follows from the product identity at the beginning of the proof and the elementary telescoping bound for finite products whose factors are bounded by one.
Moreover, following the constants in the estimates above gives the explicit growth bound
\begin{equation}\label{eq:cf-constant-growth}
    \mathrm{M}_{\mathrm{cf}}(\mathrm{R},\mathbf{x})
    \leq
    \mathrm{C}_{\mathrm{cf},\mathbf{x}}\exp\{\mathrm{C}_{\mathrm{cf},\mathbf{x}}\mathrm{R}^{2}\}(1+\mathrm{R})^{\mathrm{C}_{\mathrm{cf},\mathbf{x}}},
    \qquad \mathrm{R}\geq1.
\end{equation}
The exponential factor comes from the square-integrability estimate \eqref{eq:Z-L2-bound}, where the auxiliary real spectral value is proportional to $\mathrm{R}^{2}/\Sigma(t)^2$; all remaining factors in \eqref{eq:E2-bound} are polynomial in $\mathrm{R}$.

 We now pass from characteristic functions to the observable $H$. Since $H$ is allowed to grow polynomially, we first localize it in physical space and only then apply Fourier smoothing. This order is needed because the Fourier inversion step below is applied to a compactly supported function, while the finite first moment of the smoothing kernel converts the Lipschitz regularity of the cutoff observable into an explicit $\mathrm{R}^{-1}$ smoothing error.

 Choose
\begin{equation*}
    \mathrm{R}_{\mathrm{F}}(t)\overset{\textnormal{def}}{=}\left(\frac{\log(\gamma(t)^{-1})}{8\mathrm{C}_{\mathrm{cf},\mathbf{x}}}\right)^{1/2},\qquad \mathrm{A}(t)\overset{\textnormal{def}}{=}\left(1+\log(\gamma(t)^{-1})\right)^{1/(4(q+1))},
\end{equation*}
after increasing the lower time so that $\gamma(t)<\mathrm{e}^{-8\mathrm{C}_{\mathrm{cf},\mathbf{x}}}$ and $\mathrm{R}_{\mathrm{F}}(t)\geq1$. Let $\chi_{\mathrm{A}}$ be a Lipschitz cutoff which equals one on $|\mathbf{v}|\leq \mathrm{A}$, vanishes on $|\mathbf{v}|\geq2\mathrm{A}$, takes values in $[0,1]$, and has Lipschitz constant at most $2/\mathrm{A}$. Set $H_{\mathrm{A}}=H\chi_{\mathrm{A}}$. The growth and Lipschitz assumptions on $H$ imply, for $\mathrm{A}\geq1$,
\begin{equation}\label{eq:HA-bounds}
    \|H_{\mathrm{A}}\|_{\infty}\leq \mathrm{C}_{N,q}\mathrm{C}_{H}\mathrm{A}^{q},\qquad \operatorname{Lip}(H_{\mathrm{A}})\leq \mathrm{C}_{N,q}\mathrm{C}_{H}\mathrm{A}^{q+1},\qquad \operatorname{supp}(H_{\mathrm{A}})\subseteq\{\mathbf{v}:|\mathbf{v}|\leq2\mathrm{A}\}.
\end{equation}
Indeed, the supremum bound follows from $|H(\mathbf{v})|\leq\mathrm{C}_{H}(1+|\mathbf{v}|^{q})$ on the support of $\chi_{\mathrm{A}}$. For the Lipschitz bound, the product rule gives two contributions: the first is bounded by $\operatorname{Lip}(H)$ on $|\mathbf{v}|\leq2\mathrm{A}$, and the second is bounded by $\|H\|_{L^{\infty}(|\mathbf{v}|\leq2\mathrm{A})}\operatorname{Lip}(\chi_{\mathrm{A}})$. Both are bounded by the middle expression in \eqref{eq:HA-bounds}, after absorbing constants depending only on $N$ and $q$.

 We next record the smoothing estimate in the precise form used here. Let $\kappa:\mathbb{R}\to\mathbb{R}_{+}$ be the normalized fourth-power sinc kernel
\begin{equation*}
    \kappa(r)\overset{\textnormal{def}}{=}\frac{3}{8\pi}\left(\frac{\sin(r/4)}{r/4}\right)^4,
\end{equation*}
with the value at $r=0$ defined by continuity. With the Fourier convention $\widehat g(\xi)=\int_{\mathbb{R}}\mathrm{e}^{-ir\xi}g(r)\,\mathrm{d}r$, we have $\int_{\mathbb{R}}\kappa(r)\,\mathrm{d}r=1$, the Fourier transform of $\kappa$ is supported in $[-1,1]$, and $\int_{\mathbb{R}}|r|\kappa(r)\,\mathrm{d}r<\infty$, because $\kappa(r)$ decays like $|r|^{-4}$. For $\mathrm{R}\geq1$ define $\kappa_{\mathrm{R}}(r)=\mathrm{R}\kappa(\mathrm{R}r)$ and $\mathcal{K}_{\mathrm{R}}(\mathbf{v})=\prod_{j=1}^{N}\kappa_{\mathrm{R}}(v_j)$. Then $\mathcal{K}_{\mathrm{R}}$ has total mass one, its Fourier transform is supported in $[-\mathrm{R},\mathrm{R}]^{N}$, and
\begin{equation}\label{eq:kernel-first-moment}
    \int_{\mathbb{R}^{N}}|\mathbf{v}|\mathcal{K}_{\mathrm{R}}(\mathbf{v})\,\mathrm{d}\mathbf{v}\leq \mathrm{C}_{N}\mathrm{R}^{-1}.
\end{equation}
The last estimate follows by the triangle inequality $|\mathbf{v}|\leq |v_1|+\cdots+|v_N|$ and the one-dimensional scaling identity $\int_{\mathbb{R}}|r|\kappa_{\mathrm{R}}(r)\,\mathrm{d}r=\mathrm{R}^{-1}\int_{\mathbb{R}}|r|\kappa(r)\,\mathrm{d}r$.

Let $\mu_t$ be the law of $\mathbf{U}(t)$ and let $\nu$ be the standard Gaussian law on $\mathbb{R}^{N}$. If $f$ is bounded, Lipschitz and supported in $\{|\mathbf{v}|\leq \mathrm{B}\}$, then \eqref{eq:kernel-first-moment} gives
\begin{equation}\label{eq:smoothing-approximation}
    \left|\int f\,\mathrm{d}\mu_t-\int f*\mathcal{K}_{\mathrm{R}}\,\mathrm{d}\mu_t\right|+\left|\int f\,\mathrm{d}\nu-\int f*\mathcal{K}_{\mathrm{R}}\,\mathrm{d}\nu\right|\leq \mathrm{C}_{N}\operatorname{Lip}(f)\mathrm{R}^{-1}.
\end{equation}
For the smoothed difference, Fourier inversion applies because $f$ is compactly supported and $\widehat{\mathcal{K}_{\mathrm{R}}}$ is compactly supported. Writing $\phi_t$ and $\phi_{\nu}$ for the characteristic functions of $\mu_t$ and $\nu$, respectively, we obtain
\begin{equation}\label{eq:smoothing-fourier}
    \left|\int f*\mathcal{K}_{\mathrm{R}}\,\mathrm{d}\mu_t-\int f*\mathcal{K}_{\mathrm{R}}\,\mathrm{d}\nu\right|\leq \mathrm{C}_{N}(1+\mathrm{B})^{N}\|f\|_{\infty}\mathrm{R}^{N}\sup_{\xi\in[-\mathrm{R},\mathrm{R}]^{N}}|\phi_t(\xi)-\phi_{\nu}(\xi)|.
\end{equation}
Indeed, $|\widehat f(\xi)|\leq\|f\|_{1}\leq \mathrm{C}_{N}(1+\mathrm{B})^{N}\|f\|_{\infty}$, $\widehat{\mathcal{K}_{\mathrm{R}}}$ is bounded by a constant depending only on $N$, and the integration region has volume $(2\mathrm{R})^{N}$. Combining \eqref{eq:smoothing-approximation} and \eqref{eq:smoothing-fourier}, then using \eqref{eq:multi-CF-rate}, gives the compact-support smoothing bound
\begin{align}\label{eq:smoothing-BL}
    &\left|\mathbb{E}_{\mathbf{x}}^N[f(\mathbf{U}(t))]-\frac{1}{(2\pi)^{N/2}\!}\int_{\mathbb{R}^{N}}f(\mathbf{v})\exp\Big(-\tfrac12|\mathbf{v}|^2\Big)\,\mathrm{d}\mathbf{v}\right|\\
    &\qquad\leq \mathrm{C}_{\mathrm{sm},N}\left(\operatorname{Lip}(f)\mathrm{R}^{-1}+(1+\mathrm{B})^{N}\|f\|_{\infty}\mathrm{R}^{N}\mathrm{M}_{\mathrm{cf}}(\mathrm{R},\mathbf{x})\gamma(t)^{1/2}\right),\notag
\end{align}
for every $\mathrm{R}\geq1$, every compactly supported bounded Lipschitz $f$ with support contained in $\{|\mathbf{v}|\leq\mathrm{B}\}$, and every $t>\mathrm{T}_{\mathrm{cf}}(\mathrm{R},\mathbf{x})$.

Although \eqref{eq:multi-CF-rate} was first proved for fixed $\mathrm{R}$, the preceding estimates also apply with $\mathrm{R}=\mathrm{R}_{\mathrm{F}}(t)$ after increasing the deterministic lower time once more. Indeed, Definition~\ref{gammadef} gives $\gamma(t)\geq\Sigma(t)^{-1/2}$, and hence, after increasing the deterministic lower time, $\mathrm{R}_{\mathrm{F}}(t)\leq \mathrm{C}_{\mathbf{x}}\sqrt{1+\log\Sigma(t)}$. This logarithmic growth is dominated by every positive power of $\Sigma(t)$: for each fixed $a\geq0$, $b>0$ and $c>0$ appearing in Steps~2--6, we increase the deterministic lower time so that $\mathrm{R}_{\mathrm{F}}(t)^a\Sigma(t)^{-b}\leq c$ for every later $t$.

 Applying \eqref{eq:smoothing-BL} to $f=H_{\mathrm{A}(t)}$, with $\mathrm{B}=2\mathrm{A}(t)$ and $\mathrm{R}=\mathrm{R}_{\mathrm{F}}(t)$, gives
\begin{align}\label{eq:truncated-H-raw}
    &\left|\mathbb{E}_{\mathbf{x}}^N[H_{\mathrm{A}(t)}(\mathbf{U}(t))]-\frac{1}{(2\pi)^{N/2}\!}\int_{\mathbb{R}^{N}}H_{\mathrm{A}(t)}(\mathbf{v})\exp\Big(-\tfrac12|\mathbf{v}|^2\Big)\,\mathrm{d}\mathbf{v}\right|\\
    &\qquad\leq \mathrm{C}_{q,\mathbf{x}}\mathrm{C}_{H}\left(\mathrm{A}(t)^{q+1}\mathrm{R}_{\mathrm{F}}(t)^{-1}+(1+\mathrm{A}(t))^{N}\mathrm{A}(t)^q\mathrm{R}_{\mathrm{F}}(t)^N\mathrm{M}_{\mathrm{cf}}(\mathrm{R}_{\mathrm{F}}(t),\mathbf{x})\gamma(t)^{1/2}\right).\notag
\end{align}
The first term on the right-hand side of \eqref{eq:truncated-H-raw} has the required logarithmic decay. Indeed, by the definitions of $\mathrm{A}(t)$ and $\mathrm{R}_{\mathrm{F}}(t)$, after increasing the lower time if necessary,
\begin{equation*}
    \mathrm{A}(t)^{q+1}\mathrm{R}_{\mathrm{F}}(t)^{-1}\leq \mathrm{C}_{\mathbf{x}}\left(1+\log(\gamma(t)^{-1})\right)^{-1/4}.
\end{equation*}
For the second term, \eqref{eq:cf-constant-growth} and the definition of $\mathrm{R}_{\mathrm{F}}(t)$ give
\begin{equation*}
    \gamma(t)^{1/2}\exp\{\mathrm{C}_{\mathrm{cf},\mathbf{x}}\mathrm{R}_{\mathrm{F}}(t)^2\}\leq \gamma(t)^{3/8}.
\end{equation*}
All remaining factors in the second term of \eqref{eq:truncated-H-raw} are powers of $1+\log(\gamma(t)^{-1})$. Hence there is a constant $\mathrm{m}_{N,q,\mathbf{x}}>0$ such that this second term is bounded, up to the constant $\mathrm{C}_{q,\mathbf{x}}\mathrm{C}_{H}$, by $\left(1+\log(\gamma(t)^{-1})\right)^{\mathrm{m}_{N,q,\mathbf{x}}}\gamma(t)^{3/8}$. The elementary exponential-polynomial comparison gives a number $\mathrm{r}_{N,q,\mathbf{x}}>0$ such that $r^{\mathrm{m}_{N,q,\mathbf{x}}}\mathrm{e}^{-3r/8}\leq r^{-1/4}$ for every $r\geq\mathrm{r}_{N,q,\mathbf{x}}$. Increasing the lower time so that $1+\log(\gamma(t)^{-1})\geq\mathrm{r}_{N,q,\mathbf{x}}$, we get
\begin{equation*}
    (1+\mathrm{A}(t))^{N}\mathrm{A}(t)^q\mathrm{R}_{\mathrm{F}}(t)^N(1+\mathrm{R}_{\mathrm{F}}(t))^{\mathrm{C}_{\mathrm{cf},\mathbf{x}}}\gamma(t)^{3/8}\leq \left(1+\log(\gamma(t)^{-1})\right)^{-1/4}.
\end{equation*}
Substituting these two estimates into \eqref{eq:truncated-H-raw} yields
\begin{equation}\label{eq:truncated-H-rate}
    \left|\mathbb{E}_{\mathbf{x}}^N[H_{\mathrm{A}(t)}(\mathbf{U}(t))]-\frac{1}{(2\pi)^{N/2}\!}\int_{\mathbb{R}^{N}}H_{\mathrm{A}(t)}(\mathbf{v})\exp\Big(-\tfrac12|\mathbf{v}|^2\Big)\,\mathrm{d}\mathbf{v}\right|\leq \mathrm{C}_{q,\mathbf{x}}\mathrm{C}_{H}\left(1+\log(\gamma(t)^{-1})\right)^{-1/4}.
\end{equation}

It remains to remove the cutoff and insert the endpoint event. Choose $p=4(q+1)$ in Proposition~\ref{boundedmoments}. The endpoint moment bound \eqref{eq:Up-moment} and the Gaussian $p$-moment give
\begin{equation}\label{eq:H-tail-rate}
    \mathbb{E}_{\mathbf{x}}^N\left[
    |H(\mathbf{U}(t))|\mathbf{1}_{\{|\mathbf{U}(t)|>\mathrm{A}(t)\}}
    \mathbf{1}_{\{\mathbf{X}(t)\in\mathscr{C}_{t}^{N}\}}
    \right]
    +
    \int_{|\mathbf{v}|>\mathrm{A}(t)}|H(\mathbf{v})|\frac{\mathrm{e}^{-|\mathbf{v}|^2/2}}{(2\pi)^{N/2}}\,\mathrm{d}\mathbf{v}
    \leq
    \mathrm{C}_{q,\mathbf{x}}\mathrm{C}_{H}\left(1+\log(\gamma(t)^{-1})\right)^{-1/4}.
\end{equation}
Finally, on the endpoint complement, applying Proposition~\ref{localizationtail} to each coordinate and using the finite union bound $\{\mathbf{X}(t)\notin\mathscr{C}_{t}^{N}\}\subseteq\bigcup_{i=1}^{N}\{\mathscr{X}_i(t)\notin\mathscr{C}_{t}\}$ gives constants $\mathrm{M}_{\mathrm L,\mathbf{x}}>0$ and $\mathrm{c}_{\mathrm L}>0$ such that
\begin{equation*}
    \mathbb{P}_{\mathbf{x}}^N\left[\mathbf{X}(t)\notin\mathscr{C}_{t}^{N}\right]\leq\mathrm{M}_{\mathrm L,\mathbf{x}}\mathrm{e}^{-\mathrm{c}_{\mathrm L}\Sigma(t)}
\end{equation*}
after increasing the lower time. Combining this with the cutoff bound gives
\begin{equation}\label{eq:H-endpoint-complement-rate}
    \mathbb{E}_{\mathbf{x}}^N\left[
    |H_{\mathrm{A}(t)}(\mathbf{U}(t))|\mathbf{1}_{\{\mathbf{X}(t)\notin\mathscr{C}_{t}^{N}\}}
    \right]
    \leq
    \mathrm{C}_{q,\mathbf{x}}\mathrm{C}_{H}\left(1+\log(\gamma(t)^{-1})\right)^{q/(4(q+1))}\mathrm{e}^{-\mathrm{c}_{\mathrm L}\Sigma(t)}.
\end{equation}
The last estimate uses the coordinatewise endpoint-window bound from Proposition~\ref{localizationtail} (combined with the finite union bound above) and the cutoff bound $|H_{\mathrm{A}(t)}|\leq \mathrm{C}_{q}\mathrm{C}_{H}\mathrm{A}(t)^q$. By Definition~\ref{gammadef}, $\gamma(t)\geq\Sigma(t)^{-1/2}$, so $1+\log(\gamma(t)^{-1})\leq 1+\frac{1}{2}\log\Sigma(t)$ after increasing the lower time. The exponential factor $\mathrm{e}^{-\mathrm{c}_{\mathrm L}\Sigma(t)}$ therefore dominates the displayed logarithmic power, and we increase the lower time so that
\begin{equation*}
    \left(1+\log(\gamma(t)^{-1})\right)^{q/(4(q+1))}\mathrm{e}^{-\mathrm{c}_{\mathrm L}\Sigma(t)}
    \leq
    \left(1+\log(\gamma(t)^{-1})\right)^{-1/4}
\end{equation*}
for every later $t$. Writing $\mathcal{D}_t(H)$ for the difference inside the absolute value in \eqref{eq:free-U-Gaussian-limit}, the triangle inequality combined with \eqref{eq:truncated-H-rate}, \eqref{eq:H-tail-rate}, and \eqref{eq:H-endpoint-complement-rate} gives
\begin{equation*}
    \left|\mathcal{D}_t(H)\right|
    \leq
    \underbrace{\eqref{eq:truncated-H-rate}}_{\text{cutoff main term}}
    +\underbrace{\eqref{eq:H-tail-rate}}_{\text{tail $|\mathbf{U}|>\mathrm{A}(t)$}}
    +\underbrace{\eqref{eq:H-endpoint-complement-rate}}_{\text{endpoint complement}}
    \leq
    3\mathrm{C}_{q,\mathbf{x}}\mathrm{C}_{H}\left(1+\log(\gamma(t)^{-1})\right)^{-1/4},
\end{equation*}
which is \eqref{eq:free-U-Gaussian-limit} with the explicit error in \eqref{eq:RH-explicit} after setting $\mathrm{C}_{H,\mathbf{x}}=3\mathrm{C}_{q,\mathbf{x}}\mathrm{C}_{H}$ and $\mathrm{T}_{H,\mathbf{x}}$ equal to the maximum of all deterministic lower times introduced above.
\end{proof}
\noindent Recall the continuous Weyl chamber $\mathbf{W}_{\mathrm c}^{N}$ from Definition~\ref{def:continuous-weyl}.
\begin{proposition}[Endpoint ordering on $\mathscr{C}_{t}$]\label{twobadprobs-fixed}
Assume \textbf{(A1)} and \textbf{(A2)}. There exists $\mathrm{T}_{\mathrm{W}}\in\mathscr{T}$ such that, for all $t>\mathrm{T}_{\mathrm{W}}$,
\begin{equation}\label{eq:twobad-equality}
\left\{\mathbf{U}(t)\in\mathbf{W}_{\mathrm c}^{N},\ \mathscr{X}_i(t)\in\mathscr{C}_{t}\textnormal{ for every }i\in\llbracket N\rrbracket\right\}
=
\left\{\mathbf{X}(t)\in\mathbf{W}_{\mathrm d}^{N},\ \mathscr{X}_i(t)\in\mathscr{C}_{t}\textnormal{ for every }i\in\llbracket N\rrbracket\right\}.
\end{equation}
\end{proposition}

\begin{proof}
On the endpoint window $\mathscr{C}_{t}$, the ordering of the stationary-point coordinates is the same as the ordering of the terminal particle positions. We prove this directly from the stationary-point equation, rather than by comparing two asymptotic expansions. By Definition~\ref{form}, the phase $\mathfrak{F}_{y}(\cdot,t)$ depends on $y$ only through the terminal value $y(t)$, and hence so does the stationary point $\mathrm{u}_{y}(t)$ via Definition~\ref{definitionstationary}.

Choose $r\in(0,\mathrm{R}_2)$ small enough that $\mathrm{C}_{\mathscr{C}}r\leq 1/4$, where $\mathrm{R}_2$ and $\mathrm{C}_{\mathscr{C}}$ are the constants in Proposition~\ref{specificcasesub}. By the calibration in \eqref{axpro-small-a}, the quadratic comparison \eqref{secondconditiontosat}, and the Newton-scale estimate \eqref{axpro-a0-small},
\begin{equation*}
\begin{aligned}
    \mathrm{C}_{\mathrm{A}}\left|\frac{\mathfrak{F}_{y}^{(1)}(0,t)}{\mathfrak{F}_{y}^{(2)}(0,t)}\right| &\leq 1,\quad\text{and}\quad
    |\mathrm{u}_{y}(t)| &\leq 2\left|\frac{\mathfrak{F}_{y}^{(1)}(0,t)}{\mathfrak{F}_{y}^{(2)}(0,t)}\right|\leq 2\mathrm{B}_{\mathrm{A}}\gamma(t)^{1/2},
\end{aligned}
\end{equation*}
uniformly over $y\in\mathscr{P}$ with $y(t)\in\mathscr{C}_{t}$ and $t>\mathrm{T}_{\mathrm{A}}$; here the factor $2$ on the second line uses the calibration on the first line to absorb the quadratic remainder in \eqref{secondconditiontosat}. The threshold $\mathrm{T}_{\mathrm{A}}$ already absorbs the auxiliary times $\mathrm{T}_{2,\mathscr{C}}$, $\mathrm{T}_{\mathscr{C},1}$, and $\mathrm{T}_{4,\mathscr{C}}$ from Proposition~\ref{specificcasesub}. Since $\gamma(t)\to 0$ by \textbf{(A2)} and \textbf{(iv)} of Proposition~\ref{specificcasesub}, there exists $\mathrm{T}_{\mathrm{W},1}\in\mathscr{T}$ such that $2\mathrm{B}_{\mathrm{A}}\gamma(t)^{1/2}<r$, and hence $\mathrm{u}_{y}(t)\in(-r,r)$, for every $t>\mathrm{T}_{\mathrm{W},1}$ and every such $y$. For the second derivative we use the strengthened estimate in \textbf{(iv)} of Proposition~\ref{specificcasesub}: for $y(t)\in\mathscr{C}_{t}$ and $w\in\mathcal{B}_{\mathrm{R}_2}(0)$,
\begin{equation*}
    \left|\mathfrak{F}_{y}^{(2)}(w,t)-\Sigma(t)^2\right|\leq\mathrm{C}_{\mathscr{C}}\left(\gamma(t)^{1/2}+|w|\right)\Sigma(t)^2.
\end{equation*}
Choose $\mathrm{T}_{\mathrm{W},2}\in\mathscr{T}$ so that $\mathrm{C}_{\mathscr{C}}\gamma(t)^{1/2}\leq 1/4$ for every $t>\mathrm{T}_{\mathrm{W},2}$. With our choice of $r$, this gives $\mathrm{C}_{\mathscr{C}}(\gamma(t)^{1/2}+|w|)\leq 1/2$ whenever $|w|\leq r$, and therefore
\begin{equation}\label{eq:local-convexity-window}
    \mathfrak{F}_{y}^{(2)}(w,t)\geq\frac{1}{2}\Sigma(t)^2>0,
    \qquad y(t)\in\mathscr{C}_{t},\quad w\in[-r,r],
\end{equation}
for every $t>\mathrm{T}_{\mathrm{W}}\overset{\textnormal{def}}{=}\max(\mathrm{T}_{\mathrm{A}},\mathrm{T}_{\mathrm{W},1},\mathrm{T}_{\mathrm{W},2})$. Thus, for each path $y\in\mathscr{P}$ with $y(t)\in\mathscr{C}_{t}$, the map $w\mapsto\mathfrak{F}_{y}^{(1)}(w,t)$ is strictly increasing on the interval $[-r,r]$, and its unique zero $\mathrm{u}_{y}(t)$ lies in $(-r,r)$.

Take two paths $y,y'\in\mathscr{P}$ with $y(t),y'(t)\in\mathscr{C}_{t}$ and $y(t)<y'(t)$. For every real $w\in[-r,r]$, the explicit formula for the phase in Definition~\ref{form} gives
\begin{equation}\label{eq:endpoint-monotonicity}
    \mathfrak{F}_{y'}^{(1)}(w,t)-\mathfrak{F}_{y}^{(1)}(w,t)
    =
    \sum_{\ell=y(t)+1}^{y'(t)}\frac{1}{\lambda_{\ell}-w}.
\end{equation}
The right-hand side is strictly positive: for real $w\in[-r,r]$, the choice of $r$ together with $\mathrm{R}_2=\mathrm{R}_{-}/2=\min(\mathfrak{L},\mathfrak{U}^{-1})/2$ from Proposition~\ref{specificcasesub}\,\textbf{(i)} and \textbf{(iii)} gives
\begin{equation*}
    \lambda_{\ell}-w\geq\mathfrak{L}-r>\mathfrak{L}-\mathrm{R}_2\geq\mathfrak{L}/2>0,
\end{equation*}
using $\lambda_{\ell}\geq\mathfrak{L}$ from \textbf{(A1)}. Substituting $w=\mathrm{u}_{y'}(t)$ into \eqref{eq:endpoint-monotonicity} and using the stationary-point identity $\mathfrak{F}_{y'}^{(1)}(\mathrm{u}_{y'}(t),t)=0$ gives
\begin{equation*}
    \mathfrak{F}_{y}^{(1)}(\mathrm{u}_{y'}(t),t)<0.
\end{equation*}
Since $\mathfrak{F}_{y}^{(1)}(\cdot,t)$ is strictly increasing on $[-r,r]$ and vanishes at $\mathrm{u}_{y}(t)$, this strict inequality forces $\mathrm{u}_{y'}(t)<\mathrm{u}_{y}(t)$. Since $\Sigma(t)>0$, multiplying by $-\Sigma(t)$ reverses the inequality, giving
\begin{equation}\label{eq:stationary-order}
    y(t)<y'(t)\qquad\Longrightarrow\qquad \mathscr{U}_{y}(t)<\mathscr{U}_{y'}(t).
\end{equation}
If $y(t)=y'(t)$, then by Definition~\ref{form} the two phases coincide, and hence so do the two stationary points. Combining the strict implication \eqref{eq:stationary-order}, the equality case $y(t)=y'(t)\Rightarrow\mathscr{U}_{y}(t)=\mathscr{U}_{y'}(t)$, and the symmetric case obtained by swapping $y$ and $y'$, we conclude that the strict order of $y(t),y'(t)\in\mathscr{C}_{t}$ is equivalent to that of $\mathscr{U}_{y}(t),\mathscr{U}_{y'}(t)$. Applying this pairwise to $\mathscr{X}_1(t),\ldots,\mathscr{X}_N(t)$ proves the displayed equality of events.
\end{proof}

    \subsection{Asymptotic expansion of the Karlin--McGregor semigroup}

 We first recall two lemmas. The following is implicit in \cite{AssiotisDeterminantal,assiotis2023integrablemodelsinhomogeneousspace} as a consequence of probabilistic intertwining considerations. It can also be proven directly using the LGV formula \cite{gessel1985binomial}.
    
    \begin{lemma}[Positivity of $\mathfrak{h}_N$]\label{prop:h-positive}
        Assume \textbf{(A1)}. Then for every $N\in\mathbb{N}$ and every $\mathbf{x}\in\mathbf{W}_{\mathrm d}^{N}$,
        \begin{equation*}
            \mathfrak{h}_N(\mathbf{x})>0.
        \end{equation*}

    \end{lemma}

The following, in continuous time, is simply the well-known Karlin-McGregor formula \cite{karlin1959coincidence} while in discrete time is a straightforward consequence of the LGV formula \cite{gessel1985binomial} along with the composition of the Markov kernels.

\begin{lemma}[Karlin-McGregor semigroup]\label{lem:killing-convention}
Fix $r\in\mathscr{T}$ and let $\mathbf{x},\mathbf{y}\in\mathbf{W}_{\mathrm d}^{N}$. For every $s\in\mathscr{T}$ with $s\geq r$, the time-inhomogeneous killed transition kernel satisfies
\begin{equation}\label{KMlocal}
    \mathbb{P}_{r,\mathbf{x}}^{N}\!\left(\mathbf{X}(s)=\mathbf{y},\,\mathcal{T}^{(r)}>s\right)
    =
    \det\left(\mathbb{P}_{r,x_i}\!\left(\mathscr{X}(s)=y_j\right)\right)_{i,j=1}^{N},
\end{equation}
where the single-particle laws $\mathbb{P}_{r,x_i}$ are obtained from Definitions~\ref{contprocess}--\ref{mixprocess} by starting from $x_i$ at absolute time $r$.
\end{lemma}

    \noindent We will now convert the one-particle expansion into an $N$-particle determinant expansion. The next lemma is a Cauchy-Binet-type decomposition tailored to our coefficients and isolates the index structure that drives the leading asymptotics.
    \begin{lemma}\label{Gcauchybinet}
        Assume \hyperref[A1]{\textbf{(A1)}} and \hyperref[A2]{\textbf{(A2)}}. Let $\mathbf{x}=(x_1,\ldots,x_N)\in\mathbb{Z}_+^N$, let $t>\mathrm{T}_{\mathrm{dec}}$, and let $y_1,\ldots,y_N\in\mathscr{P}$ satisfy $y_j(t)\in\mathscr{C}_{t}$ for every $j\in\llbracket N\rrbracket$. Then
	        \begin{equation}\label{doublecauchybinet}
	            \begin{aligned}
		            \det \left(\sum_{n=0}^{x_i}\frac{\rho_{x_i}^{(n)}(\mathrm{u}_{y_j}(t))}{n!}\eta_{n,y_j}(t)\right)_{i,j=1}^{N}
	            &=
	            \frac{1}{N!}\sum_{\mathbf{k},\mathbf{n}\in\mathbb{Z}_+^{N}}
	            \det\left((-1)^{k_i}\frac{\rho_{x_i}^{(n_j+k_j)}(0)}{n_j!k_j!}\right)_{i,j=1}^{N}\\
	            &\quad\times
		            \det\left(\mathscr{U}_{y_{j}}(t)^{k_i}\eta_{n_i,y_j}(t)\right)_{i,j=1}^{N}
	            \Sigma(t)^{-\sum_{i=1}^{N}k_i},
	            \end{aligned}
	        \end{equation}
        \begin{proof}
            We start from the determinant on the left-hand side of \eqref{doublecauchybinet} and expand it by the Leibniz formula. This gives
            \begin{equation}\label{expansionused}
                \det \left(\sum_{n=0}^{x_i}\frac{\rho_{x_i}^{(n)}(\mathrm{u}_{y_j}(t))}{n!}\eta_{n,y_j}(t)\right)_{i,j=1}^{N}
                =
                \sum_{\sigma\in S_N}\textnormal{sgn}(\sigma)
                \prod_{i=1}^{N}
                \left(\sum_{n=0}^{x_i}\frac{\rho_{x_i}^{(n)}(\mathrm{u}_{y_{\sigma(i)}}(t))}{n!}\eta_{n,y_{\sigma(i)}}(t)\right).
            \end{equation}
            All rearrangements below are finite. Indeed, for each fixed $x_i$, the polynomial $\rho_{x_i}$ has degree $x_i$, so $\rho_{x_i}^{(r)}\equiv0$ whenever $r>x_i$. Thus extending sums over $n$ and $k$ to $\mathbb{Z}_+$ only adds zero terms. Expanding the product in \eqref{expansionused} therefore gives
            \begin{equation}\label{orderofoperation}
                \sum_{\sigma\in S_N}\textnormal{sgn}(\sigma)
                \sum_{\mathbf{n}\in\mathbb{Z}_+^{N}}
                \prod_{i=1}^{N}
                \frac{\rho_{x_i}^{(n_i)}(\mathrm{u}_{y_{\sigma(i)}}(t))}{n_i!}
                \eta_{n_i,y_{\sigma(i)}}(t).
            \end{equation}
            Next we Taylor expand each derivative of $\rho_{x_i}$ at the origin. Since $\rho_{x_i}$ is a polynomial, Taylor's formula is exact and finite:
            \begin{equation*}
                \rho_{x_i}^{(n_i)}(\mathrm{u}_{y_{\sigma(i)}}(t))
                =
                \sum_{k_i=0}^{\infty}
                \frac{\rho_{x_i}^{(n_i+k_i)}(0)}{k_i!}
                \mathrm{u}_{y_{\sigma(i)}}(t)^{k_i}.
            \end{equation*}
            Substituting this identity into \eqref{orderofoperation} gives the following finite expression:
            \begin{equation*}
                \sum_{\sigma\in S_N}\textnormal{sgn}(\sigma)
                \sum_{\mathbf{n},\mathbf{k}\in\mathbb{Z}_+^{N}}
                \prod_{i=1}^{N}
                \frac{\rho_{x_i}^{(n_i+k_i)}(0)}{n_i!k_i!}
                \mathrm{u}_{y_{\sigma(i)}}(t)^{k_i}
                \eta_{n_i,y_{\sigma(i)}}(t).
            \end{equation*}
            For fixed $(\mathbf{n},\mathbf{k})$, the sum over $\sigma$ is exactly the determinant expansion in the terminal index $j$. Hence the preceding display equals
            \begin{equation}\label{eq:cauchy-before-sym}
                \sum_{\mathbf{n},\mathbf{k}\in\mathbb{Z}_+^{N}}
                \prod_{i=1}^{N}
                \frac{\rho_{x_i}^{(n_i+k_i)}(0)}{n_i!k_i!}
                \det\left(\mathrm{u}_{y_j}(t)^{k_i}\eta_{n_i,y_j}(t)\right)_{i,j=1}^{N}.
            \end{equation}
            The remaining step is a symmetrisation in the pair-index $(n_i,k_i)$. Because $\mathbb{Z}_+^N\times\mathbb{Z}_+^N$ is invariant under simultaneous permutations of the pairs $(n_i,k_i)$, the value of \eqref{eq:cauchy-before-sym} is unchanged if we average over all such permutations. Thus
            \begin{align}
                &\sum_{\mathbf{n},\mathbf{k}\in\mathbb{Z}_+^{N}}
                \prod_{i=1}^{N}
                \frac{\rho_{x_i}^{(n_i+k_i)}(0)}{n_i!k_i!}
                \det\left(\mathrm{u}_{y_j}(t)^{k_i}\eta_{n_i,y_j}(t)\right)_{i,j=1}^{N}\notag\\
                &\quad\label{eq:symmetrised-cauchy-binet-sum}
                =
                \frac{1}{N!}
                \sum_{\mathbf{n},\mathbf{k}\in\mathbb{Z}_+^{N}}
                \sum_{\mu\in S_N}
                \prod_{i=1}^{N}
                \frac{\rho_{x_i}^{(n_{\mu(i)}+k_{\mu(i)})}(0)}{n_{\mu(i)}!k_{\mu(i)}!}
                \det\left(\mathrm{u}_{y_j}(t)^{k_{\mu(i)}}\eta_{n_{\mu(i)},y_j}(t)\right)_{i,j=1}^{N}.
            \end{align}
            In the determinant in \eqref{eq:symmetrised-cauchy-binet-sum}, the row indexed by $i$ is the row associated with the pair $(n_{\mu(i)},k_{\mu(i)})$. Reordering these rows back to the order $(n_i,k_i)$ contributes the factor $\textnormal{sgn}(\mu)$. Therefore the inner sum over $\mu$ is the Leibniz expansion of the determinant
            \begin{equation*}
                \det\left(\frac{\rho_{x_i}^{(n_j+k_j)}(0)}{n_j!k_j!}\right)_{i,j=1}^{N},
            \end{equation*}
            and \eqref{eq:symmetrised-cauchy-binet-sum} becomes
            \begin{equation}\label{eq:cauchy-u-form}
                \frac{1}{N!}
                \sum_{\mathbf{n},\mathbf{k}\in\mathbb{Z}_+^{N}}
                \det\left(\frac{\rho_{x_i}^{(n_j+k_j)}(0)}{n_j!k_j!}\right)_{i,j=1}^{N}
                \det\left(\mathrm{u}_{y_j}(t)^{k_i}\eta_{n_i,y_j}(t)\right)_{i,j=1}^{N}.
            \end{equation}
            Finally, by the definition of the stationary-point coordinate, $\mathscr{U}_{y_j}(t)=-\Sigma(t)\mathrm{u}_{y_j}(t)$. Hence
            \begin{equation*}
                \mathrm{u}_{y_j}(t)^{k_i}=(-1)^{k_i}\Sigma(t)^{-k_i}\mathscr{U}_{y_j}(t)^{k_i}.
            \end{equation*}
            In the second determinant in \eqref{eq:cauchy-u-form}, this factor is constant along row $i$. Pulling the factors $\Sigma(t)^{-k_i}$ out of the rows gives the global factor $\Sigma(t)^{-\sum_i k_i}$, while the product of the signs $\prod_i(-1)^{k_i}$ may equivalently be inserted into the rows of the first determinant. This gives exactly \eqref{doublecauchybinet}.
        \end{proof}
    \end{lemma}

 The next proposition lifts the one-particle estimates from Proposition~\ref{Thm2} to the killed $N$-particle semigroup. The starting point is the killed determinant formula of Lemma~\ref{lem:killing-convention} into its entries we insert the one-particle expansion of Section~\ref{decouplingsection}. We then use Lemma~\ref{Gcauchybinet} to reorganize the result into a finite sum indexed by pairs $(\mathbf{n},\mathbf{k})\in\mathbb{Z}_+^N\times\mathbb{Z}_+^N$.

 The main issue is to identify which index families contribute at order $\Sigma(t)^{-N(N-1)/2}$. The cancellation mechanism comes from the antisymmetry of the first determinant, while the parity information is supplied by \eqref{leadingordereta}. The final expectation is evaluated by combining the moment bounds of Proposition~\ref{boundedmoments} with the quantitative Gaussian estimate from Proposition~\ref{NormalisedProcess}. The moment bound is the ingredient that places the polynomially growing observable $\Delta(\mathbf{U}(t))F(\mathbf{U}(t))$ within the admissible test-function class of that proposition.

For the quantitative statement below define
\begin{equation}\label{eq:Reta-def}
    \mathrm{R}_{\eta}(t)
    \overset{\textnormal{def}}{=}
    \mathrm{C}_{\mathscr{C}}\left(\gamma(t)^{1/2}+4\left(\mathfrak{L}^{-1}\mathrm{C}_{1}^{-1}+2\right)\mathrm{C}_2\gamma(t)^{1/2}+6\Lambda\Sigma(t)^{-2}\right)
    +\sqrt{2}\,\Sigma(t)^{-1}
    +\mathrm{e}^{-\Sigma(t)^{2/3}/(8\,2^{1/3})}.
\end{equation}
We compare this coefficient-error scale with the single Gaussian error scale $\psi_{\mathrm{G}}$ from \eqref{eq:psiG-def}. Set
\begin{equation}\label{eq:CetaG-def}
    \mathrm{C}_{\eta,\mathrm{G}}\overset{\textnormal{def}}{=}\mathrm{C}_{\mathscr{C}}\left(1+4\left(\mathfrak{L}^{-1}\mathrm{C}_{1}^{-1}+2\right)\mathrm{C}_{2}+6\Lambda\right)+\sqrt{2}+1.
\end{equation}
By Definition~\ref{gammadef}, $\gamma(t)\to0$ and $\gamma(t)\geq\Sigma(t)^{-1/2}$. Hence, after increasing the lower time so that $\gamma(t)<1$, $\gamma(t)^{1/2}\leq\psi_{\mathrm{G}}(t)$, $\Sigma(t)^{-1}\leq\psi_{\mathrm{G}}(t)$, $\Sigma(t)^{-2}\leq\psi_{\mathrm{G}}(t)$, and $\exp(-\Sigma(t)^{2/3}/(8\,2^{1/3}))\leq\psi_{\mathrm{G}}(t)$. Substituting these four comparisons into \eqref{eq:Reta-def} gives
\begin{equation}\label{eq:Reta-psi-bound}
    \mathrm{R}_{\eta}(t)\leq \mathrm{C}_{\eta,\mathrm{G}}\psi_{\mathrm{G}}(t)
\end{equation}
for every $t$ beyond this enlarged lower time.

\begin{proposition}[Karlin--McGregor asymptotics]\label{ThmKM}
    Assume \textbf{(A1)} and \textbf{(A2)}. Let $\mathbf{x}\in\mathbf{W}_{\mathrm d}^{N}$ and let $F:\mathbb{R}^{N}\to\mathbb{R}$ be bounded and Lipschitz. Then there exist $\mathrm{T}_{\mathrm{KM}}(F,\mathbf{x})\in\mathscr{T}$ and $\mathrm{C}_{\mathrm{KM}}(F,\mathbf{x})>0$ such that, for every $t>\mathrm{T}_{\mathrm{KM}}(F,\mathbf{x})$,
	    \begin{align}\label{keyexpectation}
	        \left|
	        \frac{\mathbb{E}_{\mathbf{x}}^N\!\left[F(\mathbf{U}(t))\mathbf{1}_{\{\mathcal{T}>t\}}\right]}
	             {\mathfrak{h}_N(\mathbf{x})\,\Sigma(t)^{-\frac{N}{2}(N-1)}}
	        -
	        \frac{1}{(2\pi)^{N/2}}
	        \int_{\mathbf{W}_{\mathrm c}^{N}}F(\mathbf{u})\,\Delta(\mathbf{u})\,\exp\!\Big(-\tfrac12|\mathbf{u}|^2\Big)\,\mathrm{d}\mathbf{u}
	        \right|
	        \leq
		        \mathrm{C}_{\mathrm{KM}}(F,\mathbf{x})
		        \psi_{\mathrm{G}}(t).
	    \end{align}
		        
    \begin{proof}
        \noindent\textbf{Outline of strategy:} First we split the terminal sum into endpoints lying in $\mathscr{C}_{t}$ and its complement; the complement is controlled by the endpoint-window estimate in Proposition~\ref{localizationtail} applied coordinatewise. On the endpoint window, Proposition~\ref{Thm2} gives a uniform one-particle expansion, so the Karlin--McGregor determinant can be expanded by Lemma~\ref{Gcauchybinet}. The determinant cancellation and parity estimates then identify the unique leading family. The remaining expectation is evaluated using the Gaussian limit in Proposition~\ref{NormalisedProcess} and the moment bounds in Proposition~\ref{boundedmoments}.
      
        \textbf{Step 1: Restriction to the endpoint window and determinant expansion.} Set
        $\mathrm{c}_{\varphi}=1/(8\,2^{1/3})$. Recall that $\mathscr{C}_{t}^{N}$ denotes the coordinatewise endpoint window from the beginning of this section. Applying Proposition~\ref{localizationtail} to each coordinate and using the finite union bound $\{\mathbf{X}(t)\notin\mathscr{C}_{t}^{N}\}\subseteq\bigcup_{i=1}^{N}\{\mathscr{X}_i(t)\notin\mathscr{C}_{t}\}$, there exist constants $\mathrm{T}_{\mathrm{loc}}\in\mathscr{T}$, $\mathrm{c}_{\mathrm{loc}}>0$ and $\mathrm{M}_{\mathbf{x}}>0$ such that, for $t>\mathrm{T}_{\mathrm{loc}}$,
        \begin{equation}\label{eq:KM-tail-window}
            \mathbb{P}_{\mathbf{x}}^N\left[\mathbf{X}(t)\notin\mathscr{C}_{t}^{N}\right]
            \leq
            \mathrm{M}_{\mathbf{x}}\mathrm{e}^{-\mathrm{c}_{\mathrm{loc}}\Sigma(t)}.
        \end{equation}
        Since $F$ is bounded, the part of the expectation in \eqref{keyexpectation} coming from $\mathbf{X}(t)\notin\mathscr{C}_{t}^{N}$ is bounded in absolute value by $\|F\|_\infty$ times the right-hand side of \eqref{eq:KM-tail-window}. Choose $\mathrm{T}_{\mathrm{tail}}\geq\mathrm{T}_{\mathrm{loc}}$ so that $\mathrm{e}^{-\mathrm{c}_{\mathrm{loc}}\Sigma(t)}\leq\Sigma(t)^{-\mathrm{d}_N-1}$ for every $t>\mathrm{T}_{\mathrm{tail}}$, which is possible because $\Sigma(t)\to\infty$. Then
        \begin{equation}\label{eq:KM-tail-power}
            \|F\|_\infty\mathrm{M}_{\mathbf{x}}\mathrm{e}^{-\mathrm{c}_{\mathrm{loc}}\Sigma(t)}
            \leq
            \mathrm{C}_{\mathrm{tail}}\Sigma(t)^{-\mathrm{d}_N-1},
        \end{equation}
        for a constant $\mathrm{C}_{\mathrm{tail}}>0$ depending only on $N$, $\mathbf{x}$, $F$, and the constants in Proposition~\ref{localizationtail}.
        
         Now fix $\mathbf{m}\in\mathbf{W}_{\mathrm d}^{N}\cap\mathscr{C}_{t}^{N}$. For each component $m_j$, choose any path $y_j\in\mathscr{P}$ with $y_j(t)=m_j$. Since the phase and all coefficients in Proposition~\ref{Thm2} depend on this auxiliary path only through its terminal value at time $t$, we use the following local shorthand in this proof: for any endpoint $m\in\mathscr{C}_{t}$, $\mathfrak{F}_{m}(\cdot,t)$, $\mathrm{u}_{m}(t)$ and $\eta_{n,m}(t)$ denote the corresponding objects obtained from any path with terminal value $m$ at time $t$. By Lemma~\ref{lem:killing-convention} (applied at $r=0$),
        \begin{equation}\label{KarlinMM}
            \mathbb{P}_{\mathbf{x}}^N\left[\mathbf{X}(t)=\mathbf{m},\mathcal{T}>t\right]
            =
            \det\left(\mathbb{P}_{x_i}\left[\mathscr{X}(t)=m_j\right]\right)_{i,j=1}^{N}.
        \end{equation}
        Applying \eqref{ointc} from Proposition~\ref{Thm2} to each matrix entry, with the auxiliary path $y_j$ attached to the terminal endpoint $m_j$, gives
        \begin{equation*}
            \mathbb{P}_{x_i}\left[\mathscr{X}(t)=m_j\right]
            =
            \left(
            \sum_{n=0}^{x_i}\frac{\rho_{x_i}^{(n)}(\mathrm{u}_{m_j}(t))}{n!}\eta_{n,m_j}(t)
            +\mathrm{E}_{i,j}(t)
            \right)
            \mathbb{P}_{0}\left[\mathscr{X}(t)=m_j\right],
        \end{equation*}
        where
        \begin{equation*}
            |\mathrm{E}_{i,j}(t)|\leq \mathrm{C}_{x_i}\mathrm{e}^{-\mathrm{c}_{\varphi}\Sigma(t)^{2/3}}.
        \end{equation*}
        This is exactly the uniform complementary-contour error in \eqref{ointc}. The stationary points $\mathrm{u}_{m_j}(t)$ lie in a fixed compact subset of the holomorphic neighbourhood, the polynomial derivatives $\rho_{x_i}^{(n)}$ involve only finitely many $n$ because $x_i$ is fixed, and the coefficient bounds in \eqref{leadingordereta} bound the finitely many $\eta_{n,m_j}(t)$ uniformly over $m_j\in\mathscr{C}_{t}$ after the lower time is chosen to dominate the corresponding finite set of coefficient thresholds. Here and below, $\mathbb{P}_{\underline{0}}^N$ denotes the independent $N$-particle product law with initial vector $\underline{0}=(0,\ldots,0)$; it is not a killed or conditioned law. Hence determinant multilinearity gives a constant $\mathrm{C}_{N,\mathbf{x}}>0$ such that
        \begin{equation}\label{decoupedexp}
            \left|
            \frac{\mathbb{P}_{\mathbf{x}}^N[\mathbf{X}(t)=\mathbf{m},\mathcal{T}>t]}
                 {\mathbb{P}_{\underline{0}}^N[\mathbf{X}(t)=\mathbf{m}]}
            -
            \frac{1}{N!}\sum_{\mathbf{n},\mathbf{k}\in\mathbb{Z}_+^{N}}\Psi_{\mathbf{n},\mathbf{k},\mathbf{m}}(t)
            \right|
            \leq
            \mathrm{C}_{N,\mathbf{x}}\mathrm{e}^{-\mathrm{c}_{\varphi}\Sigma(t)^{2/3}},
        \end{equation}
        where
        \begin{equation*}
            \Psi_{\mathbf{n},\mathbf{k},\mathbf{m}}(t)
            \overset{\textnormal{def}}{=}
            \det\left((-1)^{k_i}\frac{\rho_{x_i}^{(n_j+k_j)}(0)}{n_j!k_j!}\right)_{i,j=1}^{N}
            \det\left(\mathscr{U}_{m_j}(t)^{k_i}\eta_{n_i,m_j}(t)\right)_{i,j=1}^{N}
            \Sigma(t)^{-\sum_{i=1}^{N}k_i}.
        \end{equation*}
        The sum is finite: if $n_j+k_j>x_i$ for every row index $i$ in some column, then the corresponding column of the first determinant is zero. Therefore a non-zero term must satisfy $n_j+k_j\leq\max_i x_i$ for every column index $j$.
        Choose $\mathrm{T}_{\mathrm{err}}\in\mathscr{T}$ so that, for every $t>\mathrm{T}_{\mathrm{err}}$,
        \begin{equation*}
            \mathrm{e}^{-\mathrm{c}_{\varphi}\Sigma(t)^{2/3}}
            \leq
            \Sigma(t)^{-\mathrm{d}_N-1}.
        \end{equation*}
        This is possible because $\Sigma(t)\to\infty$, and it ensures that the accumulated determinant-expansion error in \eqref{decoupedexp}, after multiplication by the bounded function $F$ and summation against $\mathbb{P}_{\underline{0}}^N$, is absorbed into the same remainder scale as \eqref{eq:KM-tail-power}.
      
         For $\mathbf{m}\in\mathscr{C}_{t}^{N}$ write $\mathbf{U}_{\mathbf{m}}(t)\overset{\textnormal{def}}{=}(\mathscr{U}_{m_1}(t),\ldots,\mathscr{U}_{m_N}(t))$. By Proposition~\ref{twobadprobs-fixed}, there exists $\mathrm{T}_{\mathrm{ord}}\in\mathscr{T}$ such that, for all $t>\mathrm{T}_{\mathrm{ord}}$, the ordering event for endpoints in $\mathscr{C}_{t}^{N}$ is equivalent to the ordering event for the corresponding stationary-point coordinates. Summing \eqref{decoupedexp} against $F(\mathbf{U}_{\mathbf{m}}(t))\mathbb{P}_{\underline{0}}^N[\mathbf{X}(t)=\mathbf{m}]$ over $\mathbf{m}\in\mathbf{W}_{\mathrm d}^{N}\cap\mathscr{C}_{t}^{N}$, and using \eqref{eq:KM-tail-power}, gives
        \begin{equation}\label{sumofkarlin}
            \begin{aligned}
                \mathbb{E}_{\mathbf{x}}^N\left[F(\mathbf{U}(t))\mathbf{1}_{\{\mathcal{T}>t\}}\right]
                &=
                \mathbb{E}_{\underline{0}}^N\left[
                \frac{1}{N!}\sum_{\mathbf{n},\mathbf{k}\in\mathbb{Z}_+^{N}}
                \Psi_{\mathbf{n},\mathbf{k},\mathbf{X}(t)}(t)
                F(\mathbf{U}(t))
                \mathbf{1}_{\{\mathbf{U}(t)\in\mathbf{W}_{\mathrm c}^{N},\,\mathbf{X}(t)\in\mathscr{C}_{t}^{N}\}}
                \right]\\
                &\quad+
                \mathrm{R}_{\mathrm{KM},1}(t),
            \end{aligned}
        \end{equation}
        with
        \begin{equation}\label{eq:RKM1-bound}
            |\mathrm{R}_{\mathrm{KM},1}(t)|
            \leq
            \mathrm{C}_{\mathrm{KM},1}\Sigma(t)^{-\mathrm{d}_N-1},
        \end{equation}
        for every $t>\max(\mathrm{T}_{\mathrm{tail}},\mathrm{T}_{\mathrm{ord}},\mathrm{T}_{\mathrm{dec}},\mathrm{T}_{\mathrm{err}})$, for a constant $\mathrm{C}_{\mathrm{KM},1}>0$ depending only on $N$, $\mathbf{x}$, $F$, and the constants already displayed above.
        
        \textbf{Step 2: Uniform bounds for each index family.} Recall the explicit coefficient-error scale $\mathrm{R}_{\eta}(t)$ from \eqref{eq:Reta-def}; after increasing the lower time, \eqref{eq:Reta-psi-bound} gives the simpler bound $\mathrm{R}_{\eta}(t)\leq\mathrm{C}_{\eta,\mathrm{G}}\psi_{\mathrm{G}}(t)$. Since the sum in \eqref{sumofkarlin} is finite and every non-zero term has $n_j+k_j\leq\max_i x_i$, we choose this lower time to dominate $\mathrm{T}_{\mathrm{dec},n}$ for every $0\leq n\leq\max_i x_i$. The expression inside the first parentheses in \eqref{eq:Reta-def} is exactly the endpoint error $\gamma(t)^{1/2}$ together with the local stationary-point bound
        \begin{equation*}
            |\mathrm{u}_m(t)|\leq4\left(\mathfrak{L}^{-1}\mathrm{C}_{1}^{-1}+2\right)\mathrm{C}_2\gamma(t)^{1/2}+6\Lambda\Sigma(t)^{-2},
        \end{equation*}
        which follows by combining the endpoint first-derivative estimate in \textbf{(ii)} of Proposition~\ref{specificcasesub}, the scale comparison $\mathscr{Z}(t)\leq\mathrm{C}_2\Sigma(t)^2$, the denominator lower bound $\mathfrak{F}_{m}^{(2)}(0,t)\geq\Sigma(t)^2/2$, and the quadratic estimate \eqref{secondconditiontosat}. Therefore, through \textbf{(iv)} of Proposition~\ref{specificcasesub}, the first term controls
        \begin{equation*}
            \left|\frac{\mathfrak{F}_{m}^{(2)}(\mathrm{u}_{m}(t),t)}{\Sigma(t)^2}-1\right|.
        \end{equation*}
        The middle term is supplied by the lower curvature bound \eqref{eq:SD-second-lower}; the final term is the complementary-contour contribution. Therefore \eqref{leadingordereta} implies that, for each fixed even $n$ and every $m\in\mathscr{C}_{t}$,
        \begin{equation}\label{eq:eta-even-explicit-KM}
            \left|\eta_{n,m}(t)-c_n\Sigma(t)^{-n}\right|
            \leq
            \mathrm{C}_{n}\mathrm{C}_{\eta,\mathrm{G}}\psi_{\mathrm{G}}(t)\Sigma(t)^{-n},
            \qquad
            c_n\overset{\textnormal{def}}{=}
            \begin{cases}
                1,&n=0,\\
                (-1)^{n/2}(n-1)!!,&n\geq2,
            \end{cases}
        \end{equation}
        while for each fixed odd $n$,
        \begin{equation}\label{eq:eta-odd-explicit-KM}
            |\eta_{n,m}(t)|\leq \mathrm{C}_{n}\Sigma(t)^{-n-1}.
        \end{equation}
        Applying the Leibniz formula to the second determinant in $\Psi_{\mathbf{n},\mathbf{k},\mathbf{X}(t)}(t)$ and using \eqref{eq:eta-even-explicit-KM}--\eqref{eq:eta-odd-explicit-KM} row by row, we obtain, for each fixed pair $(\mathbf{n},\mathbf{k})$, a polynomial $\mathrm{P}_{\mathbf{n},\mathbf{k}}$ with non-negative coefficients and a constant $\mathrm{C}_{\mathbf{n},\mathbf{k},\mathbf{x}}>0$ such that on $\{\mathbf{X}(t)\in\mathscr{C}_{t}^{N}\}$,
        \begin{equation}\label{eq:generic-bound}
            \left|\Psi_{\mathbf{n},\mathbf{k},\mathbf{X}(t)}(t)\right|
            \leq
            \mathrm{C}_{\mathbf{n},\mathbf{k},\mathbf{x}}
            \mathrm{P}_{\mathbf{n},\mathbf{k}}(|\mathbf{U}(t)|)
            \Sigma(t)^{-\alpha(\mathbf{n},\mathbf{k})},
        \end{equation}
        where
        \begin{equation*}
            \alpha(\mathbf{n},\mathbf{k})
            \overset{\textnormal{def}}{=}
            \sum_{i=1}^{N}\left(k_i+n_i+\mathbf{1}_{\{\text{$n_i$ is odd}\}}\right).
        \end{equation*}
        Choose $p$ larger than the degrees of the finitely many polynomials that occur. Then Proposition~\ref{boundedmoments} gives constants $\mathrm{T}_{p,N}\in\mathscr{T}$ and $\mathrm{C}_{p,N,\mathbf{x}}>0$ such that the expectations of these polynomial factors, multiplied by $\mathbf{1}_{\{\mathbf{X}(t)\in\mathscr{C}_{t}^{N}\}}$, are bounded by $\mathrm{C}_{p,N,\mathbf{x}}$ for every $t>\mathrm{T}_{p,N}$. This is the form needed here, because the summation in \eqref{sumofkarlin} is already restricted to $\mathbf{X}(t)\in\mathscr{C}_{t}^{N}$.
     
        \textbf{Step 3: Identification of the leading index family.}
 If two of the integers $n_i+k_i$ coincide, then two columns of the first determinant in $\Psi_{\mathbf{n},\mathbf{k},\mathbf{m}}(t)$ coincide, and hence that term vanishes. Thus only pairs with distinct non-negative integers $n_i+k_i$ can contribute. For such pairs,
        \begin{equation*}
            \sum_{i=1}^{N}(n_i+k_i)
            \geq
            0+1+\cdots+(N-1)
            =
            \mathrm{d}_N.
        \end{equation*}
        Since the parity contribution in $\alpha(\mathbf{n},\mathbf{k})$ is non-negative, every non-zero index family satisfies $\alpha(\mathbf{n},\mathbf{k})\geq\mathrm{d}_N$. If $\alpha(\mathbf{n},\mathbf{k})\geq\mathrm{d}_N+1$, then \eqref{eq:generic-bound} and the moment bound imply that the total contribution of all such pairs is bounded by
        \begin{equation}\label{eq:nonleading-preliminary}
            \mathrm{C}_{N,F,\mathbf{x}}\Sigma(t)^{-\mathrm{d}_N-1}.
        \end{equation}
        Equality $\alpha(\mathbf{n},\mathbf{k})=\mathrm{d}_N$ can occur only if $\{n_i+k_i:i\in\llbracket N\rrbracket\}=\{0,1,\ldots,N-1\}$ and every $n_i$ is even. Among these equality cases, the determinant with the leading part of \eqref{eq:eta-even-explicit-KM} inserted has rows proportional to $\mathscr{U}_{m_j}(t)^{k_i}$. Hence, if two of the $k_i$ coincide, this leading determinant is zero; all remaining terms in the row-by-row expansion contain at least one coefficient error from \eqref{eq:eta-even-explicit-KM} and are therefore bounded, after taking expectation and using Proposition~\ref{boundedmoments}, by $\mathrm{C}_{N,F,\mathbf{x}}\mathrm{C}_{\eta,\mathrm{G}}\psi_{\mathrm{G}}(t)\Sigma(t)^{-\mathrm{d}_N}$. Thus the only equality cases that can contribute to the limiting constant must have distinct $k_i$. If some $n_i$ were a positive even integer in such a case, then
        \begin{equation*}
            \sum_{i=1}^{N}k_i
            =
            \sum_{i=1}^{N}(n_i+k_i)-\sum_{i=1}^{N}n_i
            <
            \mathrm{d}_N.
        \end{equation*}
        But $N$ distinct non-negative integers have sum at least $\mathrm{d}_N$. Therefore, if the $k_i$ are distinct, all $n_i$ must be zero. Hence the only index family that contributes to the limiting constant is
        \begin{equation*}
            \mathcal{M}_N
            \overset{\textnormal{def}}{=}
            \left\{(\mathbf{0},\mathbf{k}_{\sigma}):\sigma\in S_N,\ k_{\sigma,i}=\sigma(i)-1\right\}.
        \end{equation*}
        For the leading family, \eqref{eq:eta-even-explicit-KM} with $n=0$ gives $|\eta_{0,m}(t)-1|\leq\mathrm{C}_0\mathrm{C}_{\eta,\mathrm{G}}\psi_{\mathrm{G}}(t)$ uniformly over $m\in\mathscr{C}_{t}$. Expanding the determinant by the Leibniz formula therefore yields, for every $\sigma\in S_N$,
        \begin{equation}\label{eq:leading-family-explicit}
            \Psi_{\mathbf{0},\mathbf{k}_{\sigma},\mathbf{X}(t)}(t)
            =
            A_{\sigma}(\mathbf{x})
            \det\left(\mathscr{U}_{j}(t)^{\sigma(i)-1}\right)_{i,j=1}^{N}
            \Sigma(t)^{-\mathrm{d}_N}
            +\mathrm{R}_{\sigma}(t),
        \end{equation}
        Here $\mathscr{U}_{j}(t)$ denotes the $j$th component of the random vector $\mathbf{U}(t)$, and the remainder satisfies
        \begin{equation*}
            \mathbb{E}_{\underline{0}}^N\left[
            |\mathrm{R}_{\sigma}(t)|\,|F(\mathbf{U}(t))|
            \mathbf{1}_{\{\mathbf{U}(t)\in\mathbf{W}_{\mathrm c}^{N},\,\mathbf{X}(t)\in\mathscr{C}_{t}^{N}\}}
            \right]
            \leq
            \mathrm{C}_{N,F,\mathbf{x}}\mathrm{C}_{\eta,\mathrm{G}}\psi_{\mathrm{G}}(t)\Sigma(t)^{-\mathrm{d}_N}.
        \end{equation*}
        Combining this bound with \eqref{eq:nonleading-preliminary}, and using $\Sigma(t)^{-1}\leq\psi_{\mathrm{G}}(t)$ after increasing the lower time, shows that every pair outside $\mathcal{M}_N$ contributes at most
        \begin{equation}\label{eq:KM-nonleading-explicit}
            \mathrm{C}_{N,F,\mathbf{x}}\left(1+\mathrm{C}_{\eta,\mathrm{G}}\right)\psi_{\mathrm{G}}(t)
            \Sigma(t)^{-\mathrm{d}_N}
        \end{equation}
        to \eqref{sumofkarlin}.
      
        \textbf{Step 4: Evaluation of the leading family.}
        \noindent For $\sigma\in S_N$, pulling out row signs and then permuting columns gives
        \begin{equation*}
            A_{\sigma}(\mathbf{x})
            =
            \det\left((-1)^{\sigma(i)-1}\frac{\rho_{x_i}^{(\sigma(j)-1)}(0)}{(\sigma(j)-1)!}\right)_{i,j=1}^{N}
            =
            \textnormal{sgn}(\sigma)\,\mathfrak{h}_N(\mathbf{x}).
        \end{equation*}
        Similarly,
        \begin{equation*}
            \det\left(\mathscr{U}_{j}(t)^{\sigma(i)-1}\right)_{i,j=1}^{N}
            =
            \textnormal{sgn}(\sigma)\,\Delta(\mathbf{U}(t)).
        \end{equation*}
        The signs cancel. Since $|\mathcal{M}_N|=N!$, substituting \eqref{eq:leading-family-explicit} and \eqref{eq:KM-nonleading-explicit} into \eqref{sumofkarlin}, and then using \eqref{eq:Reta-psi-bound} together with $\Sigma(t)^{-1}\leq\psi_{\mathrm{G}}(t)$ after increasing the lower time, gives
        \begin{equation}\label{expectform}
            \begin{aligned}
            &\left|
            \mathbb{E}_{\mathbf{x}}^N\left[F(\mathbf{U}(t))\mathbf{1}_{\{\mathcal{T}>t\}}\right]
            -
            \mathfrak{h}_N(\mathbf{x})\Sigma(t)^{-\mathrm{d}_N}
            \mathbb{E}_{\underline{0}}^N\left[
            \Delta(\mathbf{U}(t))F(\mathbf{U}(t))
            \mathbf{1}_{\{\mathbf{U}(t)\in\mathbf{W}_{\mathrm c}^{N},\,\mathbf{X}(t)\in\mathscr{C}_{t}^{N}\}}
            \right]
            \right|\\
            &\qquad\leq
            \mathrm{C}_{N,F,\mathbf{x}}
            \left(1+\mathrm{C}_{\eta,\mathrm{G}}\right)\psi_{\mathrm{G}}(t)
            \Sigma(t)^{-\mathrm{d}_N}.
            \end{aligned}
        \end{equation}
       
        \textbf{Step 5: Identification of the remaining expectation.}
        \noindent Define the observable
        \begin{equation*}
            H_F(\mathbf{u})\overset{\textnormal{def}}{=}\Delta(\mathbf{u})F(\mathbf{u})\mathbf{1}_{\{\mathbf{u}\in\mathbf{W}_{\mathrm c}^{N}\}}.
        \end{equation*}
        This is exactly the function to which we apply Proposition~\ref{NormalisedProcess}. Since $F$ is bounded and Lipschitz, $H_F$ satisfies the hypotheses of that proposition. We spell this out because the indicator of the Weyl chamber appears in the definition of $H_F$. The Vandermonde determinant is a polynomial of degree $\mathrm{d}_N$, and therefore, on every ball of radius $\mathrm{R}\geq1$, both $\Delta$ and its gradient are bounded by constants depending only on $N$ times $\mathrm{R}^{\mathrm{d}_N}$. If both points lie in $\mathbf{W}_{\mathrm c}^{N}$, the local Lipschitz bound for $H_F$ follows from the product rule and the Lipschitz bound for $F$; if both points lie outside $\mathbf{W}_{\mathrm c}^{N}$, both values of $H_F$ are zero. If one point lies in $\mathbf{W}_{\mathrm c}^{N}$ and the other does not, choose a point on the line segment between them that lies on the boundary; at this point $\Delta=0$, and applying the same gradient bound for $\Delta$ on each of the two subsegments gives the required Lipschitz estimate across the boundary. Thus there are constants $\mathrm{C}_{H_F}>0$ and $q_F=\mathrm{d}_N+1$ such that $H_F$ satisfies the growth and local Lipschitz assumptions of Proposition~\ref{NormalisedProcess}. Therefore, for every $t$ beyond the lower time supplied by that proposition,
        \begin{align}\label{general-limit}
            &\left|
            \mathbb{E}_{\underline{0}}^N\left[
            \Delta(\mathbf{U}(t))F(\mathbf{U}(t))
            \mathbf{1}_{\{\mathbf{U}(t)\in\mathbf{W}_{\mathrm c}^{N},\,\mathbf{X}(t)\in\mathscr{C}_{t}^{N}\}}
            \right]\right.\\
            &\qquad\left.
            -
            \frac{1}{(2\pi)^{N/2}}
            \int_{\mathbf{W}_{\mathrm c}^{N}}F(\mathbf{u})\,\Delta(\mathbf{u})\,\exp\!\Big(-\tfrac12|\mathbf{u}|^2\Big)\,\mathrm{d}\mathbf{u}
            \right|
            \leq
            \mathrm{C}_{H_F,\underline{0}}\psi_{\mathrm{G}}(t).
            \notag
        \end{align}
        Dividing \eqref{expectform} by $\mathfrak{h}_N(\mathbf{x})\Sigma(t)^{-\mathrm{d}_N}$ and combining the result with \eqref{general-limit} gives \eqref{keyexpectation} after increasing $\mathrm{C}_{\mathrm{KM}}(F,\mathbf{x})$ so that it dominates $\mathfrak{h}_N(\mathbf{x})^{-1}\mathrm{C}_{N,F,\mathbf{x}}(1+\mathrm{C}_{\eta,\mathrm{G}})+\mathrm{C}_{H_F,\underline{0}}$.
    \end{proof}
\end{proposition}

\section{Proof of main results}

    \subsection{Proof of Theorem~\ref{collisionThm}}

        \begin{proof}[\textbf{Proof of Theorem~\ref{collisionThm}}]
        Applying Proposition~\ref{ThmKM} with $F\equiv1$ gives a quantitative estimate. In this case the observable in Proposition~\ref{NormalisedProcess} is $H_1(\mathbf{u})=\Delta(\mathbf{u})\mathbf{1}_{\{\mathbf{u}\in\mathbf{W}_{\mathrm c}^{N}\}}$, and therefore, for every $t>\mathrm{T}_{\mathrm{KM}}(1,\mathbf{x})$,
        \begin{align*}
            &\left|
            \frac{\mathbb{P}_{\mathbf{x}}^N[\mathcal{T}>t]}
            {\mathfrak{h}_N(\mathbf{x})\Sigma(t)^{-\mathrm{d}_N}}
            -
            \frac{1}{(2\pi)^{N/2}}
            \int_{\mathbf{W}_{\mathrm c}^{N}}\Delta(\mathbf{u})
            \exp\!\Big(-\tfrac12|\mathbf{u}|^2\Big)\,\mathrm{d}\mathbf{u}
            \right|\\
            &\qquad\leq
            \mathrm{C}_{\mathrm{KM}}(1,\mathbf{x})
            \psi_{\mathrm{G}}(t).
        \end{align*}
        It remains only to identify the deterministic integral. By the normalization of the ordered $\mathrm{GOE}_N$ eigenvalue density in Definition~\ref{GOEeigenvalue},
        \begin{equation*}
            \frac{N!}{(2\pi)^{N/2}}
            \prod_{i=1}^{N}\frac{\Gamma(3/2)}{\Gamma(1+i/2)}
            \int_{\mathbf{W}_{\mathrm c}^{N}}\Delta(\mathbf{u})
            \exp\!\Big(-\tfrac12|\mathbf{u}|^2\Big)\,\mathrm{d}\mathbf{u}
            =
            1,
        \end{equation*}
  and from this we readily obtain that  Theorem~\ref{collisionThm} holds with $\mathrm{T}_{\mathrm{col},\mathbf{x}}=\mathrm{T}_{\mathrm{KM}}(1,\mathbf{x})$ and $\mathrm{C}_{\mathrm{col},\mathbf{x}}=\mathfrak{h}_N(\mathbf{x})\mathrm{C}_{\mathrm{KM}}(1,\mathbf{x})$.
        \end{proof}

    \subsection{Proof of Corollary~\ref{Thm1}}
        The argument combines the finite-horizon conditioning identity \eqref{construction} with the collision-tail asymptotics from Theorem~\ref{collisionThm} to construct the limiting law and identify the Doob $h$-transform. The only deterministic point requiring care is that, in the discrete-time regime, the survival factor after time $t$ is computed in the environment shifted by $t$. We isolate this comparison first.

\begin{lemma}[Shifted deterministic scales]\label{lem:shifted-deterministic-scales}
Fix $t\in\mathscr{T}$. In the discrete-time regime, set
\begin{equation*}
\widetilde{\alpha}_n\overset{\textnormal{def}}{=}\alpha_{t+n},\qquad
\widetilde{\beta}_n\overset{\textnormal{def}}{=}\beta_{t+n},\qquad
\widetilde{\tau}_n\overset{\textnormal{def}}{=}\tau_{t+n},
\qquad n\in\mathbb{Z}_+,
\end{equation*}
while in the continuous-time regime there is no shift in the parameters. Define the shifted centre and shifted fluctuation scale by
\begin{equation*}
\mathrm{z}_t(s)\overset{\textnormal{def}}{=}
\mathbf{1}_{\{\mathscr{T}=\mathbb{R}_+\}}\frac{s-t}{\Lambda}
+\mathbf{1}_{\{\mathscr{T}=\mathbb{Z}_+\}}\frac{1}{\Lambda}\sum_{n=t}^{s-1}\big(\alpha_n(1-\tau_n)+\beta_n\tau_n\big),
\qquad
\mathscr{Z}_t(s)\overset{\textnormal{def}}{=}\lfloor \mathrm{z}_t(s)\rceil,
\end{equation*}
and
\begin{equation*}
\Sigma_t(s)^2\overset{\textnormal{def}}{=}
\sum_{n=0}^{\mathscr{Z}_t(s)}\frac{1}{\lambda_n^2}
-\mathbf{1}_{\{\mathscr{T}=\mathbb{Z}_+\}}\sum_{n=t}^{s-1}\big(\alpha_n^2(1-\tau_n)-\beta_n^2\tau_n\big).
\end{equation*}
The corresponding shifted error scale is
\begin{equation*}
\gamma_t(s)\overset{\textnormal{def}}{=}\max\left(\sup_{r\geq s}\left|\frac{1}{\mathscr{Z}_t(r)+1}\sum_{n=0}^{\mathscr{Z}_t(r)}\frac{1}{\lambda_n}-\Lambda\right|,\sup_{r\geq s}\Sigma_t(r)^{-\frac{1}{2}}\right),
\qquad
\psi_{\mathrm{G},t}(s)\overset{\textnormal{def}}{=}\left(1+\log(\gamma_t(s)^{-1})\right)^{-\frac14},
\end{equation*}
whenever the lower time has been chosen so that $\gamma_t(s)<1$. Then the shifted model satisfies \textbf{(A1)} and \textbf{(A2)} with the same structural constants as the original model. Moreover, there exist a constant $\mathrm{C}_{\mathrm{A1}}>0$, depending only on the constants in \textbf{(A1)} and on $\Lambda$, and a threshold $\mathrm{T}_{\mathrm{sh}}(t)\in\mathscr{T}$, depending on the fixed shift $t$, such that, for every $s>\mathrm{T}_{\mathrm{sh}}(t)$,
\begin{equation}\label{eq:shifted-square-scale-bound}
\left|\Sigma_t(s)^2-\Sigma(s)^2\right|\leq \mathrm{C}_{\mathrm{A1}}(1+t).
\end{equation}
Consequently, after increasing the lower time depending on $t$,
\begin{equation}\label{eq:shifted-scale-ratio}
\left|\frac{\Sigma_t(s)}{\Sigma(s)}-1\right|\leq\frac{\mathrm{C}_{\mathrm{A1}}(1+t)}{\Sigma(s)^2}.
\end{equation}
In particular, $\Sigma_t(s)/\Sigma(s)\rightarrow1$ and $\gamma_t(s)\rightarrow0$ as $s\rightarrow\infty$.
\end{lemma}

\begin{proof}
The shifted model has the same spatial rates $(\lambda_n)_{n\in\mathbb{Z}_+}$. In discrete time, \textbf{(A1)} is preserved because it is uniform in the temporal index, and \textbf{(A2)} is unchanged because it involves only the spatial reciprocal-rate averages. The continuous-time case is immediate, since there is no temporal environment to shift.

 We compare the deterministic centres before rounding. Set
\begin{equation*}
\mathrm{C}_{\mathrm{sh}}\overset{\textnormal{def}}{=}
\frac{\max(1,\mathfrak{B}_2^{-1},\mathfrak{U})}{\Lambda}.
\end{equation*}
By the definitions of $\mathrm{z}$ and $\mathrm{z}_t$, and by the temporal bounds in \textbf{(A1)},
\begin{equation*}
\left|\mathrm{z}(s)-\mathrm{z}_t(s)\right|
\leq
\mathbf{1}_{\{\mathscr{T}=\mathbb{R}_+\}}\frac{t}{\Lambda}
+\mathbf{1}_{\{\mathscr{T}=\mathbb{Z}_+\}}\frac{1}{\Lambda}\sum_{n=0}^{t-1}\big(\alpha_n(1-\tau_n)+\beta_n\tau_n\big)
\leq
\mathrm{C}_{\mathrm{sh}}t.
\end{equation*}
Since both $\mathscr{Z}(s)$ and $\mathscr{Z}_t(s)$ are nearest-integer roundings, the two rounding errors contribute at most one in total, and hence
\begin{equation*}
\left|\mathscr{Z}(s)-\mathscr{Z}_t(s)\right|
\leq
\mathrm{C}_{\mathrm{sh}}t+1.
\end{equation*}
Therefore the spatial parts of $\Sigma(s)^2$ and $\Sigma_t(s)^2$ differ by at most
\begin{equation*}
\left|
\sum_{n=0}^{\mathscr{Z}(s)}\frac{1}{\lambda_n^2}
-\sum_{n=0}^{\mathscr{Z}_t(s)}\frac{1}{\lambda_n^2}
\right|
\leq
\mathfrak{L}^{-2}\left|\mathscr{Z}(s)-\mathscr{Z}_t(s)\right|
\leq
\mathfrak{L}^{-2}\big(\mathrm{C}_{\mathrm{sh}}t+1\big),
\end{equation*}
where the first inequality uses $\lambda_n^{-2}\leq\mathfrak{L}^{-2}$ from \textbf{(A1)}. Similarly, in the discrete-time regime,
\begin{equation*}
\left|
\sum_{n=0}^{s-1}\big(\alpha_n^2(1-\tau_n)-\beta_n^2\tau_n\big)
-\sum_{n=t}^{s-1}\big(\alpha_n^2(1-\tau_n)-\beta_n^2\tau_n\big)
\right|
\leq
\sum_{n=0}^{t-1}\left(\mathfrak{B}_2^{-2}+\mathfrak{U}^2\right)
\leq \bigl(\mathfrak{B}_2^{-2}+\mathfrak{U}^2\bigr)t,
\end{equation*}
again by \textbf{(A1)}. Taking
\begin{equation*}
\mathrm{C}_{\mathrm{A1}}\overset{\textnormal{def}}{=}
\mathfrak{L}^{-2}(\mathrm{C}_{\mathrm{sh}}+1)
+\mathfrak{B}_2^{-2}+\mathfrak{U}^2
\end{equation*}
gives \eqref{eq:shifted-square-scale-bound}.

 Dividing \eqref{eq:shifted-square-scale-bound} by $\Sigma(s)^2$ gives the corresponding bound for the squared scales. Increase the lower time so that $\mathrm{C}_{\mathrm{A1}}(1+t)\Sigma(s)^{-2}\leq1/2$. Then both scales are positive and the elementary inequality $|\sqrt{1+r}-1|\leq |r|$, valid for $r\geq-1/2$, gives \eqref{eq:shifted-scale-ratio}. Since $\Sigma(s)\to\infty$ by \textbf{(iv)} of Proposition~\ref{specificcasesub}, this proves $\Sigma_t(s)/\Sigma(s)\to1$.

 It remains only to check the shifted logarithmic scale. The first term in $\gamma_t(s)$ tends to zero by \textbf{(A2)}, since $\mathscr{Z}_t(r)\to\infty$ uniformly as $r\to\infty$ with $r\geq s$. The second term tends to zero because \eqref{eq:shifted-scale-ratio} gives $\Sigma_t(s)\geq\Sigma(s)/2$ after increasing the lower time, while $\Sigma(s)\to\infty$. Hence $\gamma_t(s)\to0$.
\end{proof}
        \begin{proof}[Proof of Corollary~\ref{Thm1}]\label{proofofcor}
        We first identify the ratio of survival probabilities which appears in the finite-horizon conditioning formula \eqref{construction}. Use the shifted notation from Lemma~\ref{lem:shifted-deterministic-scales}. The proof of Theorem~\ref{collisionThm} uses only the structural constants in \textbf{(A1)}--\textbf{(A2)}, the convergence control of the reciprocal-rate averages, and the initial point. By Lemma~\ref{lem:shifted-deterministic-scales}, the shifted model has the same structural constants, the same reciprocal-rate convergence, and $\gamma_t(s)\to0$. Thus Theorem~\ref{collisionThm} applies to the shifted model; in this shifted application the deterministic scale and Gaussian error scale are precisely $\Sigma_t(s)$ and $\psi_{\mathrm{G},t}(s)$. 
        Applying Theorem~\ref{collisionThm} to the shifted model gives constants $\mathrm{S}_{\mathrm{col},t,\mathbf{y}}\in\mathscr{T}$ and $\mathrm{C}_{\mathrm{col},t,\mathbf{y}}>0$ such that, for every $s>\mathrm{S}_{\mathrm{col},t,\mathbf{y}}$,
        \begin{equation}\label{shiftedcollisionasymptotic}
        \left|\frac{H_{t,s}(\mathbf{y})}{\mathrm{A}_N\mathfrak{h}_N(\mathbf{y})\,\Sigma_t(s)^{-\mathrm{d}_N}}-1\right|\leq\frac{\mathrm{C}_{\mathrm{col},t,\mathbf{y}}}{\mathrm{A}_N\mathfrak{h}_N(\mathbf{y})}\psi_{\mathrm{G},t}(s).
        \end{equation}
        Applying Theorem~\ref{collisionThm} at the initial time $0$ gives constants $\mathrm{S}_{\mathrm{col},\mathbf{x}}\in\mathscr{T}$ and $\mathrm{C}_{\mathrm{col},\mathbf{x}}>0$ such that, for every $s>\mathrm{S}_{\mathrm{col},\mathbf{x}}$,
        \begin{equation*}
        \left|\frac{H_{0,s}(\mathbf{x})}{\mathrm{A}_N\mathfrak{h}_N(\mathbf{x})\,\Sigma(s)^{-\mathrm{d}_N}}-1\right|\leq\frac{\mathrm{C}_{\mathrm{col},\mathbf{x}}}{\mathrm{A}_N\mathfrak{h}_N(\mathbf{x})}\psi_{\mathrm{G}}(s).
        \end{equation*}

        Set $\mathrm{C}_{N,\mathrm{pow}}\overset{\textnormal{def}}{=}\max\{1,\mathrm{d}_N2^{\mathrm{d}_N+1}\}$. By Lemma~\ref{lem:shifted-deterministic-scales}, after increasing the lower time so that the right-hand side of \eqref{eq:shifted-scale-ratio} is at most $1/2$, the mean-value theorem on the interval $[1/2,3/2]$ gives
        \begin{equation*}
        \left|\left(\frac{\Sigma(s)}{\Sigma_t(s)}\right)^{\mathrm{d}_N}-1\right|\leq \mathrm{C}_{N,\mathrm{pow}}\frac{\mathrm{C}_{\mathrm{A1}}(1+t)}{\Sigma(s)^2}.
        \end{equation*}
        Combining the two collision-tail estimates with this scale comparison, and increasing the lower time so that the two normalized collision-tail errors are at most $1/2$, gives, for every $s$ larger than the maximum of the deterministic thresholds introduced in this paragraph,
        \begin{align}\label{eq:survival-ratio-quant}
        &\left|\frac{H_{t,s}(\mathbf{y})}{H_{0,s}(\mathbf{x})}\frac{\mathfrak{h}_N(\mathbf{x})}{\mathfrak{h}_N(\mathbf{y})}-1\right|\\
        &\qquad\leq3\frac{\mathrm{C}_{\mathrm{col},t,\mathbf{y}}}{\mathrm{A}_N\mathfrak{h}_N(\mathbf{y})}\psi_{\mathrm{G},t}(s)+3\frac{\mathrm{C}_{\mathrm{col},\mathbf{x}}}{\mathrm{A}_N\mathfrak{h}_N(\mathbf{x})}\psi_{\mathrm{G}}(s)+\mathrm{C}_{N,\mathrm{pow}}\frac{\mathrm{C}_{\mathrm{A1}}(1+t)}{\Sigma(s)^2}.\notag
        \end{align}
        The right-hand side in \eqref{eq:survival-ratio-quant} tends to zero as $s\to\infty$.

      Now fix $T\in\mathscr{T}$. It follows from \cite{AssiotisDeterminantal,assiotis2023integrablemodelsinhomogeneousspace}, and can also be readily shown from the results of this paper, that,
\begin{equation*}
\mathbb{E}_\mathbf{x}^N\left[\mathfrak{h}_N(\mathbf{X}(T))\mathbf{1}_{\mathcal{T}>T}\right]=\mathfrak{h}_N(\mathbf{x}).
\end{equation*}
Hence, the following measure, with $\mathscr{B}\in \mathcal{F}_{0,T}$, 
        \begin{equation}\label{DoobTransMeasureInProof}
        \mathscr{B}\mapsto
        \frac{1}{\mathfrak{h}_N(\mathbf{x})}
        \mathbb{E}_{\mathbf{x}}^N\left[
        \mathfrak{h}_N\!\left(\mathbf{X}(T)\right)\mathbf{1}_{\mathscr{B}}\mathbf{1}_{\{\mathcal{T}>T\}}
        \right],
        \qquad \mathscr{B}\in\mathcal{F}_{0,T},
        \end{equation}
is a probability measure. To take the limit, we first work on cylinder events whose final endpoint is prescribed:
        \begin{equation*}
        \mathscr{A}\overset{\textnormal{def}}{=}\left\{\mathbf{X}(t_1)=\mathbf{y}_1,\ldots,\mathbf{X}(t_m)=\mathbf{y}_m\right\}\in\mathcal{F}_{0,T},
        \qquad 0<t_1<\cdots<t_m=T.
        \end{equation*}
        If $\mathbf{y}_m\notin\mathbf{W}_{\mathrm d}^{N}$, then both the finite-horizon conditioned probability and the Doob-transform expression below are zero because the factor $\mathbf{1}_{\{\mathcal{T}>T\}}$ forces the terminal state to remain ordered. We therefore assume $\mathbf{y}_m\in\mathbf{W}_{\mathrm d}^{N}$ in the calculation below. By \eqref{construction},
        \begin{equation*}
        \mathsf{P}_{0,\mathbf{x}}^{N,(s)}(\mathscr{A})
        =
        \mathbb{E}_{\mathbf{x}}^N\left[
        \frac{H_{T,s}(\mathbf{X}(T))}{H_{0,s}(\mathbf{x})}
        \mathbf{1}_{\mathscr{A}}\mathbf{1}_{\{\mathcal{T}>T\}}
        \right].
        \end{equation*}
        On the event $\mathscr{A}$, we have $\mathbf{X}(T)=\mathbf{y}_m$, and therefore \eqref{eq:survival-ratio-quant}, with $t=T$ and $\mathbf{y}=\mathbf{y}_m$, gives
        \begin{equation}\label{eq:terminal-cylinder-limit}
        \lim_{s\to\infty}\mathsf{P}_{0,\mathbf{x}}^{N,(s)}(\mathscr{A})
        =
        \frac{\mathfrak{h}_N(\mathbf{y}_m)}{\mathfrak{h}_N(\mathbf{x})}
        \mathbb{P}_{\mathbf{x}}^N\left[\mathscr{A},\,\mathcal{T}>T\right].
        \end{equation}
        The cylinder events of this form are a $\pi$-system generating the finite-time path sigma-field $\mathcal{F}_{0,T}$. Their limiting probabilities are consistent as $T$ varies, because they are pointwise limits of the consistent finite-horizon conditioned laws. Hence they define a limiting measure on each $\mathcal{F}_{0,T}$, and these finite-time measures are themselves consistent.
     
        The right-hand side of \eqref{eq:terminal-cylinder-limit} is the value on $\mathscr{A}$ of the probability measure in \eqref{DoobTransMeasureInProof}. Hence the two probability measures agree on the generating $\pi$-system, and the $\pi$-$\lambda$ theorem gives, for every $\mathscr{B}\in\mathcal{F}_{0,T}$,
        \begin{equation*}
        \mathsf{P}_{\mathbf{x}}^N(\mathscr{B})
        =
        \frac{1}{\mathfrak{h}_N(\mathbf{x})}
        \mathbb{E}_{\mathbf{x}}^N\left[
        \mathfrak{h}_N\!\left(\mathbf{X}(T)\right)\mathbf{1}_{\mathscr{B}}\mathbf{1}_{\{\mathcal{T}>T\}}
        \right].
        \end{equation*}
        Extending from indicators to bounded $\mathcal{F}_{0,T}$-measurable functionals by linearity and bounded monotone approximation yields
        \begin{equation*}
        \mathsf{E}^N_{\mathbf{x}}\left[F\left(\hat{\mathbf{X}}(t):t\leq T\right)\right]
        =
        \frac{1}{\mathfrak{h}_N(\mathbf{x})}
        \mathbb{E}_{\mathbf{x}}^N\left[
        \mathfrak{h}_N\!\left(\mathbf{X}(T)\right)
        F\left(\mathbf{X}(t):t\leq T\right)\mathbf{1}_{\{\mathcal{T}>T\}}
        \right],
        \end{equation*}
        thus completing the proof.
        \end{proof}

	    \subsection{Proof of Theorem~\ref{thm:stationarypoints}}

\begin{proposition}[Survival-conditioned stationary-point limit]\label{lastcorcor}
		Assume \textbf{(A1)} and \textbf{(A2)}. Let $\mathbf{x}\in\mathbf{W}_{\mathrm d}^{N}$, and let $F:\mathbb{R}^{N}\to\mathbb{R}$ be bounded and Lipschitz. Define
		\begin{equation*}
			\mathrm{L}_{F}\overset{\textnormal{def}}{=}
			\frac{1}{(2\pi)^{N/2}}
			\int_{\mathbf{W}_{\mathrm c}^{N}}F(\mathbf{u})\,\Delta(\mathbf{u})\exp\!\Big(-\tfrac12|\mathbf{u}|^2\Big)\,\mathrm{d}\mathbf{u},
			\qquad
			\mathrm{L}_{0}\overset{\textnormal{def}}{=}
			\frac{1}{N!}\prod_{i=1}^{N}\frac{\Gamma(1+i/2)}{\Gamma(3/2)}.
		\end{equation*}
		Then there exist constants $\mathrm{T}_{\mathrm{GOE},F,\mathbf{x}}\in\mathscr{T}$ and $\mathrm{C}_{\mathrm{GOE},F,\mathbf{x}}>0$ such that, for every $t>\mathrm{T}_{\mathrm{GOE},F,\mathbf{x}}$,
			\begin{align}\label{karlinMcGregor}
				&\left|
				\mathbb{E}_{\mathbf{x}}^N\!\left[F\!\left(\mathbf{U}(t)\right)\middle|\mathcal{T}>t\right]
				-
				\frac{N!}{(2\pi)^{N/2}}
				\prod_{i=1}^{N}\frac{\Gamma(3/2)}{\Gamma(1+i/2)}
				\int_{\mathbf{W}_{\mathrm c}^{N}}F(\mathbf{u})\,\Delta(\mathbf{u})
				\exp\!\Big(-\tfrac12|\mathbf{u}|^2\Big)\,\mathrm{d}\mathbf{u}
				\right|\\
				&\qquad\leq \mathrm{C}_{\mathrm{GOE},F,\mathbf{x}}\psi_{\mathrm{G}}(t).
				\notag
			\end{align}
		In particular, the conditional law of $\mathbf{U}(t)$ given $\{\mathcal{T}>t\}$ is controlled quantitatively by the ordered eigenvalue law of $\mathrm{GOE}_N$ with the single explicit error scale in \eqref{karlinMcGregor}.
\end{proposition}
\begin{proof}

		\noindent Fix bounded Lipschitz $F:\mathbb{R}^{N}\to\mathbb{R}$. Since $\{\mathcal{T}>t\}$ has strictly positive probability for every finite $t$, we may write
		\begin{equation}\label{lastcor-step1}
			\mathbb{E}_{\mathbf{x}}^N\!\left[F\!\left(\mathbf{U}(t)\right)\middle|\mathcal{T}>t\right]
			=
			\frac{\mathbb{E}_{\mathbf{x}}^N\!\left[F\!\left(\mathbf{U}(t)\right)\mathbf{1}_{\{\mathcal{T}>t\}}\right]}
			{\mathbb{P}_{\mathbf{x}}^N\!\left[\mathcal{T}>t\right]}.
		\end{equation}
		We now analyze the numerator and denominator after dividing by their common scale.
			The constant $\mathrm{L}_{0}$ is strictly positive, since every factor in the finite product defining it is positive. Applying Proposition~\ref{ThmKM} to the test function $F$ gives, for every $t>\mathrm{T}_{\mathrm{KM}}(F,\mathbf{x})$,
			\begin{equation}\label{lastcor-step2}
				\left|
				\frac{\mathbb{E}_{\mathbf{x}}^N\!\left[F\!\left(\mathbf{U}(t)\right)\mathbf{1}_{\{\mathcal{T}>t\}}\right]}
				{\mathfrak{h}_N(\mathbf{x})\Sigma(t)^{-\mathrm{d}_N}}
				-\mathrm{L}_{F}
				\right|
				\leq \mathrm{C}_{\mathrm{KM}}(F,\mathbf{x})\psi_{\mathrm{G}}(t).
			\end{equation}
			Applying the same proposition to the constant test function $1$ and then using the evaluation of the Gaussian integral from Theorem~\ref{collisionThm} gives, for every $t>\mathrm{T}_{\mathrm{KM}}(1,\mathbf{x})$,
			\begin{equation}\label{lastcor-step3}
				\left|
				\frac{\mathbb{P}_{\mathbf{x}}^N\!\left[\mathcal{T}>t\right]}
				{\mathfrak{h}_N(\mathbf{x})\Sigma(t)^{-\mathrm{d}_N}}
				-\mathrm{L}_{0}
				\right|
				\leq \mathrm{C}_{\mathrm{KM}}(1,\mathbf{x})\psi_{\mathrm{G}}(t).
			\end{equation}
			By \eqref{eq:psiG-def}, $\psi_{\mathrm{G}}(t)\to0$. Hence we may choose $\mathrm{T}_{\mathrm{GOE},F,\mathbf{x}}\in\mathscr{T}$ larger than the lower times in \eqref{lastcor-step2} and \eqref{lastcor-step3}, and satisfying $\mathrm{C}_{\mathrm{KM}}(1,\mathbf{x})\psi_{\mathrm{G}}(t)\leq \mathrm{L}_{0}/2$ for every $t>\mathrm{T}_{\mathrm{GOE},F,\mathbf{x}}$. For those $t$, the normalized denominator in \eqref{lastcor-step1} is at least $\mathrm{L}_{0}/2$, and the elementary identity for ratios gives
			\begin{multline*}
				\left|
				\mathbb{E}_{\mathbf{x}}^N\!\left[F\!\left(\mathbf{U}(t)\right)\middle|\mathcal{T}>t\right]
				-\frac{\mathrm{L}_{F}}{\mathrm{L}_{0}}
				\right|\\
				\leq
				\frac{2}{\mathrm{L}_{0}}
				\left|
				\frac{\mathbb{E}_{\mathbf{x}}^N\!\left[F\!\left(\mathbf{U}(t)\right)\mathbf{1}_{\{\mathcal{T}>t\}}\right]}
				{\mathfrak{h}_N(\mathbf{x})\Sigma(t)^{-\mathrm{d}_N}}
				-\mathrm{L}_{F}
				\right|
				+
				\frac{2|\mathrm{L}_{F}|}{\mathrm{L}_{0}^{2}}
				\left|
				\frac{\mathbb{P}_{\mathbf{x}}^N\!\left[\mathcal{T}>t\right]}
				{\mathfrak{h}_N(\mathbf{x})\Sigma(t)^{-\mathrm{d}_N}}
				-\mathrm{L}_{0}
				\right|.
			\end{multline*}
			Substituting the bounds \eqref{lastcor-step2} and \eqref{lastcor-step3} into this ratio estimate gives
			\begin{equation*}
				\left|\mathbb{E}_{\mathbf{x}}^N\!\left[F\!\left(\mathbf{U}(t)\right)\middle|\mathcal{T}>t\right]-\frac{\mathrm{L}_{F}}{\mathrm{L}_{0}}\right|\leq\left(\frac{2\mathrm{C}_{\mathrm{KM}}(F,\mathbf{x})}{\mathrm{L}_{0}}+\frac{2|\mathrm{L}_{F}|\mathrm{C}_{\mathrm{KM}}(1,\mathbf{x})}{\mathrm{L}_{0}^{2}}\right)\psi_{\mathrm{G}}(t).
			\end{equation*}
			Thus \eqref{karlinMcGregor} holds with $\mathrm{C}_{\mathrm{GOE},F,\mathbf{x}}=2\mathrm{C}_{\mathrm{KM}}(F,\mathbf{x})/\mathrm{L}_{0}+2|\mathrm{L}_{F}|\mathrm{C}_{\mathrm{KM}}(1,\mathbf{x})/\mathrm{L}_{0}^{2}$. Finally, the identity $\mathrm{L}_{F}/\mathrm{L}_{0}$ equals the displayed $\mathrm{GOE}_{N}$ integral in \eqref{karlinMcGregor}, because $\mathrm{L}_{0}^{-1}=N!\prod_{i=1}^{N}\Gamma(3/2)/\Gamma(1+i/2)$. This proves the proposition.
	    \end{proof}
       
	    \begin{proof}[Proof of Theorem~\ref{thm:stationarypoints}]
		   Fix $i\in\llbracket N\rrbracket$. By Proposition~\ref{moments}, applied to the one-particle path $\mathscr{X}_i$, there is an almost surely finite random time after which $\mathscr{X}_i(t)\in\mathscr{C}_{t}$ for all later $t$. Increasing this random time if necessary so that $t>\mathrm{T}_{\mathrm{A}}$, Proposition~\ref{axpro} applies directly to the realized terminal endpoint $\mathscr{X}_i(t)$ at every later time. Hence, after this almost surely finite random time, the minimizer in Definition~\ref{definitionstationary} is unique and equals $\mathrm{u}_{\mathscr{X}_i}(t)$, so $\mathscr{U}_i(t)=-\Sigma(t)\mathrm{u}_{\mathscr{X}_i}(t)$. This proves the finite-random-time uniqueness assertion for each coordinate.
		
         Now let $F:\mathbb{R}^{N}\to\mathbb{R}$ be bounded and Lipschitz. The function $F$ satisfies the growth and local Lipschitz hypotheses of Proposition~\ref{NormalisedProcess}, with polynomial exponent $q=0$ after increasing the local Lipschitz constant to dominate the uniform bound. Therefore \eqref{eq:free-U-Gaussian-limit} gives a bound of the form $\mathrm{C}_{F,\mathbf{x}}\psi_{\mathrm{G}}(t)$ for the expectation restricted to $\{\mathbf{X}(t)\in\mathscr{C}_{t}^{N}\}$. The omitted endpoint complement contributes at most $\|F\|_{\infty}\mathbb{P}_{\mathbf{x}}^N[\mathbf{X}(t)\notin\mathscr{C}_{t}^{N}]$, and the one-time endpoint-window estimate from Proposition~\ref{localizationtail} (applied coordinatewise and combined with the finite union bound, as in the final step of Proposition~\ref{NormalisedProcess}) bounds this probability by another constant multiple of $\psi_{\mathrm{G}}(t)$ after increasing the lower time. Absorbing these two constants into $\mathrm{C}_{\mathrm{free},F,\mathbf{x}}$ proves the free quantitative bound. For real-valued bounded Lipschitz $F$, the survival-conditioned statement is precisely Proposition~\ref{lastcorcor}; its limiting integral is the ordered $\mathrm{GOE}_N$ eigenvalue density from Definition~\ref{GOEeigenvalue}. Combining these conclusions proves the theorem.
	    \end{proof}
	    \subsection{Proof of Corollary~\ref{cor:particlecoordinates}}
	        \begin{proof}[Proof of Corollary~\ref{cor:particlecoordinates}]
	        We first prove the quantitative comparison between the stationary-point coordinates and the centered inverse-rate coordinates on the endpoint window. This is the only additional estimate needed to transfer Theorem~\ref{thm:stationarypoints} from $\mathbf{U}(t)$ to $\mathbf{V}(t)$. Since the comparison is also used inside the leading term of the Karlin--McGregor expansion, we prove it for an arbitrary deterministic initial vector $\mathbf{a}\in\mathbb{Z}_{+}^{N}$ and take expectations with respect to $\mathbb{P}_{\mathbf{a}}^{N}$ throughout this proof.
	        
            For one coordinate, applying \eqref{expressionforthefirst} to the random endpoint path $y=\mathscr{X}_i$ gives the exact identity
	        \begin{equation}\label{eq:coord-V-F1-identity}
	            \frac{1}{\Sigma(t)}\left(\sum_{n=0}^{\mathscr{X}_i(t)}\frac{1}{\lambda_n}-\Lambda\mathscr{Z}(t)\right)=\frac{\mathfrak{F}_{\mathscr{X}_i}^{(1)}(0,t)}{\Sigma(t)}+\frac{\Lambda(\mathrm{z}(t)-\mathscr{Z}(t))}{\Sigma(t)}.
	        \end{equation}
	        Since $\mathscr{Z}(t)=\lfloor\mathrm{z}(t)\rceil$, the last term in \eqref{eq:coord-V-F1-identity} has absolute value at most $\Lambda/(2\Sigma(t))$. On the event $\{\mathbf{X}(t)\in\mathscr{C}_{t}^{N}\}$, each coordinate endpoint belongs to $\mathscr{C}_{t}$, so the comparison estimate \eqref{eq:Uidentify} applies to every coordinate after increasing the deterministic lower time. Combining \eqref{eq:coord-V-F1-identity}, the rounding bound, and \eqref{eq:Uidentify}, we obtain, for every $i\in\llbracket N\rrbracket$ and every $t$ beyond this deterministic lower time,
	        \begin{equation}\label{eq:coord-pointwise-VU}
	            \left|V_i(t)-\mathscr{U}_i(t)\right|\mathbf{1}_{\{\mathbf{X}(t)\in\mathscr{C}_{t}^{N}\}}\leq\left(\frac{\Lambda}{2\Sigma(t)}+\frac{4\mathrm{C}_{\mathrm{A}}\left|\mathfrak{F}_{\mathscr{X}_i}^{(1)}(0,t)\right|^2}{\Sigma(t)^3}+\frac{2\mathrm{C}_{\mathscr{C}}\gamma(t)^{1/2}\left|\mathfrak{F}_{\mathscr{X}_i}^{(1)}(0,t)\right|}{\Sigma(t)}\right)\mathbf{1}_{\{\mathbf{X}(t)\in\mathscr{C}_{t}^{N}\}}.
	        \end{equation}
	        This pointwise estimate becomes an expectation estimate by using the moment input already stated in Proposition~\ref{boundedmoments}. More precisely, \eqref{eq:first-derivative-moment} gives first-derivative moments. We also use the deterministic lower comparison between the fluctuation scale and time: the linear lower bound on $\mathscr{Z}(t)$ from \eqref{Zlinearintbounds}, together with the upper comparison between $\mathscr{Z}(t)$ and $\Sigma(t)^2$ in \textbf{(iv)} of Proposition~\ref{specificcasesub}, gives a constant $\mathrm{c}_{\Sigma}>0$ and a deterministic lower time such that
	        \begin{equation*}
	            \Sigma(t)^2\geq \mathrm{c}_{\Sigma}(1+t)
	        \end{equation*}
	        for every later $t$. Hence, for each coordinate,
	        \begin{equation*}
	            \mathbb{E}_{\mathbf{a}}^{N}\left[\frac{\left|\mathfrak{F}_{\mathscr{X}_i}^{(1)}(0,t)\right|^2}{\Sigma(t)^3}\right]\leq\frac{\mathrm{C}_{2,\mathbf{a}}(1+t)}{\Sigma(t)^3}\leq\mathrm{C}_{2,\mathbf{a}}\mathrm{c}_{\Sigma}^{-1}\Sigma(t)^{-1},
	        \end{equation*}
	        and, by the Cauchy--Schwarz inequality followed by the same scale comparison,
	        \begin{equation*}
	            \mathbb{E}_{\mathbf{a}}^{N}\left[\frac{\left|\mathfrak{F}_{\mathscr{X}_i}^{(1)}(0,t)\right|}{\Sigma(t)}\right]\leq\frac{\mathrm{C}_{2,\mathbf{a}}^{1/2}(1+t)^{1/2}}{\Sigma(t)}\leq\mathrm{C}_{2,\mathbf{a}}^{1/2}\mathrm{c}_{\Sigma}^{-1/2}.
	        \end{equation*}
	        Summing \eqref{eq:coord-pointwise-VU} over $i$ and absorbing the finite coordinate sum into the constant gives constants $\mathrm{T}_{\mathrm{VU},1,\mathbf{a}}\in\mathscr{T}$ and $\mathrm{C}_{\mathrm{VU},1,\mathbf{a}}>0$ such that, for every $t>\mathrm{T}_{\mathrm{VU},1,\mathbf{a}}$,
	        \begin{equation}\label{eq:VU-explicit-bound}
	            \mathbb{E}_{\mathbf{a}}^{N}\left[|\mathbf{V}(t)-\mathbf{U}(t)|\mathbf{1}_{\{\mathbf{X}(t)\in\mathscr{C}_{t}^{N}\}}\right]\leq \mathrm{C}_{\mathrm{VU},1,\mathbf{a}}\left(\Sigma(t)^{-1}+\gamma(t)^{1/2}\right).
	        \end{equation}
	        The same pointwise estimate gives the square-integrable form needed in the determinant argument. We square \eqref{eq:coord-pointwise-VU}, use $(a+b+c)^2\leq3(a^2+b^2+c^2)$ and $|\mathbf{v}|^2\leq N\sum_{i=1}^{N}|v_i|^2$, and then apply \eqref{eq:first-derivative-moment} with exponents $2$ and $4$. This gives a constant $\mathrm{C}_{\mathrm{VU},2,\mathbf{a}}>0$ such that, for every $t$ beyond the deterministic lower time used in the preceding estimates,
	        \begin{equation*}
	            \mathbb{E}_{\mathbf{a}}^{N}\left[|\mathbf{V}(t)-\mathbf{U}(t)|^2\mathbf{1}_{\{\mathbf{X}(t)\in\mathscr{C}_{t}^{N}\}}\right]\leq \mathrm{C}_{\mathrm{VU},2,\mathbf{a}}^{2}\left(\Sigma(t)^{-2}+\gamma(t)\right).
	        \end{equation*}
	        Therefore,
	        \begin{equation}\label{eq:VU-L2-bound}
	            \left(\mathbb{E}_{\mathbf{a}}^{N}\left[|\mathbf{V}(t)-\mathbf{U}(t)|^2\mathbf{1}_{\{\mathbf{X}(t)\in\mathscr{C}_{t}^{N}\}}\right]\right)^{1/2}\leq \mathrm{C}_{\mathrm{VU},2,\mathbf{a}}\left(\Sigma(t)^{-1}+\gamma(t)^{1/2}\right).
	        \end{equation}
	        Finally, Definition~\ref{gammadef} gives $\gamma(t)\geq\Sigma(t)^{-1/2}$, while \eqref{eq:psiG-def} and the fact that $\gamma(t)\to0$ imply, after increasing the lower time, that both $\Sigma(t)^{-1}$ and $\gamma(t)^{1/2}$ are bounded by $\psi_{\mathrm{G}}(t)$. Hence, for the original initial state $\mathbf{x}$, \eqref{eq:VU-explicit-bound} gives
	        \begin{equation}\label{eq:VU-explicit-bound-psi}
	            \mathbb{E}_{\mathbf{x}}^{N}\left[|\mathbf{V}(t)-\mathbf{U}(t)|\mathbf{1}_{\{\mathbf{X}(t)\in\mathscr{C}_{t}^{N}\}}\right]\leq2\mathrm{C}_{\mathrm{VU},1,\mathbf{x}}\psi_{\mathrm{G}}(t),
	        \end{equation}
	        and, for every fixed initial vector $\mathbf{a}\in\mathbb{Z}_{+}^{N}$, \eqref{eq:VU-L2-bound} gives
	        \begin{equation}\label{eq:VU-L2-bound-psi}
	            \left(\mathbb{E}_{\mathbf{a}}^{N}\left[|\mathbf{V}(t)-\mathbf{U}(t)|^2\mathbf{1}_{\{\mathbf{X}(t)\in\mathscr{C}_{t}^{N}\}}\right]\right)^{1/2}\leq2\mathrm{C}_{\mathrm{VU},2,\mathbf{a}}\psi_{\mathrm{G}}(t).
	        \end{equation}
	        We also record the endpoint-complement estimate in the same place. Proposition~\ref{localizationtail}, applied coordinatewise and combined with the finite union bound, gives a constant $\mathrm{C}_{\mathrm{tail},\mathbf{a}}>0$ such that, after increasing the lower time if necessary,
	        \begin{equation}\label{eq:coord-tail-bound}
	            \mathbb{P}_{\mathbf{a}}^{N}[\mathbf{X}(t)\notin\mathscr{C}_{t}^{N}]\leq\mathrm{C}_{\mathrm{tail},\mathbf{a}}\psi_{\mathrm{G}}(t)
	        \end{equation}
	        for every later $t$.
	      
             We can now finish the free coordinate limit. Taking $\mathbf{a}=\mathbf{x}$ in \eqref{eq:VU-explicit-bound-psi} and \eqref{eq:coord-tail-bound}, the Lipschitz and boundedness assumptions on $F$ give
	        \begin{equation*}
	            \left|\mathbb{E}_{\mathbf{x}}^{N}\left[F(\mathbf{V}(t))\right]-\mathbb{E}_{\mathbf{x}}^{N}\left[F(\mathbf{U}(t))\right]\right|\leq2\mathrm{Lip}(F)\mathrm{C}_{\mathrm{VU},1,\mathbf{x}}\psi_{\mathrm{G}}(t)+2\|F\|_{\infty}\mathrm{C}_{\mathrm{tail},\mathbf{x}}\psi_{\mathrm{G}}(t).
	        \end{equation*}
	        Combining this estimate with the free part of Theorem~\ref{thm:stationarypoints} proves the first assertion of Corollary~\ref{cor:particlecoordinates}, after absorbing the displayed constants into $\mathrm{C}_{\mathrm{coord},F,\mathbf{x}}$.
	     
          It remains to identify the survival-conditioned limit of $\mathbf{V}(t)$. We apply the same determinant expansion as in Proposition~\ref{ThmKM}, but with the bounded Lipschitz factor $F(\mathbf{U}_{\mathbf{m}}(t))$ replaced by $F(\mathbf{V}_{\mathbf{m}}(t))$, where
	        \begin{equation*}
	            \mathbf{V}_{\mathbf{m}}(t)\overset{\textnormal{def}}{=}\frac{1}{\Sigma(t)}\left(\Lambda_N(\mathbf{m})-\Lambda\mathscr{Z}(t)\,\underline{1}\right).
	        \end{equation*}
	        This substitution does not alter the algebraic part of the proof. The endpoint-complement estimate in Step~1 of Proposition~\ref{ThmKM}, and the non-leading bounds in Steps~2--4, use the test function only through $\|F\|_{\infty}$. The same bounds therefore apply to $F(\mathbf{V}_{\mathbf{m}}(t))$. Thus the expansion \eqref{expectform} is replaced by the identical formula with $F(\mathbf{V}(t))$ in place of $F(\mathbf{U}(t))$, except that the final expectation becomes
	        \begin{equation*}
	            \mathbb{E}_{\underline{0}}^N\!\left[\Delta(\mathbf{U}(t))F(\mathbf{V}(t))\mathbf{1}_{\{\mathbf{U}(t)\in\mathbf{W}_{\mathrm c}^{N},\,\mathbf{X}(t)\in\mathscr{C}_{t}^{N}\}}\right].
	        \end{equation*}
	        We compare this with the leading expectation already evaluated in Proposition~\ref{ThmKM}. The Lipschitz bound gives
	        \begin{equation*}
	            \left|\Delta(\mathbf{U}(t))F(\mathbf{V}(t))-\Delta(\mathbf{U}(t))F(\mathbf{U}(t))\right|\leq\mathrm{Lip}(F)|\Delta(\mathbf{U}(t))|\,|\mathbf{V}(t)-\mathbf{U}(t)|.
	        \end{equation*}
	        Since $\Delta$ is a polynomial of degree $\mathrm{d}_{N}$, Proposition~\ref{boundedmoments}, applied with any $p\geq\max(2,2\mathrm{d}_{N})$, gives a uniform bound for the second moment of $|\Delta(\mathbf{U}(t))|$ on $\{\mathbf{X}(t)\in\mathscr{C}_{t}^{N}\}$. Hence the Cauchy--Schwarz inequality and \eqref{eq:VU-L2-bound-psi}, used with $\mathbf{a}=\underline{0}$, imply that
	        \begin{equation*}
	            \mathbb{E}_{\underline{0}}^{N}\left[|\Delta(\mathbf{U}(t))|\,|\mathbf{V}(t)-\mathbf{U}(t)|\mathbf{1}_{\{\mathbf{X}(t)\in\mathscr{C}_{t}^{N}\}}\right]\leq\mathrm{C}_{\Delta,\mathrm{VU}}\psi_{\mathrm{G}}(t)
	        \end{equation*}
	        for a constant $\mathrm{C}_{\Delta,\mathrm{VU}}>0$ and every $t$ beyond the deterministic lower time fixed above. 
	        Combining this estimate with Proposition~\ref{ThmKM} gives
	        \begin{align}\label{eq:coord-killed-numerator}
	            &\left|
	            \frac{\mathbb{E}_{\mathbf{x}}^N\!\left[F(\mathbf{V}(t))\mathbf{1}_{\{\mathcal{T}>t\}}\right]}
	            {\mathfrak{h}_N(\mathbf{x})\Sigma(t)^{-\frac{N}{2}(N-1)}}
	            -
	            \frac{1}{(2\pi)^{N/2}}
	            \int_{\mathbf{W}_{\mathrm c}^{N}}F(\mathbf{u})\Delta(\mathbf{u})
	            \exp\!\Big(-\tfrac12|\mathbf{u}|^2\Big)\,\mathrm{d}\mathbf{u}
	            \right|\\
	            &\qquad\leq
	            \mathrm{C}_{\mathrm{coordKM},F,\mathbf{x}}\psi_{\mathrm{G}}(t)
	            \notag
	        \end{align}
	        for a constant $\mathrm{C}_{\mathrm{coordKM},F,\mathbf{x}}>0$ and all $t$ beyond a deterministic lower time.
	        The denominator is controlled by Proposition~\ref{ThmKM} with the constant test function $1$, and the limiting constant is
	        \begin{equation*}
	            \mathrm{L}_{0}=\frac{1}{N!}\prod_{i=1}^{N}\frac{\Gamma(1+i/2)}{\Gamma(3/2)}.
	        \end{equation*}
	        Thus, for every $t$ beyond the lower time in Proposition~\ref{ThmKM},
	        \begin{equation*}
	            \left|\frac{\mathbb{P}_{\mathbf{x}}^N[\mathcal{T}>t]}{\mathfrak{h}_N(\mathbf{x})\Sigma(t)^{-\frac{N}{2}(N-1)}}-\mathrm{L}_{0}\right|\leq\mathrm{C}_{\mathrm{KM}}(1,\mathbf{x})\psi_{\mathrm{G}}(t).
	        \end{equation*}
	        Increase the lower time so that $\mathrm{C}_{\mathrm{KM}}(1,\mathbf{x})\psi_{\mathrm{G}}(t)\leq\mathrm{L}_{0}/2$. The elementary ratio estimate displayed in the proof of Proposition~\ref{lastcorcor}, applied with the numerator bound \eqref{eq:coord-killed-numerator}, then gives the conditional coordinate estimate with
	        \begin{equation*}
	            \mathrm{C}_{\mathrm{coordGOE},F,\mathbf{x}}\overset{\textnormal{def}}{=}\frac{2\mathrm{C}_{\mathrm{coordKM},F,\mathbf{x}}}{\mathrm{L}_{0}}+\frac{2|\mathrm{L}_{F}|\mathrm{C}_{\mathrm{KM}}(1,\mathbf{x})}{\mathrm{L}_{0}^{2}}.
	        \end{equation*}
	        The quotient is therefore controlled quantitatively by the ordered $\mathrm{GOE}_N$ eigenvalue law from Definition~\ref{GOEeigenvalue}. This proves the conditional coordinate estimate and completes the proof.
		\end{proof}

\bibliographystyle{plainurl}
\bibliography{References}
\end{document}